\documentclass[a4paper, 11pt]{article}

\usepackage{longtable}
\usepackage{geometry}
\usepackage{graphicx}
\usepackage{subfig}
\usepackage{color}
\usepackage{xcolor}
\usepackage{bm}
\usepackage{fixmath}
\usepackage{mathrsfs}
\usepackage{ulem}
\usepackage{multirow}
\usepackage{pdflscape}
\usepackage{rotating}
\usepackage{lscape}
\usepackage{pdflscape}
\usepackage{titlesec}
\usepackage{float}
\usepackage{hyperref}

\titleclass{\subsubsubsection}{straight}[\subsection]

\newcounter{subsubsubsection}[subsubsection]
\renewcommand\thesubsubsubsection{\thesubsubsection.\arabic{subsubsubsection}}

\titleformat{\subsubsubsection}
  {\normalfont\normalsize\bfseries}{\thesubsubsubsection}{1em}{}
\titlespacing*{\subsubsubsection}
{0pt}{3.25ex plus 1ex minus .2ex}{1.5ex plus .2ex}

\makeatletter

\renewcommand\paragraph{\@startsection{paragraph}{5}{\z@}%
  {3.25ex \@plus1ex \@minus.2ex}%
  {-1em}%
  {\normalfont\normalsize\bfseries}}

\renewcommand\subparagraph{\@startsection{subparagraph}{6}{\parindent}%
  {3.25ex \@plus1ex \@minus .2ex}%
  {-1em}%
  {\normalfont\normalsize\bfseries}}

\def\toclevel@subsubsubsection{4}
\def\toclevel@paragraph{5}
\def\toclevel@paragraph{6}
\def\l@subsubsubsection{\@dottedtocline{4}{7em}{4em}}
\def\l@paragraph{\@dottedtocline{5}{10em}{5em}}
\def\l@subparagraph{\@dottedtocline{6}{14em}{6em}}

\makeatother

\def\b{\begin{eqnarray}}
\def\e{\end{eqnarray}}

\makeatletter
\newcommand*\bigdot{\mathpalette\bigdot@{.5}}
\newcommand*\bigdot@[2]{\mathbin{\vcenter{\hbox{\scalebox{#2}{$\m@th#1\bullet$}}}}}
\makeatother

\begin{document}

\begin{center}
{\huge \textbf{Classification of the Real Roots of the \vskip.1cm Quartic Equation and their Pythagorean \vskip.2cm Tunes}}

\vspace{9mm}
\noindent
{\large \bf Emil M. Prodanov} \vskip.4cm
{\it School of Mathematical Sciences, Technological University Dublin,
\vskip.1cm
City Campus, Kevin Street, Dublin, D08 NF82, Ireland,}
\vskip.1cm
{\it E-Mail: emil.prodanov@tudublin.ie} \\
\vskip.5cm
\end{center}

\begin{abstract}
\noindent
Presented is a two-tier analysis of the location of the real roots of the general quartic equation $x^4 \, + \, ax^3 \, + \, bx^2 \, + \, cx \, + \, d \, = \, 0$ with real coefficients and the classification of the roots in terms of $a$, $b$, $c$, and $d$, without using any numerical approximations. Associated with the general quartic, there is a number of subsidiary quadratic equations ({\it resolvent quadratic equations}) whose roots allow this systematization as well as the determination of the bounds of the individual roots of the quartic. In many cases the root isolation intervals are found. The second tier of the analysis uses two subsidiary cubic equations ({\it auxiliary cubic equations}) and solving these, together with some of the {\it resolvent quadratic equations}, allows the full classification of the roots of the general quartic and also the determination of the isolation interval of each root. These isolation intervals involve the stationary points of the quartic (among others) and, by solving some of the {\it resolvent quadratic equations}, the isolation intervals of the stationary points of the quartic are also determined. The presented classification  of the roots of the quartic equation is particularly useful in situations in which the equation stems from a model the coefficients of which are (functions of) the model parameters and solving cubic equations, let alone using the explicit quartic formulas, is a daunting task. The only benefit in such cases would be to gain insight into the location of the roots and the proposed method provides this. Each possible case has been carefully studied and illustrated with a detailed figure containing a description of its specific characteristics, analysis based on solving cubic equations and analysis based on solving quadratic equations only. As the analysis of the roots of the quartic equation is done by studying the intersection points of the ``sub-quartic" $x^4 + ax^3 + bx^2$ with a set of suitable parallel lines, a beautiful Pythagorean analogy can be found between these intersection points and the set of parallel lines on one hand and the musical notes and the staves representing different musical pitches on the other: each particular case of the quartic equation has its own short tune.
\end{abstract}


\newpage

\section{Introduction}

\noindent
The process of solving the general quartic equation
\b
\label{q}
x^4 + ax^3 + bx^2 + cx + d = 0
\e
involves the removal of the cubic term by using the substitution $x = y - a/4$. This results in the depressed quartic $y^4 + p y^2 + q y + r = 0$, where $p = b - (3/2)a^2/4, \,\, q  =  c - (a/2) (b - a^2/4),$ and $r =  d - (3/256) a^4 + (1/16) a^2 b - (1/4) ac$. There are several algorithms for solving the depressed quartic equation and each of these involves solving a cubic equation called {\it resolvent}. The resolvents are different for the different algorithms. Finding the roots of the original quartic equation, once the roots of the resolvent cubic are known, is a straightforward procedure in each algorithm. \\
Using the methods presented in \cite{1}, this paper presents a two-tier analysis of the location of the real roots of the general quartic equation (\ref{q}) with real coefficients and their classification.  \\
Associated with each quartic, there is a number of subsidiary quadratic equations, referred to in this text as {\it resolvent quadratic equations}, with the help of which the roots of the general quartic are systematized in terms of its coefficients $a, b, c,$ and $d$. Additionally, the bounds on individual roots are determined. In many cases the root isolation intervals are found, but there is some residual ambiguity as some intervals may contain either two roots of the quartic or no roots at all. The individual root bounds themselves are associated with the roots of $x^2 + ax + b = 0$ and $cx + d = 0$, or with the non-zero stationary points of $x^2(x^2 + ax + b)$, if real roots are absent, or with the points of curvature change of $x^2(x^2 + ax + b)$, if real roots and non-zero stationary points are absent, or with the point of vanishing third derivative of $x^2(x^2 + ax + b)$, if real roots and non-zero stationary points, and curvature change points are absent. \\
Using the proposed classification of the roots of the quartic in terms of the coefficients of the quartic and based on solving quadratic equations only would be particularly useful when the quartic results from the study of some model and the coefficients of the quartic are functions of the model parameters. Even the presence of a single parameter in the equation makes the application of the cubic formulas, let alone the quartic formulas, practically impossible and the only benefit would be to get insight into the location of the roots. \\
The other tier of the presented analysis is the full classification of the roots of the quartic in terms of its coefficients and, also, the determination of the isolation interval of each root of the quartic by solving two cubic equations and some of the resolvent quadratic equations. These two cubic equations are subsidiary to each quartic and are referred to in this paper as {\it auxiliary cubic equations} --- in order to make a distinction from the resolvent cubic equation. \\
It can be argued why it is necessary to address not one, but two cubic equations (the {\it auxiliary cubic equations}) which do not yield the roots of the quartic, but only reveal their isolation intervals, and what can be gained by doing so. Indeed, it suffices to solve a single cubic equation (the resolvent) in order to find the actual roots of any quartic. The answer to this is in the search for rules and patterns through abstraction to gain insight on how different coefficient vectors affect the roots. Systematization and having predictive powers are always a bonus and merit investigation. If analysis of the roots of the quartic could be done with equations of degree two, one would expect that using equations of degree three would be more informative, despite getting into the realm where the procedure of finding the roots through the explicit formulas is applied. The analysis based on solving cubic equations complements the picture revealed by analysis based on quadratic equations. The resulting systematization and classification are not possible if one addresses the explicit formulas for the roots of the quartic --- these are rather unwieldy and one cannot trace how the variation of a particular coefficient of the given quartic affects the roots, i.e. the coefficients of the quartic enter the root formulas in an intricate combination involving the root(s) of the resolvent cubic (which, in turn, depend on the coefficients of the quartic) and one cannot discern the individual contribution of each coefficient of the quartic to the location of its roots. The proposed analysis based on cubic equations also has heuristic potential: just by observation of the coefficients and whether they fall into specific ranges, one can predict the number of real roots and also find their isolation intervals. For example, for any $a, c,$ and $d$, when $b > (3/2) a^2$, the quartic cannot have four real roots. If, additionally, the free term $d$ is negative and $c$ is also negative, then the quartic has one negative root smaller than $-d/c$ and one positive root greater than the only real non-zero root $\lambda > 0$ of the equation obtained from the quartic after removing the free term $d$. If $\mu$ denotes the only stationary point of quartic ($\mu > 0$ in this case) and $d$ is positive and greater than $\mu^4 + a \mu^3 + b \mu^2 + c \mu$ (which itself is positive), then the quartic has no real roots. If $0 < d  < \mu^4 + a \mu^3 + b \mu^2 + c \mu$, then the quartic has two positive roots --- one bigger than $-d/c$ and smaller than $\mu$, the other -- bigger than $\mu$ and smaller than $\lambda$. \\
In the presented analysis, the isolation intervals of the stationary points of the general quartic are also determined --- with the help of the {\it resolvent quadratic equations}. \\
The analysis of the roots of the quartic equation is done after studying the intersection points of the ``sub-quartic" $x^4 + ax^3 + bx^2$ with a set of suitable parallel lines. A beautiful Pythagorean analogy exists between these intersection points and the set of parallel lines on one hand and the musical notes and the staves representing different musical pitches on the other. Even more, each case of the quartic equation has its own tune. \\
Every possible situation has been individually studied and illustrated with a detailed figure containing a description of its specific characteristics, analysis based on solving cubic equations and analysis based on solving quadratic equations only.

\section{Subsidiary Cubic and Quadratic and Linear Equations}

\noindent
Each quartic is associated with a number of subsidiary cubic, quadratic, and linear equations whose roots can be used for the classification of the roots of the quartic. \\
In parallel with the presentation of the analysis, Figures 1 to 5 illustrate, through a particular example with the quartic $x^4 + x^3 - 3x^2 - x + 1$, the full procedure of finding the isolation interval of each root of the quartic, based on solving cubic equations, and the localization of the roots by determination of the individual root bounds, based on solving quadratic equations only. The isolation intervals of the stationary points of this quartic are also found. \\
Figures 6 to 12 illustrate some patterns associated with the general quartic and these are used for the classification of the roots of the quartic.
\vskip.5cm
\noindent
Taking the free term $d$ of the quartic and varying it, yields a one-parameter congruence of quartics, all having the same set of stationary points (which are either three --- two local minima and a local maximum or a saddle point and a local minimum, or just one --- a minimum). Let $\mu_i$ denote the stationary point(s). For each $\mu_i$, there is a ``special" quartic within this congruence: $x^4 + ax^3 + bx^2 + cx + \delta_i$ --- the one whose graph is tangent to the abscissa at that particular stationary point. The derivative of the ``special" quartic is also zero at this stationary point, namely, for $x^4 + ax^3 + bx^2 + cx + \delta_i$ the stationary point $\mu_i$ is also a double root [or a triple root, if the original quartic has a saddle at $\mu_i$, or a quadruple root $-a/4$ in the case of the quartic $(x - a/4)^4 = x^4 + ax^3 + (3/8)a^2 x^2 + (1/16) a^3 x + (1/256) a^4$, which coincides with its only ``special" quartic]. Setting the derivative of the quartic equal to zero yields the set of its stationary points and the resulting equation,
\b
\label{1st}
4x^3 + 3ax^2 + 2bx + c = 0,
\e
is referred to in this text as {\it first auxiliary cubic equation}. \\
Substituting each real root $\mu_i$ of this equation into the corresponding ``special" quartic equation $x^4 + ax^3 + bx^2 + cx + \delta_i = 0$, immediately gives:
\b
\label{delta}
\delta_i = - \mu_i^4 - a \mu_i^3 - b \mu_i^2 - c \mu_i.
\e
Thus the ``special" quartics are given by $x^4 - \mu_i^4 + a(x^3 - \mu_i^3) + b(x^2 - \mu_i^2)+ c(x - \mu_i)$.  \\
The discriminant of the first auxiliary cubic equation is $\Delta_1 = - 432c^2 - 432a(a^2/4 - b)c + 128b^2 [(9/8)a^2/4 - b]$. It can be viewed as a quadratic in $c$, treated as unknown, with $a$ and $b$ treated as parameters. The {\it first resolvent quadratic equation} is obtained by setting $\Delta_1 = 0$:
\b
\label{qr1}
c^2 + a \left( \frac{a^2}{4} - b \right) c  - \frac{1}{4} b^2 \left( \frac{9}{8} \frac{a^2}{4} - b \right) = 0.
\e
The roots of this equation,
\b
\label{c12}
c_{1,2}(a,b) = c_0 \pm \frac{2\sqrt{6}}{9} \, \sqrt{ \left( \frac{3}{2}\frac{a^2}{4} - b \right)^3}, \quad \mbox{with} \quad c_0(a,b) = \frac{1}{2} a \left( b - \frac{a^2}{4} \right),
\e
play a very important role in the analysis. For any given quartic, one has to see first whether the coefficient $c$ of the linear term falls between the roots $c_{1,2}(a,b)$ or outside them. If $c_2(a,b) \le c \le c_1(a, b)$, then the discriminant $\Delta_1$ is positive and the quartic has three stationary points: $\mu_{1,2,3}$. Otherwise it has just one: $\mu_1$. In the first case, the quartic can have either 0, or 2, or 4 real roots; in the second case it can have either 0 or 2 real roots. It is immediately obvious that, for any $a, c, $ and $d$, when $b > (3/2)(a^2/4)$, the quartic can have either 0 or 2 real roots only (the first resolvent quadratic equation (\ref{qr1}) has negative discriminant). \\
When the discriminant $\Delta_1$ of the first auxiliary cubic equation (\ref{1st}) is equal to zero, that is, when $c = c_{1,2}(a,b)$, the original quartic with $c$ replaced by $c_{1,2}$ has a saddle point at $\eta_{1,2}$ and a local minimum at $\theta_{1,2}$. The corresponding ``special" quartic at $\eta_{1,2}$ is $x^4 + ax^3 + bx^2 + c_{1,2}x + d_{1,2}$, where $d_{1,2} = - \eta_{1,2}^4 - a \eta_{1,2}^3 - b \eta_{1,2}^2 - c_{1,2} \eta_{1,2}$, and for the ``special" quartic, $\eta_{1,2}$ is a triple root. The points $\eta_{1,2}$ and $\theta_{1,2}$ can be easily found as the ``special" quartic $x^4 + ax^3 + bx^2 + c_{1,2}x + d_{1,2}$, its first derivative, and its second derivative are all zero at $\eta_{1,2}$. In other words, one has to start with solving the {\it second resolvent quadratic equation},
\b
\label{qr2}
6x^2 + 3ax + b = 0,
\e
the roots of which are
\b
\label{eta}
\eta_{1,2} = -\frac{1}{4} a \pm \frac{\sqrt{6}}{6} \, \sqrt{\frac{3}{2}\frac{a^2}{4} - b},
\e
then write down the vanishing first derivative of the (``special") quartic as $4x^3 + 3ax^2 + 2bx + c_{1,2} = 4(x-\eta_{1,2})^2(x - \theta_{1,2})$, and then compare the coefficients of the quadratic terms. This will give $\theta_{1,2} = -(3/4) a - 2 \eta_{1,2}$ and hence:
\b
\label{teta}
\theta_{1,2} = -\frac{1}{4} a \pm \frac{\sqrt{6}}{3} \, \sqrt{\frac{3}{2}\frac{a^2}{4} - b}.
\e
One could observe that $-a/4$ is the quadruple root of the quartic $x^4 + ax^3 + (3/8)a^2 x^2 + (1/16) a^3 x + (1/256) a^4 = 0$. \\
With the help of $\eta_{1,2}$ and $\theta_{1,2}$, the isolation intervals of the stationary points $\mu_i$ of the general quartic can be easily found (see the example on Figure 6). In the regime of increasing $c$ and starting with $c < c_2$, there is only one stationary point (local minimum) at $\widehat{\mu}_1 > \theta_2$. When $c = c_2$, the quartic has a saddle point at $\eta_2$ and a local minimum at $\theta_2$. As soon as $c$ gets bigger than $c_2$, the saddle point $\eta_2$ bifurcates into two stationary points $\mu_2$ and $\mu_3$ on either side of $\eta_2$: a local maximum at $\mu_2$ such that $\eta_2 < \mu_2 < \eta_1$ and a local minimum at $\mu_3$ such that $\theta_1 < \mu_3 < \eta_2$. The local minimum $\widetilde{\mu}_1$ remains as $\mu_1$ and is such that $\eta_1 < \mu_1 < \theta_2$. With the further increase of $c$, the local maximum at $\mu_2$ and the right local minimum (at $\mu_1$) get closer to each other and coalesce at $\eta_1$ when $c = c_1$. The left local minimum is then at $\theta_1$. This corresponds to a saddle point $\eta_1$ and a local minimum at $\theta_1$ for the quartic with $c = c_1$. When $c$ becomes bigger than $c_1$, the quartic will have only one stationary point --- the local left local minimum $\mu_3$ remains as $\widetilde{\mu}_1 < \theta_1$. \\
To summarize, the isolation intervals of the stationary points of the general quartic are as follows (dropping the tilde and the hat):
\vskip.1cm
{\bf (i)} If $c < c_2$, the quartic has a single local minimum $\mu_1 > \theta_2$.
\vskip.1cm
{\bf (ii)} If $c = c_2$, the quartic has a saddle point at $\eta_2$ and a local minimum at $\theta_2$.
\vskip.1cm
{\bf (iii)} If $c_2 < c < c_1$, the quartic has a local minimum at $\mu_3$ where $\theta_1 < \mu_3 < \eta_2$, a local maximum at $\mu_2$ where $\eta_2 < \mu_2 < \eta_1$, and local minimum at $\mu_1$ where $\eta_1 < \mu_1 < \theta_2$.
\vskip.1cm
{\bf (iv)} If $c = c_1$, the quartic has a saddle point at $\eta_1$ and a local minimum at $\theta_1$.
\vskip.1cm
{\bf (v)} If $c > c_1$, the quartic has a single local minimum $\mu_1 < \theta_1$.
\vskip.1cm
\noindent
For the analysis further, one also needs to determine the other two roots $\xi^{(i)}_{1,2}$ of the ``special" quartics $x^4 + ax^3 + bx^2 + cx + \delta_i$ (recall that $\mu_i$ is at least a double root for them). One has:
\b
\label{jkl}
x^4 + ax^3 + bx^2 + cx + \delta_i = (x - \mu_i)^2(x - \xi^{(i)}_1)(x-\xi^{(i)}_2) = 0.
\e
Vi\`ete formulas give: $\xi^{(i)}_1 + \xi^{(i)}_2 = - a - 2\mu_i$ and $2 \mu_i [\xi^{(i)}_1 + \xi^{(i)}_2] + \mu_i^2 + \xi^{(i)}_1 \xi^{(i)}_2 = b$. From these one finds that $\xi^{(i)}_1 \xi^{(i)}_2 = b + 3 \mu_i^2 + 2 a \mu_i$. If $x \ne \mu_i$ then (\ref{jkl}) reduces to
\b
\label{qr3}
x^2 + (a + 2 \mu_i) x + b + 3 \mu_i^2 + 2 a \mu_i = 0.
\e
This is the {\it third resolvent quadratic equation}. The roots of this equation are:
\b
\label{ksi}
\xi^{(i)}_{1,2} = -\frac{1}{2}a - \mu_i \pm \frac{1}{2} \sqrt{a^2 - 4 a \mu_i - 8 \mu_i^2 - 4b}.
\e
Note that the roots of the {\it third resolvent quadratic equation} (\ref{qr3}) depend on the roots $\mu_i$ of the first auxiliary cubic equation (\ref{1st}). That is, to find the $\xi$'s, one needs to find at least one of the stationary points of the quartic. Thus, the third resolvent quadratic equation should be used in the analysis based on solving cubic equations. \\
Separately, for one of the three ``special" quartics $x^4 + ax^3 + bx^2 + cx + \delta_i$, the roots $\xi^{(i)}_{1,2}$ of the third resolvent quadratic equation (\ref{qr3}) are not real, while for the remaining two they are real (unless the original quartic has the same value at its two local minima, in which case two of the ``special" quartics $x^4 + ax^3 + bx^2 + cx + \delta_i$ coincide and so all ``special" quartics will have four real roots) --- see Figure 2. \\
\vskip.3cm
\noindent
In the congruence of quartics, there is another significant quartic --- the one that passes through the origin --- i.e. this is a privileged quartic as it is the only one that has zero as a root. It is obtained from the original quartic by removing the free term $d$. The remaining three roots of this privileged quartic are found by solving the {\it second auxiliary cubic equation}:
\b
\label{2nd}
x^3 + ax^2 + bx + c = 0.
\e
If the discriminant $\Delta_2 = -27 c^2+(-4 a^3+18 a b) c + a^2 b^2 - 4b^3$ of this equation is negative, there is only one real root: $\lambda_1$. If it is not negative, there are three real roots: $\lambda_{0,1,2}$. To determine which of these occurs, set $\Delta_2 = 0$ to obtain the {\it fourth resolvent quadratic equation}:
\b
\label{qr4}
c^2 + \frac{2}{3} a \left( \frac{8}{9}  \frac{a^2}{4} - b \right)c - \frac{4}{27} b^2 \left(\frac{1}{4} a^2 - b \right) = 0.
\e
Once again, setting a discriminant of a cubic equal to zero is viewed as a quadratic in the unknown $c$ with $a$ and $b$ treated as parameters. The roots of this equation,
\b
\label{g12}
\gamma_{1,2}(a,b) = \frac{1}{3} a \left( b - \frac{8}{9} \frac{a^2}{4} \right)  \pm \frac{2 \sqrt{3}}{9} \, \sqrt{\left( \frac{4}{3}\frac{a^2}{4}- b \right)^3}.
\e
also play a very important role in the analysis. If the given $c$ is such that $\gamma_2(a,b) \le c \le \gamma_1(a,b)$, then the quartic $x^4 + a x^3 +b x^2 + c x$ has four real roots: 0 and $\lambda_{0,1,2}$ (there may be zeros among the $\lambda$'s). Otherwise, the $x^4 + a x^3 +b x^2 + c x$ has only two real roots: zero and $\lambda_1$ (which may also be zero).
\vskip.3cm
\noindent
Following the ideas of \cite{1}, the four ``degrees of freedom" of the general quartic $x^4 + ax^3 + bx^2 + cx +d$ are split equally between two separate polynomials, $x^4 + ax^3 + bx^2$ and $-cx - d$, the difference of which comprises the given quartic and the ``interaction" between which gives the roots of the quartic:
\b
\label{split}
x^2 (x^2 + ax + b) = - cx - d
\e
--- see Figure 3 (and also Figures 4 and 5) where this is illustrated with an example. \\
It may be tempting to depress the quartic and analyse only its ``three-dimensional projection" $y^4 + p y^2 + q y + r = 0$, but by doing so, study of how the coefficients of the original quartic affect its roots would not be possible as these coefficients would be ``dissolved" into $p$, $q$, and $r$. \\
It is quite easy to analyse the two parts of the ``split" quartic and, hence, the quartic itself --- one of the ``components" is a straight line, while the other is a quadratic ``in disguise" --- in the sense that it is a quartic having zero as a double root and allowing analysis not more difficult than that of a genuine quadratic (with the possible addition of a pair of stationary points and/or a pair of curvature change points, and the addition of a point where the third derivative vanishes). \\
For any given $a$ and $b$, i.e. for any ``sub-quartic" $x^2 (x^2 + ax + b)$, one can find a straight line $-c^\dagger x - d^\dagger$ such that it will be tangent to $x^2 (x^2 + ax + b)$ at two points, say $\alpha$ and $\beta$. This means that the obtained in this manner quartic, $x^4 + ax^3 + bx^2 + c^\dagger x + d^\dagger$, has two double roots: $\alpha$ and $\beta$. That is, $x^4 + ax^3 + bx^2 + c^\dagger x + d^\dagger = (x - \alpha)^2 (x- \beta)^2 = x^4 - 2(\alpha + \beta) x^3 + (\alpha^2 + \beta^2 + 4\alpha \beta) x^2 - 2 \alpha \beta (\alpha + \beta) x + \alpha^2 \beta^2.$ Comparing the corresponding coefficients yields: $c^\dagger = (1/2) a ( b - a^2/4)$. This is exactly equal to $c_0$ --- see the roots (\ref{c12}) of the {\it first resolvent quadratic equation} (\ref{qr1}). One also gets $d^\dagger = (1/4) (b - a^2/4)^2 > 0$. From now on, $c_0$ and $d_0$ will be used instead of $c^\dagger$ and $d^\dagger$.\\
The quartic $x^4 + ax^3 + bx^2 + c_0 x + d_0$ is very important for the classification of the roots of the general quartic. This is better visualized on the tablecloth of the ``split" quartic (\ref{split}). The determination of whether $c$ is smaller, equal, or greater than $c_0$ will determine the ordering of the $\delta$'s and will also determine the number of intersection points between the ``sub-quartic" $x^2 (x^2 + ax + b)$ and $- cx - d$ for any value of $d$. For example, on Figure 5, one has $0 < c_0 = -1.63 < c = - 1$. Thus, if one studies the intersection points of $x^2 (x^2 + ax + b)$ with $- cx - d$ in the regime of increasing $(-d)$ starting from $-\infty$, i.e. ``sliding" a straight line with fixed slope $(-c)$ upwards, intersections of this straight line with $x^2 (x^2 + ax + b)$ will occur first in the third quadrant before they occur in the forth. \\
From the roots (\ref{c12}) of the {\it first resolvent quadratic equation} (\ref{qr1}), it is immediately obvious that $c_2(a,b) \le c_0(a,b) \le c_1(a,b)$ for any $a$ and $b$. The graph of $c_0$ as a function of $a$ and $b$ is shown on Figure 7. Figures 8 to 11 show that, for any $a$ and $b$, the following holds for the general quartic:
\b
\label{seq}
c_2(a,b) \le \gamma_2(a,b) \le c_0(a,b) \le \gamma_1(a,b) \le c_1(a,b).
\e
Depending on $a$ and $b$, the place of zero in the above chain of inequalities could be anywhere. For the analysis of the general quartic, the very first step is the determination of these numbers and the following step is to find the place of 0 and the given $c$ in the above. Figure 12 shows this chain when $a$ and $b$ are both negative, in which case one has $c_2 < \gamma_2 < 0 < c_0 < \gamma_1 < c_1$. On Figure 12, the ``separator" straight lines $- c_{1,2} x - d_{1,2}$, $-\gamma_{1,2} x$, and $-c_o x - d_0$ are drawn and this clearly demonstrates that (including also the coordinate axes) there are seven ranges in which $c$ may fall. Each of these is individually studied. It has its own peculiarities that reflect on the number of roots and their localization. \\
Next, for the classification of the roots based on solving cubic equations, one needs to find the place of $(-d)$ among the $(-\delta)$'s and zero --- see Figure 4 where this is illustrated with the quartic equation $x^4 + x^3 - 3x^2 - x + 1 = 0$. The role of the {\it second auxiliary cubic equation} (\ref{2nd}) and its roots $\lambda_{0,1,2}$ now becomes clear. Keeping $a$ and $b$ fixed [i.e. not changing the ``sub-quartic" $x^2 (x^2 + ax + b)$], and only varying $c$ would ``move" the stationary points $\mu_i$ along the ``sub-quartic" $x^2 (x^2 + ax + b)$. For the example with $a = 1$ and $b = -3$ on Figures 1 to 5, for as long as $c < 0$, one always has $-\delta_2 > 0$ for the ``sub-quartic" $x^2 (x^2 + ax + b)$. On the other hand however, $-\delta_1$ and $-\delta_3$ could be anywhere. Because the {\it second auxiliary cubic equation} (\ref{2nd}) has three real roots $\lambda_{0,1,2}$,
$-\delta_1$ and $-\delta_3$ are both negative --- these are on the ``other side" (opposite side of $-\delta_2$) of the straight line $-cx - 0$ (the linear ``part" of the privileged quartic $x^4 + ax^3 + bx^2 + cx$). If the {\it second auxiliary cubic equation} (\ref{2nd}) had just one root $\lambda_{1}$, then $-\delta_1$ and $-\delta_3$ would both be positive. And because $c > c_0$, one has $-\delta_3 < -\delta_1$.

\vskip.3cm
\noindent
For the other tier of the analysis --- based on solving quadratic equations only --- one does not have the $\mu$'s, the $\delta$'s, and the $\lambda's$ explicitly. The isolation intervals of the $\delta$'s and the $\lambda$'s can be found in a manner similar to the one used for the determination of the isolation intervals of the $\mu$'s or one can see \cite{1} for the full classification of the roots of the cubic equation and the determination of the isolation intervals of its roots. The ``separator" line $-cx$ can still be used without knowing the loci of its point(s) of intersection with the ``sub-quartic" $x^2 (x^2 + ax + b)$, but knowing if these are three or just one. The ``separator" lines $- c x - \delta_i$ can no longer be used for analysis based on solving quadratic equations only. There is a way however, to find ``replacements". \\
The ``sub-quartic" $x^2 (x^2 + ax + b)$ has zero as a double root and two more roots which are the roots of the {\it fifth resolvent quadratic equation}
\b
\label{qr5}
x^2 + ax + b = 0,
\e
the roots of which are
\b
\label{r12}
\rho_{1,2} = -\frac{1}{2} a \pm \sqrt{\frac{a^2}{4} - b}.
\e
These are real for $b \le a^2/4$. \\
Then, one takes the given $c$ and draws the two parallel lines with equations $-c (x - \rho_{1,2})$. These straight lines are the sought ``replacements" of the ``separator" lines $- c x - \delta_i$ and their intersections with the ordinate --- the ``marker" points $c \rho_{1,2}$ --- are the ``replacements" of the $(-\delta)$'s. This allows the analysis based on solving quadratic equations only to be performed in manner fully analogous to that of the analysis based on solving cubic equations --- see Figure 5. \\
All possibilities for this analysis are shown on Figures 1.1 to 1.14 and 2.1 to 2.14.\\
If $b > a^2/4$ ($\rho_{1,2}$ not being real), a different pair of characteristic points of the ``sub-quartic" $x^2 (x^2 + ax + b)$ should be chosen as ``marker" points --- the two non-zero critical points $\sigma_{1,2}$ of $x^2 (x^2 + ax + b)$. Clearly, recourse to these can be made for  $a^2/4 < b \le  (9/8)a^2/4$ --- as can be easily seen from the {\it sixth resolvent quadratic equation}
\b
\label{qr6}
4x^2 + 3ax + 2b = 0,
\e
the roots of which are the non-zero critical points of $x^2 (x^2 + ax + b)$ given by
\b
\label{s12}
\sigma_{h,H} = -\frac{3}{8} a \pm \frac{\sqrt{2}}{2} \sqrt{\frac{9}{8}\frac{a^2}{4} - b}.
\e
One then calculates the values $H$ and $h$ of $x^2 (x^2 + ax + b)$ at $\sigma_H$ and $\sigma_h$ respectively and draws the parallel lines $-c (x - \sigma_{H}) + H$ and $-c (x - \sigma_{h}) + h$ to serve as ``separators". These intersect the ordinate at the ``marker" points $c \sigma_{H} + H$ and $c \sigma_{h} + h$. One has to be careful because, depending on $c$, one can have zero, $c \sigma_{H} + H$, and $c \sigma_{h} + h$ in any order. \\
All possibilities for this analysis are shown on Figures 3.1 to 3.14.\\
Should $b$ be greater than $(9/8)(a^2/4)$, then the ``sub-quartic" $x^2 (x^2 + ax + b)$ will not have critical points. One should then use the points of curvature change [non-zero first derivative, but vanishing second derivative of the ``sub-quartic" $x^2 (x^2 + ax + b)$]. These are real for $b \le  (3/2)(a^2/4)$. They are the roots of the {\it second resolvent quadratic equation} and on Figures 4.1 to 4.14 and, also 5.1 to 5.10 (where the relevant analysis is), they are denoted by
\b
\label{t12}
\tau_{h,H} = -\frac{1}{4} a \pm \frac{\sqrt{6}}{6} \, \sqrt{\frac{3}{2}\frac{a^2}{4} - b}.
\e
As with the critical points $\sigma_{1,2}$, one then calculates the values $H$ and $h$ of $x^2 (x^2 + ax + b)$ at $\tau_H$ and $\tau_h$ respectively and draws the parallel lines $-c (x - \tau_{H}) + H$ and $-c (x - \tau_{h}) + h$ to serve as ``separators". These intersect the ordinate at the ``marker" points $c \tau_{H} + H$ and $c \tau_{h} + h$. The slope of the straight line joining the two points of curvature change is equal to $\pm c_0$. Thus, which of $c \tau_{H} + H$ and $c \tau_{h} + h$ is bigger depends on whether $c$ is bigger or smaller than $c_0$. Again, care should be exercised as zero, $c \tau_{H} + H$ and $c \tau_{h} + h$ could be in any order --- it is the number of real roots of the {\it second auxiliary cubic equation} that determines this order.  \\
Finally, when $b > (3/2)a^2/4$, not one of the {\it resolvent quadratic equations} has real roots. One still needs to find identifiable points the ``sub-quartic" $x^2 (x^2 + ax + b)$ from which ``separator" lines can be drawn. There is just one such point --- where the third derivative of $x^2 (x^2 + ax + b)$ vanishes. This point is the only root $\phi = -a/4$ of the {\it resolvent linear equation}: $4x + a = 0$. One then draws the only available ``separator" --- the line $- c(x + a/4) + (1/16)a^2 [b - (3/4)(a^2/4)]$. It intersects the ordinate at the ``marker" point $t = -(1/4)ac + (1/16)a^2 [b - (3/4)(a^2/4)]$. If the signs of $a$ and $c$ are opposite, a second ``separator" line, $- cx + (1/16)a^2 [b - (3/4)(a^2/4)]$ (it is parallel to the first), can provide sharper bounds on the roots (this line is not useful if $a$ and $c$ are with the same sign). Care should be exercised as one could have $t < 0, \,\, t = 0$,  or $t > 0$ --- see Figures 6.1 to 6.4 where the corresponding analysis can be found.

\section{Classification of the Roots of the Quartic Equation}

\noindent
Every possible case for non-zero $a$, $b$, and $c$ has been thoroughly analyzed. The cases of $c$ equal to $c_{1,2}$, or $\gamma_{1,2}$, or $c_0$ do not get special attention either --- should one or more of $a$, $b$, and $c$ be zero or should $c$ be equal to one of the above, the analysis (not presented here) follows trivially. \\
The results of the investigation are presented in figures with labels $i.j$, where $i$ and $j$ are positive integers . The figures can be grouped into a (rather large) table. For ease of reference, an effort has been made to keep the individual figures independent from each other and for this purpose, each Figure $i. j$ contains a short description of the situation, analysis based on solving cubic equations, and analysis based on solving quadratic equations only. \\
The index $i$ in Figure $i.j$ runs from 1 to 6 and labels the rows of the table:
\vskip.1cm
$\bm{i = 1}\!: \quad$ This is the case of $b < 0$. The roots $c_{1,2}, \,\, \gamma_{1,2},$ and $\rho_{1,2}$ of the {\it first, fourth, and fifth resolvent quadratic equation}, respectively, are real. The roots $\rho_{1,2}$ have opposite signs. When $i = 1$, the index $j$ runs from 1 to 14 (there are fourteen columns). The first seven of these correspond to the seven possible ranges for $c$ when $a < 0$; the remaining seven are the seven possible ranges for $c$ when $a > 0$.
\vskip.1cm
$\bm{i = 2}\!: \quad$ This is the case of $0 < b \le a^2/4$. The roots $c_{1,2}, \,\, \gamma_{1,2},$ and $\rho_{1,2}$ of the {\it first, fourth, and fifth resolvent quadratic equation}, respectively, are again real. This time the roots $\rho_{1,2}$ have the same sign (opposite to the sign of $a$). When $i = 2$, the index $j$ again runs from 1 to 14 (there are fourteen columns) with the first seven of these corresponding to the seven possible ranges for $c$ when $a < 0$ and the remaining seven corresponding to the seven possible ranges for $c$ when $a > 0$.
\vskip.1cm
$\bm{i = 3}\!: \quad$ This is the case of $a^2/4 < b \le (9/8) a^2/4$. The roots $c_{1,2}$ and $\gamma_{1,2}$ of the {\it first and fourth resolvent quadratic equation}, respectively, are real, but the roots $\rho_{1,2}$ of the {\it fifth resolvent quadratic equation}, are not real. The roots $\sigma_{h,H}$ of the {\it sixth resolvent quadratic equation} are real and these are used in the analysis. When $i = 3$, the index $j$ again runs from 1 to 14 (there are fourteen columns). The first seven of these correspond to the seven possible ranges for $c$ when $a < 0$; the remaining seven are the seven possible ranges for $c$ when $a > 0$.
\vskip.1cm
$\bm{i = 4}\!: \quad$ This is the case of $(9/8) a^2/4 < b \le (4/3) a^2/4$. The roots $c_{1,2}$ and $\gamma_{1,2}$ of the {\it first and fourth resolvent quadratic equation}, respectively, are real, but the roots $\rho_{1,2}$ and $\sigma_{h,H}$ of the {\it fifth and sixth resolvent quadratic equation}, respectively, are not real. The roots $\tau_{h,H}$ of the {\it second resolvent quadratic equation} are real and these are used in the analysis. When $i = 4$, the index $j$ again runs from 1 to 14 (there are fourteen columns). The first seven of these correspond to the seven possible ranges for $c$ when $a < 0$; the remaining seven are the seven possible ranges for $c$ when $a > 0$.
\vskip.1cm
$\bm{i = 5}\!: \quad$ This is case of $(4/3) a^2/4 < b \le (3/2) a^2/4$. The roots $c_{1,2}$ of the {\it first resolvent quadratic equation}  are real, but the roots $\gamma_{1,2}$, $\rho_{1,2}$, and $\sigma_{h,H}$,  of the {\it fourth, fifth, and sixth resolvent quadratic equation}, respectively, are not real. The roots $\tau_{h,H}$ of the {\it second resolvent quadratic equation} are real and these are again used in the analysis. When $i = 5$, the index $j$ runs from 1 to 10 (there are ten columns). The first five of these correspond to the five possible ranges for $c$ when $a < 0$; the remaining five are the five possible ranges for $c$ when $a > 0$ --- there are no $\gamma$'s anymore.
\vskip.1cm
$\bm{i = 6}\!: \quad$ This is the case of $(3/2) a^2/4 < b$. Not one of the {\it resolvent quadratic equations} has real roots. But the graph of the ``sub-quartic" $x^2(ax^2 + bx + c)$  still has a ``blemish" that can be used and extract one or two ``separator" lines for the analysis. This is the point $\phi = -a/4$ where the third derivative of $x^2(ax^2 + bx + c)$ is zero. i.e. $\phi$ is the root of the {\it resolvent linear equation} $4x + a = 0$. When $i = 6$, the index $j$ runs from 1 to 4 (there are only four columns): $a < 0$ with $c <0, \,\,$ $a< 0$ with $c > 0, \,\,$ $a > 0$ with $c < 0$, and $a > 0$ with $c >0$.
\vskip.1cm
\noindent
To find the relevant case, one has to see first in which of the above ranges the coefficient $b$ of the quadratic term falls into. The relevant row of the table is then selected by the particular value of $b$. \\
Next, if the sign of the coefficient $a$ of the cubic term is negative, one should look at columns 1 to 7 for the first four rows of the table, columns 1 to 5 for the fifth row and columns 1 and 2 for the sixth row of the table. If $a$ is positive, one should look at the other columns of the relevant row. \\
It is the place of the coefficient $c$ of the linear term within the chain of inequalities (\ref{seq}) that determines which particular column applies, namely, $c$ selects the individual cell in the table that is relevant.  \\
Finally, the value of the coefficient $d$ of the free term determines, within that cell, which one of the cases labeled by lower case Roman numerals applies (either for the analysis based on cubic equations or for the analysis based on solving quadratic equations only). \\
All figures are on pages 12--61.


\newgeometry{top=30mm, bottom=30mm,left=15mm, right=15mm}

\begin{center}

\begin{tabular}{|c||c|}
\hline
{\footnotesize \bf Figure 1} &  {\footnotesize \bf Figure 2} \\
\hline
\hline
\includegraphics[width=77mm]{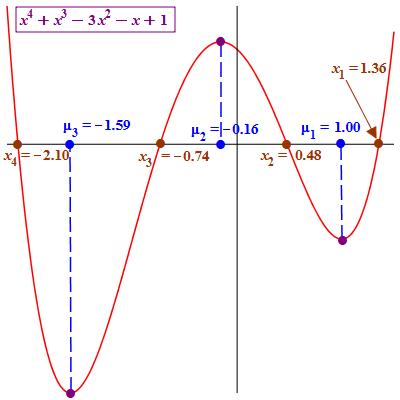} & \includegraphics[width=77mm]{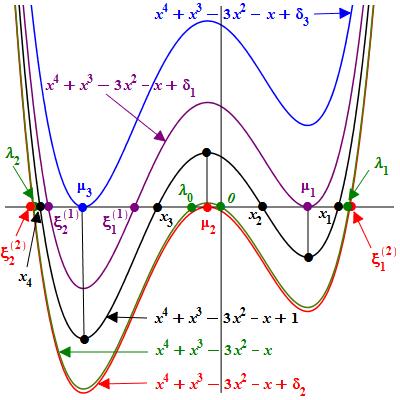} \\
\hline
\multicolumn{1}{|c||}{\begin{minipage}{20.7em}
\scriptsize
\vskip.15cm
The quartic polynomial $x^4 + x^3 - 3x^2 - x + 1$ has four real roots $x_i$ and three stationary points: a local minimum at $\mu_{1} = 1.00$, a local maximum at $\mu_2 = -0.16$, and another local minimum at $\mu_3 = -1.59$  --- roots of the first auxiliary cubic equation: $4x^3 + 3x^2 - 6x - 1 =0.$ \\
(The second auxiliary cubic equation is $x^3 + x^2 - 3x - 1 = 0$ and it has three real roots: $\lambda_{0,1,2}$ --- see also Figure 2 and Figure 4.) \\
\end{minipage}}
& \multicolumn{1}{|c|}{\begin{minipage}{21.9em}
\scriptsize
\vskip-.5cm
Shown are the significant quartics of the congruence obtained by varying the free term of the original quartic $x^4 + x^3 - 3x^2 - x + 1$ (also pictured):  the three ``special" quartics: $x^4 + x^3 - 3x^2 - x + 2.00, \,\, $ $x^4 + x^3 - 3x^2 - x - 0.08, \,\, $ and $x^4 + x^3 - 3x^2 - x + 3.62,$ together with the quartic without a free term, $x^4 + x^3 - 3x^2 - x.$ \\
\end{minipage}}  \\
\hline
\end{tabular}
\end{center}

\vskip2cm

\begin{center}
\begin{tabular}{|c|c|}
\hline
\multicolumn{2}{|c|}{\footnotesize \bf Figure 3}  \\
\hline
\hline
& \\
\includegraphics[width=80mm,height=77mm]{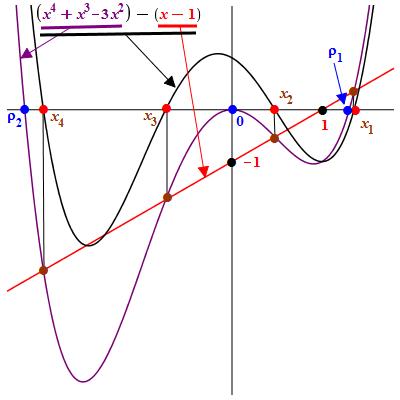} &
\multicolumn{1}{|c|}{\begin{minipage}{22em}
\scriptsize
\vskip-11.2cm
The quartic polynomial $x^4 + x^3 - 3x^2 - x + 1$ is viewed as the difference of two polynomials: $f(x) = x^2 (x^2 + x - 3)$ and $g(x) = x - 1$. The intersection points of $f(x)$ and $g(x)$ are the roots of the given quartic equation $x^4 + x^3 - 3x^2 - x + 1 = 0$. The quartic $f(x)$ has a double root at zero and two more roots: $\rho_2 = -2.30$ and $\rho_1 = 1.30$. The straight line $g(x)$ has crosses the abscissa at $x = 1$. \\
The roots of the quartic equation $x^4 + x^3 - 3x^2 - x + 1 = 0$ are studied in a two-tier analysis --- by solving cubic equations (see Figure 4) and
by solving quadratic equations only (see Figure 5). This equation (with $a = 1$, $b = -3$, and $c = -1$) belongs to the sub-case shown on Figure 1.11.
\end{minipage}}
\\
\hline
\end{tabular}
\end{center}

\newgeometry{top=30mm, bottom=30mm,left=15mm, right=15mm}

\begin{center}
\begin{tabular}{|c|c|}
\hline
\multicolumn{2}{|c|}{\footnotesize \bf Figure 4}  \\
\hline
\hline
& \\
\includegraphics[width=80mm,height=80mm]{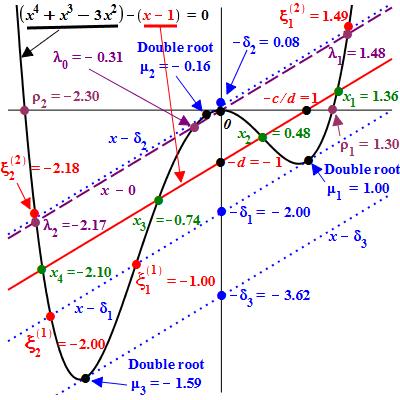} &
\multicolumn{1}{|c|}{\begin{minipage}{22em}
\scriptsize
\vskip-7.2cm
\begin{center}
\underline{{\color{red}{\it Analysis based on solving cubic equations}}}
\end{center}
\vskip-.3cm
The equation $x^4 + x^3 - 3x^2 - x + 1 = 0$ comes with $c = -1$ which is between the roots $c_2 = - 5.00$ and $c_1 = 1.75$ of the {\it first resolvent quadratic equation} (\ref{qr1}). Thus, the {\it first auxiliary cubic equation} (\ref{1st}) has three real roots, i.e. the quartic has three stationary points ($\mu_{1} = 1.00, \, \mu_2 = -0.16, \, \mu_3 = -1.59$) and can have 0, 2 or 4 real roots.
\vskip.1cm
As $c = -1$ is also between the roots $\gamma_2 = -3.42$ and $\gamma_1 = 1.27$ of the {\it fourth resolvent quadratic equation} (\ref{qr4}), the {\it second auxiliary cubic equation} (\ref{2nd}) has three real roots ($\lambda_1 = 1.48, \, \lambda_0 = -0.31, \, \lambda_2 = -2.17$), i.e. there are three intersection points of the straight line $x-1$ with the quartic $x^2(x^2 + x - 3)$. This is why $- \delta_1 < 0$ and $-\delta_3 < 0$. On the other hand, $-\delta_2 > 0$ for all $c \ne 0$.
\vskip.1cm
As $0 < c_0 = -1.63 < c$ (see Figure 5 where the line $-c_0 x - d_0$ is shown), one has $-\delta_3 < -\delta_1$.
\vskip.1cm
The real roots $\xi^{(1,2)}_{1,2}$ of the {\it third resolvent quadratic equation} (\ref{qr3}) are those corresponding to $\mu_{1,2}$.
\vskip.1cm
The line $x -1$ intersects the abscissa at $1$ --- to the left of $\rho_1 = 1.30$.
\vskip.1cm
As the free term $d = 1$ satisfies $-\delta_3 = -3.62 < -\delta_1 = -2.00 < -d < 0 < -\delta_2 = 0.08$, the number of real roots is exactly 4.
\vskip.1cm
Therefore, there is a negative root $x_4$ between $\lambda_2$ and $\xi^{(1)}_2$, a negative root $x_3$ between $\xi^{(1)}_1$ and $\lambda_0$, a positive root $x_2$ smaller than min $\{\mu_1, -d/c$\}, and a positive root $x_1$ between $\mu_1$ and $\lambda_1$ (a sharper bound for $x_1$ is $\rho_1 < x_1 < \lambda_1$). See also Figure 1.11. \\
\end{minipage}}
\\
\hline
\end{tabular}
\end{center}

\vskip2cm

\begin{center}
\begin{tabular}{|c|c|}
\hline
\multicolumn{2}{|c|}{\footnotesize \bf Figure 5}  \\
\hline
\hline
& \\
\includegraphics[width=80mm,height=80mm]{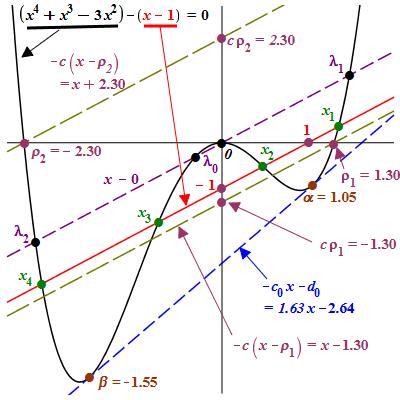} &
\multicolumn{1}{|c|}{\begin{minipage}{22em}
\scriptsize
\vskip-8.2cm
\begin{center}
\underline{{\color{red}{\it Analysis based on solving quadratic equations only}}}
\end{center}
\vskip-.2cm
To perform this analysis, it is necessary first to determine, by solving the {\it first and third resolvent quadratic equations}, all numbers in the chain of inequalities $c_2 \le \gamma_2 \le c_0 \le \gamma_1 \le c_1$ (see also Figures 7 to 12) and then find in this chain the places of 0 and the given $c = -1$.
\vskip.1cm
For the given $a = 1$ and $b = -3$, one finds: $c_2 = -5 < \gamma_2 = -3.42 < c_0 = - 1.63 < c = -1 < 0 < \gamma_1 = 1.27 < c_1 = 1.75$.
\vskip.1cm
One also that $d_0 = 2.64$ and, by solving another quadratic equation, $\alpha = 1.05$ and $\beta = -1.55$.
\vskip.1cm
Due to $c$ being between $c_2$ and $c_1$, the given quartic equation could have 0, 2, or 4 real roots.
\vskip.1cm
Due to $c$ being between $\gamma_2$ and $\gamma_1$, there are three intersection points ($\lambda_{0,1,2}$) of $x^4 + ax^3 + bx^2$ with $-cx$.
\vskip.1cm
Due to $c > c_0$, it is guaranteed that there are two intersection points of $x^4 + ax^3 + bx^2$ with $-cx -d$ in the third quadrant, provided that $ - d > c \rho_1$.
\vskip.1cm
For the given $a$, $b$, and $c$, the analysis shown on Figure 1.11 applies with {\bf (ii)} being the relevant sub-case, namely, $c \rho_1 = -1.30 \le - d = -1 < 0$. Thus there are four real roots: the two negative roots $x_{3,4}  > \rho_2 = -2.30$, the positive root $x_2 < -d/c = 1$ and the positive root $x_1 > \rho_1 = 1.30$.
\end{minipage}}
\\
\hline
\end{tabular}
\end{center}

\newgeometry{top=10mm, bottom=15mm,left=15mm, right=15mm}

\begin{center}
\begin{tabular}{|c||c|}
\hline
\multicolumn{2}{|c|}{\footnotesize \bf Figure 6}  \\
\hline
\hline
\multirow{2}{*}{{\begin{minipage}{22em}
\vskip-3.75cm
\includegraphics[width=80mm]{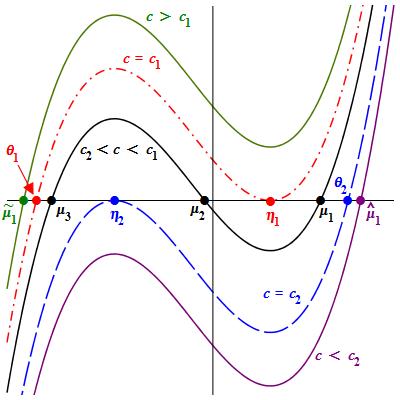}\end{minipage}}} &
\multicolumn{1}{|c|}{\begin{minipage}{22em}
\includegraphics[width=80mm]{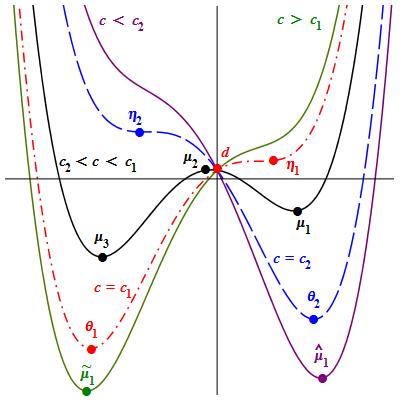}
\end{minipage}} \\
\hline
\multicolumn{2}{|c|}{\begin{minipage}{45em}
\scriptsize
\vskip.3cm
\begin{center}
\underline{{\color{red}{\it Isolation Intervals of the Stationary Points of the General Quartic}}}
\end{center}
\vskip-.2cm
The example chosen here for illustration is a quartic with coefficients: $a = 1$, $b = -5$, $c = -1$, and $d = 1$.
\vskip.1cm
The {\it first auxiliary cubic} polynomial $4x^3 + 3ax^2 + 2bx + c$ is plotted on the left pane with fixed $a = 1$ and $b= -5$ and with $c$ varying in different ranges, determined by the roots  $c_2 = -9.41$ and $c_1 = 4.16$ of the {\it first resolvent quadratic equation} (\ref{qr1}). The quartic $x^4 + ax^3 + bx^2 + cx + d$ is plotted on the right pane, again with fixed $a = 1$ and $b = -5$, also with fixed $d = 1$, and with $c$ varying over the same ranges as on the right pane. One can immediately see how the coefficient $c$ in the linear term affects the quartic and determine the number and type of stationary points of the quartic.
\vskip.1cm
By solving the {\it second resolvent quadratic equation}, one immediately finds that $\eta_1 = 0.70$ and $\eta_2 = -1.20$. By using (\ref{teta}), one finds that $\theta_1 = -2.14$ and $\theta_2 = 1.64$.
\vskip.2cm
{\bf (i)} If $c < c_2$, the quartic has a single local minimum $\widehat{\mu}_1 > \theta_2$ (purple solid curve).
\vskip.1cm
{\bf (ii)} If $c = c_2$, the quartic has a saddle point at $\eta_2$ and a local minimum at $\theta_2$ (blue dashed curve).
\vskip.1cm
{\bf (iii)} If $c_2 < c < c_1$, the quartic has a local minimum at $\mu_3$ with $\theta_1 < \mu_3 < \eta_2$, a local maximum at $\mu_2$ with $\eta_2 < \mu_2 < \eta_1$, and local minimum at $\mu_1$ with $\eta_1 < \mu_1 < \theta_2$ (black solid curve).
\vskip.1cm
{\bf (iv)} If $c = c_1$, the quartic has a saddle point at $\eta_1$ and a local minimum at $\theta_1$ (red dash-dotted curve).
\vskip.1cm
{\bf (v)} If $c > c_1$, the quartic has a single local minimum $\widetilde{\mu}_1 < \theta_1$ (green solid curve).
\vskip.3cm
Thus $-2.14 < \mu_3 < -1.20 < \mu_2 < 0.70 < \mu_3 < 1.64$. The actual values are: $\mu_3 = -1.96$, \, $\mu_2 = -0.10,$ and $\mu_1 = 1.31$. \\
\end{minipage}} \\
\hline
\end{tabular}
\end{center}


\vskip1cm

\begin{center}
\begin{tabular}{|c||c|}
\hline
{\footnotesize \bf Figure 7} &  {\footnotesize \bf Figure 8} \\
\hline
\hline
\includegraphics[width=80mm]{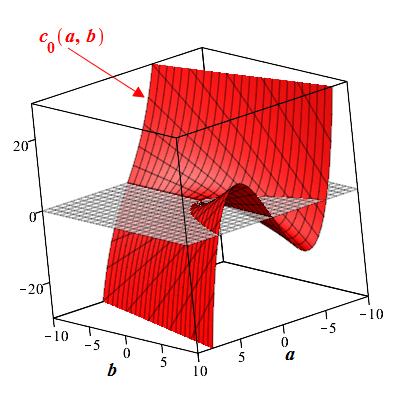} & \includegraphics[width=80mm]{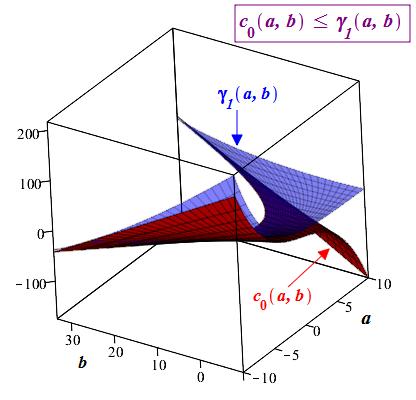} \\
\hline
\multicolumn{1}{|c||}{\begin{minipage}{22em}
\scriptsize
\vskip-.4cm
\begin{center}
$c_0(a,b) = \frac{1}{2} a \left( b - \frac{a^2}{4} \right)$.
\end{center}
\end{minipage}}
& \multicolumn{1}{|c|}{\begin{minipage}{22em}
\scriptsize
\vskip.1cm
\begin{center}
$\gamma_{1}(a,b) = \frac{1}{3} a \left( b - \frac{8}{9} \frac{a^2}{4} \right)  + \frac{2 \sqrt{3}}{9} \, \sqrt{\left( \frac{4}{3}\frac{a^2}{4}- b \right)^3} \linebreak \ge c_0(a,b) = \frac{1}{2} a \left( b - \frac{a^2}{4} \right).$
\end{center}
\end{minipage}}  \\
\hline
\end{tabular}
\end{center}

\newgeometry{top=20mm, bottom=20mm,left=15mm, right=15mm}

\begin{center}
\begin{tabular}{|c||c|}
\hline
{\footnotesize \bf Figure 9} &  {\footnotesize \bf Figure 10} \\
\hline
\hline
\includegraphics[width=80mm]{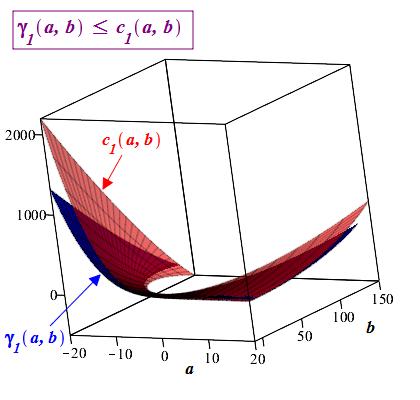} & \includegraphics[width=80mm]{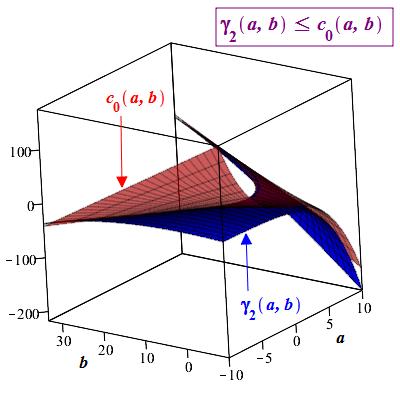} \\
\hline
\multicolumn{1}{|c||}{\begin{minipage}{22em}
\scriptsize
\vskip.1cm
\begin{center}
$\gamma_{1}(a,b) = \frac{1}{3} a \left( b - \frac{8}{9} \frac{a^2}{4} \right)  + \frac{2 \sqrt{3}}{9} \, \sqrt{\left( \frac{4}{3}\frac{a^2}{4}- b \right)^3} \linebreak \le c_{1}(a,b) = c_0(a,b) + \frac{2\sqrt{6}}{9} \, \sqrt{ \left( \frac{3}{2}\frac{a^2}{4} - b \right)^3}$. \\
\end{center}
\end{minipage}}
& \multicolumn{1}{|c|}{\begin{minipage}{22em}
\scriptsize
\vskip.1cm
\begin{center}
$\gamma_{2}(a,b) = \frac{1}{3} a \left( b - \frac{8}{9} \frac{a^2}{4} \right)  - \frac{2 \sqrt{3}}{9} \, \sqrt{\left( \frac{4}{3}\frac{a^2}{4}- b \right)^3} \linebreak \le c_0(a,b) = \frac{1}{2} a \left( b - \frac{a^2}{4} \right).$ \\
\end{center}
\end{minipage}}  \\
\hline
\end{tabular}
\end{center}

\vskip2cm

\begin{center}
\begin{tabular}{|c||c|}
\hline
{\footnotesize \bf Figure 11} &  {\footnotesize \bf Figure 12} \\
\hline
\hline
\includegraphics[width=80mm]{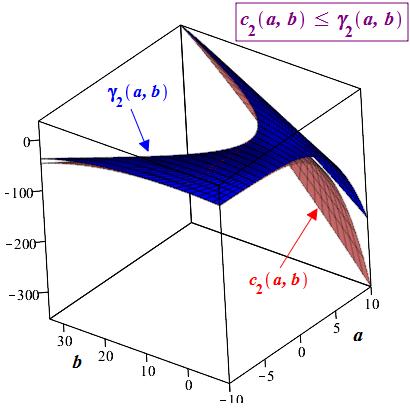} & \includegraphics[width=80mm]{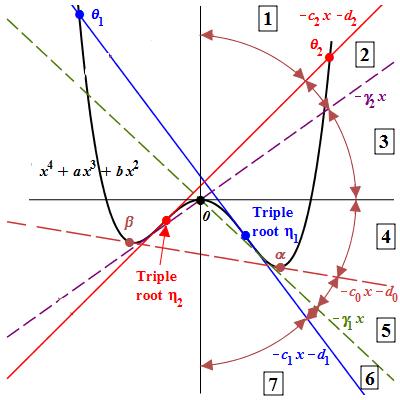} \\
\hline
\multicolumn{1}{|c||}{\begin{minipage}{22em}
\scriptsize
\vskip-2.2cm
\begin{center}
$c_{2}(a,b) = c_0(a,b) - \frac{2\sqrt{6}}{9} \, \sqrt{ \left( \frac{3}{2}\frac{a^2}{4} - b \right)^3} \linebreak \le
\gamma_{2}(a,b) = \frac{1}{3} a \left( b - \frac{8}{9} \frac{a^2}{4} \right)  - \frac{2 \sqrt{3}}{9} \, \sqrt{\left( \frac{4}{3}\frac{a^2}{4}- b \right)^3}.$
\end{center}
\end{minipage}}
& \multicolumn{1}{|c|}{\begin{minipage}{22em}
\scriptsize
\vskip.2cm
Pictured here are the ``sub-quartic" $x^2(x^2 + ax + b)$ and the ``separator" straight lines $- c_{1,2} \, x - d_{1,2}$, $-\gamma_{1,2} \, x$, and $-c_0 \, x - d_0$ when $a$ and $b$ are both negative (in which case $c_0$ is positive). \\
On the diagram, $\eta_{1,2}$ are the roots of the {\it second resolvent quadratic equation}, i.e. the triple roots of the ``special" quartics $x^4 + ax^3 + bx^2  + c_{1,2}x + d_{1,2}$. \\
For any $a$ and $b$, one always has $c_2 \le \gamma_2 \le c_0 \le \gamma_1 \le c_1$. \\
For any particular pair $(a, b)$, the place of zero has to be found in this chain of inequalities and then one has to determine in which of the seven ranges (determined by the ``separator" lines and the coordinate axes) the coefficient $c$ of the linear term falls. \\
\end{minipage}}  \\
\hline
\end{tabular}
\end{center}


\newgeometry{top=5mm, bottom=14mm,left=15mm, right=15mm}

\begin{center}
\begin{tabular}{|c||c|}
\hline
{\footnotesize \bf Figure 1.1} &  {\footnotesize \bf Figure 1.2} \\
\hline
\hline
{\color{violet} $\bm{b < 0, \quad a <  0, \quad c < 0,$}} &  {\color{violet}$\bm{b < 0, \quad a <  0, \quad c < 0,}$} \\
{\color{violet} $\bm{{\color{red}c} < c_2 < \gamma_2 < 0 < c_0 < \gamma_1 < c_1$}} & {\color{violet} $\bm{c_2 < {\color{red}c} < \gamma_2 < 0 < c_0 < \gamma_1 < c_1$}} \\
\hline
\hline
\includegraphics[width=80mm]{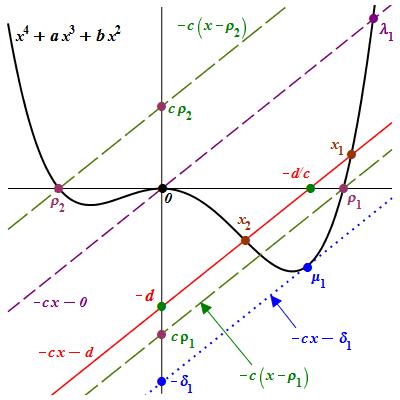} & \includegraphics[width=80mm]{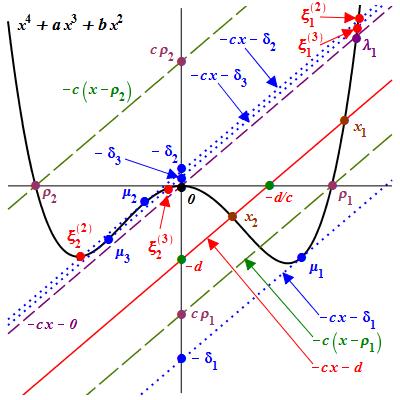} \\
\hline
\multicolumn{1}{|c||}{\begin{minipage}{22em}
\scriptsize
\vskip-2.9cm
\vskip.3cm
\hskip3.89cm \underline{{\color{red}{\it Notes}}}
\vskip.3cm
As $c < c_2 < c_1$, the quartic has a single stationary point $\mu_1$ and the number of real roots can be either 0 or 2. There is only one tangent $- c x - \delta_1$ to $x^4 + ax^3 + bx^2$ --- the one at $\mu_1$.
\vskip.1cm
As $c < \gamma_2 < \gamma_1$, the straight line $-cx$ intersects $x^4 + ax^3 + bx^2$ at one point only ($\lambda_1$).
\vskip.1cm
One could have $\mu_1 < \rho_1$ (pictured) or $\mu_1 \ge \rho_1$.
\vskip.1cm
Obviously, $ -\delta_1 \le c \rho_1 < 0 < c \rho_2$.
\vskip.1cm
The minimum of $x^4 + ax^3 + bx^2$ at negative $x$ is greater than the minimum at positive $x$ by the amount of $-(1/256)a(9a^2-32b)^{3/2} > 0$.
\vskip.1cm
Consideration of whether $- d < c \rho_1$ or $-d > c \rho_1$ would yield sharper bounds in the analysis based on solving cubic equations. \\
\hrule
\vskip.3cm
\hskip1.54cm \underline{{\color{red}{\it Analysis based on solving cubic equations}}}
\vskip.3cm
{\bf (i)} If $-d < -\delta_1$, then there are no real roots.
\vskip.1cm
{\bf (ii)} If $-\delta_1 \le -d < 0$ (pictured), then there is one positive root $x_2 \le $ min$\{\mu_1, -d/c\}$ and another positive root $x_1$ such that $\mu_1 \le x_1 < \lambda_1$.
\vskip.1cm
{\bf (iii)}  If $0 \le - d$, then there is a non-positive root $x_2$ such that min$\{\rho_2, -d/c \} \le x_2 \le $ max$\{\rho_2, -d/c \}$  and a positive root $x_1 \ge \lambda_1$. \\
\hrule
\vskip.3cm
\hskip.88cm \underline{{\color{red}{\it Analysis based on solving quadratic equations only}}}
\vskip.3cm
{\bf (i)} If $-d < c \rho_1$, then there are either no real roots or there are two positive roots (these are greater than $\rho_1$ if $\mu_1 > \rho_1$ or smaller than or equal to $\rho_1$ if $\mu_1 \le \rho_1$).
\vskip.1cm
{\bf (ii)} If $c \rho_1 \le - d < 0$ (pictured), then there is a positive root $x_2 < -d/c$ and another positive root $x_1 \ge \rho_1$.
\vskip.1cm
{\bf (iii)} If $0 \le - d < c \rho_2$, then there is one negative root $x_2$ such that  $\rho_2 < x_2 \le -d/c$ and a positive root $x_1 > \rho_1$.
\vskip.1cm
{\bf (iv)} If $c \rho_2 \le -d$, then there is one negative root $x_2$ such that $-d/c < x_2 \le \rho_2$ and a positive root $x_1 > \rho_1$.
\end{minipage}}
& \multicolumn{1}{|c|}{\begin{minipage}{22em}
\scriptsize
\scriptsize
\vskip.3cm
\hskip3.89cm \underline{{\color{red}{\it Notes}}}
\vskip.3cm
As $c_2 < c < c_1$, the quartic has three stationary points $\mu_i$ and the number of real roots can be either 0, or 2, or 4. There are three tangents $- c x - \delta_i$ to $x^4 + ax^3 + bx^2$ --- each at $\mu_i$.
\vskip.1cm
As $c < \gamma_2 < \gamma_1$, the straight line $-cx$ intersects $x^4 + ax^3 + bx^2$ at one point only ($\lambda_1$). Thus $-\delta_3 > 0$. Also: $-\delta_3 > 0$ and $-\delta_1 < 0$.
\vskip.1cm
Obviously, $ -\delta_1 < c \rho_1 < 0 < -\delta_3 < - \delta_2 < c \rho_2$.
\vskip.1cm
The minimum of $x^4 + ax^3 + bx^2$ at negative $x$ is greater than the minimum at positive $x$ by the amount of $-(1/256)a(9a^2-32b)^{3/2} > 0$.
\vskip.1cm
Consideration of whether $- d < c \rho_1$ or $-d > c \rho_1$ would yield sharper bounds in the analysis based on solving cubic equations. \\
\hrule
\vskip.3cm
\hskip1.54cm \underline{{\color{red}{\it Analysis based on solving cubic equations}}}
\vskip.3cm
{\bf (i)} If $-d < -\delta_1$, then there are no real roots.
\vskip.1cm
{\bf (ii)} If $-\delta_1 \le -d < 0$ (pictured), then there is one positive root $x_2 \le $ min$\{\mu_1, -d/c\}$ and another positive root $x_1$ such that $\mu_1 \le x_1 < \lambda_1$.
\vskip.1cm
{\bf (iii)}  If $0 \le -d < -\delta_3$, then there is a positive root $x_1$ such that $\lambda_1 \le x_1 < \xi^{(3)}_1$, and a negative root $x_2$ such that $\xi^{(3)}_2 < x_2  \le -d/c$ .
\vskip.1cm
{\bf (iv)} If $-\delta_3 \le -d < - \delta_2$, then there is a positive root $x_1$ such that $\xi^{(3)}_1 \le x_1 < \xi^{(2)}_1$, a negative root $x_2$ such that $\mu_2 < x_2 \le $ min$\{\xi^{(3)}_2,-d/c\}$, another negative root $x_3$ such that $\mu_3 \le x_3 < \mu_2$, and a third negative root $x_4$ such that $\mu_3 \le x_4 < \xi^{(2)}_2$.
\vskip.1cm
{\bf (v)} If $-\delta_2 \le -d$, then there is a negative root $x_2$ such that min$\{\rho_2, -d/c \} \le x_2 \le $ min$\{$max$\{\rho_2, -d/c \}, \xi^{(2)}_2 \}$  and a positive root $x_1 \ge \xi^{(2)}_1$. \\
\hrule
\vskip.3cm
\hskip.88cm \underline{{\color{red}{\it Analysis based on solving quadratic equations only}}}
\vskip.3cm
{\bf (i)} If $-d < c \rho_1$, then there are either no real roots or there are two positive roots $x_{1,2} < \rho_1$.
\vskip.1cm
{\bf (ii)} If $c \rho_1 \le - d < 0$ (pictured), then there is a positive root $x_2 < -d/c$ and another positive root $x_1 \ge \rho_1$.
\vskip.1cm
{\bf (iii)} If $0 \le -d < c \rho_2$, then there are either three or one negative roots greater than $\rho_2$ and smaller than or equal to $-d/c$ and there also is one positive root greater than $\rho_1$.
\vskip.1cm
{\bf (iv)} If $c \rho_2 \le -d$, then there is one negative root $x_2$ such that $-d/c < x_2 \le \rho_2$ and a positive root $x_1 > \rho_1$. \\
\end{minipage}}
\\
\hline
\end{tabular}
\end{center}



\begin{center}
\begin{tabular}{|c||c|}
\hline
{\footnotesize \bf Figure 1.3} &  {\footnotesize \bf Figure 1.4} \\
\hline
\hline
{\color{violet} $\bm{b < 0, \quad a <  0, \quad c < 0,$}} &  {\color{violet} $\bm{b < 0, \quad a <  0, \quad c > 0,$}} \\
{\color{violet} $\bm{c_2 < \gamma_2 < {\color{red}c} < 0 < c_0 < \gamma_1 < c_1$}}  & {\color{violet} $\bm{c_2 < \gamma_2 < 0 < {\color{red}c} < c_0 < \gamma_1 < c_1$}} \\
\hline
\hline
\includegraphics[width=80mm]{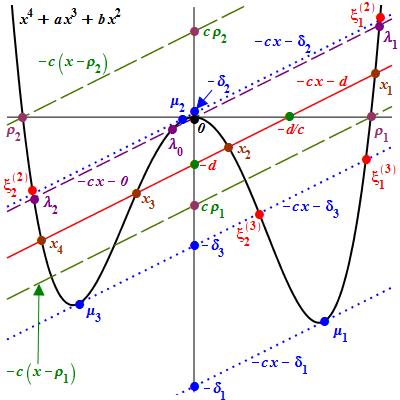} & \includegraphics[width=80mm]{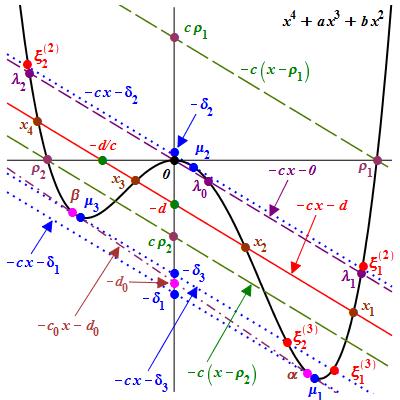} \\
\hline
\multicolumn{1}{|c||}{\begin{minipage}{22em}
\scriptsize
\vskip-.1cm
\hskip3.89cm \underline{{\color{red}{\it Notes}}}
\vskip.01cm
As $c_2 < c < c_1$, the quartic has three stationary points $\mu_i$ and the number of real roots can be either 0, or 2, or 4. There are three tangents $- c x - \delta_i$ to $x^4 + ax^3 + bx^2$ --- each at $\mu_i$.
\vskip.1cm
As $\gamma_2 < c < \gamma_1$, the straight line $-cx$ intersects $x^4 + ax^3 + bx^2$ at three points: $\lambda_{0,1,2}$. Thus $-\delta_3 < 0$. Also: $-\delta_1 < 0$ and $-\delta_2 > 0$.
\vskip.1cm
One can have $ -\delta_1 < c \rho_1 \le -\delta_3 < 0 < - \delta_2 < c \rho_2$ or $ -\delta_1 < -\delta_3 < c \rho_1 < 0 < - \delta_2 < c \rho_2$ (pictured).
\vskip.1cm
The minimum of $x^4 + ax^3 + bx^2$ at negative $x$ is greater than the minimum at positive $x$ by the amount of $-(1/256)a(9a^2-32b)^{3/2} > 0$.
\vskip.1cm
Consideration of whether $- d < c \rho_1$ or $-d > c \rho_1$ would yield sharper bounds in the analysis based on solving cubic equations.
\\
\hrule
\vskip.1cm
\hskip1.54cm \underline{{\color{red}{\it Analysis based on solving cubic equations}}}
\vskip.1cm
{\bf (i)} If $-d < -\delta_1$, then there are no real roots.
\vskip.1cm
{\bf (ii)} If $-\delta_1 \le -d < -\delta_3$, then there is a positive root $x_2$ such that $\xi^{(3)}_2 < x_2 \le \mu_1$ and another positive root $x_1$ such that $\mu_1 \le x_1 < \xi^{(3)}_1$.
\vskip.1cm
{\bf (iii)} If $-\delta_3 \le -d < 0$ (pictured), then there is a negative root $x_4$ such that $\lambda_2 < x_4 \le \mu_3$, another negative root $x_3$ such that $\mu_3 \le x_4 < \lambda_0$, a positive root $x_2 < $ min$ \{ \xi^{(3)}_2, -d/c \}$ and another positive root $x_1$ such that $\xi^{(3)}_1 \le x_1 < \lambda_1$.
\vskip.1cm
{\bf (iv)}  If $0 \le -d < -\delta_2$, then there is a negative root $x_4$ such that $\xi^{(2)}_2 < x_4 \le \lambda_2$, another negative root $x_3$ such that $\lambda_0 \le x_3 < \mu_2$, a third negative root $x_2$ such that $\mu_2 < x_2 < -d/c$, and a positive root $x_1$ such that $\lambda_1 \le x_1 <  \xi^{(2)}_1$.
\vskip.1cm
{\bf (v)} If $-\delta_2 \le -d$, then there is a negative root $x_2$ such that min$\{\rho_2, -d/c \} \le x_2 \le $ min$\{$max$\{\rho_2, -d/c \}, \xi^{(2)}_2 \}$  and a positive root $x_1 \ge \xi^{(2)}_1$.
\vskip.1cm
\hrule
\vskip.1cm
\hskip.88cm \underline{{\color{red}{\it Analysis based on solving quadratic equations only}}}
\vskip.1cm
{\bf (i)} If $-d < c \rho_1$, then there are either no real roots, or there are two positive roots $x_{1,2} < \rho_1$, or there are two two positive roots $x_{1,2} < \rho_1$ together with two negative roots $x_{3,4} > \rho_2$ (the latter appear when $c \rho_1 \ge -\delta_3$).
\vskip.1cm
{\bf (ii)} If $c \rho_1 \le - d < 0$ (pictured), then there are either two positive roots: $x_2 < -d/c$ and $x_1 \ge \rho_1$, or there are two negative roots $x_{3,4} > \rho_2$ together with two positive roots: $x_2 < -d/c$ and $x_1 \ge \rho_1$.
\vskip.1cm
{\bf (iii)} If $0 \le -d < c \rho_2$, then there are either three negative roots $x_{2,3,4}$ such that $\rho_2 < x_{2,3,4} \le -d/c$ together with a positive root $x_1 > \rho_1$, or there is one negative root $x_2$ such that $\rho_2 < x_2 \le -d/c$ and a positive root $x_1 > \rho_1$.
\vskip.1cm
{\bf (iv)} If $c \rho_2 \le -d$, then there is one negative root $x_2$ such that $-d/c < x_2 \le \rho_2$ and one positive root $x_1 > \rho_1$. \\
\end{minipage}}
& \multicolumn{1}{|c|}{\begin{minipage}{22em}
\scriptsize
\vskip.1cm
\hskip3.89cm \underline{{\color{red}{\it Notes}}}\\
\vskip-.3cm
As $c_2 < c < c_1$, the quartic has three stationary points $\mu_i$ and the number of real roots can be either 0, or 2, or 4. There are three tangents $- c x - \delta_i$ to $x^4 + ax^3 + bx^2$ --- each at $\mu_i$.
\vskip.1cm
As $\gamma_2 < c < \gamma_1$, the straight line $-cx$ intersects $x^4 + ax^3 + bx^2$ at three points: $\lambda_{0,1,2}$. Thus $-\delta_1 < 0$. Also $-\delta_3 < 0$ and $-\delta_2 > 0$.
\vskip.1cm
The straight line $ - c_0 x - d_0$ with $c_0 = (1/2)a(b-a^2/4) > 0$ and $d_0 = (1/4)(b-a^2/4)^2 > 0$ is tangent to $x^4 + ax^3 + bx^2$ at two points: $\alpha$ and $\beta = -\alpha - a/2$, where $\alpha$ and $\beta$ are the simultaneous double roots of the varied quartic $x^4 + ax^3 + bx^2 + c_0 x + d_0$.
\vskip.1cm
As $0 < c < c_0$, one has $- \delta_1 < -\delta_3$. Thus $ -\delta_1 < -\delta_3 < c \rho_2 < 0 < - \delta_2 < c \rho_1$.
\vskip.1cm
The minimum of $x^4 + ax^3 + bx^2$ at negative $x$ is greater than the minimum at positive $x$ by the amount of $-(1/256)a(9a^2-32b)^{3/2} > 0$.
\vskip.1cm
Consideration of whether $- d < c \rho_2$ or $-d > c \rho_2$ would yield sharper bounds in the analysis based on solving cubic equations.
\\
\hrule
\vskip.05cm
\hskip1.54cm \underline{{\color{red}{\it Analysis based on solving cubic equations}}} \\
\vskip-.3cm
{\bf (i)} If $-d < -\delta_1$, then there are no real roots.
\vskip.1cm
{\bf (ii)} If $-\delta_1 \le -d < -\delta_3$, then there is one positive root $x_2$ such that $\xi^{(3)}_2 < x_2 \le \mu_1$ and another positive root $x_1$ such that $\mu_1 \le x_2 < \xi^{(3)}_1$.
\vskip.1cm
{\bf (iii)} If $-\delta_3 \le -d < 0$ (pictured), then there is a negative root $x_4$ such that $\lambda_2 < x_4 \le \mu_3$, another negative root $x_3 \ge $ max$\{-d/c,\mu_3\}$, a positive root $x_2$ such that $\lambda_0 < x_2 \le \xi^{(3)}_2$ and another positive root $x_1$ such that $\xi^{(3)}_1 \le x_1 < \lambda_1$.
\vskip.1cm
{\bf (iv)}  If $0 \le -d < -\delta_2$, then there is a negative root $x_4$ such that $\xi^{(2)}_2 < x_4 \le \lambda_2$ and three positive roots: $x_3$ such that $-d/c \le x_3 < \mu_2$, \,\, $x_2$ such that $\mu_2 < x_2 \le \lambda_0$, and  $x_1$ such that $\lambda_1 \le x_1 <  \xi^{(2)}_1$.
\vskip.1cm
{\bf (v)}  If $-\delta_2 \le - d $, then there is a negative root $x_2 \le \xi^{(2)}_2$ and a positive root $x_1$ such that max$\{ \xi^{(2)}_1,$ min$\{\rho_1, -d/c \} \le x_1 < $ max$\{ \rho_1, -d/c \}$.
\\
\hrule
\vskip.05cm
\hskip.88cm \underline{{\color{red}{\it Analysis based on solving quadratic equations only}}} \\
\vskip-.3cm
{\bf (i)} If $- d < c \rho_2$, then there are either no real roots, or there are two positive roots $x_{1,2} < \rho_1$, or there are two negative roots $x_{3,4} > \rho_2$ together with two positive roots $x_{1,2} < \rho_1$.
\vskip.1cm
{\bf (ii)} If $c \rho_2 \le - d < 0$ (pictured), then there is one negative root $x_4 \le \rho_2$, another negative root $x_3 > -d/c$, and two positive roots $x_{1,2}  < \rho_1$.
\vskip.1cm
{\bf (iii)} If $0 \le - d < c \rho_1$, then there is one negative root smaller than $\rho_2$ and either one or three positive roots greater than or equal to $-d/c$ and smaller than $\rho_1$.
\vskip.1cm
{\bf (iv)} If $c \rho_1 \le -d$, then there is one negative root $x_2 < \rho_2$ and a positive root $x_1$ such that $\rho_1 \le x_1 < -d/c$. \\
\end{minipage}}  \\
\hline
\end{tabular}
\end{center}


\begin{center}
\begin{tabular}{|c||c|}
\hline
{\footnotesize \bf Figure 1.5} &  {\footnotesize \bf Figure 1.6} \\
\hline
\hline
{\color{violet} $\bm{b < 0, \quad a <  0, \quad c > 0,$}} &  {\color{violet}$\bm{b < 0, \quad a <  0, \quad c > 0,}$} \\
{\color{violet} $\bm{c_2 < \gamma_2  < 0 < c_0 < {\color{red}c} < \gamma_1 < c_1$}} & {\color{violet} $\bm{c_2 < \gamma_2 < 0 < c_0 < \gamma_1 < {\color{red}c} < c_1$}} \\
\hline
\hline
\includegraphics[width=80mm]{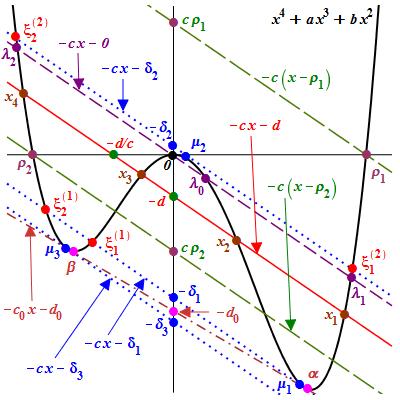} & \includegraphics[width=80mm]{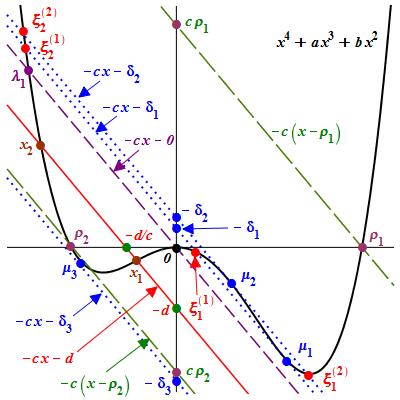} \\
\hline
\multicolumn{1}{|c||}{\begin{minipage}{22em}
\scriptsize
\vskip.1cm
\hskip3.89cm \underline{{\color{red}{\it Notes}}}\\
\vskip-.2cm
As $c_2 < c < c_1$, the quartic has three stationary points $\mu_i$ and the number of real roots can be either 0, or 2, or 4. There are three tangents $- c x - \delta_i$ to $x^4 + ax^3 + bx^2$ --- each at $\mu_i$.
\vskip.1cm
As $\gamma_2 < c < \gamma_1$, the straight line $-cx$ intersects $x^4 + ax^3 + bx^2$ at three points: $\lambda_{0,1,2}$. Thus $-\delta_1 < 0$. Also: $-\delta_3 < 0$ and $-\delta_2 > 0$.
\vskip.1cm
The straight line $ - c_0 x - d_0$ with $c_0 = (1/2)a(b-a^2/4) > 0$ and $d_0 = (1/4)(b-a^2/4)^2 > 0$ is tangent to $x^4 + ax^3 + bx^2$ at two points: $\alpha$ and $\beta = -\alpha - a/2$, where $\alpha$ and $\beta$ are the simultaneous double roots of the varied quartic $x^4 + ax^3 + bx^2 + c_0 x + d_0$.
\vskip.1cm
As $c > c_0$, one has $- \delta_3 < -\delta_1$. Thus one can have $ -\delta_3 < c \rho_2 \le -\delta_1 < 0 < -\delta_2 < c \rho_1$  or $ -\delta_3 < -\delta_1 < c \rho_2 < 0 < -\delta_2 < c \rho_1$ (pictured).
\vskip.1cm
The minimum of $x^4 + ax^3 + bx^2$ at negative $x$ is greater than the minimum at positive $x$ by the amount of $-(1/256)a(9a^2-32b)^{3/2} > 0$.
\vskip.1cm
Consideration of whether $- d < c \rho_2$ or $-d > c \rho_2$ would yield sharper bounds in the analysis based on solving cubic equations. \\
\hrule
\vskip.1cm
\hskip1.54cm \underline{{\color{red}{\it Analysis based on solving cubic equations}}} \\
\vskip-.3cm
{\bf (i)} If $-d < -\delta_3$, then there are no real roots.
\vskip.1cm
{\bf (ii)} If $-\delta_3 \le -d < -\delta_1$, then there is a negative root $x_2$ such that $\xi^{(1)}_2 < x_2 \le \mu_3$ and another negative root $x_1$ such that $\mu_3 \le x_1 < \xi^{(1)}_1$.
\vskip.1cm
{\bf (iii)} If $-\delta_1 \le -d < 0$ (pictured), then there is a negative root $x_4$ such that $\lambda_2 < x_4 \le \xi^{(1)}_2$, another negative root $x_3 \ge $ max$\{-d/c, \xi^{(1)}_1 \}$ a positive root $x_2$ such that $\lambda_0 < x_2 \le \mu_1$ and another positive root $x_1$ such that  $\mu_1 \le x_1 < \lambda_1$.
\vskip.1cm
{\bf (iv)}  If $0 \le -d < -\delta_2$, then there is a negative root $x_4$ such that $\xi^{(2)}_2 < x_4 \le \lambda_2$ and three positive roots $x_{1,2,3}$ such that: $-d/c \le x_3 < \mu_2$, \,\, $\mu_2 < x_2 \le \lambda_0$, and $\lambda_1 \le x_1 <  \xi^{(2)}_1$.
\vskip.1cm
{\bf (v)}  If $-\delta_2 \le - d$, then there is a negative root $x_2 \le \xi^{(2)}_2$ and a positive root $x_1$ such that max$\{ \xi^{(2)}_1,$ min$\{\rho_1, -d/c \} \le x_1 < $ max$\{ \rho_1, -d/c \}$. \\
\hrule
\vskip.1cm
\hskip.88cm \underline{{\color{red}{\it Analysis based on solving quadratic equations only}}} \\
\vskip-.3cm
{\bf (i)} If $-d < c \rho_2$, then there are either no real roots, or there are two negative roots $x_{1,2} > \rho_2$, or there are two negative roots $x_{3,4} > \rho_2$ and two positive roots $x_{1,2} < \rho_1$ (the latter appear when $c \rho_2 \ge -\delta_1$).
\vskip.1cm
{\bf (ii)} If $c \rho_2 \le - d < 0$ (pictured), then there are either two negative roots: $x_2 \le \rho_2$ and $x_1 > -d/c$ , or there are two negative roots: $x_4 \le \rho_2$ and $x_3 > -d/c$ and two positive roots: $x_{1,2} < \rho_1$.
\vskip.1cm
{\bf (iii)} If $0 \le - d < c \rho_1$, then there is one negative root smaller than $\rho_2$ together with either one or three positive roots greater than or equal to $-d/c$ and smaller than $\rho_1$.
\vskip.1cm
{\bf (iv)} If $c \rho_1 \le -d$, then there is one negative root $x_2 < \rho_2$ and one positive root $x_1$ such that $\rho_1 \le x_1 < -d/c$.
\\
\end{minipage}}
& \multicolumn{1}{|c|}{\begin{minipage}{22em}
\scriptsize
\vskip-2.9cm
\hskip3.89cm \underline{{\color{red}{\it Notes}}}\\
\vskip-.25cm
As $c_2 < c < c_1$, the quartic has three stationary points $\mu_i$ and the number of real roots can be either 0, or 2, or 4. There are three tangents $- c x - \delta_i$ to $x^4 + ax^3 + bx^2$ --- each at $\mu_i$.
\vskip.1cm
As $\gamma_2 < \gamma_1 < c$, the straight line $-cx$ intersects $x^4 + ax^3 + bx^2$ at one point only ($\lambda_1$). Thus $-\delta_1 > 0$. Also: $-\delta_3 < 0$ and $-\delta_2 > 0$.
\vskip.1cm
Obviously, $ -\delta_3 < c \rho_2 < 0 < -\delta_1 < -\delta_2 < c \rho_1$.
\vskip.1cm
The minimum of $x^4 + ax^3 + bx^2$ at negative $x$ is greater than the minimum at positive $x$ by the amount of $-(1/256)a(9a^2-32b)^{3/2} > 0$.
\vskip.1cm
Consideration of whether $- d < c \rho_2$ or $-d > c \rho_2$ would yield sharper bounds in the analysis based on solving cubic equations.
\vskip.1cm
\hrule
\vskip.3cm
\hskip1.54cm \underline{{\color{red}{\it Analysis based on solving cubic equations}}} \\
\vskip-.1cm
{\bf (i)} If $-d < -\delta_3$, then there are no real roots.
\vskip.1cm
{\bf (ii)} If $-\delta_3 \le -d < 0$ (pictured), then there is a negative root $x_2$ such that $\lambda_1 < x_2 \le \mu_3$ and another negative root $x_1 \ge $ max$ \{ \mu_3, -d/c \}$.
\vskip.1cm
{\bf (iii)} If $0 \le -d < -\delta_1$, then there is a negative root $x_2$ such that  $\xi^{(1)}_2 < x_2 \le \lambda_1$ and a non-negative root $x_1 \le -d/c$.
\vskip.1cm
{\bf (iv)}  If $-\delta_1 \le - d < -\delta_2$, then there is a negative root $x_4$ such that $\xi^{(2)}_2 < x_4  \le \xi^{(1)}_2$  and three positive roots: $x_3$ such that max$\{ -d/c, \xi^{(1)}_1 \} \le x_3 < \mu_2$, \,\, $x_2$ such that $\mu_2 < x_2 \le \mu_1$ and $x_1$ such that $\mu_1 \le x_1 < \xi^{(2)}_1$.
\vskip.1cm
{\bf (v)}  If $-\delta_2 \le - d$, then there is a negative root $x_2 \le \xi^{(2)}_2$ and a positive root $x_1$ such that max$\{ \xi^{(2)}_1,$ min$\{\rho_1, -d/c \} \le x_1 < $ max$\{ \rho_1, -d/c \}$.
\vskip.1cm
\hrule
\vskip.3cm
\hskip.88cm \underline{{\color{red}{\it Analysis based on solving quadratic equations only}}}
\vskip.3cm
{\bf (i)} If $- d < c \rho_2$, then there are either no real roots or there are two negative roots $x_{1,2} > \rho_2$.
\vskip.1cm
{\bf (ii)} If $c \rho_2 \le - d < 0$ (pictured), then there is a negative root $x_2 \le \rho_2$ and another negative root $x_1 > -d/c$.
\vskip.1cm
{\bf (iii)} If $0 \le - d < c \rho_1$, then there is one negative root smaller than $\rho_2$ together with either one or three positive roots greater than $-d/c$ and smaller than or equal to $\rho_1$.
\vskip.1cm
{\bf (iv)} If $c \rho_1 \le -d$, then there is one negative root $x_2 < \rho_2$ and a positive root $x_1$ such that $\rho_1 \le x_1 < -d/c$. \\
\end{minipage}}  \\
\hline
\end{tabular}
\end{center}



\begin{center}
\begin{tabular}{|c||c|}
\hline
{\footnotesize \bf Figure 1.7} &  {\footnotesize \bf Figure 1.8} \\
\hline
\hline
{\color{violet} $\bm{b < 0, \quad a < 0, \quad c > 0,$}} &  {\color{violet}$\bm{b < 0, \quad a > 0, \quad c < 0,}$} \\
{\color{violet} $\bm{c_2 < \gamma_2 < 0  < c_0 < \gamma_1 < c_1 < {\color{red}c} $}} &  {\color{violet} $\bm{{\color{red}c} < c_2 < \gamma_2  < c_0 < 0 < \gamma_1 < c_1 $}} \\
\hline
\hline
\includegraphics[width=80mm]{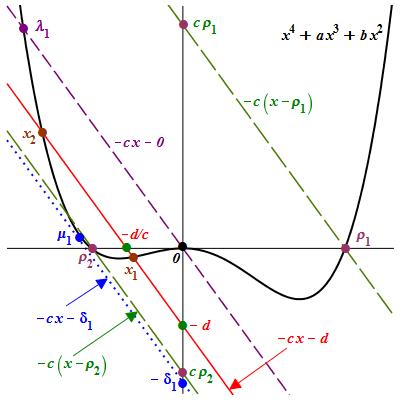} & \includegraphics[width=80mm]{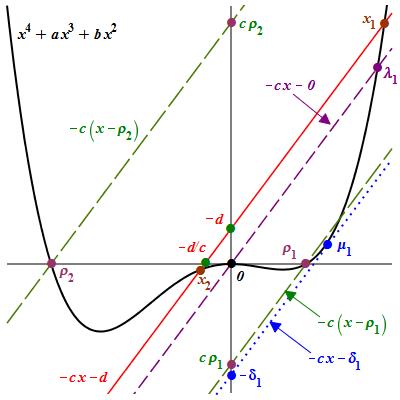} \\
\hline
\multicolumn{1}{|c||}{\begin{minipage}{22em}
\scriptsize
\vskip0.3cm
\hskip3.89cm \underline{{\color{red}{\it Notes}}}\\
\vskip.05cm
As $c_2 < c_1 < c$, the quartic has a single stationary point $\mu_1$ and the number of real roots can be either 0 or 2. There is only one tangent $- c x - \delta_1$ to $x^4 + ax^3 + bx^2$ --- the one at $\mu_1$.
\vskip.1cm
As $\gamma_2 < \gamma_1 < c$, the straight line $-cx$ intersects $x^4 + ax^3 + bx^2$ at one point only ($\lambda_1$).
\vskip.1cm
One could have $\mu_1 < \rho_2$ (pictured) or $\mu_1 \ge \rho_2$.
\vskip.1cm
Obviously, $ -\delta_1 < c \rho_2 < 0 < c \rho_1$.
\vskip.1cm
The minimum of $x^4 + ax^3 + bx^2$ at negative $x$ is greater than the minimum at positive $x$ by the amount of $-(1/256)a(9a^2-32b)^{3/2} > 0$.
\vskip.1cm
Consideration of whether $- d < c \rho_2$ or $-d > c \rho_2$ would yield sharper bounds in the analysis based on solving cubic equations. \\
\hrule
\vskip.3cm
\hskip1.54cm \underline{{\color{red}{\it Analysis based on solving cubic equations}}} \\
\vskip.1cm
{\bf (i)} If $-d < -\delta_1$, then there are no real roots.
\vskip.1cm
{\bf (ii)} If $-\delta_1 \le -d < 0$ (pictured), then there is a negative root $x_2$ such that $\lambda_1 < x_2 \le \mu_1$ and another negative root $x_1$ such that $\mu_1 \le x_1 < -d/c$.
\vskip.1cm
{\bf (iii)} If $0 \le -d$, then there is a negative root $x_2 \le \lambda_1$ and a positive root $x_1$ such that
min$\{\rho_1, -d/c \} \le x_1 \le $ max$\{\rho_1, -d/c \}$. \\
\hrule
\vskip.3cm
\hskip.88cm \underline{{\color{red}{\it Analysis based on solving quadratic equations only}}}
\vskip.3cm
{\bf (i)} If $- d < c \rho_2$, then there are either no real roots or there are two negative roots (these are smaller than $\rho_2$ if $\mu_1 < \rho_2$ or greater than or equal to $\rho_2$ if $\mu_1 \ge \rho_2$).
\vskip.1cm
{\bf (ii)} If $c \rho_2 \le - d < 0$ (pictured), then there is a negative root $x_2 \le \rho_2$ and another negative root $x_1 \ge -d/c$.
\vskip.1cm
{\bf (iii)} If $0 \le - d < c \rho_1$, then there is one negative root $x_2 < \rho_2$ and a non-negative root $x_1$ such that $-d/c \le x_1 < \rho_1$.
\vskip.1cm
{\bf (iv)} If $c \rho_1 \le -d$, then there is one negative root $x_2 < \rho_2$ and a positive root $x_1$ such that $\rho_1 \le x_1 \le -d/c$. \\
\end{minipage}}
& \multicolumn{1}{|c|}{\begin{minipage}{22em}
\scriptsize
\vskip0.1cm
\hskip3.89cm \underline{{\color{red}{\it Notes}}}\\
\vskip.05cm
As $c < c_2 < c_1$, the quartic has a single stationary point $\mu_1$ and the number of real roots can be either 0 or 2. There is only one tangent $- c x - \delta_1$ to $x^4 + ax^3 + bx^2$ --- the one at $\mu_1$.
\vskip.1cm
As $c < \gamma_2 < \gamma_1$, the straight line $-cx$ intersects $x^4 + ax^3 + bx^2$ at one point only ($\lambda_1$).
\vskip.1cm
One could have $\mu_1 > \rho_1$ (pictured) or $\mu_1 \le \rho_1$.
\vskip.1cm
Obviously, $ -\delta_1 < c \rho_1 < 0 < c \rho_2$.
\vskip.1cm
The minimum of $x^4 + ax^3 + bx^2$ at negative $x$ is smaller than the minimum at positive $x$ by the amount of $(1/256)a(9a^2-32b)^{3/2}$.
\vskip.1cm
Consideration of whether $- d < c \rho_1$ or $-d > c \rho_1$ would yield sharper bounds in the analysis based on solving cubic equations. \\
\vskip.1cm
\hrule
\vskip.3cm
\hskip1.54cm \underline{{\color{red}{\it Analysis based on solving cubic equations}}} \\
\vskip.1cm
{\bf (i)} If $-d < -\delta_1$, then there are no real roots.
\vskip.1cm
{\bf (ii)} If $-\delta_1 \le -d < 0$, then there is a positive root $x_2 \le \mu_1$ and another positive root $x_1$ such that $\mu_1 \le x_1 < \lambda_1$.
\vskip.1cm
{\bf (iii)} If $0 \le -d$ (pictured), then there is a negative root $x_2$ such that min$\{\rho_2, -d/c \} \le x_2 \le $ max$\{\rho_2, -d/c \}$  and a positive root $x_1 \ge \lambda_1$. \\
\hrule
\vskip.3cm
\hskip.88cm \underline{{\color{red}{\it Analysis based on solving quadratic equations only}}}
\vskip.3cm
{\bf (i)} If $- d < c \rho_1$, then there are either no real roots or there are two positive roots (these are greater than $\rho_1$ if $\mu_1 > \rho_1$ or smaller than or equal to $\rho_1$ if $\mu_1 \le \rho_1$).
\vskip.1cm
{\bf (ii)} If $c \rho_1 \le - d < 0$, then there is a positive root $x_2 < -d/c$ and another positive root $x_1 \ge \rho_1$.
\vskip.1cm
{\bf (iii)} If $0 \le - d < c \rho_2$ (pictured), then there is one non-positive root $x_2$ such that $\rho_2 < x_2 \le -d/c$ and a positive root $x_1 > \rho_1$.
\vskip.1cm
{\bf (iv)} If $c \rho_2 \le -d$, then there is one negative root $x_2$ such that $-d/c \le x_2 \le \rho_2$ and a positive root $x_1 > \rho_1$. \\
\end{minipage}}  \\
\hline
\end{tabular}
\end{center}



\begin{center}
\begin{tabular}{|c||c|}
\hline
{\footnotesize \bf Figure 1.9} &  {\footnotesize \bf Figure 1.10} \\
\hline
\hline
{\color{violet} $\bm{b < 0, \quad a > 0, \quad c < 0,$}} &  {\color{violet}$\bm{b < 0, \quad a > 0, \quad c < 0,}$} \\
{\color{violet} $\bm{c_2 < {\color{red}c} < \gamma_2  < c_0 < 0 < \gamma_1 < c_1$}} & {\color{violet} $\bm{c_2 < \gamma_2  < {\color{red}c} < c_0 < 0 < \gamma_1 < c_1$}} \\
\hline
\hline
\includegraphics[width=80mm]{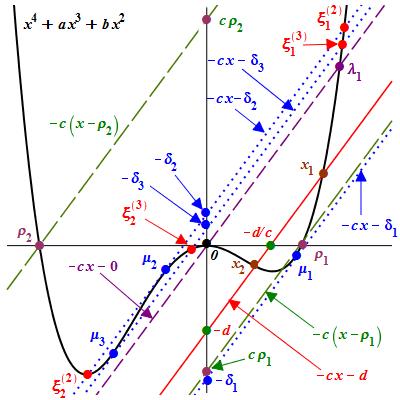} & \includegraphics[width=80mm]{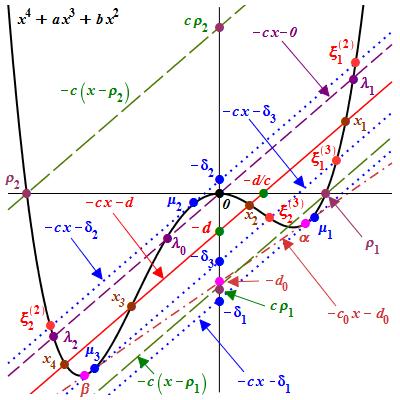} \\
\hline
\multicolumn{1}{|c||}{\begin{minipage}{22em}
\scriptsize
\vskip-2.6cm
\hskip3.89cm \underline{{\color{red}{\it Notes}}}\\
\vskip-.3cm
As $c_2 < c < c_1$, the quartic has three stationary points $\mu_i$ and the number of real roots can be either 0, or 2, or 4. There are three tangents $- c x - \delta_i$ to $x^4 + ax^3 + bx^2$ --- each at $\mu_i$.
\vskip.1cm
As $\gamma_2 < \gamma_1 < c$, the straight line $-cx$ intersects $x^4 + ax^3 + bx^2$ at one point only ($\lambda_1$). Thus $-\delta_3 > 0$. Also: $-\delta_1 < 0$ and $-\delta_2 > 0$.
\vskip.1cm
Obviously, $ -\delta_1 < c \rho_1 < 0 < -\delta_3 < -\delta_2 < c \rho_2$.
\vskip.1cm
The minimum of $x^4 + ax^3 + bx^2$ at negative $x$ is smaller than the minimum at positive $x$ by the amount of $(1/256)a(9a^2-32b)^{3/2}$.
\vskip.1cm
Consideration of whether $- d < c \rho_1$ or $-d > c \rho_1$ would yield sharper bounds in the analysis based on solving cubic equations.
\\
\hrule
\vskip.2cm
\hskip1.54cm \underline{{\color{red}{\it Analysis based on solving cubic equations}}} \\
\vskip-.1cm
{\bf (i)} If $-d < -\delta_1$, then there are no real roots.
\vskip.1cm
{\bf (ii)} If $-\delta_1 \le -d < 0$ (pictured), then there is a positive root $x_2 \le $ min$ \{ -d/c, \mu_1 \}$ and another positive root $x_1$ such that $\mu_1 \le x_1 <  \lambda_1 $.
\vskip.1cm
{\bf (iii)} If $0 \le -d < -\delta_3$, then there is a positive root $x_1$ such that $\lambda_1 \le x_1 < \xi^{(3)}_1$, and a non-positive root $x_2$ such that $\xi^{(3)}_2 < x_2  \le -d/c$ .
\vskip.1cm
{\bf (iv)} If $-\delta_3 \le -d < - \delta_2$, then there is a positive root $x_1$ such that $\xi^{(3)}_1 \le x_1 < \xi^{(2)}_1$, a negative root $x_2$ such that $\mu_2 < x_2 \le $ min$\{\xi^{(3)}_2,-d/c\}$, another negative root $x_3$ such that $\mu_3 \le x_3 < \mu_2$, and a third negative root $x_4$ such that $\xi^{(2)}_2 < x_4 \le \mu_3$.
\vskip.1cm
{\bf (v)} If $-\delta_2 \le -d$, then there is a negative root $x_2$ such that min$\{\rho_2, -d/c \} \le x_2 \le $ min$\{$max$\{\rho_2, -d/c \}, \xi^{(2)}_2 \}$  and a positive root $x_1 \ge \xi^{(2)}_1$.
\vskip.1cm
\hrule
\vskip.3cm
\hskip.88cm \underline{{\color{red}{\it Analysis based on solving quadratic equations only}}}
\vskip.3cm
{\bf (i)} If $-d < c \rho_1$, then there are either no real roots or there are two positive roots $x_{1,2} < \rho_1$.
\vskip.1cm
{\bf (ii)} If $c \rho_1 \le - d < 0$ (pictured), then there is one positive root $x_2 < -d/c$ and another positive root $x_1 \ge \rho_1$.
\vskip.1cm
{\bf (iii)} If $0 \le -d < c \rho_2$, then there is one positive root $x_1 > \rho_1$ together with either three negative roots $x_{1,2,3}$ such that $\rho_2 < x_{1,2,3} \le -d/c$, or with one negative root $x_2$ such that $\rho_2 < x_2 \le -d/c$.
\vskip.1cm
{\bf (iv)} If $c \rho_2 \le -d$, then there is one negative root $x_2$ such that $-d/c < x_2 \le \rho_2$ and one positive root $x_1 > \rho_1.$ \\
\end{minipage}}
& \multicolumn{1}{|c|}{\begin{minipage}{22em}
\scriptsize
\vskip.05cm
\hskip3.89cm \underline{{\color{red}{\it Notes}}}\\
\vskip-.35cm
As $c_2 < c < c_1$, the quartic has three stationary points $\mu_i$ and the number of real roots can be either 0, or 2, or 4. There are three tangents $- c x - \delta_i$ to $x^4 + ax^3 + bx^2$ --- each at $\mu_i$.
\vskip.1cm
As $\gamma_2 < c < \gamma_1$, the straight line $-cx$ intersects $x^4 + ax^3 + bx^2$ at three points: $\lambda_{0,1,2}$. Thus $-\delta_3 < 0$. Also: $-\delta_1 < 0$ and $-\delta_2 > 0$.
\vskip.1cm
The straight line $ - c_0 x - d_0$ with $c_0 = (1/2)a(b-a^2/4) < 0$ and $d_0 = (1/4)(b-a^2/4)^2 > 0$ is tangent to $x^4 + ax^3 + bx^2$ at two points: $\alpha$ and $\beta = -\alpha - a/2$, where $\alpha$ and $\beta$ are the simultaneous double roots of the varied quartic $x^4 + ax^3 + bx^2 + c_0 x + d_0$.
\vskip.1cm
As $c < c_0$, one has $-\delta_1 < - \delta_3$. Thus one can have $ -\delta_1 < c \rho_1 < -\delta_3 < 0 < -\delta_2 < c \rho_2$ (pictured) or $ -\delta_1 < -\delta_3 \le c \rho_1 < 0 < -\delta_2 < c \rho_2$.
\vskip.1cm
The minimum of $x^4 + ax^3 + bx^2$ at negative $x$ is smaller than the minimum at positive $x$ by the amount of $(1/256)a(9a^2-32b)^{3/2}$.
\vskip.1cm
Consideration of whether $- d < c \rho_1$ or $-d > c \rho_1$ would yield sharper bounds in the analysis based on solving cubic equations.
\hrule
\vskip.1cm
\hskip1.54cm \underline{{\color{red}{\it Analysis based on solving cubic equations}}} \\
\vskip-.3cm
{\bf (i)} If $-d < -\delta_1$, then there are no real roots.
\vskip.1cm
{\bf (ii)} If $-\delta_1 \le -d < -\delta_3$, then there is a positive root $x_2$ such that $\xi^{(3)}_2 < x_2 \le $ min$ \{ -d/c, \mu_1 \}$ and another positive root $x_1$ such that $\mu_1 \le x_1 <  \xi^{(3)}_1 $.
\vskip.1cm
{\bf (iii)}  If $-\delta_3 \le -d < 0$ (pictured), then there is a negative root $x_4$ such that $\lambda_2 < x_4 \le \mu_3$, another negative root $x_3$ such that $\mu_3 \le x_3 < \lambda_0$, a positive root $x_2 < $ min$ \{  \xi^{(3)}_2, -d/c \}$ and another positive root $x_1$ such that  $\xi^{(3)}_1 \le x_1 < \lambda_1$.
\vskip.1cm
{\bf (iv)}  If $0 \le -d < -\delta_2$, then there is a negative root $x_4$ such that $\xi^{(2)}_2 < x_4 \le \lambda_2$, a negative root $x_3$ such that $\lambda_0 \le x_3 < \mu_2$, a non-positive root $x_2$ such that $\mu_2 < x_2 \le -d/c$ and a positive root $x_1$ such that $\lambda_1 \le x_1 <  \xi^{(2)}_1$.
\vskip.1cm
{\bf (v)} If $-\delta_2 \le -d$, then there is a negative root $x_2$ such that min$\{\rho_2, -d/c \} \le x_2 \le $ min$\{$max$\{\rho_2, -d/c \}, \xi^{(2)}_2 \}$  and a positive root $x_1 \ge \xi^{(2)}_1$.
\\
\hrule
\vskip.1cm
\hskip.88cm \underline{{\color{red}{\it Analysis based on solving quadratic equations only}}}
\vskip.1cm
{\bf (i)} If $-d < c \rho_1$, then there are either no real roots, or there are two positive roots $x_{1,2} < \rho_1$, or there are two positive roots $x_{1,2} < \rho_1$ and two negative roots $x_{3,4} > \rho_2$ (the latter appear when $c \rho_1 \ge -\delta_3$).
\vskip.1cm
{\bf (ii)} If $c \rho_1 \le - d < 0$ (pictured), then there are either two positive roots $x_{1,2}$ such that $x_1 < -d/c$ and $x_2 \ge \rho_1$ or there are two positive roots  $x_{1,2}$ such that $x_1 < -d/c$ and $x_2 \ge \rho_1$ and two negative roots $x_{3,2} > \rho_2$.
\vskip.1cm
{\bf (iii)} If $0 \le -d < c \rho_2$, then there is one positive root $x_1 > \rho_1$ together with either one or three negative roots greater than $\rho_2$ and smaller than or equal to $-d/c$.
\vskip.1cm
{\bf (iv)} If $c \rho_2 \le -d$, then there is one negative root $x_2$ such that $-d/c < x_2 \le \rho_2$ and one positive root $x_1 > \rho_1.$
\\
\end{minipage}}  \\
\hline
\end{tabular}
\end{center}


\begin{center}
\begin{tabular}{|c||c|}
\hline
{\footnotesize \bf Figure 1.11} &  {\footnotesize \bf Figure 1.12} \\
\hline
\hline
{\color{violet} $\bm{b < 0, \quad a > 0, \quad c < 0,$}} &  {\color{violet}$\bm{b < 0, \quad a > 0, \quad c > 0,}$} \\
{\color{violet} $\bm{c_2 < \gamma_2 < c_0 < {\color{red}c} < 0 < \gamma_1 < c_1$}}  & {\color{violet} $\bm{c_2 < \gamma_2 < c_0 < 0 < {\color{red}c} < \gamma_1 < c_1$}} \\
\hline
\hline
\includegraphics[width=80mm]{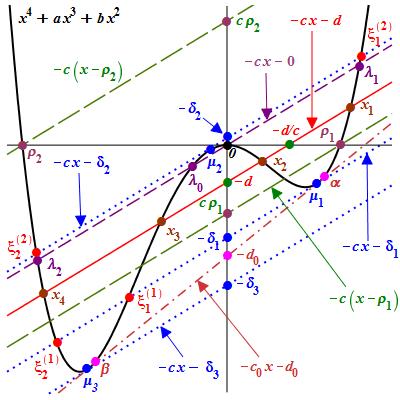} & \includegraphics[width=80mm]{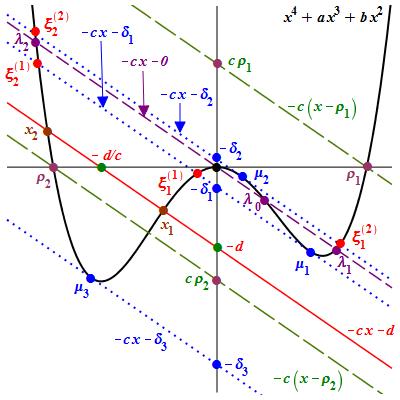} \\
\hline
\multicolumn{1}{|c||}{\begin{minipage}{22em}
\scriptsize
\vskip.1cm
\hskip3.89cm \underline{{\color{red}{\it Notes}}}\\
\vskip-.1cm
As $c_2 < c < c_1$, the quartic has three stationary points $\mu_i$ and the number of real roots can be either 0, or 2, or 4. There are three tangents $- c x - \delta_i$ to $x^4 + ax^3 + bx^2$ --- each at $\mu_i$.
\vskip.1cm
As $\gamma_2 < c < \gamma_1$, the straight line $-cx$ intersects $x^4 + ax^3 + bx^2$ at three points: $\lambda_{0,1,2}$. Thus $-\delta_3 < 0$. Also: $-\delta_1 < 0$ and $-\delta_2 > 0$.
\vskip.1cm
The straight line $ - c_0 x - d_0$ with $c_0 = (1/2)a(b-a^2/4) < 0$ and $d_0 = (1/4)(b-a^2/4)^2 > 0$ is tangent to $x^4 + ax^3 + bx^2$ at two points: $\alpha$ and $\beta = -\alpha - a/2$, where $\alpha$ and $\beta$ are the simultaneous double roots of the varied quartic $x^4 + ax^3 + bx^2 + c_0 x + d_0$.
\vskip.1cm
As $c_0 < c$, one has $- \delta_3 < -\delta_1$.
\vskip.1cm
Obviously, $ -\delta_3 < -\delta_1 < c \rho_1 < 0 < -\delta_2 < c \rho_2$.
\vskip.1cm
The minimum of $x^4 + ax^3 + bx^2$ at negative $x$ is smaller than the minimum at positive $x$ by the amount of $(1/256)a(9a^2-32b)^{3/2}$.
\vskip.1cm
Consideration of whether $- d < c \rho_1$ or $-d > c \rho_1$ would yield sharper bounds in the analysis based on solving cubic equations.
\\
\hrule
\vskip.1cm
\hskip1.54cm \underline{{\color{red}{\it Analysis based on solving cubic equations}}} \\
\vskip-.2cm
{\bf (i)} If $-d < -\delta_3$, then there are no real roots.
\vskip.1cm
{\bf (ii)} If $-\delta_3 \le -d < -\delta_1$, then there is a negative root $x_2$ such that $\xi^{(1)}_2 < x_2 \le \mu_3$ and another negative root $x_1$ such that $\mu_3 \le x_2 < \xi^{(1)}_1$.
\vskip.1cm
{\bf (iii)} If $-\delta_1 \le -d < 0$ (pictured), then there is a negative root $x_4$ such that $\lambda_2 < x_4 \le \xi^{(1)}_2$, another negative root $x_3$ such that $\xi^{(1)}_1 \le x_3 < \lambda_0$, a positive root $x_2 \le$ min$ \{ -d/c, \mu_1 \}$ and another positive root $x_1$ such that $\mu_1 \le x_1 < \lambda_1$.
\vskip.1cm
{\bf (iv)}  If $0 \le -d < -\delta_2$, then there is a negative root $x_4$ such that $\xi^{(2)}_2 < x_4 \le \lambda_2$, a negative root $x_3$ such that $\lambda_0 \le x_3 < \mu_2$, a non-positive root $x_2$ such that $\mu_2 < x_2 \le -d/c$ and a positive root $x_1$ such that $\lambda_1 \le x_1 <  \xi^{(2)}_1$.
\vskip.1cm
{\bf (v)} If $-\delta_2 \le -d$, then there is a negative root $x_2$ such that min$\{\rho_2, -d/c \} \le x_2 \le $ min$\{$max$\{\rho_2, -d/c \}, \xi^{(2)}_2 \}$  and a positive root $x_1 \ge \xi^{(2)}_1$.
\\
\hrule
\vskip.1cm
\hskip.88cm \underline{{\color{red}{\it Analysis based on solving quadratic equations only}}}
\vskip.1cm
{\bf (i)} If $- d < c \rho_1$, then there are either no real roots, or there are two negative roots $x_{1,2} > \rho_2$, or there are two negative roots $x_{3,4} > \rho_2$ together with two positive roots $x_{1,2} < \rho_1$.
\vskip.1cm
{\bf (ii)} If $c \rho_1 \le - d < 0$ (pictured), then there are two negative roots $x_{3,4}  > \rho_2$, a positive root $x_2 < -d/c$ and another positive root $x_1 \ge \rho_1$.
\vskip.1cm
{\bf (iii)} If $0 \le - d < c \rho_2$, then there is one positive root bigger than $\rho_1$ and either one or three negative roots greater than $\rho_2$ and smaller than or equal to $-d/c$.
\vskip.1cm
{\bf (iv)} If $c \rho_2 \le -d$, then there is a negative root $x_2$ such that $-d/c < x_2 \le \rho_2$ and a positive root $x_1 > \rho_1$.
\\
\end{minipage}}
& \multicolumn{1}{|c|}{\begin{minipage}{22em}
\scriptsize
\vskip-0.85cm
\hskip3.89cm \underline{{\color{red}{\it Notes}}}
\vskip.2cm
As $c_2 < c < c_1$, the quartic has three stationary points $\mu_i$ and the number of real roots can be either 0, or 2, or 4. There are three tangents $- c x - \delta_i$ to $x^4 + ax^3 + bx^2$ --- each at $\mu_i$.
\vskip.1cm
As $\gamma_2 < c < \gamma_1$, the straight line $-cx$ intersects $x^4 + ax^3 + bx^2$ at three points: $\lambda_{0,1,2}$. Thus $-\delta_1 < 0$. Also: $-\delta_3 < 0$ and $-\delta_2 > 0$.
\vskip.1cm
One can have $ -\delta_3 <  c \rho_2 < -\delta_1 < 0 < -\delta_2 < c \rho_1$ (pictured) or $ -\delta_3 < -\delta_1 \le c \rho_2 < 0 < -\delta_2 < c \rho_1$.
\vskip.1cm
The minimum of $x^4 + ax^3 + bx^2$ at negative $x$ is smaller than the minimum at positive $x$ by the amount of $(1/256)a(9a^2-32b)^{3/2}$.
\vskip.1cm
Consideration of whether $- d < c \rho_2$ or $-d > c \rho_2$ would yield sharper bounds in the analysis based on solving cubic equations.
\\
\hrule
\vskip.1cm
\hskip1.54cm \underline{{\color{red}{\it Analysis based on solving cubic equations}}} \\
\vskip-.2cm
{\bf (i)} If $-d < -\delta_3$, then there are no real roots.
\vskip.1cm
{\bf (ii)} If $-\delta_3 \le - d < -\delta_1$ (pictured), then there is a negative root $x_2$ such that $\xi^{(1)}_2 < x_2 \le \mu_3$ and another negative root $x_1$ such that max$ \{-d/c, \mu_3 \} \le x_1 < \xi^{(1)}_1 $.
\vskip.1cm
{\bf (iii)}  If $-\delta_1 \le -d < 0$, then there is a negative root $x_4$ such that $\lambda_2 < x_4 \le \xi^{(1)}_2$, another negative root $x_3 \ge$ max$\{-d/c, \xi^{(1)}_1 \}$, a positive root $x_2$ such that $\lambda_0 < x_2 \le \mu_1$ and another positive root $x_1$ such that  $\mu_1 \le x_1 < \lambda_1$.
\vskip.1cm
{\bf (iv)}  If $0 \le -d < -\delta_2$, then there is a negative root $x_4$ such that $\xi^{(2)}_2 < x_4 \le \lambda_2$, a non-negative root $x_3$ such that $-d/c \le x_3 < \mu_2$, a positive root $x_2$ such that $\mu_2 < x_2 \le \lambda_0$ and another positive root $x_1$ such that $\lambda_1 \le x_1 <  \xi^{(2)}_1$.
\vskip.1cm
{\bf (v)} If $-\delta_2 \le -d$, then there is a negative root $x_2 \le \xi^{(2)}_2$ and a positive root $x_1$ such that max$\{$min$\{\rho_1, -d/c \}, \xi^{(2)}_1 \} \le x_1 < $ max$\{\rho_1, -d/c \}$ .
\\
\hrule
\vskip.1cm
\hskip.88cm \underline{{\color{red}{\it Analysis based on solving quadratic equations only}}}
\vskip.1cm
{\bf (i)} If $- d < c \rho_2$, then there are either no real roots, or there are two negative roots $x_{1,2} > \rho_2$, or there are two negative roots $x_{3,4} > \rho_2$ together with two positive roots $x_{1,2} < \rho_1$ (the latter appear when $c \rho_2 \ge -\delta_1$).
\vskip.1cm
{\bf (ii)} If $c \rho_2 \le - d < 0$ (pictured), then there are either two negative roots $x_{1,2}$ such that $x_2 \le \rho_2$ and $x_1 > -d/c$, or there are two negative roots $x_{3,4}$ such that $x_4 \le \rho_2$ and $x_3 > -d/c$, together with two positive roots $x_{1,2} < \rho_1$.
\vskip.1cm
{\bf (iii)} If $0 \le - d < c \rho_1$, then there is one negative root smaller than $\rho_2$ together with either one or three positive roots greater than or equal to $-d/c$ and smaller than $\rho_1$.
\vskip.1cm
{\bf (iv)} If $c \rho_1 \le -d$, then there is a negative root $x_2 < \rho_2$ and a positive root $x_1$ such that $\rho_1 \le x_1 < -d/c$. \\
\end{minipage}}  \\
\hline
\end{tabular}
\end{center}



\begin{center}
\begin{tabular}{|c||c|}
\hline
{\footnotesize \bf Figure 1.13} &  {\footnotesize \bf Figure 1.14} \\
\hline
\hline
{\color{violet} $\bm{b < 0, \quad a > 0, \quad c > 0,$}} &  {\color{violet}$\bm{b < 0, \quad a <  0, \quad c > 0,}$} \\
{\color{violet} $\bm{c_2 < \gamma_2 < c_0 < 0 < \gamma_1 < {\color{red}c} < c_1$}} & {\color{violet} $\bm{c_2 < \gamma_2 < c_0 < 0 < \gamma_1 < c_1 <{\color{red}c}$}} \\
\hline
\hline
\includegraphics[width=80mm]{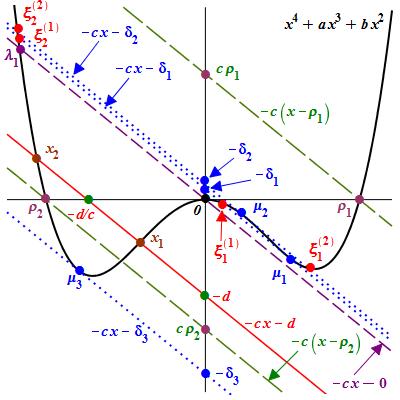} & \includegraphics[width=80mm]{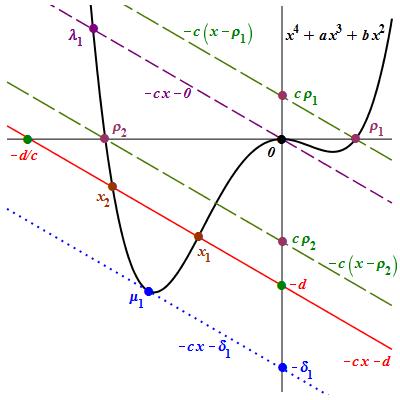} \\
\hline
\multicolumn{1}{|c||}{\begin{minipage}{22em}
\scriptsize
\vskip.3cm
\hskip3.89cm \underline{{\color{red}{\it Notes}}}\\
\vskip-.2cm
As $c_2 < c < c_1$, the quartic has three stationary points $\mu_i$ and the number of real roots can be either 0, or 2, or 4. There are three tangents $- c x - \delta_i$ to $x^4 + ax^3 + bx^2$ --- each at $\mu_i$.
\vskip.1cm
As $\gamma_2 < \gamma_1 < c$, the straight line $-cx$ intersects $x^4 + ax^3 + bx^2$ at one point only ($\lambda_1$). Thus $-\delta_1 > 0$. Also: $-\delta_3 < 0$ and $-\delta_2 > 0$.
\vskip.1cm
Obviously, $ -\delta_3 < c \rho_2 < 0 < -\delta_1 < -\delta_2 < c \rho_2$.
\vskip.1cm
The minimum of $x^4 + ax^3 + bx^2$ at negative $x$ is smaller than the minimum at positive $x$ by the amount of $(1/256)a(9a^2-32b)^{3/2}$.
\vskip.1cm
Consideration of whether $- d < c \rho_2$ or $-d > c \rho_2$ would yield sharper bounds in the analysis based on solving cubic equations.\\
\hrule
\vskip.3cm
\hskip1.54cm \underline{{\color{red}{\it Analysis based on solving cubic equations}}} \\
\vskip-.1cm
{\bf (i)} If $-d < -\delta_3$, then there are no real roots.
\vskip.1cm
{\bf (ii)} If $-\delta_3 \le -d < 0$ (pictured), then there is a negative root $x_2$ such that $\lambda_1 < x_2 \le \mu_3$ and another negative root $x_1 \ge $ max$ \{ \mu_3, -d/c \}$.
\vskip.1cm
{\bf (iii)} If $0 \le -d < -\delta_1$, then there is a negative root $x_2$ such that  $\xi^{(1)}_2 < x_2 \le \lambda_1$ and a non-negative root $x_1$ such that $- d/c \le x_1 < \xi^{(1)}_1$.
\vskip.1cm
{\bf (iv)}  If $-\delta_1 \le -d < -\delta_2$, then there is a negative root $x_4$ such that $\xi^{(2)}_2 < x_4  \le \xi^{(1)}_2$  and three positive roots: $x_3$ such that max$\{\xi^{(1)}_1, -d/c \} \le x_3 < \mu_2$, \,\, $x_2$ such that $\mu_2 < x_2 \le \mu_1$ and $x_1$ such that $\mu_1 \le x_1 < \xi^{(2)}_1$.
\vskip.1cm
{\bf (v)}  If $-\delta_2 \le -d$, then there is a negative root $x_2 \le \xi^{(2)}_2$ and a positive root $x_1$ such that max$\{ \xi^{(2)}_1,$ min$\{\rho_1, -d/c \} \le x_1 < $ max$\{ \rho_1, -d/c \}$.\\
\hrule
\vskip.3cm
\hskip.88cm \underline{{\color{red}{\it Analysis based on solving quadratic equations only}}}
\vskip.3cm
{\bf (i)} If $- d < c \rho_2$, then there are either no real roots or there are two negative roots $x_{1,2} > \rho_2$.
\vskip.1cm
{\bf (ii)} If $c \rho_2 \le - d < 0$ (pictured), then there is a negative root $x_2 \le \rho_2$ and another negative root $x_1 > -d/c$.
\vskip.1cm
{\bf (iii)} If $0 \le - d < c \rho_1$, then there is one negative root smaller than $\rho_2$ and either one or three positive roots greater than or equal to $-d/c$ and smaller than  $\rho_1$.
\vskip.1cm
{\bf (iv)} If $c \rho_1 \le -d$, then there is one negative root $x_2 < \rho_2$ and a positive root $x_1$ such that $\rho_1 \le x_1 < -d/c$. \\
\end{minipage}}
& \multicolumn{1}{|c|}{\begin{minipage}{22em}
\scriptsize
\vskip-2.8cm
\hskip3.89cm \underline{{\color{red}{\it Notes}}}
\vskip.1cm
As $c_2 < c_1 < c$, the quartic has a single stationary point $\mu_1$ and the number of real roots can be either 0 or 2. There is only one tangent $- c x - \delta_1$ to $x^4 + ax^3 + bx^2$ --- the one at $\mu_1$.
\vskip.1cm
As $\gamma_2 < \gamma_1 < c$, the straight line $-cx$ intersects $x^4 + ax^3 + bx^2$ at one point only ($\lambda_1$).
\vskip.1cm
One could have $\mu_1 < \rho_2$ (pictured) or $\mu_1 \ge \rho_2$.
\vskip.1cm
Obviously, $-\delta_1 \le c \rho_2 < 0 < c \rho_1$.
\vskip.1cm
The minimum of $x^4 + ax^3 + bx^2$ at negative $x$ is smaller than the minimum at positive $x$ by the amount of $(1/256)a(9a^2-32b)^{3/2}$.
\vskip.1cm
Consideration of whether $- d < c \rho_2$ or $-d > c \rho_2$ would yield sharper bounds in the analysis based on solving cubic equations. \\
\hrule
\vskip.3cm
\hskip1.54cm \underline{{\color{red}{\it Analysis based on solving cubic equations}}} \\
\vskip-.1cm
{\bf (i)} If $-d < -\delta_1$, then there are no real roots.
\vskip.1cm
{\bf (ii)} If $-\delta_1 \le -d < 0$ (pictured), then there is a negative root $x_2$ such that $\lambda_1 < x_2 \le \mu_1$ and another negative root $x_1 \ge $ max$ \{ \mu_1, -d/c \}$.
\vskip.1cm
{\bf (iii)} If $0 \le -d$, then there is a negative root $x_2 \le \lambda_1$ and a positive root $x_1$ such that min$\{ \rho_1, -d/c \} \le x_1 \le $ max$\{ \rho_1, -d/c \}$. \\
\hrule
\vskip.3cm
\hskip.88cm \underline{{\color{red}{\it Analysis based on solving quadratic equations only}}}
\vskip.3cm
{\bf (i)} If $- d < c \rho_2$ (pictured), then there are either no real roots or there are two negative roots (these are smaller than $\rho_2$ if $\mu_1 < \rho_2$ or greater than or equal to $\rho_2$ if $\mu_1 \ge \rho_2$).
\vskip.1cm
{\bf (ii)} If $c \rho_2 \le - d < 0$, then there is a negative root $x_2 \le \rho_2$ and another negative root $x_1 > -d/c$.
\vskip.1cm
{\bf (iii)} If $0 \le - d < c \rho_1$, then there is one negative root $x_2 < \rho_2$ and one positive root $x_1$ such that $-d/c \le x_1 < \rho_1$.
\vskip.1cm
{\bf (iv)} If $c \rho_1 \le -d$, then there is one negative root $x_2 < \rho_2$ and one positive root $x_1$ such that $\rho_1 \le x_1 < -d/c$. \\
\end{minipage}}  \\
\hline
\end{tabular}
\end{center}



\begin{center}
\begin{tabular}{|c||c|}
\hline
{\footnotesize \bf Figure 2.1} &  {\footnotesize \bf Figure 2.2} \\
\hline
\hline
{\color{violet} $\bm{0 < b \le \frac{a^2}{4}\, , \quad a < 0, \quad c < 0,$}} &  {\color{violet}$\bm{0 < b \le \frac{a^2}{4}\, , \quad a < 0, \quad c < 0,$}} \\
{\color{violet} $\bm{{\color{red}c} < c_2 < \gamma_2 < 0 < c_0 < \gamma_1 <  c_1$}} & {\color{violet} $\bm{c_2 < {\color{red}c} < \gamma_2 < 0 < c_0 < \gamma_1 < c_1 $}} \\
\hline
\hline
\includegraphics[width=80mm]{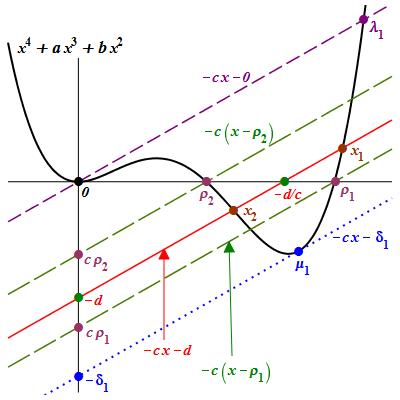} & \includegraphics[width=80mm]{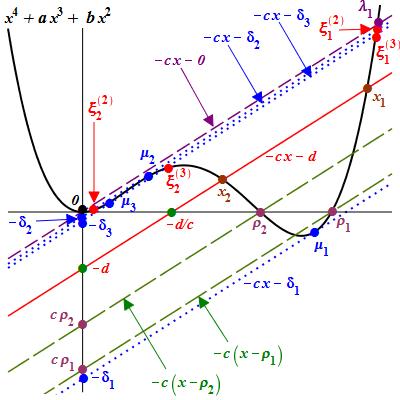} \\
\hline
\multicolumn{1}{|c||}{\begin{minipage}{22em}
\scriptsize
\vskip-2.7cm
\hskip3.89cm \underline{{\color{red}{\it Notes}}}
\vskip.3cm
As $c < c_2 < c_1$, the quartic has a single stationary point $\mu_1$ and the number of real roots can be either 0 or 2. There is only one tangent $- c x - \delta_1$ to $x^4 + ax^3 + bx^2$ --- the one at $\mu_1$.
\vskip.1cm
As $c < \gamma_2 < \gamma_1$, the straight line $-cx$ intersects $x^4 + ax^3 + bx^2$ at one point only ($\lambda_1$).
\vskip.1cm
One could have $\mu_1 < \rho_1$ (pictured) or $\mu_1 \ge \rho_1$.
\vskip.1cm
Obviously, $-\delta_1 \le c \rho_1 < c \rho_2 < 0$.
\vskip.1cm
Consideration of whether $- d < c \rho_1$ or $-d > c \rho_1$ would yield sharper bounds in the analysis based on solving cubic equations. \\
\hrule
\vskip.3cm
\hskip1.54cm \underline{{\color{red}{\it Analysis based on solving cubic equations}}}
\vskip.3cm
{\bf (i)} If $-d < -\delta_1$, then there are no real roots.
\vskip.1cm
{\bf (ii)} If $-\delta_1 \le -d < 0$ (pictured), then there is one positive root $x_2$ such that min$\{\rho_2, - d/c\} \le x_2  \le $
min$\{$max$\{\rho_2, -d/c \}, \mu_1 \}$ and another positive root $x_1$ such that $\mu_1 \le x_1 < \lambda_1$.
\vskip.1cm
{\bf (iii)}  If $0 \le - d$, then there is a non-positive root $x_2 \ge -d/c$ and a positive root $x_1 \ge \lambda_1$. \\
\hrule
\vskip.3cm
\hskip.88cm \underline{{\color{red}{\it Analysis based on solving quadratic equations only}}}
\vskip.3cm
{\bf (i)} If $-d < c \rho_1$, then there are either no real roots or there are two positive roots (these are greater than $\rho_1$ if $\mu_1 > \rho_1$ or smaller than $\rho_1$ if $\mu_1 \le \rho_1$).
\vskip.1cm
{\bf (ii)} If $c \rho_1 \le - d < c \rho_2$ (pictured), then there is a positive root $x_2$ such that $\rho_2 < x_2 \le -d/c$ and another positive root $x_1 \ge \rho_1$.
\vskip.1cm
{\bf (iii)} If $c \rho_2 \le - d < 0$, then there is a positive root $x_2$ such that $-d/c < x_2 \le \rho_2$ and another positive root $x_1 > \rho_1$.
\vskip.1cm
{\bf (iv)} If $0 \le -d$, then there is one non-positive root $x_2 \ge -d/c$ and a positive root $x_1 > \rho_1$.\\
\end{minipage}}
& \multicolumn{1}{|c|}{\begin{minipage}{22em}
\scriptsize
\scriptsize
\vskip.3cm
\hskip3.89cm \underline{{\color{red}{\it Notes}}}
\vskip.3cm
As $c_2 < c < c_1$, the quartic has three stationary points $\mu_i$ and the number of real roots can be either 0, or 2, or 4. There are three tangents $- c x - \delta_i$ to $x^4 + ax^3 + bx^2$ --- each at $\mu_i$.
\vskip.1cm
As $c < \gamma_2 < \gamma_1$, the straight line $-cx$ intersects $x^4 + ax^3 + bx^2$ at one point only ($\lambda_1$). Thus $-\delta_2 < 0$. Also: $-\delta_3 < 0$ and $-\delta_1 < 0$.
\vskip.1cm
Obviously, $-\delta_1 < c \rho_1 <     c \rho_2 < -\delta_3 < -\delta_2 < 0$.
\vskip.1cm
Consideration of whether $- d < c \rho_1$ or $-d > c \rho_1$ would yield sharper bounds in the analysis based on solving cubic equations. \\
\hrule
\vskip.3cm
\hskip1.54cm \underline{{\color{red}{\it Analysis based on solving cubic equations}}}
\vskip.3cm
{\bf (i)} If $-d < -\delta_1$, then there are no real roots.
\vskip.1cm
{\bf (ii)} If $-\delta_1 \le -d < -\delta_3$ (pictured), then there is one positive root $x_2$ such that max$ \{ \xi^{(3)}_2, $ min$\{\rho_2, -d/c \}\} < x_2 \le $ min$\{$max$\{\rho_2, -d/c \}, \mu_1 \}$ and another positive root $x_1$ such that $\mu_1 \le x_2 < \xi^{(3)}_1$.
\vskip.1cm
{\bf (iii)}  If $-\delta_3 \le -d < -\delta_2$, then there are four positive roots: $x_4$ such that max$\{\xi^{2}_2, -d/c\} < x_4 \le \mu_3$, $\,\, x_3$ such that $\mu_3 \le x_3 < \mu_2, \,\,$ $x_2$ such that $\mu_2 < x_2 \le \xi^{(3)}_2,$ and $x_1$ such that $\xi^{(3)}_1 \le x_1 < \xi^{(2)}_1$.
\vskip.1cm
{\bf (iv)} If $-\delta_2 \le -d < 0$, then there is a positive root $x_2$ such that $-d/c < x_2 \le \xi^{(2)}_2$ and another positive root $x_1$ such that $\xi^{(2)}_1 \le x_1 < \lambda_1$.
\vskip.1cm
{\bf (v)} If $0 \le -d$, then there is a non-positive root $x_2 \ge -d/c$ and a positive root $x_1 \ge \lambda_1$. \\
\hrule
\vskip.3cm
\hskip.88cm \underline{{\color{red}{\it Analysis based on solving quadratic equations only}}}
\vskip.3cm
{\bf (i)} If $-d < c \rho_1$, then there are either no real roots or there are two positive roots $x_{1,2} < \rho_1$.
\vskip.1cm
{\bf (ii)} If $c \rho_1 \le - d < c \rho_2$, then there is a positive root $x_2$ such that $\rho_2 < x_2 \le -d/c$ and another positive root $x_1 \ge \rho_1$.
\vskip.1cm
{\bf (iii)} If $c \rho_2 \le - d < 0$ (pictured), then there is a positive root $x_1 > \rho_1$ together with either one or three positive roots greater than $-d/c$ and smaller than or equal to $\rho_2$.
\vskip.1cm
{\bf (iv)} If $0 \le -d$, then there is one non-positive root $x_2 \ge -d/c$ and a positive root $x_1 > \rho_1$.\\
\end{minipage}}
\\
\hline
\end{tabular}
\end{center}



\begin{center}
\begin{table}[H]
\begin{tabular}{|c||c|}
\hline
{\footnotesize \bf Figure 2.3} &  {\footnotesize \bf Figure 2.4} \\
\hline
\hline
{\color{violet} $\bm{0 < b \le \frac{a^2}{4}\, , \quad a < 0, \quad c < 0,$}} &  {\color{violet} $\bm{0 < b \le \frac{a^2}{4}\, , \quad a < 0, \quad c > 0,$}} \\
{\color{violet} $\bm{c_2 < \gamma_2 < {\color{red}c} < 0 < c_0 < \gamma_1 < c_1$}} & {\color{violet} $\bm{c_2 < \gamma_2 < 0 <{\color{red}c} < c_0 < \gamma_1 < c_1 $}} \\
\hline
\hline
\includegraphics[width=80mm]{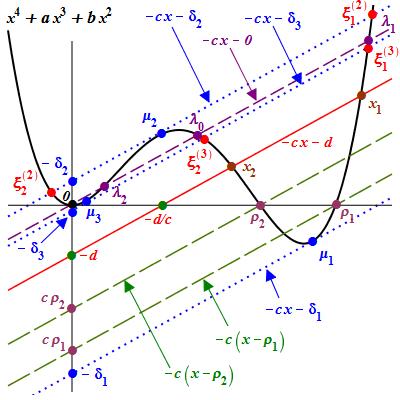} & \includegraphics[width=80mm]{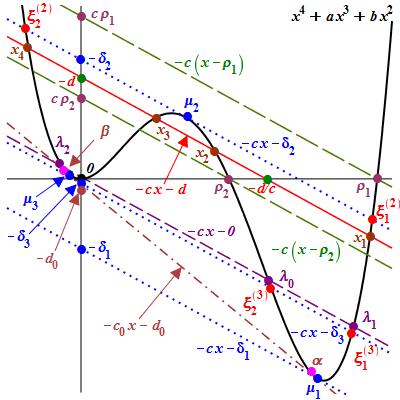} \\
\hline
\multicolumn{1}{|c||}{\begin{minipage}{22em}
\scriptsize
\vskip-3.6cm
\hskip3.89cm \underline{{\color{red}{\it Notes}}}
\vskip.01cm
As $c_2 < c < c_1$, the quartic has three stationary points $\mu_i$ and the number of real roots can be either 0, or 2, or 4. There are three tangents $- c x - \delta_i$ to $x^4 + ax^3 + bx^2$ --- each at $\mu_i$.
\vskip.1cm
As $\gamma_2 < c < \gamma_1$, the straight line $-cx$ intersects $x^4 + ax^3 + bx^2$ at three points: $\lambda_{0,1,2}$.  Thus $-\delta_2 > 0$. Also: $-\delta_3 < 0$ and $-\delta_1 < 0$.
\vskip.1cm
Obviously, $-\delta_1 < c \rho_1 < c \rho_2 < -\delta_3 < 0  < -\delta_2 $.
\vskip.1cm
Consideration of whether $- d < c \rho_1$ or $-d > c \rho_1$ would yield sharper bounds in the analysis based on solving cubic equations.
\\
\hrule
\vskip.1cm
\hskip1.54cm \underline{{\color{red}{\it Analysis based on solving cubic equations}}}
\vskip.1cm
{\bf (i)} If $-d < -\delta_1$, then there are no real roots.
\vskip.1cm
{\bf (ii)} If $-\delta_1 \le -d < -\delta_3$ (pictured), then there is one positive root $x_2$ such that max$ \{ \xi^{(3)}_2, $ min$\{\rho_2, -d/c \}\} < x_2 \le $ min$\{$max$\{\rho_2, -d/c \}, \mu_1 \}$ and another positive root $x_1$ such that $\mu_1 \le x_2 < \xi^{(3)}_1$.
\vskip.1cm
\vskip.1cm
{\bf (iii)} If $-\delta_3 \le -d < 0$, then there are four positive roots: $x_4$ such that $-d/c < x_4 \le \mu_3, \, \, $
$x_3$ such that $\mu_3 \le x_3 < \lambda_2, \,\,$  $x_2$ such that $\lambda_0 < x_2 \le \xi^{(3)}_2,$ and $x_1$ such that $\xi^{(3)}_1 \le x_1 < \lambda_1$.
\vskip.1cm
{\bf (iv)}  If $0 \le -d < -\delta_2$, then there is a non-positive root $x_4 \ge $ max$\{\xi^{(2)}_2, - d/c \}$ and three positive roots: $x_3$ such that $\lambda_2 \le x_3 < \mu_2, \,\,$ $x_2$ such that $\mu_2 < x_2 \le \lambda_0$, and $x_1$ such that $\lambda_1 \le x_1 <  \xi^{(2)}_1$.
\vskip.1cm
{\bf (v)} If $-\delta_2 \le -d$, then there is a negative root $x_2$ such that $-d/c < x_2 \le \xi^{(2)}_2 $  and a positive root $x_1 \ge \xi^{(2)}_1$.
\\
\hrule
\vskip.1cm
\hskip.88cm \underline{{\color{red}{\it Analysis based on solving quadratic equations only}}}
\vskip.1cm
{\bf (i)} If $-d < c \rho_1$, then there are either no real roots or there are two positive roots $x_{1,2} < \rho_1$.
\vskip.1cm
{\bf (ii)} If $c \rho_1 \le - d < c \rho_2$, then there is a positive root $x_2$ such that $\rho_2 < x_2 \le -d/c$ and another positive root $x_1 \ge \rho_1$.
\vskip.1cm
{\bf (iii)} If $c \rho_2 \le -d < 0$ (pictured), then there are either three negative roots $x_{2,3,4}$ such that $-d/c < x_{2,3,4} \le \rho_2$ together with a positive root $x_1 > \rho_1$, or there is one negative root $x_2$ such that $-d/c < x_2 \le \rho_2$ and a positive root $x_1 > \rho_1$.
\vskip.1cm
{\bf (iv)} If $0 \le -d$, then there is one non-positive root greater than or equal to $-d/c$, a positive root greater than $\rho_1$ and either zero or two positive roots smaller than $\rho_2$. \\
\end{minipage}}
& \multicolumn{1}{|c|}{\begin{minipage}{22em}
\scriptsize
\vskip.1cm
\hskip3.89cm \underline{{\color{red}{\it Notes}}}\\
As $c_2 < c < c_1$, the quartic has three stationary points $\mu_i$ and the number of real roots can be either 0, or 2, or 4. There are three tangents $- c x - \delta_i$ to $x^4 + ax^3 + bx^2$ --- each at $\mu_i$.
\vskip.1cm
As $\gamma_2 < c < \gamma_1$, the straight line $-cx$ intersects $x^4 + ax^3 + bx^2$ at three points: $\lambda_{0,1,2}$. Thus $-\delta_1 < 0$. Also: $-\delta_2 > 0$ and $-\delta_3 < 0$.
\vskip.1cm
The straight line $ - c_0 x - d_0$ with $c_0 = (1/2)a(b-a^2/4) > 0$ and $d_0 = (1/4)(b-a^2/4)^2 > 0$ is tangent to $x^4 + ax^3 + bx^2$ at two points: $\alpha$ and $\beta = -\alpha - a/2$, where $\alpha$ and $\beta$ are the simultaneous double roots of the varied quartic $x^4 + ax^3 + bx^2 + c_0 x + d_0$.
\vskip.1cm
As $0 < c < c_0$, one has $- \delta_1 < -\delta_3$.
\vskip.1cm
One could have $-\delta_1 < - \delta_3 < 0 < c \rho_2 < -\delta_2 < c \rho_1$ (pictured) or $-\delta_1 < -\delta_3 < 0 < c \rho_2 < c \rho_1 \le -\delta_2 $.
\vskip.1cm
Consideration of whether $- d < c \rho_2$ or $-d > c \rho_2$ would yield sharper bounds in the analysis based on solving cubic equations.
\hrule
\vskip.1cm
\hskip1.54cm \underline{{\color{red}{\it Analysis based on solving cubic equations}}} \\
\vskip-.2cm
{\bf (i)} If $-d < -\delta_1$, then there are no real roots.
\vskip.1cm
{\bf (ii)} If $-\delta_1 \le -d < -\delta_3$, then there is one positive root $x_2$ such that $\xi^{(3)}_2 < x_2 \le \mu_1$ and another positive root $x_1$ such that $\mu_1 \le x_2 < \xi^{(3)}_1$.
\vskip.1cm
{\bf (iii)} If $-\delta_3 \le -d < 0$, then there is a negative root $x_4$ such that $\lambda_2 < x_4 \le \mu_3$, another negative root $x_3$ such that $\mu_3 \le x_3 < -d/c$, a positive root $x_2$ such that $\lambda_0 < x_2 \le \xi^{(3)}_2$ and another positive root $x_1$ such that $\xi^{(3)}_1 \le x_1 < \lambda_1$.
\vskip.1cm
{\bf (iv)}  If $0 \le -d < -\delta_2$ (pictured), then there is a negative root $x_4$ such that $\xi^{(2)}_2 < x_4 \le \lambda_2$, a non-negative root $x_3 \le $ min$\{-d/c,\mu_2\}$, a positive root $x_2$ such that $\mu_2 < x_2 \le \lambda_0$ and another positive $x_1$ such that $\lambda_1 \le x_1 <  \xi^{(2)}_1$.
\vskip.1cm
{\bf (v)}  If $-\delta_2 \le - d $, then there is a negative root $x_2 \le \xi^{(2)}_2$ and a positive root $x_1$ such that max$\{ \xi^{(2)}_1,$ min$\{\rho_1, -d/c \} \le x_1 < $ max$\{ \rho_1, -d/c \}$.\\
\hrule
\vskip.1cm
\hskip.88cm \underline{{\color{red}{\it Analysis based on solving quadratic equations only}}} \\
\vskip-.2cm
{\bf (i)} If $- d < 0$, then there are either no real roots, or there are two positive roots smaller than $\rho_1$, or there are two positive roots smaller than $\rho_1$ and two negative roots smaller than or equal to $-d/c$.
\vskip.1cm
{\bf (ii)} If $0 \le - d < c \rho_2$, then there is one negative root $x_4$, a non-negative root $x_3 \le -d/c$ and two positive roots $x_{1,2}$ such that $\rho_2 < x_{1,2} < \rho_1$.
\vskip.1cm
{\bf (iii)} If $c \rho_2 \le - d < c \rho_1$ (pictured), then there is one negative root, either zero or two positive roots smaller than or equal to $\rho_2$, and one positive root greater than $-d/c$ and smaller than $\rho_1$.
\vskip.1cm
{\bf (iv)} If $c \rho_1 \le -d$, then there is one negative root $x_2$ and a positive root $x_1$ such that $\rho_1 \le x_1 < -d/c$ and either zero or two positive roots smaller than or equal to $\rho_2$ (the latter appear when $c \rho_1 \le -\delta_2$). \\
\end{minipage}}  \\
\hline
\end{tabular}
\end{table}
\end{center}



\begin{center}
\begin{tabular}{|c||c|}
\hline
{\footnotesize \bf Figure 2.5} &  {\footnotesize \bf Figure 2.6} \\
\hline
\hline
{\color{violet} $\bm{0 < b \le \frac{a^2}{4}\, , \quad a < 0, \quad c > 0,$}} & {\color{violet} $\bm{0 < b \le \frac{a^2}{4}\, , \quad a < 0, \quad c > 0,$}} \\
{\color{violet} $\bm{c_2 < \gamma_2 < 0 < c_0 < {\color{red}c} < \gamma_1 < c_1$}} & {\color{violet} $\bm{c_2 < \gamma_2 < 0 < c_0 < \gamma_1 < {\color{red}c} < c_1 $}} \\
\hline
\hline
\includegraphics[width=80mm]{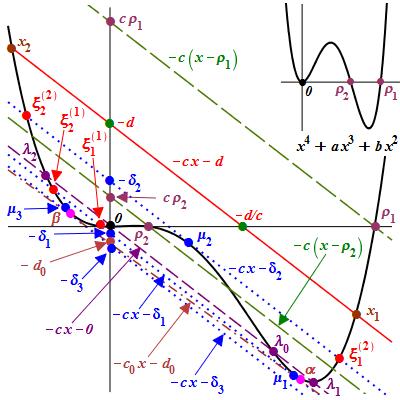} & \includegraphics[width=80mm]{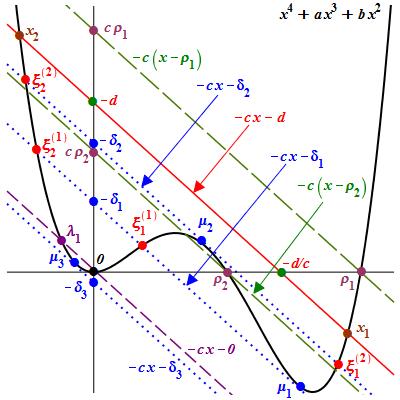} \\
\hline
\multicolumn{1}{|c||}{\begin{minipage}{22em}
\scriptsize
\vskip.05cm
\hskip3.89cm \underline{{\color{red}{\it Notes}}}\\
\vskip-.35cm
As $c_2 < c < c_1$, the quartic has three stationary points $\mu_i$ and the number of real roots can be either 0, or 2, or 4. There are three tangents $- c x - \delta_i$ to $x^4 + ax^3 + bx^2$ --- each at $\mu_i$.
\vskip.1cm
As $\gamma_2 < c < \gamma_1$, the straight line $-cx$ intersects $x^4 + ax^3 + bx^2$ at three points: $\lambda_{0,1,2}$. Thus $-\delta_1 < 0$. Also: $-\delta_2 > 0$ and $-\delta_3 < 0$.
\vskip.1cm
The straight line $ - c_0 x - d_0$ with $c_0 = (1/2)a(b-a^2/4) > 0$ and $d_0 = (1/4)(b-a^2/4)^2 > 0$ is tangent to $x^4 + ax^3 + bx^2$ at two points: $\alpha$ and $\beta = -\alpha - a/2$, where $\alpha$ and $\beta$ are the simultaneous double roots of the varied quartic $x^4 + ax^3 + bx^2 + c_0 x + d_0$.
\vskip.1cm
As $c > c_0$, one has $- \delta_3 < -\delta_1$.
\vskip.1cm
Obviously, $-\delta_3  < -\delta_1 < 0 < c \rho_2 < -\delta_2 < c \rho_1$.
\vskip.1cm
The graph of $x^4 + ax^3 + bx^2$ is shown on its own in the top-right corner.
\vskip.1cm
Point $\mu_2$ could be in the first quadrant or in the fourth (pictured).
\vskip.1cm
Consideration of whether $- d < c \rho_2$ or $-d > c \rho_2$ would yield sharper bounds in the analysis based on solving cubic equations.
\hrule
\vskip.05cm
\hskip1.54cm \underline{{\color{red}{\it Analysis based on solving cubic equations}}} \\
\vskip-.3cm
{\bf (i)} If $-d < -\delta_3$, then there are no real roots.
\vskip.1cm
{\bf (ii)} If $-\delta_3 \le -d < -\delta_1$, then there is a negative root $x_2$ such that $\xi^{(1)}_2 < x_2 \le \mu_3$ and another negative root $x_1$ such that $\mu_3 \le x_1 < $ min$\{-d/c,  \xi^{(1)}_1\}$.
\vskip.1cm
{\bf (iii)} If $-\delta_1 \le -d < 0$, then there is a negative root $x_4$ such that $\lambda_2 < x_4 \le \xi^{(1)}_2$, another negative root $x_3$ such that $\xi^{(1)}_1 \le x_3 < -d/c,$ a positive root $x_2$ such that $\lambda_0 < x_2 \le \mu_1$ and another positive root $x_1$ such that  $\mu_1 \le x_1 < \lambda_1$.
\vskip.1cm
{\bf (iv)}  If $0 \le -d < -\delta_2$, then there is a negative root $x_4$ such that $\xi^{(2)}_2 < x_4 \le \lambda_2$, a non-negative root $x_3 \le -d/c$, a positive root $x_2$ such that $\mu_2 < x_2 \le \lambda_0$, and another positive root $x_1$ such that $\lambda_1 \le x_1 <  \xi^{(2)}_1$.
\vskip.1cm
{\bf (v)}  If $-\delta_2 \le - d$ (pictured), then there is a negative root $x_2 \le \xi^{(2)}_2$ and a positive root $x_1$ such that max$\{ \xi^{(2)}_1,$ min$\{\rho_1, -d/c \} \le x_1 < $ max$\{ \rho_1, -d/c \}$.
\hrule
\vskip.05cm
\hskip.88cm \underline{{\color{red}{\it Analysis based on solving quadratic equations only}}} \\
\vskip-.3cm
{\bf (i)} If $- d < 0$, then there are either no real roots, or there are two negative roots smaller than $-d/c$, or there are two negative roots smaller than $-d/c$ and two positive roots greater than $\rho_2$ and smaller than $\rho_1$.
\vskip.1cm
{\bf (ii)} If $0 \le - d < c \rho_2$, then there is one negative root $x_4$, a non-negative root $x_3 \le -d/c$ and two positive roots greater than $\rho_2$ and smaller than $\rho_1$.
\vskip.1cm
{\bf (iii)} If $c \rho_2 \le - d < c \rho_1$ (pictured), then there is one negative root and either one or three positive roots greater than or equal to $-d/c$ and smaller than $\rho_1$.
\vskip.1cm
{\bf (iv)} If $c \rho_1 \le -d$, then there is one negative root $x_2$ and a positive root $x_1$ such that $\rho_1 \le x_1 < -d/c$.
\end{minipage}}
& \multicolumn{1}{|c|}{\begin{minipage}{22em}
\scriptsize
\vskip-2.cm
\hskip3.89cm \underline{{\color{red}{\it Notes}}}\\
\vskip-.2cm
As $c_2 < c < c_1$, the quartic has three stationary points $\mu_i$ and the number of real roots can be either 0, or 2, or 4. There are three tangents $- c x - \delta_i$ to $x^4 + ax^3 + bx^2$ --- each at $\mu_i$.
\vskip.1cm
As $\gamma_2 < \gamma_1 < c$, the straight line $-cx$ intersects $x^4 + ax^3 + bx^2$ at one point only ($\lambda_1$). Thus $-\delta_1 > 0$. Also: $-\delta_2 > 0$ and $-\delta_3 < 0$.
\vskip.1cm
Obviously, $-\delta_3  < 0 < -\delta_1 < c \rho_2 < -\delta_2 < c \rho_1$.
\vskip.1cm
Point $\mu_2$ could be in the first quadrant (pictured) or in the fourth.
\vskip.1cm
Consideration of whether $- d < c \rho_2$ or $-d > c \rho_2$ would yield sharper bounds in the analysis based on solving cubic equations.
\vskip.1cm
\hrule
\vskip.3cm
\hskip1.54cm \underline{{\color{red}{\it Analysis based on solving cubic equations}}} \\
\vskip-.1cm
{\bf (i)} If $-d < -\delta_3$, then there are no real roots.
\vskip.1cm
{\bf (ii)} If $-\delta_3 \le -d < 0$, then there is a negative root $x_2$ such that $\lambda_1 < x_2 \le \mu_3$ and another negative root $x_1$ such that $\mu_3 \le x_1 < -d/c$.
\vskip.1cm
{\bf (iii)} If $0 \le -d < -\delta_1$, then there is a negative root $x_2$ such that  $\xi^{(1)}_2 < x_2 \le \lambda_1$ and a non-negative root $x_1 < $ min$\{-d/c, \xi^{(1)}_1 \}$.
\vskip.1cm
{\bf (iv)}  If $-\delta_1 \le - d < -\delta_2$, then there is a negative root $x_4$ such that $\xi^{(2)}_2 < x_4  \le \xi^{(1)}_2$  and three positive roots: $x_3$ such that $\xi^{(1)}_1 \le x_3 < \mu_2$, \,\, $x_2$ such that $\mu_2 < x_2 \le \mu_1$ and $x_1$ such that $\mu_1 \le x_1 < \xi^{(2)}_1$.
\vskip.1cm
{\bf (v)}  If $-\delta_2 \le - d$ (pictured), then there is a negative root $x_2 \le \xi^{(2)}_2$ and a positive root $x_1$ such that max$\{ \xi^{(2)}_1,$ min$\{\rho_1, -d/c \} \le x_1 < $ max$\{ \rho_1, -d/c \}$.
\vskip.1cm
\hrule
\vskip.3cm
\hskip.88cm \underline{{\color{red}{\it Analysis based on solving quadratic equations only}}}
\vskip.3cm
{\bf (i)} If $- d < 0$, then there are either no real roots or there are two negative roots smaller than $-d/c$.
\vskip.1cm
{\bf (ii)} If $0 \le - d < c \rho_2$, then there is one negative root, a non-negative root smaller than or equal to $-d/c$ and either zero or two positive roots greater than $\rho_2$ and smaller than $\rho_1$.
\vskip.1cm
{\bf (iii)} If $c \rho_2 \le - d < c \rho_1$ (pictured), then there is one negative root, either zero or two positive roots smaller than or equal to $\rho_2$, and one positive root greater than or equal to $-d/c$ and smaller than $\rho_1$.
\vskip.1cm
{\bf (iv)} If $c \rho_1 \le -d$, then there is one negative root $x_2$ and a positive root $x_1$ such that $\rho_1 \le x_1 < -d/c$. \\
\end{minipage}}  \\
\hline
\end{tabular}
\end{center}



\begin{center}
\begin{tabular}{|c||c|}
\hline
{\footnotesize \bf Figure 2.7} &  {\footnotesize \bf Figure 2.8} \\
\hline
\hline
{\color{violet} $\bm{0 < b \le \frac{a^2}{4}\, , \quad a < 0, \quad c > 0,$}} &  {\color{violet} $\bm{0 < b \le \frac{a^2}{4}\, , \quad a > 0, \quad c < 0,$}} \\
{\color{violet} $\bm{c_2 < \gamma_2 < 0 < c_0 < \gamma_1 < c_1 < {\color{red}c}$}} & {\color{violet} $\bm{{\color{red}c} < c_2 < \gamma_2 < c_0 < 0 < \gamma_1 < c_1 $}} \\
\hline
\hline
\includegraphics[width=80mm]{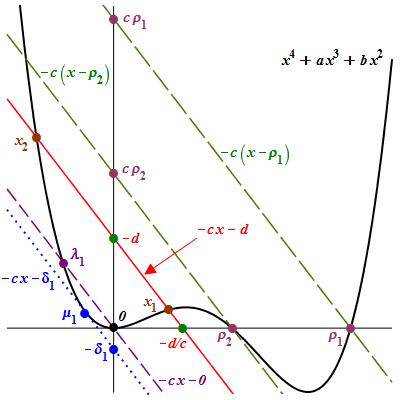} & \includegraphics[width=80mm]{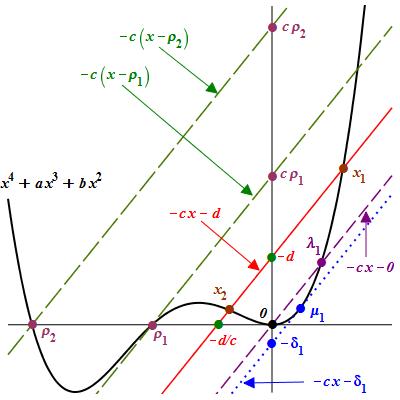} \\
\hline
\multicolumn{1}{|c||}{\begin{minipage}{22em}
\scriptsize
\vskip0.3cm
\hskip3.89cm \underline{{\color{red}{\it Notes}}}\\
\vskip.05cm
As $c_2 < c_1 < c$, the quartic has a single stationary point $\mu_1$ and the number of real roots can be either 0 or 2. There is only one tangent $- c x - \delta_1$ to $x^4 + ax^3 + bx^2$ --- the one at $\mu_1$.
\vskip.1cm
As $\gamma_2 < \gamma_1 < c$, the straight line $-cx$ intersects $x^4 + ax^3 + bx^2$ at one point only ($\lambda_1$).
\vskip.1cm
Obviously, $ -\delta_1 < 0 < c \rho_2 < c \rho_1$.
\\
\hrule
\vskip.3cm
\hskip1.54cm \underline{{\color{red}{\it Analysis based on solving cubic equations}}} \\
\vskip.1cm
{\bf (i)} If $-d < -\delta_1$, then there are no real roots.
\vskip.1cm
{\bf (ii)} If $-\delta_1 \le -d < 0$, then there is a negative root $x_2$ such that $\lambda_1 < x_2 \le \mu_1$ and another negative root $x_1$ such that $\mu_1 \le x_1 < -d/c$.
\vskip.1cm
{\bf (iii)} If $0 \le - d < c \rho_2$ (pictured), then there is a negative root $x_2 \le \lambda_1$ and a non-negative root $x_1 \le -d/c$.
\vskip.1cm
{\bf (iv)} If $c \rho_2 \le -d$, then there is a negative root $x_2 < \lambda_1$ and a positive root $x_1$ such that
min$\{\rho_1, -d/c \} \le x_1 \le $ max$\{\rho_1, -d/c \}$. \\
\hrule
\vskip.3cm
\hskip.88cm \underline{{\color{red}{\it Analysis based on solving quadratic equations only}}}
\vskip.3cm
{\bf (i)} If $- d < 0$, then there are either no real roots or there are two negative roots smaller than $-d/c$.
\vskip.1cm
{\bf (ii)} If $0 \le - d < c \rho_2$ (pictured), then there is a negative root $x_2$ and a non-negative root $x_1 \le -d/c$.
\vskip.1cm
{\bf (iii)} If $c \rho_2 \le - d < c \rho_1$, then there is one negative root $x_2$ and a positive root $x_1$ such that $-d/c \le x_1 < \rho_1$.
\vskip.1cm
{\bf (iv)} If $c \rho_1 \le -d$, then there is one negative root $x_2$ and a positive root $x_1$ such that $\rho_1 \le x_1 < -d/c$. \\
\end{minipage}}
& \multicolumn{1}{|c|}{\begin{minipage}{22em}
\scriptsize
\vskip0.3cm
\hskip3.89cm \underline{{\color{red}{\it Notes}}}\\
\vskip.05cm
As $c < c_2 < c_1$, the quartic has a single stationary point $\mu_1$ and the number of real roots can be either 0 or 2. There is only one tangent $- c x - \delta_1$ to $x^4 + ax^3 + bx^2$ --- the one at $\mu_1$.
\vskip.1cm
As $c < \gamma_2 < \gamma_1$, the straight line $-cx$ intersects $x^4 + ax^3 + bx^2$ at one point only ($\lambda_1$).
\vskip.1cm
Obviously, $ -\delta_1 < 0 < c \rho_1 < c \rho_2$.
\\
\hrule
\vskip.3cm
\hskip1.54cm \underline{{\color{red}{\it Analysis based on solving cubic equations}}} \\
\vskip.1cm
{\bf (i)} If $-d < -\delta_1$, then there are no real roots.
\vskip.1cm
{\bf (ii)} If $-\delta_1 \le -d < 0$, then there is a positive root $x_2$ such that $-d/c < x \le \mu_1$ and another positive root $x_1$ such that $\mu_1 \le x_2 < \lambda_1$.
\vskip.1cm
{\bf (iii)} If $0 \le - d < c \rho_1$ (pictured), then there is a non-positive root $x_2 \ge -d/c$ and a positive root $x_1 \ge \lambda_1$.
\vskip.1cm
{\bf (iv)} If $c \rho_1 \le -d$, then there is a positive root $x_1 > \lambda_1$ and a negative root $x_2$ such that
min$\{\rho_2, -d/c \} \le x_2 \le $ max$\{\rho_2, -d/c \}$. \\
\hrule
\vskip.3cm
\hskip.88cm \underline{{\color{red}{\it Analysis based on solving quadratic equations only}}}
\vskip.3cm
{\bf (i)} If $- d < 0$, then there are either no real roots or there are two positive roots greater than $-d/c$.
\vskip.1cm
{\bf (ii)} If $0 \le - d < c \rho_1$ (pictured), then there is a positive root $x_1$ and a non-positive root $x_2 \ge -d/c$.
\vskip.1cm
{\bf (iii)} If $c \rho_1 \le - d < c \rho_2$, then there is one positive root $x_1$ and a negative root $x_2$ such that $\rho_2 < x_2 \le -d/c$.
\vskip.1cm
{\bf (iv)} If $c \rho_2 \le -d$, then there is one positive root $x_1$ and a negative root $x_2$ such that $-d/c <  x_1 \le \rho_2$. \\
\end{minipage}}  \\
\hline
\end{tabular}
\end{center}


\newpage

\begin{center}
\begin{tabular}{|c||c|}
\hline
{\footnotesize \bf Figure 2.9} &  {\footnotesize \bf Figure 2.10} \\
\hline
\hline
{\color{violet} $\bm{0 < b \le \frac{a^2}{4}\, , \quad a > 0, \quad c < 0,$}} &  {\color{violet} $\bm{0 < b \le \frac{a^2}{4}\, , \quad a > 0, \quad c < 0,$}} \\
{\color{violet} $\bm{c_2 < {\color{red}c} < \gamma_2 < c_0 < 0 < \gamma_1 < c_1 $}} & {\color{violet} $\bm{c_2 < \gamma_2 < {\color{red}c} < c_0 < 0 < \gamma_1 < c_1 $}} \\
\hline
\hline
\includegraphics[width=80mm]{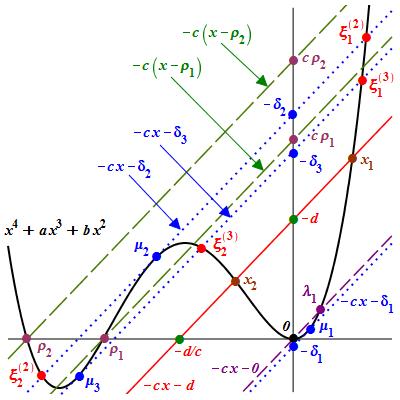} & \includegraphics[width=80mm]{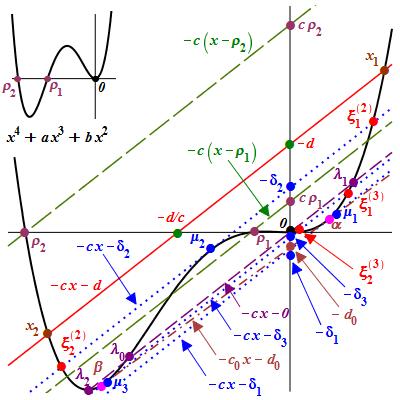} \\
\hline
\multicolumn{1}{|c||}{\begin{minipage}{22em}
\scriptsize
\vskip-2.3cm
\hskip3.89cm \underline{{\color{red}{\it Notes}}}\\
\vskip-.2cm
As $c_2 < c < c_1$, the quartic has three stationary points $\mu_i$ and the number of real roots can be either 0, or 2, or 4. There are three tangents $- c x - \delta_i$ to $x^4 + ax^3 + bx^2$ --- each at $\mu_i$.
\vskip.1cm
As $c < \gamma_2 < \gamma_1$, the straight line $-cx$ intersects $x^4 + ax^3 + bx^2$ at one point only ($\lambda_1$). Thus $-\delta_3 > 0$. Also: $-\delta_1 <  0$ and $-\delta_2 > 0$.
\vskip.1cm
Obviously, $-\delta_1  < 0 < -\delta_3 < c \rho_1 < -\delta_2 < c \rho_2$.
\vskip.1cm
Point $\mu_2$ could be in the second quadrant (pictured) or in the third.
\vskip.1cm
Consideration of whether $- d < c \rho_1$ or $-d > c \rho_1$ would yield sharper bounds in the analysis based on solving cubic equations.
\vskip.1cm
\hrule
\vskip.2cm
\hskip1.54cm \underline{{\color{red}{\it Analysis based on solving cubic equations}}} \\
\vskip-.1cm
{\bf (i)} If $-d < -\delta_1$, then there are no real roots.
\vskip.1cm
{\bf (ii)} If $-\delta_1 \le -d < 0$, then there is a positive root $x_2$ such that $ -d/c < x_2 \le \mu_1$ and another positive root $x_1$ such that $\mu_1 \le x_1 <  \lambda_1$.
\vskip.1cm
{\bf (iii)} If $0 \le -d < -\delta_3$ (pictured), then there is a negative root $x_2$ such that $x_2 \ge $ max$\{\xi^{(3)}_2, -d/c \}$, and a positive root $x_1$ such that $\lambda_1 \le x_1 < \xi^{(3)}_1$ .
\vskip.1cm
{\bf (iv)} If $-\delta_3 \le -d < - \delta_2$, then there is a negative root $x_4$ such that $\xi^{(2)}_2 < x_4 \le $ min$\{-d/c, \mu_3 \}$, a negative root $x_3$ such that max$ \{ - d/c, \mu_3 \} \le x_3 < \mu_2$, another negative root $x_2$ such that $\mu_2 < x_2 \le \xi^{(3)}_2$, and a positive root $x_1$ such that $\xi^{(3)}_1 \le x_1 < \xi^{(2)}_1$.
\vskip.1cm
{\bf (v)} If $-\delta_2 \le -d$, then there is a negative root $x_2$ such that min$\{\rho_2, -d/c \} \le x_2 \le $ min$\{$max$\{\rho_2, -d/c \}, \xi^{(2)}_2 \}$  and a positive root $x_1 \ge \xi^{(2)}_1$.
\vskip.1cm
\hrule
\vskip.3cm
\hskip.88cm \underline{{\color{red}{\it Analysis based on solving quadratic equations only}}}
\vskip.3cm
{\bf (i)} If $-d < 0$, then there are either no real roots or there are two positive roots greater that $-d/c$.
\vskip.1cm
{\bf (ii)} If $0 \le - d <  c \rho_1$ (pictured), then there is one positive root, one non-positive root greater than or equal to $-d/c$ and either zero or two negative roots greater than $\rho_2$ and smaller than or equal to $\rho_1$.
\vskip.1cm
{\bf (iii)} If $ c \rho_1 \le -d < c \rho_2$, then there is one positive root, one negative root greater than $\rho_2$ and smaller than or equal to $-d/c$, and either zero or two negative roots smaller than or equal to $\rho_1$.
\vskip.1cm
{\bf (iv)} If $c \rho_2 \le -d$, then there is one negative root $x_2$ such that $-d/c < x_2 \le \rho_2$ and one positive root $x_1$. \\
\end{minipage}}
& \multicolumn{1}{|c|}{\begin{minipage}{22em}
\scriptsize
\vskip.05cm
\hskip3.89cm \underline{{\color{red}{\it Notes}}}\\
\vskip-.35cm
As $c_2 < c < c_1$, the quartic has three stationary points $\mu_i$ and the number of real roots can be either 0, or 2, or 4. There are three tangents $- c x - \delta_i$ to $x^4 + ax^3 + bx^2$ --- each at $\mu_i$.
\vskip.1cm
As $\gamma_2 < c < \gamma_1$, the straight line $-cx$ intersects $x^4 + ax^3 + bx^2$ at three points: $\lambda_{0,1,2}$.
Thus $-\delta_3 < 0$. Also: $-\delta_1 <  0$ and $-\delta_2 > 0$.
\vskip.1cm
The straight line $ - c_0 x - d_0$ with $c_0 = (1/2)a(b-a^2/4) < 0$ and $d_0 = (1/4)(b-a^2/4)^2 > 0$ is tangent to $x^4 + ax^3 + bx^2$ at two points: $\alpha$ and $\beta = -\alpha - a/2$, where $\alpha$ and $\beta$ are the simultaneous double roots of the varied quartic $x^4 + ax^3 + bx^2 + c_0 x + d_0$.
\vskip.1cm
As $c < c_0$, one has $-\delta_1 < - \delta_3$. Obviously, $-\delta_1  < -\delta_3  < 0 < c \rho_1 < -\delta_2 < c \rho_2$.
\vskip.1cm
The graph of $x^4 + ax^3 + bx^2$ is shown on its own in the top-left corner.
\vskip.1cm
Point $\mu_2$ could be in the second quadrant or in the third (pictured).
\vskip.1cm
Consideration of whether $- d < c \rho_1$ or $-d > c \rho_1$ would yield sharper bounds in the analysis based on solving cubic equations.
\vskip.1cm
\hrule
\vskip.1cm
\hskip1.54cm \underline{{\color{red}{\it Analysis based on solving cubic equations}}} \\
\vskip-.3cm
{\bf (i)} If $-d < -\delta_1$, then there are no real roots.
\vskip.1cm
{\bf (ii)} If $-\delta_1 \le -d < -\delta_3$, then there is a positive root $x_2$ such that $-d/c < x_2 \le \mu_1$ and another positive root $x_1$ such that $\mu_1 \le x_1 <  \xi^{(3)}_1 $.
\vskip.1cm
{\bf (iii)}  If $-\delta_3 \le -d < 0$, then there is a negative root $x_4$ such that $\lambda_2 < x_4 \le \mu_3$, another negative root $x_3$ such that $\mu_3 \le x_3 < \lambda_0$, a positive root $x_2$ such that $ -d/c < x_2 \le  \xi^{(3)}_2$ and another positive root $x_1$ such that  $\xi^{(3)}_1 \le x_1 < \lambda_1$.
\vskip.1cm
{\bf (iv)}  If $0 \le -d < -\delta_2$, then there is a negative root $x_4$ such that $\xi^{(2)}_2 < x_4 \le \lambda_2$, a negative root $x_3$ such that $\lambda_0 \le x_3 < \mu_2$, a non-positive root $x_2 \ge $ max$\{\mu_2, -d/c \}$ and a positive root $x_1$ such that $\lambda_1 \le x_1 <  \xi^{(2)}_1$.
\vskip.1cm
{\bf (v)} If $-\delta_2 \le -d$ (pictured), then there is a negative root $x_2$ such that min$\{\rho_2, -d/c \} \le x_2 \le $ min$\{$max$\{\rho_2, -d/c \}, \xi^{(2)}_2 \}$  and a positive root $x_1 \ge \xi^{(2)}_1$.
\vskip.1cm
\hrule
\vskip.1cm
\hskip.88cm \underline{{\color{red}{\it Analysis based on solving quadratic equations only}}}
\vskip.1cm
{\bf (i)} If $-d < 0$, then there are either no real roots, or there are two positive roots greater than $-d/c$, or there are two positive roots greater than $-d/c$ and two negative roots greater than $\rho_2$ and smaller than $\rho_1$.
\vskip.1cm
{\bf (ii)} If $0 \le - d <  c \rho_1$, then there is one positive root, one non-positive root greater than or equal to $-d/c$ and two negative roots greater than $\rho_2$ and smaller than $\rho_1$.
\vskip.1cm
{\bf (iii)} If $ c \rho_1 \le -d < c \rho_2$ (pictured), then there is one positive root together with either one or three negative roots greater than $\rho_2$ and smaller than or equal to $-d/c$.
\vskip.1cm
{\bf (iv)} If $c \rho_2 \le -d$, then there is one negative root $x_2$ such that $-d/c < x_2 \le \rho_2$ and one positive root $x_1$.
\\
\end{minipage}}  \\
\hline
\end{tabular}
\end{center}



\begin{center}
\begin{tabular}{|c||c|}
\hline
{\footnotesize \bf Figure 2.11} &  {\footnotesize \bf Figure 2.12} \\
\hline
\hline
{\color{violet} $\bm{0 < b \le \frac{a^2}{4}\, , \quad a > 0, \quad c < 0,$}} &  {\color{violet} $\bm{0 < b \le \frac{a^2}{4}\, , \quad a > 0, \quad c >  0,$}} \\
{\color{violet} $\bm{c_2 < \gamma_2 < c_0 < {\color{red}c} < 0 < \gamma_1 < c_1 $}} & {\color{violet} $\bm{c_2 < \gamma_2 < c_0 < 0 < {\color{red}c} < \gamma_1 < c_1 $}} \\
\hline
\hline
\includegraphics[width=80mm]{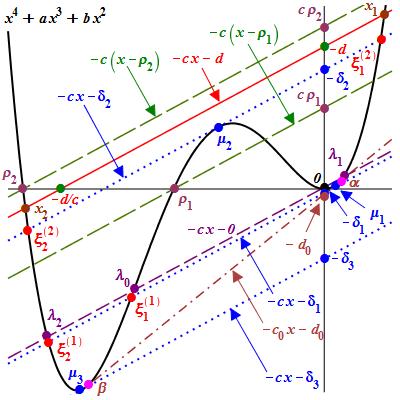} & \includegraphics[width=80mm]{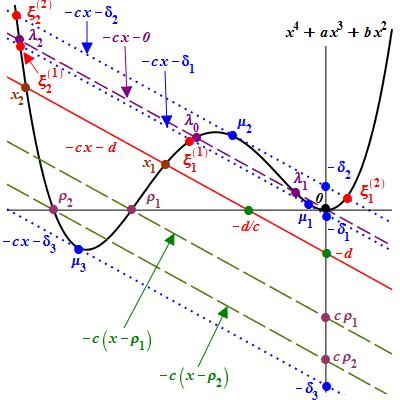} \\
\hline
\multicolumn{1}{|c||}{\begin{minipage}{22em}
\scriptsize
\vskip.1cm
\hskip3.89cm \underline{{\color{red}{\it Notes}}}\\
\vskip-.1cm
As $c_2 < c < c_1$, the quartic has three stationary points $\mu_i$ and the number of real roots can be either 0, or 2, or 4. There are three tangents $- c x - \delta_i$ to $x^4 + ax^3 + bx^2$ --- each at $\mu_i$.
\vskip.1cm
As $\gamma_2 < c < \gamma_1$, the straight line $-cx$ intersects $x^4 + ax^3 + bx^2$ at three points: $\lambda_{0,1,2}$. Thus $-\delta_3 < 0$. Also: $-\delta_1 <  0$ and $-\delta_2 > 0$.
\vskip.1cm
The straight line $ - c_0 x - d_0$ with $c_0 = (1/2)a(b-a^2/4) < 0$ and $d_0 = (1/4)(b-a^2/4)^2 > 0$ is tangent to $x^4 + ax^3 + bx^2$ at two points: $\alpha$ and $\beta = -\alpha - a/2$, where $\alpha$ and $\beta$ are the simultaneous double roots of the varied quartic $x^4 + ax^3 + bx^2 + c_0 x + d_0$.
\vskip.1cm
As $c_0 < c$, one has $- \delta_3 < -\delta_1$. One can have $-\delta_3  < -\delta_1  < 0 < c \rho_1 < -\delta_2 < c \rho_2$ (pictured) or $-\delta_3  < -\delta_1  < 0 < c \rho_1  < c \rho_2 \le -\delta_2$.
\vskip.1cm
Consideration of whether $- d < c \rho_1$ or $-d > c \rho_1$ would yield sharper bounds in the analysis based on solving cubic equations. \\
\hrule
\vskip.1cm
\hskip1.54cm \underline{{\color{red}{\it Analysis based on solving cubic equations}}} \\
\vskip-.3cm
{\bf (i)} If $-d < -\delta_3$, then there are no real roots.
\vskip.1cm
{\bf (ii)} If $-\delta_3 \le -d < -\delta_1$, then there is a negative root $x_2$ such that $\xi^{(1)}_2 < x_2 \le \mu_3$ and another negative root $x_1$ such that $\mu_3 \le x_2 < \xi^{(1)}_1$.
\vskip.1cm
{\bf (iii)} If $-\delta_1 \le -d < 0$, then there is a negative root $x_4$ such that $\lambda_2 < x_4 \le \xi^{(1)}_2$, another negative root $x_3$ such that $\xi^{(1)}_1 \le x_3 < \lambda_0$, a positive root $x_2$ such that $ -d/c < x_2 \le \mu_1$ and another positive root $x_1$ such that $\mu_1 \le x_1 < \lambda_1$.
\vskip.1cm
{\bf (iv)}  If $0 \le -d < -\delta_2$, then there is a negative root $x_4$ such that $\xi^{(2)}_2 < x_4 \le \lambda_2$, a negative root $x_3$ such that $\lambda_0 \le x_3 < \mu_2$, a non-positive root $x_2 \ge $ max$\{ \mu_2, -d/c \}$ and a positive root $x_1$ such that $\lambda_1 \le x_1 <  \xi^{(2)}_1$.
\vskip.1cm
{\bf (v)} If $-\delta_2 \le -d$ (pictured), then there is a negative root $x_2$ such that min$\{\rho_2, -d/c \} \le x_2 \le $ min$\{$max$\{\rho_2, -d/c \}, \xi^{(2)}_2 \}$  and a positive root $x_1 \ge \xi^{(2)}_1$.
\\
\hrule
\vskip.1cm
\hskip.88cm \underline{{\color{red}{\it Analysis based on solving quadratic equations only}}}\\
\vskip-.3cm
{\bf (i)} If $-d < 0$, then there are either no real roots, or there are two negative roots greater than $\rho_2$ and smaller than $\rho_1$, or there are two negative roots greater than $\rho_2$ and smaller than $\rho_1$ and two positive roots greater than $-d/c$.
\vskip.1cm
{\bf (ii)} If $0 \le - d <  c \rho_1$, then there is one positive root, one non-positive root greater than or equal to $-d/c$ and two negative roots greater than $\rho_2$ and smaller than $\rho_1$.
\vskip.1cm
{\bf (iii)} If $ c \rho_1 \le -d < c \rho_2$ (pictured), then there is one positive root, a negative root greater than $\rho_2$ and smaller than or equal to $-d/c$ and either zero or two negative roots smaller than or equal to $\rho_1$.
\vskip.1cm
{\bf (iv)} If $c \rho_2 \le -d$, then there is one negative root $x_2$ such that $-d/c < x_2 \le \rho_2$, one positive root $x_1$ and either zero or two negative roots greater than $\rho_1$ (the latter appear when $c \rho_2 \le -\delta_2$).
\\
\end{minipage}}
& \multicolumn{1}{|c|}{\begin{minipage}{22em}
\scriptsize
\vskip-3.8cm
\hskip3.89cm \underline{{\color{red}{\it Notes}}}\\
\vskip-.1cm
As $c_2 < c < c_1$, the quartic has three stationary points $\mu_i$ and the number of real roots can be either 0, or 2, or 4. There are three tangents $- c x - \delta_i$ to $x^4 + ax^3 + bx^2$ --- each at $\mu_i$.
\vskip.1cm
As $\gamma_2 < c < \gamma_1$, the straight line $-cx$ intersects $x^4 + ax^3 + bx^2$ at three points: $\lambda_{0,1,2}$. Thus $-\delta_2 > 0$. Also: $-\delta_3 < 0$ and $-\delta_1 < 0$.
\vskip.1cm
Obviously, $-\delta_3 < c \rho_2 < c \rho_1 < -\delta_1 < 0  < -\delta_2$.
\vskip.1cm
Consideration of whether $- d < c \rho_2$ or $-d > c \rho_2$ would yield sharper bounds in the analysis based on solving cubic equations.
\vskip.1cm
\hrule
\vskip.1cm
\hskip1.54cm \underline{{\color{red}{\it Analysis based on solving cubic equations}}} \\
\vskip-.2cm
{\bf (i)} If $-d < -\delta_3$, then there are no real roots.
\vskip.1cm
{\bf (ii)} If $-\delta_3 \le - d < -\delta_1$ (pictured), then there is a negative root $x_2$ such that $\xi^{(1)}_2 < x_2 \le \mu_3$ and another negative root $x_1$ such that max$\{$min$\{\rho_1, -d/c \}, \mu_3 \} \le x_1 < $min$ \{ \xi^{(1)}_1, $ max$\{\rho_1, -d/c \}\}$.
\vskip.1cm
{\bf (iii)}  If $-\delta_1 \le -d < 0$, then there is a negative root $x_4$ such that $\lambda_2 < x_4 \le \xi^{(1)}_2$, a negative root $x_3$ such that $\xi^{(1)}_1 \le x_3 < \lambda_0$, another negative root $x_2$ such that $\lambda_1 < x_2 \le \mu_1$ and a fourth negative root $x_1$ such that $\mu_1 \le x_1 < -d/c$.
\vskip.1cm
{\bf (iv)}  If $0 \le -d < -\delta_2$, then there is a negative root $x_4$ such that $\xi^{(2)}_2 < x_4 \le \lambda_2$,
a negative root $x_3$ such that $\lambda_0 \le x_3 < \mu_2$, a third negative root $x_2$ such that $\mu_2 < x_2 \le \lambda_1$ and a non-negative root $x_1 \le $ min$\{-d/c, \xi^{(2)}_1\}$.
\vskip.1cm
{\bf (v)} If $-\delta_2 \le -d$, then there is a negative root $x_2 \le \xi^{(2)}_2$ and a positive root $x_1$ such that  $\xi^{(2)}_1 \le x_1 < -d/c$.
\vskip.1cm
\hrule
\vskip.1cm
\hskip.88cm \underline{{\color{red}{\it Analysis based on solving quadratic equations only}}}
\vskip.1cm
{\bf (i)} If $- d < c \rho_2$, then there are either no real roots or there are two negative roots $x_{1,2}$ such that $\rho_2 < x_{1,2} < \rho_1$.
\vskip.1cm
{\bf (ii)} If $c \rho_2 \le - d < c \rho_1$, then there are two negative roots $x_{1,2}$ such that $x_2 \le \rho_2$ and $-d/c \le  x_1 < \rho_1$.
\vskip.1cm
{\bf (iii)} If $c \rho_1 \le - d < 0$ (pictured), then there is one negative root smaller than $\rho_2$ and either one or three negative roots greater than or equal to $\rho_1$ and smaller than $-d/c$.
\vskip.1cm
{\bf (iv)} If $0 \le -d$, then there is a negative root smaller than $\rho_2$, either zero or two negative roots greater than $\rho_1$ and one  non-negative root smaller than or equal to $-d/c$.\\
\end{minipage}}  \\
\hline
\end{tabular}
\end{center}


\begin{center}
\begin{tabular}{|c||c|}
\hline
{\footnotesize \bf Figure 2.13} &  {\footnotesize \bf Figure 2.14} \\
\hline
\hline
{\color{violet} $\bm{0 < b \le \frac{a^2}{4}\, , \quad a > 0, \quad c > 0,$}} &  {\color{violet} $\bm{0 < b \le \frac{a^2}{4}\, , \quad a > 0, \quad c >  0,$}} \\
{\color{violet} $\bm{c_2 < \gamma_2 < c_0 < 0 < \gamma_1 < {\color{red}c} < c_1 $}} & {\color{violet} $\bm{c_2 < \gamma_2 < c_0 < 0  < \gamma_1 < c_1 < {\color{red}c}$}} \\
\hline
\hline
\includegraphics[width=80mm]{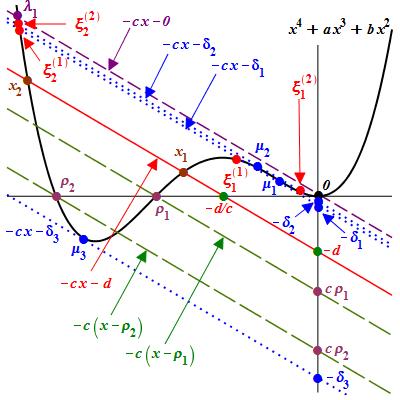} & \includegraphics[width=80mm]{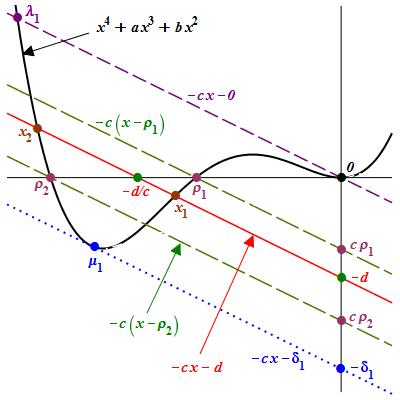} \\
\hline
\multicolumn{1}{|c||}{\begin{minipage}{22em}
\scriptsize
\vskip.3cm
\hskip3.89cm \underline{{\color{red}{\it Notes}}}\\
\vskip-.2cm
As $c_2 < c < c_1$, the quartic has three stationary points $\mu_i$ and the number of real roots can be either 0, or 2, or 4. There are three tangents $- c x - \delta_i$ to $x^4 + ax^3 + bx^2$ --- each at $\mu_i$.
\vskip.1cm
As $\gamma_2 < \gamma_1 < c$, the straight line $-cx$ intersects $x^4 + ax^3 + bx^2$ at one point only ($\lambda_1$). Thus $- \delta_2 < 0.$ Also:  $-\delta_1 < 0$ and $- \delta_3 < 0$.
\vskip.1cm
Obviously, $-\delta_3 < c \rho_2 < c \rho_1 < -\delta_1 < -\delta_2 < 0$.
\vskip.1cm
Consideration of whether $- d < c \rho_2$ or $-d > c \rho_2$ would yield sharper bounds in the analysis based on solving cubic equations.\\
\hrule
\vskip.3cm
\hskip1.54cm \underline{{\color{red}{\it Analysis based on solving cubic equations}}} \\
\vskip-.1cm
{\bf (i)} If $-d < -\delta_3$, then there are no real roots.
\vskip.1cm
{\bf (ii)} If $-\delta_3 \le - d < -\delta_1$ (pictured), then there is a negative root $x_2$ such that $\xi^{(1)}_2 < x_2 \le \mu_3$ and another negative root $x_1$ such that max$\{$min$\{\rho_1, -d/c \}, \mu_3 \} \le x_1 < $min$ \{ \xi^{(1)}_1, $ max$\{\rho_1, -d/c \}\}$.
\vskip.1cm
{\bf (iii)}  If $-\delta_1 \le -d < -\delta_2$, then there is a negative root $x_4$ such that $\xi^{(2)}_2 < x_4 \le \xi^{(1)}_2$, a negative root $x_3$ such that $\xi^{(1)}_1 \le x_3 < \mu_2$, another negative root $x_2$ such that $\mu_2 < x_2 \le \mu_1$ and a fourth negative root $x_1$ such that $\mu_1 \le x_1 < $ min$\{-d/c, \xi^{(2)}_1\}$.
\vskip.1cm
{\bf (iv)}  If $-\delta_2 \le -d < 0$, then there is a negative root $x_2$ such that $\lambda_1 < x_2 \le \xi^{(2)}_2$ and a negative root $x_1$ such that $\xi^{(2)}_1 \le x_2 < -d/c$.
\vskip.1cm
{\bf (v)} If $0 \le -d$, then there is a negative root $x_2 \le \lambda_1$ and a non-negative root $x_1 \le -d/c$.
\vskip.1cm
\hrule
\vskip.3cm
\hskip.88cm \underline{{\color{red}{\it Analysis based on solving quadratic equations only}}}
\vskip.3cm
{\bf (i)} If $- d < c \rho_2$, then there are either no real roots or there are two negative roots $x_{1,2}$ such that $\rho_2 < x_{1,2} < \rho_1$.
\vskip.1cm
{\bf (ii)} If $c \rho_2 \le - d < c \rho_1$, then there are two negative roots $x_{1,2}$ such that $x_2 \le \rho_2$ and $-d/c < x_1 < \rho_1$.
\vskip.1cm
{\bf (iii)} If $c \rho_1 \le - d < 0$ (pictured), then there is one negative root smaller than $\rho_2$ and either one or three negative roots greater than or equal to $\rho_1$ and smaller than $-d/c$.
\vskip.1cm
{\bf (iv)} If $0 \le -d$, then there is a negative root smaller than $\rho_2$ and one non-negative root smaller than or equal to $-d/c$. \\
\end{minipage}}
& \multicolumn{1}{|c|}{\begin{minipage}{22em}
\scriptsize
\vskip-1.8cm
\hskip3.89cm \underline{{\color{red}{\it Notes}}}
\vskip.1cm
As $c_2 < c_1 < c$, the quartic has a single stationary point $\mu_1$ and the number of real roots can be either 0 or 2. There is only one tangent $- c x - \delta_1$ to $x^4 + ax^3 + bx^2$ --- the one at $\mu_1$.
\vskip.1cm
As $\gamma_2 < \gamma_1 < c$, the straight line $-cx$ intersects $x^4 + ax^3 + bx^2$ at one point only ($\lambda_1$).
\vskip.1cm
One could have $\mu_1 > \rho_2$ (pictured) or $\mu_1 \le \rho_2$.
\vskip.1cm
Obviously, $-\delta_1 \le c \rho_2 < c \rho_1 < 0$.
\vskip.1cm
Consideration of whether $- d < c \rho_2$ or $-d > c \rho_2$ would yield sharper bounds in the analysis based on solving cubic equations. \\
\hrule
\vskip.3cm
\hskip1.54cm \underline{{\color{red}{\it Analysis based on solving cubic equations}}} \\
\vskip-.1cm
{\bf (i)} If $-d < -\delta_1$, then there are no real roots.
\vskip.1cm
{\bf (ii)} If $-\delta_1 \le -d < 0$ (pictured), then there is a negative root $x_2$ such that $\lambda_1 < x_2 \le \mu_1$ and another negative root $x_1$ such that max$\{$min$\{\rho_1, -d/c \}, \mu_1 \} \le x_1 < $ max$\{\rho_1, -d/c \}$.
\vskip.1cm
{\bf (iii)} If $0 \le -d$, then there is a negative root $x_2 \le \lambda_1$ and a non-negative root $x_1 \le -d/c$. \\
\hrule
\vskip.3cm
\hskip.88cm \underline{{\color{red}{\it Analysis based on solving quadratic equations only}}}
\vskip.3cm
{\bf (i)} If $- d < c \rho_2$, then there are either no real roots or there are two negative roots (these are smaller than $\rho_2$  if
$\mu_1 < \rho_2$ and greater than $\rho_2$ if $\mu_1 \ge \rho_2$.
\vskip.1cm
{\bf (ii)} If $c \rho_2 \le - d < c \rho_1$ (pictured), then there are two negative roots $x_{1,2}$ such that $x_2 \le \rho_2$ and $-d/c \le  x_1 < \rho_1$.
\vskip.1cm
{\bf (iii)} If $c \rho_1 \le - d < 0$, then there is a negative root smaller than $\rho_2$ and another negative root greater than or equal to $\rho_1$ and smaller than $-d/c$.
\vskip.1cm
{\bf (iv)} If $0 \le -d$, then there is a negative root smaller than $\rho_2$ and non-negative root smaller than or equal to $-d/c$. \\
\end{minipage}}  \\
\hline
\end{tabular}
\end{center}


\begin{center}
\begin{tabular}{|c||c|}
\hline
{\footnotesize \bf Figure 3.1} &  {\footnotesize \bf Figure 3.2} \\
\hline
\hline
{\color{violet} $\bm{\frac{a^2}{4} < b \le \frac{9}{8}\,\frac{a^2}{4}\, , \quad a < 0, \quad c < 0,$}} &  {\color{violet} $\bm{\frac{a^2}{4} < b \le \frac{9}{8}\,\frac{a^2}{4}\, , \quad a < 0, \quad c < 0,$}} \\
{\color{violet} $\bm{{\color{red}c} < c_2 < \gamma_2 < c_0 < \gamma_1 <  0 < c_1 $}} & {\color{violet} $\bm{c_2 < {\color{red}c} < \gamma_2 < c_0 < \gamma_1 <  0 < c_1 $}} \\
\hline
\hline
\includegraphics[width=80mm]{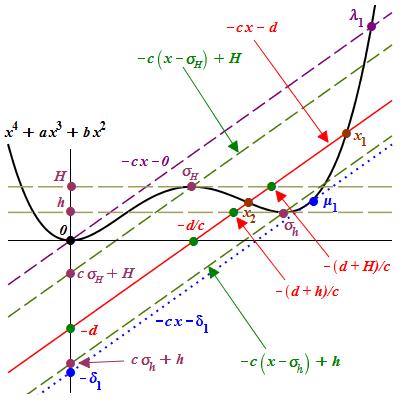} & \includegraphics[width=80mm]{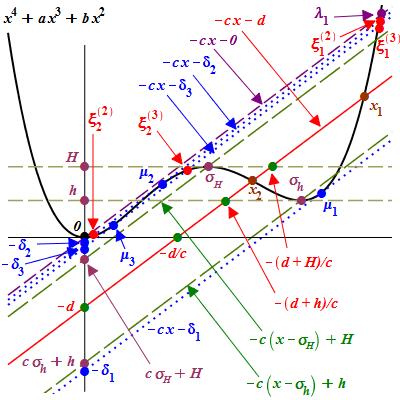} \\
\hline
\multicolumn{1}{|c||}{\begin{minipage}{22em}
\scriptsize
\vskip-3.5cm
\vskip.3cm
\hskip3.89cm \underline{{\color{red}{\it Notes}}}
\vskip.3cm
As $c < c_2 < c_1$, the quartic has a single stationary point $\mu_1$ and the number of real roots can be either 0 or 2. There is only one tangent $- c x - \delta_1$ to $x^4 + ax^3 + bx^2$ --- the one at $\mu_1$.
\vskip.1cm
As $c < \gamma_2 < \gamma_1$, the straight line $-cx$ intersects $x^4 + ax^3 + bx^2$ at one point only ($\lambda_1$).
\vskip.1cm
Obviously, $-\delta_1 < c \sigma_h + h < c \sigma_H + H < 0$.
\vskip.1cm
Consideration of whether $- d \le c \sigma_h + h$ or $-d > c \sigma_h + h$ and also whether $- d \le c \sigma_H + H$ or $-d > c \sigma_H + H$ would yield sharper bounds in the analysis based on solving cubic equations. \\
\hrule
\vskip.3cm
\hskip1.54cm \underline{{\color{red}{\it Analysis based on solving cubic equations}}}
\vskip.3cm
{\bf (i)} If $-d < -\delta_1$, then there are no real roots.
\vskip.1cm
{\bf (ii)} If $-\delta_1 \le -d < 0$ (pictured), then there is one positive root $x_2$ such that $-d/c < x_2 \le \mu_1$ and another positive root $x_1$ such that $\mu_1 \le x_2 < \lambda_1$.
\vskip.1cm
{\bf (iii)}  If $0 \le - d$, then there is a non-positive root $x_2 \ge -d/c$ and a positive root $x_1 \ge \lambda_1$. \\
\hrule
\vskip.3cm
\hskip.88cm \underline{{\color{red}{\it Analysis based on solving quadratic equations only}}}
\vskip.3cm
{\bf (i)} If $-d < c \sigma_h + h$, then there are either no real roots or there are two positive roots $x_{1,2} > \sigma_h$.
\vskip.1cm
{\bf (ii)} If $c \sigma_h + h \le - d < c \sigma_H + H$ (pictured), then there are two positive roots: $x_2$ such that max$\{ \sigma_H, -(d+h)/c \} < x_2 < $ min$\{ \sigma_h, -(d+H)/c \}$ and $x_1 > \sigma_h$.
\vskip.1cm
{\bf (iii)} If $c \sigma_H + H \le - d < 0$, then there are two positive roots: $x_2$ such that $-d/c < x_2 \le \sigma_H$ and $x_1 > \sigma_h$.
\vskip.1cm
{\bf (iv)} If $0 \le -d$, then there is one non-positive root $x_2 \ge -d/c$ and one positive root $x_1 > \sigma_h$. \\
\end{minipage}}
& \multicolumn{1}{|c|}{\begin{minipage}{22em}
\scriptsize
\vskip.3cm
\hskip3.89cm \underline{{\color{red}{\it Notes}}}
\vskip.3cm
As $c_2 < c < c_1$, the quartic has three stationary points $\mu_i$ and the number of real roots can be either 0, or 2, or 4. There are three tangents $- c x - \delta_i$ to $x^4 + ax^3 + bx^2$ --- each at $\mu_i$.
\vskip.1cm
As $c < \gamma_2 < \gamma_1$, the straight line $-cx$ intersects $x^4 + ax^3 + bx^2$ at one point only ($\lambda_1$). Thus $-\delta_2 < 0$ and $-\delta_1 < 0$. Also:  $-\delta_3 < 0$.
\vskip.1cm
Obviously, $-\delta_1 < c \sigma_h + h < c \sigma_H + H < -\delta_3 < -\delta_2 < 0$.
\vskip.1cm
Consideration of whether $- d \le c \sigma_h + h$ or $-d > c \sigma_h + h$ and also whether $- d \le c \sigma_H + H$ or $-d > c \sigma_H + H$ would yield sharper bounds in the analysis based on solving cubic equations. \\
\hrule
\vskip.3cm
\hskip1.54cm \underline{{\color{red}{\it Analysis based on solving cubic equations}}}
\vskip.3cm
{\bf (i)} If $-d < -\delta_1$, then there are no real roots.
\vskip.1cm
{\bf (ii)} If $-\delta_1 \le -d < -\delta_3$ (pictured), then there is one positive root $x_2$ such that max$\{ \xi^{(3)}_2, -d/c \} < x_2 \le \mu_1$ and another positive root $x_1$ such that $\mu_1 \le x_1 < \xi^{(3)}_1$.
\vskip.1cm
{\bf (iii)}  If $-\delta_3 \le -d < -\delta_2$, then there are four positive roots: $x_4$ such that max$\{ \xi^{(2)}_2, -d/c \} < x_4 \le \mu_3, \,\,$ $x_3$ such that $\mu_3 \le x_3 < \mu_2, \,\,$ $x_2$ such that $\mu_2 < x_2 \le \xi^{(3)}_2$, and $x_1$ such that $\xi^{(3)}_1 \le x_1 < \xi^{(2)}_1$. \vskip.1cm
{\bf (iv)} If $-\delta_2 \le -d < 0$, then there is a positive root $x_2$ such that $-d/c < x_2 \le \xi^{(2)}_2$ and another positive root $x_1$ such that $\xi^{(2)}_1 \le x_1 < \lambda_1$.
\vskip.1cm
{\bf (v)}  If $0 \le - d$, then there is a non-positive root $x_2 \ge -d/c$ and a positive root $x_1 \ge \lambda_1$. \\
\hrule
\vskip.3cm
\hskip.88cm \underline{{\color{red}{\it Analysis based on solving quadratic equations only}}}
\vskip.3cm
{\bf (i)} If $-d < c \sigma_h + h$, then there are either no real roots or there are two positive roots $x_{1,2} > \sigma_h$.
\vskip.1cm
{\bf (ii)} If $c \sigma_h + h \le - d < c \sigma_H + H$ (pictured), then there is a positive root $x_2$ such that max$\{ \sigma_H, -(d+h)/c \} < x_2 \le $ min$\{ \sigma_h, -(d+H)/c \}$ and a positive root $x_1 > \sigma_h$.
\vskip.1cm
{\bf (iii)} If $c \sigma_H + H \le - d < 0$, then there is a positive root greater than $\sigma_h$ and either one or three positive roots greater than $-d/c$ and smaller than or equal to $\sigma_H$.
\vskip.1cm
{\bf (iv)} If $0 \le -d$, then there is one non-positive root $x_2 \ge -d/c$ and one positive root $x_1 > \sigma_h$. \\
\end{minipage}}
\\
\hline
\end{tabular}
\end{center}


\begin{center}
\begin{tabular}{|c|c|}
\hline
\multicolumn{2}{|c|}{\footnotesize \bf Figure 3.3}  \\
\multicolumn{2}{|c|}{\footnotesize \bf (continues on next page)}  \\
\hline
\hline
\multicolumn{2}{|c|}{{\color{violet} $\bm{\frac{a^2}{4} < b \le \frac{9}{8}\,\frac{a^2}{4}\, , \quad a < 0, \quad c < 0,$}}} \\
\multicolumn{2}{|c|}{{\color{violet} $\bm{c_2 < \gamma_2 < {\color{red}c} < c_0 < \gamma_1 <  0 < c_1 $}}} \\
\hline
\hline
\includegraphics[width=80mm,height=80mm]{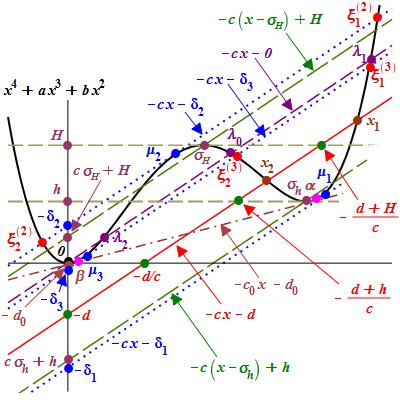} & \includegraphics[width=80mm,height=80mm]{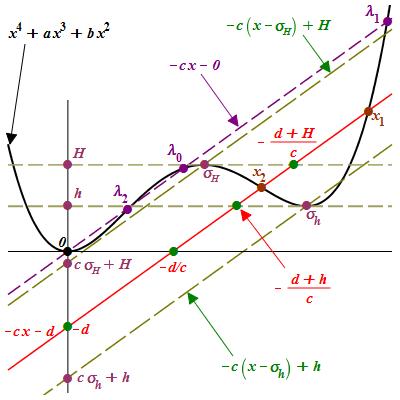} \\
\hline
\multicolumn{1}{|c|}{\begin{minipage}{22em}
\scriptsize
\begin{center}
\vskip-7.3cm
\underline{{\color{red}{\it Notes (apply to all panes)}}}
\end{center}
\vskip-.2cm
As $c_2 < c < c_1$, the quartic has three stationary points $\mu_i$ and the number of real roots can be either 0, or 2, or 4. There are three tangents $- c x - \delta_i$ to $x^4 + ax^3 + bx^2$ --- each at $\mu_i$ (only shown on the top-left pane).
\vskip.1cm
As $\gamma_2 < c < \gamma_1$, the straight line $-cx$ intersects $x^4 + ax^3 + bx^2$ at three points: $\lambda_{0,1,2}$. Thus $-\delta_2 > 0$ and $-\delta_1 < 0$. Also: $-\delta_3 < 0$.
\vskip.1cm
The straight line $ - c_0 x - d_0$ with $c_0 = (1/2)a(b-a^2/4) < 0$ and $d_0 = (1/4)(b-a^2/4)^2 > 0$ is tangent to $x^4 + ax^3 + bx^2$ at two points: $\alpha$ and $\beta = -\alpha - a/2$, where $\alpha$ and $\beta$ are the simultaneous double roots of the varied quartic $x^4 + ax^3 + bx^2 + c_0 x + d_0$.
\vskip.1cm
As $c < c_0$, one has $- \delta_1 < -\delta_3$.
\vskip.1cm
One can have $c \sigma_h + h < 0 < c \sigma_H + H$ (top-left pane), $c \sigma_h + h < c \sigma_H + H < 0$ (top-right pane; $c$ closer to $\gamma_2$), and $0 < c \sigma_h + h < c \sigma_H + H$ (bottom-right pane; $c$ closer to $c_0$). For each pane, there are further sub-cases affecting the analysis based on solving quadratic equations only (see below).
\vskip.1cm
Consideration of whether $- d \le c \sigma_h + h$ or $-d > c \sigma_h + h$ and also whether $- d \le c \sigma_H + H$ or $-d > c \sigma_H + H$ would yield sharper bounds in the analysis based on solving cubic equations.
\\
\end{minipage}}
 & \includegraphics[width=80mm,height=80mm]{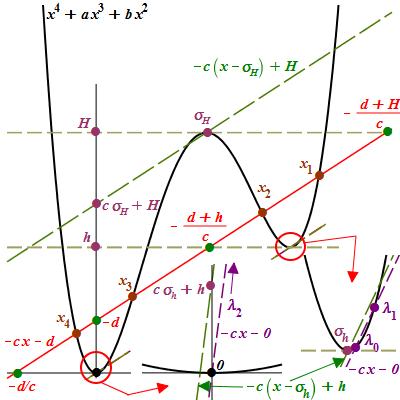} \\
\hline
\multicolumn{2}{|c|}{\begin{minipage}{45em}
\scriptsize
\vskip.1cm
\begin{center}
\underline{{\color{red}{\it Analysis based on solving cubic equations (applies to all panes)}}}
\end{center}
\vskip-.3cm
{\bf (i)} If $-d < -\delta_1$, then there are no real roots.
\vskip.1cm
{\bf (ii)} If $-\delta_1 \le -d < -\delta_3$ (pictured), then there is one positive root $x_2$ such that max$\{ \xi^{(3)}_2, -(d+h)/c \} < x_2 \le $ min$\{ \mu_1, -(d+H)/c \}$ and another positive root $x_1$ such that $\mu_1 \le x_1 < \xi^{(3)}_1$.
\vskip.1cm
{\bf (iii)} If $-\delta_3 \le -d < 0$, then there are four positive roots: $x_4$ such that $- d/c < x_4 \le \mu_3, \,\,$ $x_3$ such that $\mu_3 \le x_3 < \lambda_2, \,\,$ $x_2$ such that $\lambda_0 < x_2 \le \xi^{(3)}_2$, and $x_1$ such that $\xi^{(3)}_1 \le x_1 < \lambda_1$.
\vskip.1cm
{\bf (iv)}  If $0 \le -d < -\delta_2$, then there is one non-positive root $x_4 \ge $ max$\{ \xi^{(2)}_2, -d/c \}$ and tree positive roots: $x_3$ such that $\lambda_2 \le x_3 < \mu_2, \,\,$ $x_2$ such that $\mu_2 < x_2 \le \lambda_0$, and $x_1$ such that $\lambda_1 \le x_1 < \xi^{(2)}_1$.
\vskip.1cm
{\bf (v)}  If $-\delta_2 \le - d $, then there is a negative root $x_2$ such that $-d/c < x_2 \le \xi^{(2)}_2$ and a positive root $x_1 \ge \xi^{(2)}_1$. \\
\hrule
\begin{center}
\vskip-.3cm
\underline{{\color{red}{\it Analysis based on solving quadratic equations only (continues on next page)}}}
\end{center}
\vskip-.3cm
{\bf \underline{\color{red}{Top-left pane}}}
\vskip.2cm
One can have either $-\delta_2 > c \sigma_H + H > 0 > -\delta_3 > c \sigma_h + h > -\delta_1$ (pictured) or $-\delta_2 > c \sigma_H + H > 0 > c \sigma_h + h \ge -\delta_3 > -\delta_1$.
\vskip.2cm
{\bf (i)} If $- d < c \sigma_h + h$, then there are either no real roots, or there are two positive roots greater than $\sigma_h$, or there are two positive roots greater than $\sigma_h$ and another pair of positive roots greater than $-d/c$ and smaller than $\sigma_H$ (the latter appear when $c \sigma_h + h \ge -\delta_3$).
\vskip.1cm
{\bf (ii)} If $c \sigma_h + h \le - d < 0$ (pictured), then there is one positive root $x_2$ such that max$\{ \sigma_H, -(d+h)/c \} < x_2 \le $ min$\{ \sigma_h, -(d+H)/c \}$, another positive root $x_1 > \sigma_h$ and either zero or two positive roots $x_{3,4}$ greater than $-d/c$ and smaller than $\sigma_H$ ($x_{3,4}$ are always present if $c \sigma_h + h \ge -\delta_3$, while for the pictured $c \sigma_h + h < -\delta_3$ the roots $x_{3,4}$ may or may not be there).
\vskip.1cm
{\bf (iii)} If $0 \le - d < c \sigma_H + H$, then there is one non-positive root $x_4 \ge -d/c$, a positive root $x_3 < \sigma_H$, another positive root $x_2$ such that $\sigma_H < x_2 < -(d+H)/c$, and a third positive root $x_1 > \sigma_h$.
\vskip.1cm
{\bf (iv)} If $c \sigma_H + H \le -d$, then there is one negative root $x_2 > -d/c$, a positive root $x_1 > \sigma_h$ and either zero or two positive roots smaller than or equal to $\sigma_H$. \\
\end{minipage}}\\
\hline
\end{tabular}
\end{center}


\begin{center}
\begin{tabular}{|c|c|}
\hline
\multicolumn{2}{|c|}{\footnotesize \bf Figure 3.3}  \\
\multicolumn{2}{|c|}{\footnotesize \bf (continued from previous page)}  \\
\hline
\hline
\multicolumn{2}{|c|}{{\color{violet} $\bm{\frac{a^2}{4} < b \le \frac{9}{8}\,\frac{a^2}{4}\, , \quad a < 0, \quad c < 0,$}}} \\
\multicolumn{2}{|c|}{{\color{violet} $\bm{c_2 < \gamma_2 < {\color{red}c} < c_0 < \gamma_1 <  0 < c_1 $}}} \\
\hline
\hline
\multicolumn{2}{|c|}{\begin{minipage}{45em}
\scriptsize
\begin{center}
\vskip.2cm
\underline{{\color{red}{\it Analysis based on solving quadratic equations only --- continued from previous page}}}
\end{center}
{\bf \underline{\color{red}{Top-right pane}}}
\vskip.2cm
One can have either $-\delta_2 > 0 > -\delta_3 > c \sigma_H + H >  c \sigma_h + h > -\delta_1$, or $-\delta_2 > 0 > c \sigma_H + H \ge -\delta_3 > c \sigma_h + h > -\delta_1$, or $-\delta_2 > 0 > c \sigma_H + H > c \sigma_h + h \ge -\delta_3 > -\delta_1$ (when $\sigma_h$ and $\sigma_H$ are close).
\vskip.2cm
{\bf (i)} If $- d < c \sigma_h + h$, then there are either no real roots, or there are two positive roots greater than $\sigma_h$, or there are two positive roots greater than $\sigma_h$ and two positive roots greater than $-d/c$ and smaller than $\sigma_H$ (the latter appear when $c \sigma_h + h \ge -\delta_3$).
\vskip.1cm
{\bf (ii)} If $c \sigma_h + h \le - d < c \sigma_H + H$ (pictured), then there is a positive root $x_2$ such that max$\{ \sigma_H, -(d+h)/c \} < x_2 \le $ min$\{ \sigma_h, -(d+H)/c \}$ and a positive root $x_1 > \sigma_h$. If $-\delta_3 > c \sigma_H + H$, there are no other roots. If $c \sigma_h + h \ge -\delta_3$, there are two more positive roots $x_{3,4}$ greater than $-d/c$ and smaller than $\sigma_H$. If $c \sigma_H + H \ge -\delta_3 > c \sigma_h + h$, the roots $x_{3,4}$ may or may not be there.
\vskip.1cm
{\bf (iii)} If $c \sigma_H + H \le - d < 0$, then there is a positive root greater than $\sigma_h$ and, if $-\delta_3 > c \sigma_H + H$, there is either one or three positive roots greater than $-d/c$ and smaller than $\sigma_H$. If however, $c \sigma_H + H \ge -\delta_3$, then there will be a positive root greater than $\sigma_h$ and three positive roots greater than $-d/c$ and smaller than $\sigma_H$.
\vskip.1cm
{\bf (iv)} If $0 \le -d$, then there is a non-positive root greater than or equal to $-d/c$, either zero or two positive roots smaller than $\sigma_H$, and a positive root greater than $\sigma_h$. \\
\hrule
\vskip.3cm
{\bf \underline{\color{red}{Bottom-right pane}}}
\vskip.2cm
One can only have $-\delta_2 >  c \sigma_H + H >  c \sigma_h + h > 0 > -\delta_3 > -\delta_1$.
\vskip.2cm
{\bf (i)} If $- d < 0$, then there are either no real roots, or there are two positive roots greater than $\sigma_h$, or there are four positive roots: two greater than $\sigma_h$ and the other two --- greater than $-d/c$ and smaller than $\sigma_H$.
\vskip.1cm
{\bf (ii)} If $0 \le - d < c \sigma_h + h$, then there is one non-positive root greater than $-d/c$, one positive root smaller than $\sigma_H$, and two positive roots greater than $\sigma_h$.
\vskip.1cm
{\bf (iii)} If $c \sigma_h + h \le - d < c \sigma_H + H$ (pictured), then there is one negative root $x_4 > -d/c$, a positive root $x_3 < \sigma_H$, a positive root $x_2$ such that max$\{\sigma_H, -(d+h)/c \} < x_2 \le $ min$\{ -(d+H)/c, \sigma_h \}$ and another positive root $x_1 \ge \sigma_h$.
\vskip.1cm
{\bf (iv)} If $c \sigma_H + H \le -d$, then there is a negative root greater than $-d/c$, a positive root greater than $\sigma_h$ and either zero or two positive roots smaller than $\sigma_H$. \\
\end{minipage}}\\
\hline
\end{tabular}
\end{center}


\begin{center}
\begin{tabular}{|c||c|}
\hline
\multicolumn{2}{|c|}{\footnotesize \bf Figure 3.4}  \\
\hline
\hline
\multicolumn{2}{|c|}{{\color{violet} $\bm{\frac{a^2}{4} < b \le \frac{9}{8}\,\frac{a^2}{4}\, , \quad a < 0, \quad c < 0,$}}} \\
\multicolumn{2}{|c|}{{\color{violet} $\bm{c_2 < \gamma_2 < c_0 < {\color{red}c} < \gamma_1 <  0 < c_1 $}}}  \\
\hline
\hline
\multirow{2}{*}{{\begin{minipage}{22em}
\vskip-3.75cm
\includegraphics[width=80mm]{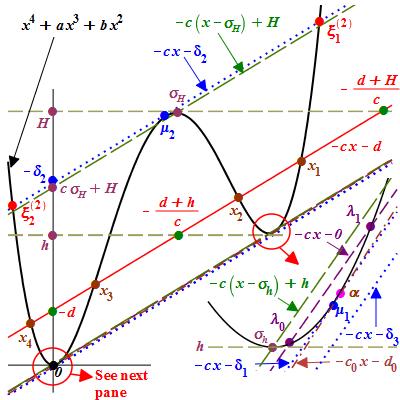}\end{minipage}}} &
\multicolumn{1}{|c|}{\begin{minipage}{22em}
\includegraphics[width=80mm]{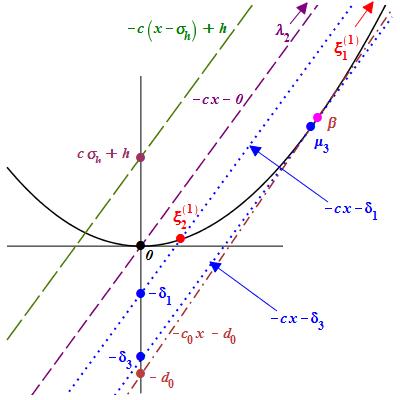}
\end{minipage}} \\
\hline
\multicolumn{2}{|c|}{\begin{minipage}{45em}
\scriptsize
\vskip.5cm
\begin{center}
\underline{{\color{red}{\it Notes}}}
\end{center}
\vskip-.3cm
As $c_2 < c < c_1$, the quartic has three stationary points $\mu_i$ and the number of real roots can be either 0, or 2, or 4. There are three tangents $- c x - \delta_i$ to $x^4 + ax^3 + bx^2$ --- each at $\mu_i$.
\vskip.1cm
As $\gamma_2 < c < \gamma_1$, the straight line $-cx$ intersects $x^4 + ax^3 + bx^2$ at three points: $\lambda_{0,1,2}$. Thus $-\delta_2 > 0$ and $-\delta_1 < 0$. Also: $-\delta_3 < 0$.
\vskip.1cm
The straight line $ - c_0 x - d_0$ with $c_0 = (1/2)a(b-a^2/4) < 0$ and $d_0 = (1/4)(b-a^2/4)^2 > 0$ is tangent to $x^4 + ax^3 + bx^2$ at two points: $\alpha$ and $\beta = -\alpha - a/2$, where $\alpha$ and $\beta$ are the simultaneous double roots of the varied quartic $x^4 + ax^3 + bx^2 + c_0 x + d_0$.
\vskip.1cm
As $c > c_0$, one has $- \delta_1 > -\delta_3$.
\vskip.1cm
Obviously, $-\delta_3 < -\delta_1 < 0 < c \sigma_h + h <  c \sigma_H + H < -\delta_2$.
\vskip.1cm
Consideration of whether $- d \le c \sigma_h + h$ or $-d > c \sigma_h + h$ and also whether $- d \le c \sigma_H + H$ or $-d > c \sigma_H + H$ would yield sharper bounds in the analysis based on solving cubic equations. \\
\hrule
\vskip.3cm
\begin{center}
\underline{{\color{red}{\it Analysis based on solving cubic equations}}}
\end{center}
\vskip-.2cm
{\bf (i)} If $-d < -\delta_3$, then there are no real roots.
\vskip.1cm
{\bf (ii)} If $-\delta_3 \le -d < -\delta_1$, then there is one positive root $x_2$ such that max$\{ \xi^{(1)}_2, -d/c \} < x_2 \le \mu_3$  and another positive root $x_1$ such that $\mu_3 \le x_1 < \xi^{(1)}_1$.
\vskip.1cm
{\bf (iii)} If $-\delta_1 \le -d < 0$, then there are four positive roots: $x_4$ such that $- d/c < x_4 \le \xi^{(1)}_2, \,\,$ $x_3$ such that $\xi^{(1)}_1 \le x_3 < \lambda_2, \,\,$ $x_2$ such that $\lambda_0 < x_2 \le \mu_1$, and $x_1$ such that $\mu_1 \le x_1 < \lambda_1$.
\vskip.1cm
{\bf (iv)}  If $0 \le -d < -\delta_2$ (pictured), then there is one non-positive root $x_4 \ge $ max $\{ \xi^{(2)}_2, -d/c \}$ and tree positive roots: $x_3$ such that $\lambda_2 \le x_3 < \mu_2, \,\,$ $x_2$ such that $\mu_2 < x_2 \le \lambda_0$, and $x_1$ such that $\lambda_1 \le x_1 < \xi^{(2)}_1$.
\vskip.1cm
{\bf (v)}  If $-\delta_2 \le - d $, then there is a negative root $x_2$ such that $-d/c < x_2 \le \xi^{(2)}_2$ and a positive root $x_1 \ge \xi^{(2)}_1$. \\
\hrule
\vskip.3cm
\begin{center}
\underline{{\color{red}{\it Analysis based on solving quadratic equations only}}}
\end{center}
\vskip-.2cm
{\bf (i)} If $- d < 0$, then there are either no real roots, or there are two positive roots greater than $-d/c$ and smaller than $\sigma_H$, or there are four positive roots, two of which greater than $-d/c$ and smaller than $\sigma_H$ and the other two --- greater than $\sigma_h$.
\vskip.1cm
{\bf (ii)} If $0 \le - d < c \sigma_h + h$, then there is one non-positive root $x_4 \ge -d/c$, a non-negative root $x_3 < \sigma_H$ and twho positive roots $x_{1,2} > \sigma_h.$
\vskip.1cm
{\bf (iii)} If $c \sigma_h + h \le - d < c \sigma_H + H$ (pictured), then there is one negative root $x_4 > -d/c$, a positive root $x_3 < \sigma_H$, another positive root $x_2$ such that max$\{ \sigma_H, -(d+h)/c \} < x_2 \le $ min$\{ \sigma_h, -(d+H)/c \}$, and a third positive root $x_1 > \sigma_h$.
\vskip.1cm
{\bf (iv)} If $c \sigma_H + H \le -d$, then there is a negative root $x_2 > -d/c$ and a positive root $x_1 > \sigma_h$ and either zero or two positive roots smaller than $\sigma_H$.
\\
\end{minipage}} \\
\hline
\end{tabular}
\end{center}


\begin{center}
\begin{tabular}{|c||c|}
\hline
\multicolumn{2}{|c|}{\footnotesize \bf Figure 3.6}  \\
\multicolumn{2}{|c|}{\footnotesize \bf (Figure 3.5 is on next page)}  \\
\hline
\hline
\multicolumn{2}{|c|}{{\color{violet} $\bm{\frac{a^2}{4} < b \le \frac{9}{8}\,\frac{a^2}{4}\, , \quad a < 0, \quad c > 0,$}}}  \\
\multicolumn{2}{|c|}{{\color{violet} $\bm{c_2 < \gamma_2  < c_0 < \gamma_1 <  0 < {\color{red}c} < c_1 $}}}  \\
\hline
\hline
\includegraphics[width=80mm]{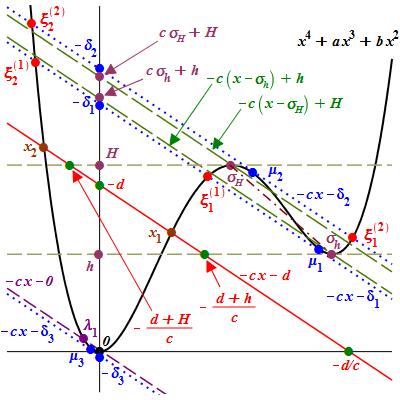} & \includegraphics[width=80mm]{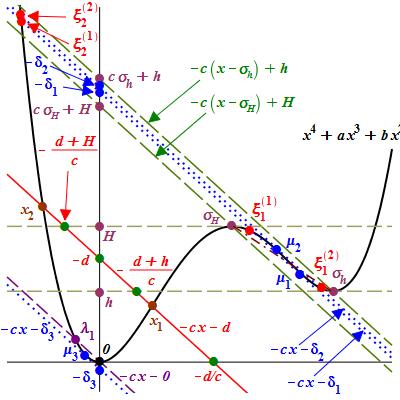} \\
\hline
\multicolumn{2}{|c|}{\scriptsize \underline{{\color{red}{\it Notes (apply to both panes)}}}} \\
\multicolumn{2}{|c|}{\begin{minipage}{45em}
\scriptsize
\vskip.1cm
As $c_2 < c < c_1$, the quartic has three stationary points $\mu_i$ and the number of real roots can be either 0, or 2, or 4. There are three tangents $- c x - \delta_i$ to $x^4 + ax^3 + bx^2$ --- each at $\mu_i$.
\vskip.1cm
As $\gamma_2 < \gamma_1 < c$, the straight line $-cx$ intersects $x^4 + ax^3 + bx^2$ at one point only ($\lambda_{1}$). Thus $-\delta_1 > 0$ and $-\delta_2 > 0$. Also, $-\delta_3 < 0$.
\vskip.1cm
For both panes, the straight line joining the two stationary points in the first quadrant has slope $-\hat{c} = (1/2)a[(9/8)(a^2/4)-b] \le 0$ such that $0 \le \hat{c} < c_1$.
\vskip.1cm
With the increase of $c$ from zero towards $c_1$, one has: $-\delta_3 < 0 < -\delta_1 < c \sigma_h + h < c \sigma_H + H < -\delta_2$ (pictured on left pane; $0 < c < \hat{c} < c_1$), then $-\delta_3 < 0 < -\delta_1 \le c \sigma_H + H \le c \sigma_h + h \le -\delta_2$ (at $c = \hat{c}$ there is a swap between $c \sigma_H + H$ and $c \sigma_h + h$), then $-\delta_3 < 0 < c \sigma_H + H  \le -\delta_1 <  c \sigma_h + h \le -\delta_2$, or $-\delta_3 < 0 < -\delta_1 \le c \sigma_H + H < -\delta_2 \le c \sigma_h + h$, and finally: $-\delta_3 < 0 < c \sigma_H + H \le -\delta_1 < -\delta_2 \le c \sigma_h + h$ (pictured on right pane; $0 < \hat{c} < c < c_1$).
\vskip.1cm
Consideration of whether $- d \le c \sigma_h + h$ or $-d > c \sigma_h + h$ and also whether $- d \le c \sigma_H + H$ or $-d > c \sigma_H + H$ would yield sharper bounds in the analysis based on solving cubic equations.
\vskip.1cm
\hrule
\begin{center}
\vskip-.3cm
\underline{{\color{red}{\it Analysis based on solving cubic equations (applies to both panes)}}} \\
\end{center}
\vskip-.3cm
{\bf (i)} If $-d < -\delta_3$, then there are no real roots.
\vskip.1cm
{\bf (ii)} If $-\delta_3 \le -d < 0$, then there is a negative root $x_2$ such that $\lambda_1 < x_2 \le \mu_3$ and a negative root $x_1$ such that $\mu_3 \le x_1 < -d/c$.
\vskip.1cm
{\bf (iii)} If $0 \le -d < -\delta_1$ (pictured), then there is a negative root $x_2$ such that $\xi^{(1)}_2 < x_2 \le \lambda_1$ and a non-negative root $x_1 \le $ min$\{\xi^{(1)}_1, -d/c \}$.
\vskip.1cm
{\bf (iv)}  If $-\delta_1 \le -d < -\delta_2$, then there is one negative root $x_4$ such that  $\xi^{(2)}_2 < x_4 \le \xi^{(1)}_2$ and tree positive roots: $x_3$ such that $\xi^{(1)}_1 \le x_3 < \mu_2, \,\,$ $x_2$ such that $\mu_2 < x_2 \le \mu_1$, and $x_1$ such that $\mu_1 \le x_1 < \xi^{(2)}_1$.
\vskip.1cm
{\bf (v)}  If $-\delta_2 \le - d $, then there is a negative root $x_2  \le \xi^{(2)}_2$ and a positive root $x_1 \ge \xi^{(2)}_1$.
\vskip.1cm
\hrule
\end{minipage}}\\
\multicolumn{1}{|c|}{\begin{minipage}{22em}
\scriptsize
\begin{center}
\vskip-0.8cm
\underline{{\color{red}{\it Analysis based on solving quadratic equations only}}}
\vskip.2cm
\underline{$0 < c < \hat{c} < c_1$}
\end{center}
\vskip-.2cm
{\bf (i)} If $- d < 0$, then there are either no real roots or there are two negative roots smaller than $-d/c$.
\vskip.1cm
{\bf (ii)} If $0 \le - d < c \sigma_h + h$ (pictured), then there is a non-positive root $x_2$, a non-negative root $x_1$ such that max$\{0, -(d+H)/c \} \le x_1 < $ min$\{ -d/c, \sigma_H \}$ and either zero or two positive roots greater than $\sigma_H$ and smaller than $\sigma_h$.
\vskip.1cm
{\bf (iii)} If $c \sigma_h + h \le - d < c \sigma_H + H$, then there is one negative root $x_4$, a positive root $x_3$ such that $-(d+H)/c < x_3 < \sigma_H$, another positive root $x_2$ such that $\sigma_H < x_2 \le \sigma_h$, and a third positive root $x_1$ such that $\sigma_h \le x_1 < -(d+h)/c$.
\vskip.1cm
{\bf (iv)} If $c \sigma_H + H \le -d$, then there is a negative root $x_2$, a positive root $x_1$ such that  $\sigma_h \le x_1 < -(d+h)/c$ and either zero or two positive roots greater than or equal to $-(d+H)/c$ and smaller than $\sigma_h$.\\
\end{minipage}}
& \multicolumn{1}{|c|}{\begin{minipage}{22em}
\scriptsize
\begin{center}
\vskip0.2cm
\underline{{\color{red}{\it Analysis based on solving quadratic equations only}}}
\vskip.2cm
\underline{$0 < \hat{c} < c < c_1$}
\end{center}
\vskip-.2cm
{\bf (i)} If $- d < 0$, then there are either no real roots or there are two negative roots smaller than $-d/c$.
\vskip.1cm
{\bf (ii)} If $0 \le - d < c \sigma_H + H$ (pictured), then there is a non-positive root $x_2$, a non-negative root $x_1$ such that max$\{0, -(d+H)/c \} \le x_1 \le $ min$\{ -d/c, \sigma_H \}$, and either zero or two positive roots greater than $\sigma_H$ and smaller than $\sigma_h$ (the latter appear when $-\delta_1 \le c \sigma_H + H$).
\vskip.1cm
{\bf (iii)} If $c \sigma_H + H \le - d < c \sigma_h + h$, then, if $-\delta_1 \le c \sigma_H + H$ and $-\delta_2 \ge c \sigma_h + h$, there will be one negative root together with three positive roots greater than or equal to $\sigma_H$ and smaller than $\sigma_h$. If however, $c \sigma_H + H < -\delta_1$ or $-\delta_2 < c \sigma_h + h$, there will be one negative root together with either one or three positive roots greater than or equal to $\sigma_H$ and smaller than $\sigma_h$.
\vskip.1cm
{\bf (iv)} If $c \sigma_h + h \le -d$, then there is a negative root $x_2$, a positive root $x_1$ such that  $\sigma_h \le x_1 < -(d+h)/c$ and either zero or two positive roots greater than or equal to $-(d+H)/c$ and smaller than $\sigma_h$ (the latter appear when $c \sigma_h + h \le -\delta_2$).
\\
\end{minipage}}  \\
\hline
\end{tabular}
\end{center}


\begin{center}
\begin{tabular}{|c||c|}
\hline
{\footnotesize \bf Figure 3.5} &  {\footnotesize \bf Figure 3.7} \\
\hline
\hline
{\color{violet} $\bm{\frac{a^2}{4} < b \le \frac{9}{8}\,\frac{a^2}{4}\, , \quad a < 0, \quad c < 0,$}} &  {\color{violet} $\bm{\frac{a^2}{4} < b \le \frac{9}{8}\,\frac{a^2}{4}\, , \quad a < 0, \quad c > 0,$}} \\
{\color{violet} $\bm{c_2 < \gamma_2 <  c_0 < \gamma_1 <  {\color{red}c} < 0 < c_1 $}} & {\color{violet} $\bm{c_2 < \gamma_2 < c_0 < \gamma_1 <  0 < c_1 < {\color{red}c}  $}} \\
\hline
\hline
\includegraphics[width=80mm]{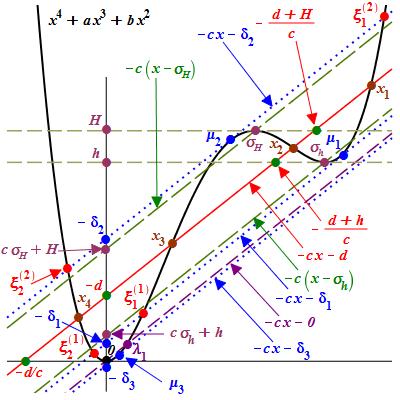} & \includegraphics[width=80mm]{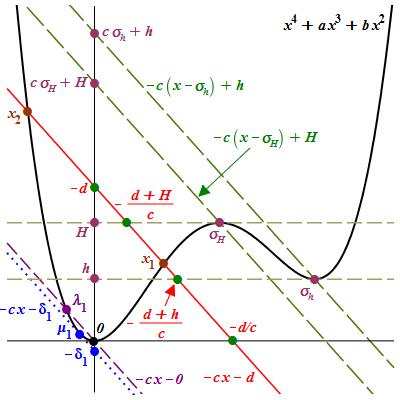} \\
\hline
\multicolumn{1}{|c||}{\begin{minipage}{22em}
\scriptsize
\vskip.2cm
\begin{center}
\underline{{\color{red}{\it Notes}}}
\end{center}
\vskip-.4cm
As $c_2 < c < c_1$, the quartic has three stationary points $\mu_i$ and the number of real roots can be either 0, or 2, or 4. There are three tangents $- c x - \delta_i$ to $x^4 + ax^3 + bx^2$ --- each at $\mu_i$.
\vskip.1cm
As $\gamma_2 < \gamma_1 < c$, the straight line $-cx$ intersects $x^4 + ax^3 + bx^2$ at one point only ($\lambda_1$).  Thus $-\delta_1 > 0$ and $-\delta_2 > 0$. Also: $-\delta_3 < 0$.
\vskip.1cm
Obviously, $-\delta_3 < 0 < -\delta_1 < c \sigma_h + h < c \sigma_H + H < -\delta_2$.
\vskip.1cm
Consideration of whether $- d \le c \sigma_h + h$ or $-d > c \sigma_h + h$ and also whether $- d \le c \sigma_H + H$ or $-d > c \sigma_H + H$ would yield sharper bounds in the analysis based on solving cubic equations. \\
\hrule
\vskip.2cm
\hskip1.54cm \underline{{\color{red}{\it Analysis based on solving cubic equations}}} \\
\vskip-.2cm
{\bf (i)} If $-d < -\delta_3$, then there are no real roots.
\vskip.1cm
{\bf (ii)} If $-\delta_3 \le -d < 0$, then there is one positive root $x_2$ such that $-d/c < x_2 \le \mu_3$ and another positive root $x_1$ such that $\mu_3 \le x_1 < \lambda_1$.
\vskip.1cm
{\bf (iii)} If $0 \le -d < -\delta_1$, then there is one non-positive root $x_2 > $ max$\{-d/c, \xi^{(1)}_2 \}$ and one positive root $x_1$ such that $\lambda_1 \le x_1 < \xi^{(1)}_1.$
\vskip.1cm
{\bf (iv)}  If $-\delta_1 \le -d < -\delta_2$ (pictured), then there is one negative root $x_4$ such that max $\{ \xi^{(2)}_2, -d/c \} \le x_4 < \xi^{(1)}_2$ and tree positive roots: $x_3$ such that $\xi^{(1)}_1 \le x_3 < \mu_2, \,\,$ $x_2$ such that $\mu_2 < x_2 \le \mu_1$, and $x_1$ such that $\mu_1 \le x_1 < \xi^{(2)}_1$.
\vskip.1cm
{\bf (v)}  If $-\delta_2 \le - d $, then there is a negative root $x_2$ such that $-d/c < x_2 \le \xi^{(2)}_2$ and a positive root $x_1 \ge \xi^{(2)}_1$.
\\
\hrule
\vskip.2cm
\hskip.88cm \underline{{\color{red}{\it Analysis based on solving quadratic equations only}}} \\
\vskip-.2cm
{\bf (i)} If $- d < 0$, then there are either no real roots or there are two positive roots greater than $-d/c$ and smaller than $\sigma_H$.
\vskip.1cm
{\bf (ii)} If $0 \le - d < c \sigma_h + h$, then there is one non-positive root $x_2 \ge -d/c$, one positive root smaller than $\sigma_H$, and either zero or two positive roots greater than $\sigma_h$.
\vskip.1cm
{\bf (iii)} If $c \sigma_h + h \le - d < c \sigma_H + H$ (pictured), then there is one negative root $x_4 > -d/c$, a positive root $x_3 < \sigma_H$, another positive root $x_2$ such that max$\{ \sigma_H, -(d+h)/c \} < x_2 \le $ min$\{ \sigma_h, -(d+H)/c \}$, and a third positive root $x_1 > \sigma_h$.
\vskip.1cm
{\bf (iv)} If $c \sigma_H + H \le -d$, then there is one negative root $x_2 > -d/c$, one positive root $x_1 > \sigma_h$ and either zero or two positive roots smaller than or equal to $\sigma_H$.\\
\end{minipage}}
& \multicolumn{1}{|c|}{\begin{minipage}{22em}
\scriptsize
\vskip-3.8cm
\begin{center}
\underline{{\color{red}{\it Notes}}}
\end{center}
\vskip-.4cm
As $c_2 < c_1 < c$, the quartic has a single stationary point $\mu_1$ and the number of real roots can be either 0 or 2. There is only one tangent $- c x - \delta_1$ to $x^4 + ax^3 + bx^2$ --- the one at $\mu_1$.
\vskip.1cm
As $\gamma_2 < \gamma_1 < c$, the straight line $-cx$ intersects $x^4 + ax^3 + bx^2$ at one point only ($\lambda_1$).
\vskip.1cm
Obviously, $-\delta_1 < 0 < c \sigma_H + H < c \sigma_h + h$.
\vskip.1cm
Consideration of whether $- d \le c \sigma_h + h$ or $-d > c \sigma_h + h$ and also whether $- d \le c \sigma_H + H$ or $-d > c \sigma_H + H$ would yield sharper bounds in the analysis based on solving cubic equations. \\
\hrule
\vskip.2cm
\hskip1.54cm \underline{{\color{red}{\it Analysis based on solving cubic equations}}}
\vskip.1cm
{\bf (i)} If $-d < -\delta_1$, then there are no real roots.
\vskip.1cm
{\bf (ii)} If $-\delta_1 \le -d < 0$, then there is one negative root $x_2$ such that $\lambda_1 < x_2 \le \mu_1$ and another negative root $x_1$ such that $\mu_1 \le x_1 < -d/c$.
\vskip.1cm
{\bf (iii)}  If $0 \le - d$ (pictured), then there is a negative root $x_2 \le \lambda_1$ and a non-negative root $x_1 \le -d/c$. \\
\hrule
\vskip.2cm
\hskip.88cm \underline{{\color{red}{\it Analysis based on solving quadratic equations only}}}
\vskip.2cm
{\bf (i)} If $- d < 0$, then there are either no real roots or there are two negative roots smaller than $-d/c$.
\vskip.1cm
{\bf (ii)} If $0 \le - d < c \sigma_H + H$ (pictured), then there is a non-positive root $x_2$ and a non-negative root $x_1$ such that max$\{0, -(d+H)/c \} \le x_1 \le $ min$\{ -d/c, \sigma_H \}$.
\vskip.1cm
{\bf (iii)} If $c \sigma_H + H \le - d < c \sigma_h + h$, then there is one negative root $x_2$ and one positive root $x_1$ such that $\sigma_H \le x_1 < \sigma_h$.
\vskip.1cm
{\bf (iv)} If $c \sigma_h + h \le -d$, then there is a negative root $x_2$ and a positive root $x_1$ such that  $\sigma_h \le x_1 < -(d+h)/c$. \\
\end{minipage}}  \\
\hline
\end{tabular}
\end{center}


\begin{center}
\begin{tabular}{|c||c|}
\hline
\multicolumn{2}{|c|}{\footnotesize \bf Figure 3.9}  \\
\multicolumn{2}{|c|}{\footnotesize \bf (Figure 3.8 is on next page)}  \\
\hline
\hline
\multicolumn{2}{|c|}{{\color{violet} $\bm{\frac{a^2}{4} < b \le \frac{9}{8}\,\frac{a^2}{4}\, , \quad a > 0, \quad c < 0,$}}}  \\
\multicolumn{2}{|c|}{{\color{violet} $\bm{c_2 < {\color{red}c} < 0 < \gamma_2  < c_0 < \gamma_1 < c_1 $}}}  \\
\hline
\hline
\includegraphics[width=80mm]{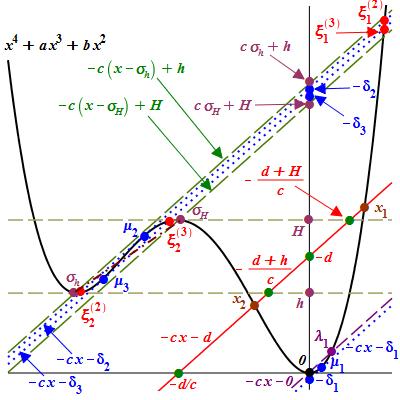} & \includegraphics[width=80mm]{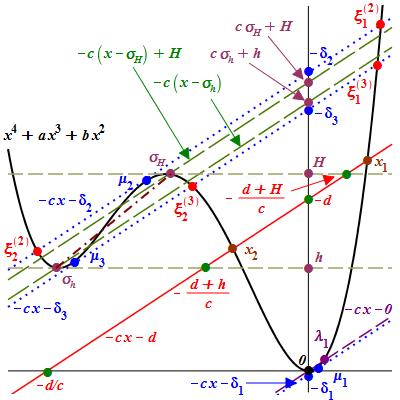} \\
\hline
\multicolumn{2}{|c|}{\scriptsize \underline{{\color{red}{\it Notes (apply to both panes)}}}} \\
\multicolumn{2}{|c|}{\begin{minipage}{45em}
\scriptsize
\vskip.2cm
As $c_2 < c < c_1$, the quartic has three stationary points $\mu_i$ and the number of real roots can be either 0, or 2, or 4. There are three tangents $- c x - \delta_i$ to $x^4 + ax^3 + bx^2$ --- each at $\mu_i$.
\vskip.1cm
As $c < \gamma_2 < \gamma_1$, the straight line $-cx$ intersects $x^4 + ax^3 + bx^2$ at one point only ($\lambda_{1}$). Thus $-\delta_2 > 0$ and $-\delta_3 > 0$. Also, $-\delta_1 < 0$.
\vskip.1cm
For both panes, the straight line joining the two stationary points in the second quadrant has slope $-\hat{c} = (1/2)a[(9/8)(a^2/4)-b] \ge 0$ such that $c_2 < \hat{c} \le 0$.
\vskip.1cm
With the increase of $c$ from $c_2$ towards zero, one has: $-\delta_1 < 0 < c \sigma_H + H \le -\delta_3 < -\delta_2 \le c \sigma_h + h$ (pictured on left pane; $c_2 < c < \hat{c} < 0$), then $-\delta_1 < 0 < c \sigma_H + H  \le -\delta_3 <  c \sigma_h + h \le -\delta_2$, or $-\delta_1 < 0 < -\delta_3 \le c \sigma_H + H < -\delta_2 \le c \sigma_h + h$, then $-\delta_1 < 0 < -\delta_3 \le c \sigma_H + H \le c \sigma_h + h \le -\delta_2$ (at $c = \hat{c}$ there is a swap between $c \sigma_H + H$ and $c \sigma_h + h$), and finally $-\delta_1 < 0 < -\delta_3 < c \sigma_h + h < c \sigma_H + H < -\delta_2$ (pictured on right pane; $c_2 < \hat{c} < c < 0$).
\vskip.1cm
Consideration of whether $- d \le c \sigma_h + h$ or $-d > c \sigma_h + h$ and also whether $- d \le c \sigma_H + H$ or $-d > c \sigma_H + H$ would yield sharper bounds in the analysis based on solving cubic equations. \\
\hrule
\begin{center}
\vskip-.3cm
\underline{{\color{red}{\it Analysis based on solving cubic equations (applies to both panes)}}} \\
\end{center}
\vskip-.3cm
{\bf (i)} If $-d < -\delta_1$, then there are no real roots.
\vskip.1cm
{\bf (ii)} If $-\delta_1 \le -d < 0$, then there is a positive root $x_2$ such that $-d/c < x_2 \le \mu_1$ and another positive root $x_1$ such that $\mu_1 \le x_1 < \lambda_1$.
\vskip.1cm
{\bf (iii)} If $0 \le -d < -\delta_3$ (pictured), then there is a non-positive root $x_2 \ge $ max$\{\xi^{(3)}_2, -d/c \}$ and a positive root $x_1$ such that $\lambda_1 \le x_1 < \xi^{(3)}_1$.
\vskip.1cm
{\bf (iv)}  If $-\delta_3 \le -d < -\delta_2$, then there is a negative root $x_4$ such that  $\xi^{(2)}_2 < x_4 \le \mu_3$, another negative root $x_3$ such that $\mu_3 \le x_3 < \mu_2,$ a third negative root $x_2$ such that $\mu_2 < x_2 \le \xi^{(3)}_2$, and a positive root $x_1$ such that $\xi^{(3)}_1 \le x_1 < \xi^{(2)}_1$.
\vskip.1cm
{\bf (v)}  If $-\delta_2 \le - d $, then there is a negative root $x_2  \le \xi^{(2)}_2$ and a positive root $x_1 \ge \xi^{(2)}_1$.
\vskip.1cm
\hrule
\end{minipage}}\\
\multicolumn{1}{|c|}{\begin{minipage}{22em}
\scriptsize
\begin{center}
\vskip0.2cm
\underline{{\color{red}{\it Analysis based on solving quadratic equations only}}}
\vskip.2cm
\underline{$c_2 < c < \hat{c} < 0$}
\end{center}
\vskip-.2cm
{\bf (i)} If $- d < 0$, then there are either no real roots or there are two positive roots greater than $-d/c$.
\vskip.1cm
{\bf (ii)} If $0 \le - d < c \sigma_H + H$ (pictured), then there is a non-negative root $x_1$, a non-positive root $x_2$ such that max$\{ -d/c, \sigma_H \} < x_2 \le $ min$\{0, -(d+H)/c \}$, and either zero or two positive roots greater than $\sigma_h$ and smaller than $\sigma_H$ (the latter appear when $-\delta_3 \le c \sigma_H + H$).
\vskip.1cm
{\bf (iii)} If $c \sigma_H + H \le - d < c \sigma_h + h$, then, if $c \sigma_H + H < -\delta_3$  or $c \sigma_h + h > -\delta_2$, there will be one positive root together with either one or three negative roots greater than $\sigma_h$ and smaller than or equal to $\sigma_H$. If however,$-\delta_3 \le c \sigma_H + H$ and $c \sigma_h + h \le -\delta_2$, there will be one positive root together with three negative roots greater than $\sigma_h$ and smaller than or equal to $\sigma_H$.
\vskip.1cm
{\bf (iv)} If $c \sigma_h + h \le -d$, then there is one negative root $x_2$ such that $-(d+h)/c < x_2 \le \sigma_h$, one positive root $x_1$ and either zero or two positive roots smaller than or equal to $-(d+H)/c$ and greater than $\sigma_h$ (the latter appear when $c \sigma_h + h < -\delta_2$).
\\
\end{minipage}}
& \multicolumn{1}{|c|}{\begin{minipage}{22em}
\scriptsize
\begin{center}
\vskip-1.15cm
\underline{{\color{red}{\it Analysis based on solving quadratic equations only}}}
\vskip.2cm
\underline{$c_2 < \hat{c} < c < 0$}
\end{center}
\vskip-.2cm
{\bf (i)} If $- d < 0$, then there are either no real roots or there are two positive roots greater than $-d/c$.
\vskip.1cm
{\bf (ii)} If $0 \le - d < c \sigma_h + h$ (pictured), then there is a non-negative root $x_1$, a non-positive root $x_2$ such that
max$\{ -d/c, \sigma_h \} < x_1 \le $ min$\{0, -(d+H)/c \}$ and either zero or two positive roots greater than or equal to $\sigma_h$ and smaller than $\sigma_H$.
\vskip.1cm
{\bf (iii)} If $c \sigma_h + h \le - d < c \sigma_H + H$, then there is one positive root $x_1$, a negative root $x_2$ such that $\sigma_H < x_2 < -(d+H)/c $, another negative root $x_3$ such that $\sigma_h \le x_3 < \sigma_H$, and a third negative root $x_4$ such that $-(d+h)/c <  x_4 \le \sigma_h$.
\vskip.1cm
{\bf (iv)} If $c \sigma_H + H \le -d$, then there is a positive root $x_1$, a negative root $x_2$ such that  $-(d+h)/c < x_1 \le \sigma_h$ and either zero or two positive roots smaller than or equal to $\sigma_H$ and greater than $\sigma_h$. \\
\end{minipage}}  \\
\hline
\end{tabular}
\end{center}


\begin{center}
\begin{tabular}{|c||c|}
\hline
{\footnotesize \bf Figure 3.8} &  {\footnotesize \bf Figure 3.10} \\
\hline
\hline
{\color{violet} $\bm{\frac{a^2}{4} < b \le \frac{9}{8}\,\frac{a^2}{4}\, , \quad a > 0, \quad c < 0,$}} &  {\color{violet} $\bm{\frac{a^2}{4} < b \le \frac{9}{8}\,\frac{a^2}{4}\, , \quad a > 0, \quad c > 0,$}} \\
{\color{violet} $\bm{{\color{red}c} < c_2 < 0 < \gamma_2  < c_0 < \gamma_1 < c_1 $}} & {\color{violet} $\bm{c_2 < 0 < {\color{red}c} < \gamma_2 < c_0 < \gamma_1 < c_1 $}} \\
\hline
\hline
\includegraphics[width=80mm]{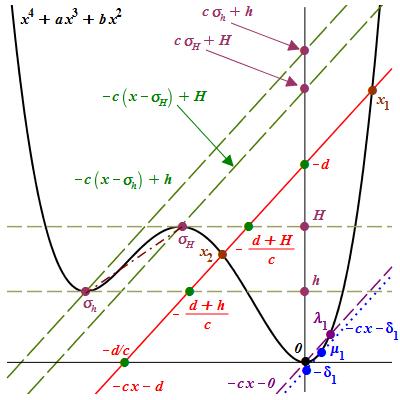} & \includegraphics[width=80mm]{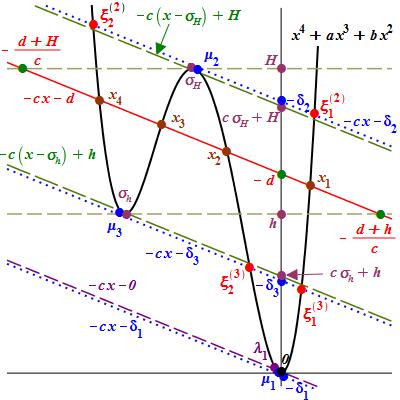} \\
\hline
\multicolumn{1}{|c||}{\begin{minipage}{22em}
\scriptsize
\vskip-3.05cm
\begin{center}
\underline{{\color{red}{\it Notes}}}
\end{center}
\vskip-.3cm
As $c < c_2 < c_1$, the quartic has a single stationary point $\mu_1$ and the number of real roots can be either 0 or 2. There is only one tangent $- c x - \delta_1$ to $x^4 + ax^3 + bx^2$ --- the one at $\mu_1$.
\vskip.1cm
As $c < \gamma_2 < \gamma_1$, the straight line $-cx$ intersects $x^4 + ax^3 + bx^2$ at one point only ($\lambda_1$).
\vskip.1cm
Obviously, $-\delta_1 < 0 < c \sigma_H + H < c \sigma_h + h$.
\vskip.1cm
Consideration of whether $- d \le c \sigma_h + h$ or $-d > c \sigma_h + h$ and also whether $- d \le c \sigma_H + H$ or $-d > c \sigma_H + H$ would yield sharper bounds in the analysis based on solving cubic equations. \\
\hrule
\vskip.3cm
\hskip1.54cm \underline{{\color{red}{\it Analysis based on solving cubic equations}}}
\vskip.3cm
{\bf (i)} If $-d < -\delta_1$, then there are no real roots.
\vskip.1cm
{\bf (ii)} If $-\delta_1 \le -d < 0$, then there is one positive root $x_2$ such that $-d/c < x_1 \le \mu_1$ and another positive root $x_1$ such that $\mu_1 \le x_2 < \lambda_1$.
\vskip.1cm
{\bf (iii)}  If $0 \le - d$ (pictured), then there is a non-positive root $x_2 \ge -d/c$ and a positive root $x_1 \ge \lambda_1$. \\
\hrule
\vskip.3cm
\hskip.88cm \underline{{\color{red}{\it Analysis based on solving quadratic equations only}}}
\vskip.3cm
{\bf (i)} If $- d < 0$, then there are either no real roots or there are two positive roots greater than $-d/c$.
\vskip.1cm
{\bf (ii)} If $0 \le - d < c \sigma_H + H$ (pictured), then there is a non-negative root $x_1$ and a non-positive root $x_2$ such that max$\{ -d/c, \sigma_H \} < x_2 \le $ min$\{0, -(d+H)/c \}$.
\vskip.1cm
{\bf (iii)} If $c \sigma_H + H \le - d < c \sigma_h + h$, then there is one positive root $x_1$ and one negative root $x_2$ such that $\sigma_h < x_1 \le \sigma_H$.
\vskip.1cm
{\bf (iv)} If $c \sigma_h + h \le -d$, then there is a positive root $x_1$ and a negative root $x_2$ such that  $-(d+h)/c < x_2 \le \sigma_h$. \\
\end{minipage}}
& \multicolumn{1}{|c|}{\begin{minipage}{22em}
\scriptsize
\vskip.3cm
\begin{center}
\underline{{\color{red}{\it Notes}}}
\end{center}
\vskip-.2cm
As $c_2 < c < c_1$, the quartic has three stationary points $\mu_i$ and the number of real roots can be either 0, or 2, or 4. There are three tangents $- c x - \delta_i$ to $x^4 + ax^3 + bx^2$ --- each at $\mu_i$.
\vskip.1cm
As $c < \gamma_2 < \gamma_1$, the straight line $-cx$ intersects $x^4 + ax^3 + bx^2$ at one point only ($\lambda_{1}$). Thus $-\delta_2 > 0$ and $-\delta_3 > 0$. Also: $-\delta_1 < 0$.
\vskip.1cm
Obviously, $-\delta_1 < 0 < -\delta_3 < c \sigma_h + h < c \sigma_H + H < -\delta_2$.
\vskip.1cm
Consideration of whether $- d \le c \sigma_h + h$ or $-d > c \sigma_h + h$ and also whether $- d \le c \sigma_H + H$ or $-d > c \sigma_H + H$ would yield sharper bounds in the analysis based on solving cubic equations.
\\
\hrule
\vskip.2cm
\hskip1.54cm \underline{{\color{red}{\it Analysis based on solving cubic equations}}} \\
\vskip-.1cm
{\bf (i)} If $-d < -\delta_1$, then there are no real roots.
\vskip.1cm
{\bf (ii)} If $-\delta_1 \le -d < 0$, then there is one negative root $x_2$ such that $-\lambda_1 < x_2 \le \mu_1$ and another negative root $x_1$ such that $\mu_1 \le x_1 < -d/c$.
\vskip.1cm
{\bf (iii)} If $0 \le -d < -\delta_3$, then there is one negative root $x_2$ such that $\xi^{(3)}_2 < x_2 \le \lambda_1$ and one non-negative root $x_1 \le $ min$\{ -d/c, \xi^{(3)}_1 \}.$
\vskip.1cm
{\bf (iv)}  If $-\delta_3 \le -d < -\delta_2$ (pictured), then there is one negative root $x_4$ such that $\xi^{(2)}_2 < x_4 \le \mu_3$, another negative root $x_3$ such that $\mu_3 \le x_3 < \mu_2$, a third negative root $x_2$ such that $\mu_2 < x_2 \le \xi^{(3)}_2$, and a positive root $x_1$ such that $\xi^{(3)}_1 \le x_1 < $ min$\{\xi^{(2)}_1, -d/c \}$.
\vskip.1cm
{\bf (v)}  If $-\delta_2 \le - d $, then there is a negative root $x_2 \le \xi^{(2)}_2$ and a positive root $x_1$ such that $\xi^{(2)}_1 \le x_2 < -d/c.$
\\
\hrule
\vskip.2cm
\hskip.88cm \underline{{\color{red}{\it Analysis based on solving quadratic equations only}}} \\
\vskip-.1cm
{\bf (i)} If $- d < 0$, then there are either no real roots or there are two negative roots smaller than $-d/c$ and greater than $\sigma_H$.
\vskip.1cm
{\bf (ii)} If $0 \le - d < c \sigma_h + h$, then there is one non-negative root $x_1 \le -d/c$, one non-positive root greater than $\sigma_H$, and either zero or two negative roots smaller than $\sigma_h$.
\vskip.1cm
{\bf (iii)} If $c \sigma_h + h \le - d < c \sigma_H + H$ (pictured), then there is one positive root $x_1 < -d/c$, a negative root $x_2 > \sigma_H$, another negative root $x_3$ such that max$\{ \sigma_h, -(d+H)/c \} \le x_3 < $ min$\{ \sigma_H, -(d+h)/c \}$, and a third negative root $x_4 \le \sigma_h$.
\vskip.1cm
{\bf (iv)} If $c \sigma_H + H \le -d$, then there is one positive root $x_1 < -d/c$, one negative root $x_2 < \sigma_h$ and either zero or two negative roots greater than or equal to $\sigma_H$. \\
\end{minipage}}  \\
\hline
\end{tabular}
\end{center}


\begin{center}
\begin{tabular}{|c||c|}
\hline
\multicolumn{2}{|c|}{\footnotesize \bf Figure 3.11}  \\
\hline
\hline
\multicolumn{2}{|c|}{{\color{violet} $\bm{\frac{a^2}{4} < b \le \frac{9}{8}\,\frac{a^2}{4}\, , \quad a > 0 \quad c > 0,$}}} \\
\multicolumn{2}{|c|}{{\color{violet} $\bm{c_2 < 0 < \gamma_2  < {\color{red}c} < c_0 < \gamma_1 < c_1$}}}  \\
\hline
\hline
\includegraphics[width=80mm]{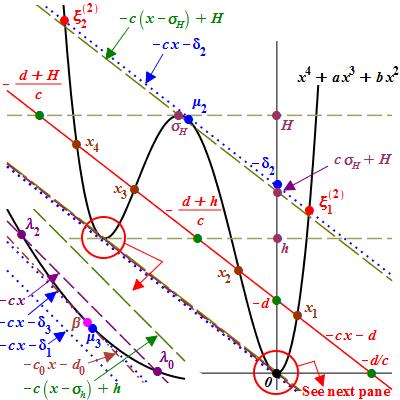} & \includegraphics[width=80mm]{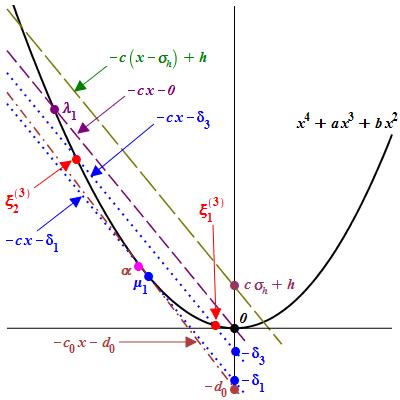} \\
\hline
\multicolumn{2}{|c|}{
\begin{minipage}{45em}
\scriptsize
\scriptsize
\vskip.3cm
\begin{center}
\underline{{\color{red}{\it Notes}}}
\end{center}
\vskip-.3cm
As $c_2 < c < c_1$, the quartic has three stationary points $\mu_i$ and the number of real roots can be either 0, or 2, or 4. There are three tangents $- c x - \delta_i$ to $x^4 + ax^3 + bx^2$ --- each at $\mu_i$ (only shown on the top-left pane).
\vskip.1cm
As $\gamma_2 < c < \gamma_1$, the straight line $-cx$ intersects $x^4 + ax^3 + bx^2$ at three points: $\lambda_{0,1,2}$. Thus $-\delta_2 > 0$ and $-\delta_3 < 0$. Also: $-\delta_1 < 0$.
\vskip.1cm
The straight line $ - c_0 x - d_0$ with $c_0 = (1/2)a(b-a^2/4) > 0$ and $d_0 = (1/4)(b-a^2/4)^2 > 0$ is tangent to $x^4 + ax^3 + bx^2$ at two points: $\alpha$ and $\beta = -\alpha - a/2$, where $\alpha$ and $\beta$ are the simultaneous double roots of the varied quartic $x^4 + ax^3 + bx^2 + c_0 x + d_0$.
\vskip.1cm
As  $c < c_0$, one has $-\delta_1 < - \delta_3$.
\vskip.1cm
Obviously, $-\delta_1 < -\delta_3 < 0 < c \sigma_h + h <  c \sigma_H + H < -\delta_2$.
\vskip.1cm
Consideration of whether $- d \le c \sigma_h + h$ or $-d > c \sigma_h + h$ and also whether $- d \le c \sigma_H + H$ or $-d > c \sigma_H + H$ would yield sharper bounds in the analysis based on solving cubic equations. \\
\hrule
\begin{center}
\vskip-.2cm
\underline{{\color{red}{\it Analysis based on solving cubic equations}}}
\end{center}
\vskip-.3cm
{\bf (i)} If $-d < -\delta_1$, then there are no real roots.
\vskip.1cm
{\bf (ii)} If $-\delta_1 \le -d < -\delta_3$, then there is one negative root $x_2$ such that $-\xi^{(3)}_2 < x_2 \le \mu_1$ and another negative root $x_1$ such that $\mu_1 \le x_1 < $ min$\{-d/c, \xi^{(3)}_1 \}$.
\vskip.1cm
{\bf (iii)} If $-\delta_3 \le -d < 0$, then there are four negative roots: $x_4$ such that $\lambda_2 < x_4 \le \mu_3, \,\, x_3$ such that $\mu_3 \le x_3 < \lambda_0, \,\, x_2$ such that $\lambda_1 < x_2 \le \xi^{(3)}_2$, and $x_1$ such that $\xi^{(3)}_1 \le x_1 < -d/c.$
\vskip.1cm
{\bf (iv)}  If $0 \le -d < -\delta_2$ (pictured), then there are three negative roots: $x_4$ such that $\xi^{(2)}_2 < x_4 \le \lambda_2, \,\, x_3$ such that $\lambda_0 \le x_3 < \mu_2$, and $x_2$ such that $\mu_2 < x_2 \le \lambda_1$, together with the non-negative root $x_1 \le $ min$\{-d/c, \xi^{(2)}_1 \}$.
\vskip.1cm
{\bf (v)}  If $-\delta_2 \le - d $, then there is a negative root $x_2 \le \xi^{(2)}_2$ and a positive root $x_1$ such that $\xi^{(2)}_1 \le x_2 < -d/c.$
\vskip.1cm
\hrule
\begin{center}
\vskip-.2cm
\underline{{\color{red}{\it Analysis based on solving quadratic equations only}}}
\end{center}
\vskip-.3cm
{\bf (i)} If $- d < 0$, then there are either no real roots, or there are two negative roots smaller than $-d/c$ and greater than $\sigma_H$, or there are four negative roots --- two smaller than $-d/c$ and greater than $\sigma_H$ and the other two --- smaller than $\sigma_h$.
\vskip.1cm
{\bf (ii)} If $0 \le - d < c \sigma_h + h$, then there is one non-negative root $x_1 \le -d/c$, one non-positive root greater than $\sigma_H$, and two negative roots smaller than $\sigma_h$.
\vskip.1cm
{\bf (iii)} If $c \sigma_h + h \le - d < c \sigma_H + H$ (pictured), then there is one positive root $x_1 < -d/c$, a negative root $x_2 > \sigma_H$, another negative root $x_3$ such that max$\{ \sigma_h, -(d+H)/c \} \le x_3 < $ min$\{ \sigma_H, -(d+h)/c \}$, and a third negative root $x_4 \le \sigma_h$.
\vskip.1cm
{\bf (iv)} If $c \sigma_H + H \le -d$, then there is one positive root $x_1 < -d/c$, one negative root $x_2 < \sigma_h$ and either zero or two negative roots greater than or equal to $\sigma_H$. \\
\end{minipage}} \\
\hline
\end{tabular}
\end{center}


\begin{center}
\begin{tabular}{|c|c|}
\hline
\multicolumn{2}{|c|}{\footnotesize \bf Figure 3.12}  \\
\multicolumn{2}{|c|}{\footnotesize \bf (continues on next page)}  \\
\hline
\hline
\multicolumn{2}{|c|}{{\color{violet} $\bm{\frac{a^2}{4} < b \le \frac{9}{8}\,\frac{a^2}{4}\, , \quad a > 0, \quad c > 0,$}}} \\
\multicolumn{2}{|c|}{{\color{violet} $\bm{c_2 < 0 < \gamma_2 < c_0 < {\color{red}c} < \gamma_1 < c_1 $}}} \\
\hline
\hline
\includegraphics[width=80mm,height=80mm]{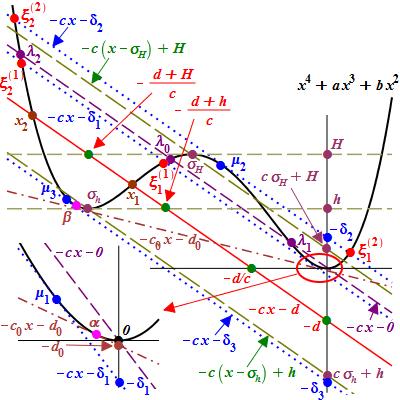} & \includegraphics[width=80mm,height=80mm]{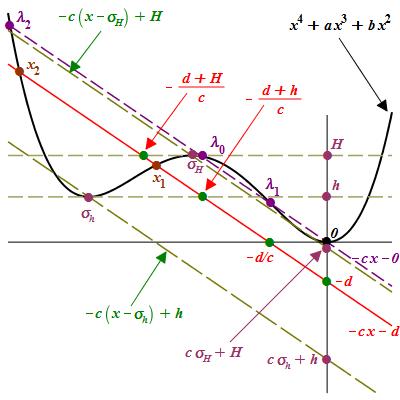} \\
\hline
\multicolumn{1}{|c|}{\begin{minipage}{22em}
\scriptsize
\begin{center}
\vskip-7.3cm
\underline{{\color{red}{\it Notes (apply to all panes)}}}
\end{center}
\vskip-.2cm
As $c_2 < c < c_1$, the quartic has three stationary points $\mu_i$ and the number of real roots can be either 0, or 2, or 4. There are three tangents $- c x - \delta_i$ to $x^4 + ax^3 + bx^2$ --- each at $\mu_i$ (only shown on the top-left pane).
\vskip.1cm
As $\gamma_2 < c < \gamma_1$, the straight line $-cx$ intersects $x^4 + ax^3 + bx^2$ at three points: $\lambda_{0,1,2}$. Thus $-\delta_2 > 0$ and $-\delta_3 < 0$. Also: $-\delta_1 < 0$.
\vskip.1cm
The straight line $ - c_0 x - d_0$ with $c_0 = (1/2)a(b-a^2/4) < 0$ and $d_0 = (1/4)(b-a^2/4)^2 > 0$ is tangent to $x^4 + ax^3 + bx^2$ at two points: $\alpha$ and $\beta = -\alpha - a/2$, where $\alpha$ and $\beta$ are the simultaneous double roots of the varied quartic $x^4 + ax^3 + bx^2 + c_0 x + d_0$.
\vskip.1cm
As $c > c_0$, one has $- \delta_3 < -\delta_1$.
\vskip.1cm
One can have $c \sigma_h + h < 0 < c \sigma_H + H$ (top-left pane), $c \sigma_h + h < c \sigma_H + H < 0$ (top-right pane; $c$ closer to $\gamma_1$), and $0 < c \sigma_h + h < c \sigma_H + H$ (bottom-right pane; $c$ closer to $c_0$). For each pane, there are further sub-cases affecting the analysis based on solving quadratic equations only (see below).
\vskip.1cm
Consideration of whether $- d \le c \sigma_h + h$ or $-d > c \sigma_h + h$ and also whether $- d \le c \sigma_H + H$ or $-d > c \sigma_H + H$ would yield sharper bounds in the analysis based on solving cubic equations.
\\
\end{minipage}}
 & \includegraphics[width=80mm,height=80mm]{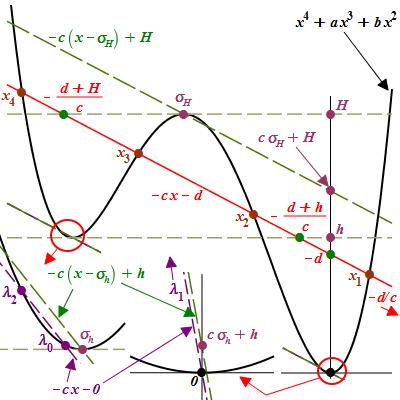} \\
\hline
\multicolumn{2}{|c|}{\begin{minipage}{45em}
\scriptsize
\vskip.1cm
\begin{center}
\underline{{\color{red}{\it Analysis based on solving cubic equations (applies to all panes)}}}
\end{center}
\vskip-.3cm
{\bf (i)} If $-d < -\delta_3$, then there are no real roots.
\vskip.1cm
{\bf (ii)} If $-\delta_3 \le -d < -\delta_1$ (pictured), then there is one negative root $x_2$ such that $\xi^{(1)}_2 < x_2 \le \mu_3$ and another negative root $x_1$ such that max$\{ \mu_3, -(d+H)/c \} \le x_2 < $ min$\{ \xi^{(1)}_1, -(d+h)/c \}$.
\vskip.1cm
{\bf (iii)} If $-\delta_1 \le -d < 0$, then there are four negative roots: $x_4$ such that $\lambda_2 < x_1  \le \xi^{(1)}_2, \,\, x_3$ such that
$\xi^{(1)}_1 \le x_2 < \lambda_0, \,\, x_2$ such that $\lambda_1 < x_2 \le \mu_1$, and $x_1$ such that $\mu_1 \le x_1 < - d/c.$
\vskip.1cm
{\bf (iv)}  If $0 \le -d < -\delta_2$, then there are tree negative roots: $x_4$ such that $\xi^{(2)}_2 < x_4 \le \lambda_2, \,\,$ $x_3$ such that $\lambda_0 \le x_3 < \mu_2$, and $x_2$ such that $\mu_2 < x_2 \le \lambda_1$, together with the non-negative root $x_1 \le $ min$\{ \xi^{(2)}_1, -d/c \}$.
\vskip.1cm
{\bf (v)}  If $-\delta_2 \le - d $, then there is a negative root $x_2 \le \xi^{(2)}_2$ and a positive root $x_1$ such that $\xi^{(2)}_1 \le x_2 -d/c$. \\
\hrule
\begin{center}
\vskip-.3cm
\underline{{\color{red}{\it Analysis based on solving quadratic equations only (continues on next page)}}}
\end{center}
\vskip-.3cm
{\bf \underline{\color{red}{Top-left pane}}}
\vskip.2cm
One can have either $-\delta_2 > c \sigma_H + H > 0 > -\delta_1 > c \sigma_h + h > -\delta_3$ (pictured) or $-\delta_2 > c \sigma_H + H > 0 > c \sigma_h + h \ge -\delta_1 > -\delta_3$.
\vskip.2cm
{\bf (i)} If $- d < c \sigma_h + h$, then there are either no real roots, or there are two negative roots smaller than $\sigma_h$, or there are two negative roots smaller than $\sigma_h$ and another pair of negative roots smaller than $-d/c$ and greater than $\sigma_H$ (the latter appear when $c \sigma_h + h \ge -\delta_1$).
\vskip.1cm
{\bf (ii)} If $c \sigma_h + h \le - d < 0$ (pictured), then there is one negative root smaller than $\sigma_h$, another negative root greater than or equal to max$\{ \sigma_h, -(d+H)/c \}$ and smaller than min$\{ \sigma_H, -(d+h)/c \}$, and either zero or two negative roots smaller than $-d/c$ and greater than $\sigma_H$ (the latter are always present if $c \sigma_h + h \ge -\delta_1$, while for the pictured $c \sigma_h + h < -\delta_1$ they may or may not be there).
\vskip.1cm
{\bf (iii)} If $0 \le - d < c \sigma_H + H$, then there is one non-negative root $x_1 \le -d/c$, a negative root $x_2 > \sigma_H$, another negative root $x_3$ such that $-(d+H)/c < x_3 < \sigma_H$ , and a third negative root $x_4 < \sigma_h$.
\vskip.1cm
{\bf (iv)} If $c \sigma_H + H \le -d$, then there is one positive root $x_1 < -d/c$, a negative root $x_2 < \sigma_h$ and either zero or two negative roots greater than or equal to $\sigma_H$.\\
\end{minipage}}\\
\hline
\end{tabular}
\end{center}


\begin{center}
\begin{tabular}{|c|c|}
\hline
\multicolumn{2}{|c|}{\footnotesize \bf Figure 3.12}  \\
\multicolumn{2}{|c|}{\footnotesize \bf (continued from previous page)}  \\
\hline
\hline
\multicolumn{2}{|c|}{{\color{violet} $\bm{\frac{a^2}{4} < b \le \frac{9}{8}\,\frac{a^2}{4}\, , \quad a > 0, \quad c > 0,$}}} \\
\multicolumn{2}{|c|}{{\color{violet} $\bm{c_2 < 0 < \gamma_2 < c_0 < {\color{red}c} < \gamma_1 < c_1 $}}} \\
\hline
\hline
\multicolumn{2}{|c|}{\begin{minipage}{45em}
\scriptsize
\begin{center}
\vskip.2cm
\underline{{\color{red}{\it Analysis based on solving quadratic equations only --- continued from previous page}}}
\end{center}
{\bf \underline{\color{red}{Top-right pane}}}
\vskip.2cm
One can have either $-\delta_2 > 0 > -\delta_1 > c \sigma_H + H >  c \sigma_h + h > -\delta_3$, or $-\delta_2 > 0 > c \sigma_H + H \ge -\delta_1 > c \sigma_h + h > -\delta_3$, or $-\delta_2 > 0 > c \sigma_H + H > c \sigma_h + h \ge -\delta_1 > -\delta_3$ (when $\sigma_h$ and $\sigma_H$ are close).
\vskip.2cm
{\bf (i)} If $- d < c \sigma_h + h$, then there are either no real roots, or there are two negative roots smaller than $\sigma_h$, or there are two negative roots smaller than $\sigma_h$ and two negative roots smaller than $-d/c$ and greater than $\sigma_H$ (the latter appear when $c \sigma_h + h \ge -\delta_1$).
\vskip.1cm
{\bf (ii)} If $c \sigma_h + h \le - d < c \sigma_H + H$ (pictured), then there is a negative root $x_1$ such that max$\{ \sigma_h, -(d+H)/c \} \le x_1 < $ min$\{ \sigma_H, -(d+h)/c \}$ and a negative root $x_2 < \sigma_h$. If $-\delta_1 > c \sigma_H + H$, there are no other roots. If $c \sigma_h + h \ge -\delta_1$, there are two more negative roots smaller than $-d/c$ and greater than $\sigma_H$. If $c \sigma_H + H \ge -\delta_1 > c \sigma_h + h$, these two roots may or may not be there.
\vskip.1cm
{\bf (iii)} If $c \sigma_H + H \le - d < 0$, then there is a negative root smaller than $\sigma_h$ and, if $-\delta_1 > c \sigma_H + H$, there is either one or three negative roots smaller than $-d/c$ and greater than $\sigma_H$. If however, $c \sigma_H + H \ge -\delta_1$, then there will be a negative root smaller than $\sigma_h$ and three negative roots smaller than $-d/c$ and greater than $\sigma_H$.
\vskip.1cm
{\bf (iv)} If $0 \le -d$, then there is a non-negaive root smaller than or equal to $-d/c$, either zero or two negative roots greater than $\sigma_H$, and a negative root smaller than $\sigma_h$. \\
\hrule
\vskip.3cm
{\bf \underline{\color{red}{Bottom-right pane}}}
\vskip.2cm
One can only have $-\delta_2 >  c \sigma_H + H >  c \sigma_h + h > 0 > -\delta_1 > -\delta_3$.
\vskip.1cm
{\bf (i)} If $- d < 0$, then there are either no real roots, or there are two negative roots smaller than $\sigma_h$, or there are four negative roots: two smaller than $\sigma_h$ and two greater than $\sigma_H$ and smaller than $-d/c$.
\vskip.1cm
{\bf (ii)} If $0 \le - d < c \sigma_h + h$, then there is one non-negative root smaller than $-d/c$, one negative root greater than $\sigma_H$, and two negative roots smaller than $\sigma_h$.
\vskip.1cm
{\bf (iii)} If $c \sigma_h + h \le - d < c \sigma_H + H$ (pictured), then there is one positive root $x_1 < -d/c$, a negative root $x_2 > \sigma_H$, a negative root $x_3$ such that max$\{ -(d+H)/c, \sigma_h \} \le x_3 < $ min$\{\sigma_H, -(d+h)/c \}$  and another negative root $x_4 \le \sigma_h$.
\vskip.1cm
{\bf (iv)} If $c \sigma_H + H \le -d$, then there is a positive root smaller than $-d/c$, a negative root smaller than $\sigma_h$ and either zero or two negative roots greater than $\sigma_H$.
\\
\end{minipage}}\\
\hline
\end{tabular}
\end{center}


\begin{center}
\begin{tabular}{|c||c|}
\hline
{\footnotesize \bf Figure 3.13} &  {\footnotesize \bf Figure 3.14} \\
\hline
\hline
{\color{violet} $\bm{\frac{a^2}{4} < b \le \frac{9}{8}\,\frac{a^2}{4}\, , \quad a > 0, \quad c > 0,$}} &  {\color{violet} $\bm{\frac{a^2}{4} < b \le \frac{9}{8}\,\frac{a^2}{4}\, , \quad a > 0, \quad c > 0,$}} \\
{\color{violet} $\bm{c_2 < 0 < \gamma_2 <  c_0 < \gamma_1 < {\color{red}c} < c_1 $}} & {\color{violet} $\bm{c_2 < 0 < \gamma_2 < c_0  < \gamma_1 < c_1 < {\color{red}c} $}} \\
\hline
\hline
\includegraphics[width=80mm]{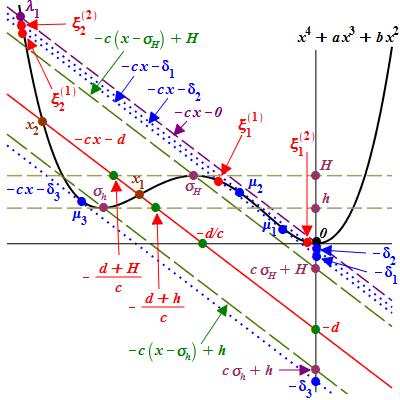} & \includegraphics[width=80mm]{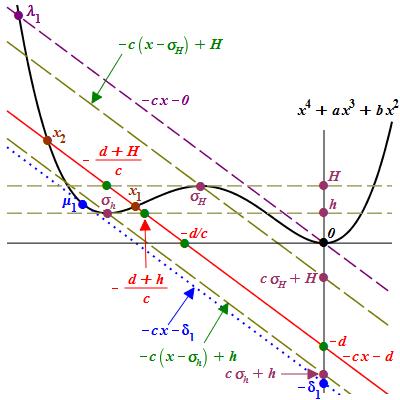} \\
\hline
\multicolumn{1}{|c||}{\begin{minipage}{22em}
\scriptsize
\vskip.3cm
\hskip3.89cm \underline{{\color{red}{\it Notes}}}
\vskip.3cm
As $c_2 < c < c_1$, the quartic has three stationary points $\mu_i$ and the number of real roots can be either 0, or 2, or 4. There are three tangents $- c x - \delta_i$ to $x^4 + ax^3 + bx^2$ --- each at $\mu_i$.
\vskip.1cm
As $\gamma_2 < \gamma_1 < c$, the straight line $-cx$ intersects $x^4 + ax^3 + bx^2$ at one point only ($\lambda_1$). Thus $-\delta_1 < 0$ and $-\delta_2 < 0$. Also:  $-\delta_3 < 0$.
\vskip.1cm
Obviously, $-\delta_3 < c \sigma_h + h < c \sigma_H + H < -\delta_1 < -\delta_2 < 0$.
\vskip.1cm
Consideration of whether $- d \le c \sigma_h + h$ or $-d > c \sigma_h + h$ and also whether $- d \le c \sigma_H + H$ or $-d > c \sigma_H + H$ would yield sharper bounds in the analysis based on solving cubic equations. \\
\hrule
\vskip.3cm
\hskip1.54cm \underline{{\color{red}{\it Analysis based on solving cubic equations}}}
\vskip.3cm
{\bf (i)} If $-d < -\delta_3$, then there are no real roots.
\vskip.1cm
{\bf (ii)} If $-\delta_3 \le -d < -\delta_1$ (pictured), then there is one negative root $x_2$ such that $\xi^{(1)}_2 < x_2 \le \mu_3$ and another negative root $x_1$ such that $\mu_3 \le x_1 < $ min$\{ \xi^{(1)}_1, -d/c \}.$
\vskip.1cm
{\bf (iii)}  If $-\delta_1 \le -d < -\delta_2$, then there are four negative roots: $x_4$ such that $\xi^{(2)}_2 < x_4 \le \xi^{(1)}_2, \,\, x_3$ such that $\xi^{(1)}_1 \le x_2 < \mu_2, \,\, x_2 $ such that $\mu_2 < x_2 \le \mu_1$, and $x_1$ such that $\mu_1 \le x_1 < $ min$\{\xi^{(2)}_1, -d/c \}$.
\vskip.1cm
{\bf (iv)} If $-\delta_2 \le -d < 0$, then there is a negative root $x_2$ such that  $\lambda_1 < x_2 \le \xi^{(2)}_2$ and another negative root $x_1$ such that $\xi^{(2)}_1 \le x_1 < -d/c$.
\vskip.1cm
{\bf (v)}  If $0 \le - d$, then there is a non-negative root $x_1 \le -d/c$ and a negative root $x_2 \le \lambda_1$. \\
\hrule
\vskip.3cm
\hskip.88cm \underline{{\color{red}{\it Analysis based on solving quadratic equations only}}}
\vskip.3cm
{\bf (i)} If $-d < c \sigma_h + h$, then there are either no real roots or there are two negative roots $x_{1,2} < \sigma_h$.
\vskip.1cm
{\bf (ii)} If $c \sigma_h + h \le - d < c \sigma_H + H$ (pictured), then there is a negative root $x_1$ such that max$\{ \sigma_h, -(d+H)/c \} \le x_1 < $ min$\{ \sigma_H, -(d+h)/c \}$ and a negative root $x_2 < \sigma_h$.
\vskip.1cm
{\bf (iii)} If $c \sigma_H + H \le - d < 0$, then there is a negative root smaller than $\sigma_h$ and either one or three negative roots smaller than $-d/c$ and greater than or equal to $\sigma_H$.
\vskip.1cm
{\bf (iv)} If $0 \le -d$, then there is one non-negative root $x_1 \le -d/c$ and one negative root $x_2 < \sigma_h$. \\
\end{minipage}}
& \multicolumn{1}{|c|}{\begin{minipage}{22em}
\scriptsize
\vskip-2.45cm
\hskip3.89cm \underline{{\color{red}{\it Notes}}}
\vskip.3cm
As $c_2 < c_1 < c$, the quartic has a single stationary point $\mu_1$ and the number of real roots can be either 0 or 2. There is only one tangent $- c x - \delta_1$ to $x^4 + ax^3 + bx^2$ --- the one at $\mu_1$.
\vskip.1cm
As $\gamma_2 < \gamma_1 < c$, the straight line $-cx$ intersects $x^4 + ax^3 + bx^2$ at one point only ($\lambda_1$).
\vskip.1cm
Obviously, $-\delta_1 < c \sigma_h + h < c \sigma_H + H < 0$.
\vskip.1cm
Consideration of whether $- d \le c \sigma_h + h$ or $-d > c \sigma_h + h$ and also whether $- d \le c \sigma_H + H$ or $-d > c \sigma_H + H$ would yield sharper bounds in the analysis based on solving cubic equations. \\
\hrule
\vskip.3cm
\hskip1.54cm \underline{{\color{red}{\it Analysis based on solving cubic equations}}}
\vskip.3cm
{\bf (i)} If $-d < -\delta_1$, then there are no real roots.
\vskip.1cm
{\bf (ii)} If $-\delta_1 \le -d < 0$ (pictured), then there is one negative root $x_2$ such that $\lambda_1 < x_2 \le \mu_1$ and another negative root $x_1$ such that $\mu_1 \le x_1 < -d/c$.
\vskip.1cm
{\bf (iii)}  If $0 \le - d$, then there is a non-negative root $x_1 \le -d/c$ and a negative root $x_2 \le \lambda_1$. \\
\hrule
\vskip.3cm
\hskip.88cm \underline{{\color{red}{\it Analysis based on solving quadratic equations only}}}
\vskip.3cm
{\bf (i)} If $-d < c \sigma_h + h$, then there are either no real roots or there are two negative roots $x_{1,2} < \sigma_h$.
\vskip.1cm
{\bf (ii)} If $c \sigma_h + h \le - d < c \sigma_H + H$ (pictured), then there is a negative root $x_1$ such that max$\{ \sigma_h, -(d+H)/c \} \le x_1 < $ min$\{ \sigma_H, -(d+h)/c \}$ and another negative root $x_2 < \sigma_h$.
\vskip.1cm
{\bf (iii)} If $c \sigma_H + H \le - d < 0$, then there is a negative root smaller than $\sigma_h$ and another negative root smaller than $-d/c$ and greater than or equal to $\sigma_H$.
\vskip.1cm
{\bf (iv)} If $0 \le -d$, then there is one non-negative root $x_1 \le -d/c$ and one negative root $x_2 < \sigma_h$. \\
\end{minipage}}  \\
\hline
\end{tabular}
\end{center}


\begin{center}
\begin{tabular}{|c||c|}
\hline
{\footnotesize \bf Figure 4.1} &  {\footnotesize \bf Figure 4.2} \\
\hline
\hline
{\color{violet} $\bm{\frac{9}{8}\,\frac{a^2}{4} < b \le \frac{4}{3}\,\frac{a^2}{4}\, , \quad a < 0, \quad c < 0,$}} &
{\color{violet} $\bm{\frac{9}{8}\,\frac{a^2}{4} < b \le \frac{4}{3}\,\frac{a^2}{4}\, , \quad a < 0, \quad c < 0,$}} \\
{\color{violet} $\bm{{\color{red}c} < c_2 < \gamma_2 < c_0 < \gamma_1 < c_1 < 0$}} & {\color{violet} $\bm{c_2 < {\color{red}c} < \gamma_2 < c_0 < \gamma_1 < c_1 < 0$}} \\
\hline
\hline
\includegraphics[width=80mm]{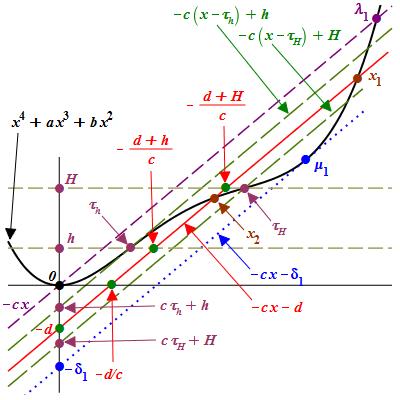} & \includegraphics[width=80mm]{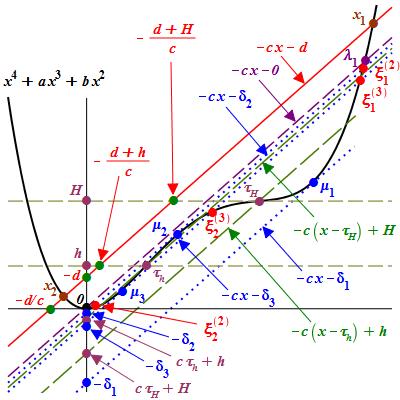} \\
\hline
\multicolumn{1}{|c||}{\begin{minipage}{22em}
\scriptsize
\vskip-3.8cm
\hskip3.89cm \underline{{\color{red}{\it Notes}}}
\vskip.1cm
As $c < c_2 < c_1$, the quartic has a single stationary point $\mu_1$ and the number of real roots can be either 0 or 2. There is only one tangent $- c x - \delta_1$ to $x^4 + ax^3 + bx^2$ --- the one at $\mu_1$.
\vskip.1cm
As $c < \gamma_2 < \gamma_1$, the straight line $-cx$ intersects $x^4 + ax^3 + bx^2$ at one point only ($\lambda_1$).
\vskip.1cm
The straight line joining the two points of curvature change ($\tau_h$ and $\tau_H$) has the same slope $-c_0 = -(1/2)a(b-a^2/4) > 0$ as that of the straight line $ - c_0 x - d_0$, tangent to $x^4 + ax^3 + bx^2$ at two points: $\alpha$ and $\beta = -\alpha - a/2$, where $\alpha$ and $\beta$ are the simultaneous double roots of the varied quartic $x^4 + ax^3 + bx^2 + c_0 x + d_0$ [with $d_0 = (1/4)(b-a^2/4)^2 > 0$].
\vskip.1cm
As $c < c_0 < 0$, one has $-\delta_1 < c \sigma_H + H < c \sigma_h + h < 0$.
\vskip.1cm
Consideration of whether $- d \le c \tau_h + h$ or $-d > c \tau_h + h$ and also whether $- d \le c \tau_H + H$ or $-d > c \tau_H + H$ would yield sharper bounds in the analysis based on solving cubic equations. \\
\hrule
\vskip.3cm
\hskip1.54cm \underline{{\color{red}{\it Analysis based on solving cubic equations}}}
\vskip.3cm
{\bf (i)} If $-d < -\delta_1$, then there are no real roots.
\vskip.1cm
{\bf (ii)} If $-\delta_1 \le -d < 0$ (pictured), then there is one positive root $x_2$ such that $-d/c < x_2 \le \mu_1$ and another positive root $x_1$ such that $\mu_1 \le x_2 < \lambda_1$.
\vskip.1cm
{\bf (iii)}  If $0 \le - d$, then there is a non-positive root $x_2 \ge -d/c$ and a positive root $x_1 \ge \lambda_1$. \\
\hrule
\vskip.3cm
\hskip.88cm \underline{{\color{red}{\it Analysis based on solving quadratic equations only}}}
\vskip.3cm
{\bf (i)} If $-d < c \tau_H + H$, then there are either no real roots or there are two positive roots $x_{1,2} > -(d+H)/c$.
\vskip.1cm
{\bf (ii)} If $c \tau_H + H \le - d < c \tau_h + h$ (pictured), then there are two positive roots: $x_2$ such that $-(d+h)/c < x_2 \le -(d+H)/c$ and $x_1 > \tau_H$.
\vskip.1cm
{\bf (iii)} If $c \tau_h + h \le - d < 0$, then there are two positive roots: $x_2$ such that $-d/c < x_2 \le -(d+h)/c$ and $x_1 > \tau_H$.
\vskip.1cm
{\bf (iv)} If $0 \le -d$, then there is one non-positive root $x_2 \ge -d/c$ and one positive root $x_1 > \tau_H$. \\
\end{minipage}}
& \multicolumn{1}{|c|}{\begin{minipage}{22em}
\scriptsize
\vskip.1cm
\hskip3.89cm \underline{{\color{red}{\it Notes}}}
\vskip.1cm
As $c_2 < c < c_1$, the quartic has three stationary points $\mu_i$ and the number of real roots can be either 0, or 2, or 4. There are three tangents $- c x - \delta_i$ to $x^4 + ax^3 + bx^2$ --- each at $\mu_i$.
\vskip.1cm
As $c < \gamma_2 < \gamma_1$, the straight line $-cx$ intersects $x^4 + ax^3 + bx^2$ at one point only ($\lambda_1$). Thus $-\delta_1 < 0$ and $-\delta_2 < 0$. Also:  $-\delta_3 < 0$.
\vskip.1cm
The straight line joining the two points of curvature change ($\tau_h$ and $\tau_H$) has the same slope $-c_0 = -(1/2)a(b-a^2/4) > 0$ as that of the straight line $ - c_0 x - d_0$, tangent to $x^4 + ax^3 + bx^2$ at two points: $\alpha$ and $\beta = -\alpha - a/2$, where $\alpha$ and $\beta$ are the simultaneous double roots of the varied quartic $x^4 + ax^3 + bx^2 + c_0 x + d_0$ [with $d_0 = (1/4)(b-a^2/4)^2 > 0$].
\vskip.1cm
As $c < c_0 < 0$, one has $-\delta_1 < c \sigma_H + H < -\delta_3 < c \sigma_h + h < -\delta_2 < 0$.
\vskip.1cm
Consideration of whether $- d \le c \tau_h + h$ or $-d > c \tau_h + h$ and also whether $- d \le c \tau_H + H$ or $-d > c \tau_H + H$ would yield sharper bounds in the analysis based on solving cubic equations. \\
\hrule
\vskip.1cm
\hskip1.54cm \underline{{\color{red}{\it Analysis based on solving cubic equations}}}
\vskip.2cm
{\bf (i)} If $-d < -\delta_1$, then there are no real roots.
\vskip.1cm
{\bf (ii)} If $-\delta_1 \le -d < -\delta_3$, then there is one positive root $x_2$ such that max$\{ \xi^{(3)}_2, -d/c \} < x_2 \le \mu_1$ and another positive root $x_1$ such that $\mu_1 \le x_1 < \xi^{(3)}_1$.
\vskip.1cm
{\bf (iii)}  If $-\delta_3 \le -d < -\delta_2$, then there are four positive roots: $x_4$ such that max$\{ \xi^{(2)}_2, -d/c \} < x_4 \le \mu_3, \,\,$ $x_3$ such that $\mu_3 \le x_3 < \mu_2, \,\,$ $x_2$ such that $\mu_2 < x_2 \le \xi^{(3)}_2$, and $x_1$ such that $\xi^{(3)}_1 \le x_1 < \xi^{(2)}_1$. \vskip.1cm
{\bf (iv)} If $-\delta_2 \le -d < 0$, then there is a positive root $x_2$ such that $-d/c < x_2 \le \xi^{(2)}_2$ and another positive root $x_1$ such that $\xi^{(2)}_1 \le x_1 < \lambda_1$.
\vskip.1cm
{\bf (v)}  If $0 \le - d$ (pictured), then there is a non-positive root $x_2 \ge -d/c$ and a positive root $x_1 \ge \lambda_1$. \\
\hrule
\vskip.1cm
\hskip.88cm \underline{{\color{red}{\it Analysis based on solving quadratic equations only}}}
\vskip.2cm
{\bf (i)} If $-d < c \tau_H + H$, then there are either no real roots or there are two positive roots $x_{1,2} > -(d+H)/c$.
\vskip.1cm
{\bf (ii)} If $c \tau_H + H \le - d < c \tau_h + h$, then there are two positive roots: $x_2$ such that $-(d+h)/c < x_2 \le -(d+H)/c$ and $x_1 > \tau_H$, together with either zero or two positive roots $x_{3,4}$ such that $-d/c < x_{3,4} < \tau_h$.
\vskip.1cm
{\bf (iii)} If $c \tau_h + h \le - d < 0$, then there are two positive roots: $x_2$ such that $-d/c < x_2 \le -(d+h)/c$ and $x_1 > \tau_H$, together with either zero or two positive roots greater than or equal to $\tau_h$ and smaller than $-(d+H)/c$.
\vskip.1cm
{\bf (iv)} If $0 \le -d$ (pictured), then there is one non-positive root $x_2 \ge -d/c$ and one positive root $x_1 > \tau_H$.  \\
\end{minipage}}
\\
\hline
\end{tabular}
\end{center}


\begin{center}
\begin{tabular}{|c|c|}
\hline
\multicolumn{2}{|c|}{\footnotesize \bf Figure 4.3}  \\
\multicolumn{2}{|c|}{\footnotesize \bf (continues on next page)}  \\
\hline
\hline
\multicolumn{2}{|c|}{{\color{violet} $\bm{\frac{9}{8}\,\frac{a^2}{4} < b \le \frac{4}{3}\,\frac{a^2}{4}\, , \quad a < 0, \quad c < 0,$}}} \\
\multicolumn{2}{|c|}{{\color{violet} $\bm{c_2 < \gamma_2 < {\color{red}c} < c_0 < \gamma_1 < c_1 < 0$}}} \\
\hline
\hline
\includegraphics[width=80mm,height=80mm]{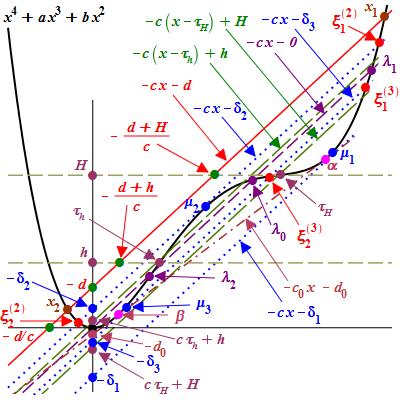} & \includegraphics[width=80mm,height=80mm]{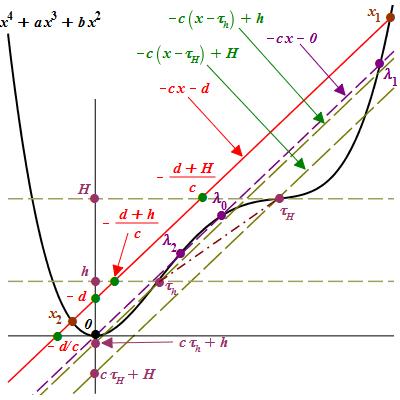} \\
\hline
& \\
& \\
\multicolumn{1}{|c|}{\begin{minipage}{22em}
\scriptsize
\begin{center}
\vskip-8.2cm
\underline{{\color{red}{\it Notes (apply to all panes)}}}
\end{center}
\vskip-.4cm
As $c_2 < c < c_1$, the quartic has three stationary points $\mu_i$ and the number of real roots can be either 0, or 2, or 4. There are three tangents $- c x - \delta_i$ to $x^4 + ax^3 + bx^2$ --- each at $\mu_i$ (only shown on the top-left pane).
\vskip.1cm
As $\gamma_2 < c < \gamma_1$, the straight line $-cx$ intersects $x^4 + ax^3 + bx^2$ at three points: $\lambda_{0,1,2}$. Thus $-\delta_2 > 0$ and $-\delta_1 < 0$. Also: $-\delta_3 < 0$.
\vskip.1cm
The straight line joining the two points of curvature change ($\tau_h$ and $\tau_H$) has the same slope $-c_0 = -(1/2)a(b-a^2/4) > 0$ as that of the straight line $ - c_0 x - d_0$, tangent to $x^4 + ax^3 + bx^2$ at two points: $\alpha$ and $\beta = -\alpha - a/2$, where $\alpha$ and $\beta$ are the simultaneous double roots of the varied quartic $x^4 + ax^3 + bx^2 + c_0 x + d_0$ [with $d_0 = (1/4)(b-a^2/4)^2 > 0$].
\vskip.1cm
As $c < c_0$, one has $- \delta_1 < -\delta_3$ and also $c \tau_H + H < c \tau_h + h$.
\vskip.1cm
\underline{Top-left pane}: one can have $-\delta_1 < c \tau_H + H < -\delta_3 < 0 < c \tau_h + h < -\delta_2$ (pictured) or $-\delta_1 < -\delta_3 \le c \tau_H + H < 0 < c \tau_h + h < -\delta_2$.
\vskip.1cm
\underline{Top-right pane} ($c$ closer to $\gamma_2$): one can have $-\delta_1 < c \tau_H + H < c \tau_h + h < -\delta_3 < 0 <  -\delta_2$, or $-\delta_1 < c \tau_H + H < -\delta_3 < c \tau_h + h < 0 < -\delta_2$, or $-\delta_1 < -\delta_3 < c \tau_H + H  < c \tau_h + h < 0 < -\delta_2$ (when $\tau_h$ and $\tau_H$ are close).
\vskip.1cm
\underline{Bottom-right pane} ($c$ closer to $c_0$): one can only have $-\delta_1 < -\delta_3 < 0 <  c \tau_H + H < c \tau_h + h < -\delta_2$.
\vskip.1cm
Consideration of whether $- d \le c \tau_h + h$ or $-d > c \tau_h + h$ and also whether $- d \le c \tau_H + H$ or $-d > c \tau_H + H$ would yield sharper bounds in the analysis based on solving cubic equations.
\\
\end{minipage}}
 & \includegraphics[width=80mm,height=80mm]{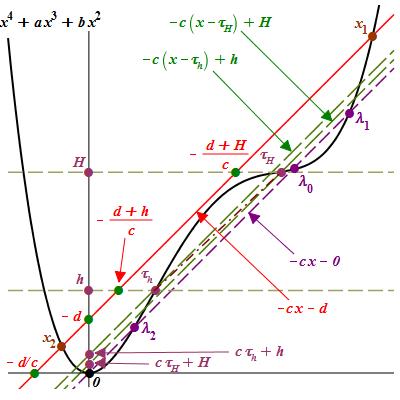} \\
\hline
\multicolumn{2}{|c|}{\begin{minipage}{45em}
\scriptsize
\vskip.1cm
\begin{center}
\underline{{\color{red}{\it Analysis based on solving cubic equations (applies to all panes)}}}
\end{center}
\vskip-.4cm
{\bf (i)} If $-d < -\delta_1$, then there are no real roots.
\vskip.1cm
{\bf (ii)} If $-\delta_1 \le -d < -\delta_3$, then there is one positive root $x_2$ such that $\xi^{(3)}_2  < x_2 \le \mu_1$ and another positive root $x_1$ such that $\mu_1 \le x_1 < \xi^{(3)}_1$.
\vskip.05cm
{\bf (iii)} If $-\delta_3 \le -d < 0$, then there are four positive roots: $x_4$ such that $- d/c < x_4 \le \mu_3, \,\,$ $x_3$ such that $\mu_3 \le x_3 < \lambda_2, \,\,$ $x_2$ such that $\lambda_0 < x_2 \le \xi^{(3)}_2$, and $x_1$ such that $\xi^{(3)}_1 \le x_1 < \lambda_1$.
\vskip.05cm
{\bf (iv)}  If $0 \le -d < -\delta_2$, then there is one non-positive root $x_4 \ge $ max$\{ \xi^{(2)}_2, -d/c \}$ and tree positive roots: $x_3$ such that $\lambda_2 \le x_3 < \mu_2, \,\,$ $x_2$ such that $\mu_2 < x_2 \le \lambda_0$, and $x_1$ such that $\lambda_1 \le x_1 < \xi^{(2)}_1$.
\vskip-.1cm
{\bf (v)}  If $-\delta_2 \le - d $ (pictured), then there is a negative root $x_2$ such that $-d/c < x_2 \le \xi^{(2)}_2$ and a positive root $x_1 \ge \xi^{(2)}_1$. \\
\hrule
\begin{center}
\vskip-.3cm
\underline{{\color{red}{\it Analysis based on solving quadratic equations only (continues on next page)}}}
\end{center}
\vskip-.4cm
{\bf \underline{\color{red}{Top-left pane}}}
\vskip.1cm

{\bf (i)} If $- d < c \tau_H + H$, then there are either no real roots, or there are two positive roots $x_{1,2} > -(d+H)/c$,  or there are two positive roots $x_{1,2} > -(d+H)/c$ together with two positive roots $x_{3,4}$ such that $-d/c < x_{3,4} < \tau_h$ ($x_{3,4}$ appear when $c \tau_H + H \ge -\delta_3$).
\vskip.1cm
{\bf (ii)} If $c \tau_H + H \le - d < 0$, then there is one positive root $x_2$ such that $-(d+h)/c < x_2 \le -(d+H)/c$, another positive root $x_1 > \tau_H$ and either zero or two positive roots $x_{3,4}$ such that $-d/c < x_{3,4} < \tau_h$ ($x_{3,4}$ are always present if $c \tau_H + H \ge -\delta_3$, while for the pictured $c \tau_H + H < -\delta_3$ the roots $x_{3,4}$ may or may not be there).
\vskip.1cm
{\bf (iii)} If $0 \le - d < c \tau_h + h$, then there is one non-positive root $x_4 \ge -d/c$, a positive root $x_3 < \tau_h$, another positive root $x_2$ such that $-(d+h)/c < x_2 < -(d+H)/c$, and a third positive root $x_1 > \tau_H$.
\vskip.1cm
{\bf (iv)} If $c \tau_h + h \le -d$ (pictured), then there is one negative root $x_2 > -d/c$, a positive root $x_1 > \tau_H$ and either zero or two positive roots greater than $\tau_h$ and smaller than or equal to $-(d+H)/c$. \\
\end{minipage}}\\
\hline
\end{tabular}
\end{center}


\begin{center}
\begin{tabular}{|c|c|}
\hline
\multicolumn{2}{|c|}{\footnotesize \bf Figure 4.3}  \\
\multicolumn{2}{|c|}{\footnotesize \bf (continued from previous page)}  \\
\hline
\hline
\multicolumn{2}{|c|}{{\color{violet} $\bm{\frac{9}{8}\,\frac{a^2}{4} < b \le \frac{4}{3}\,\frac{a^2}{4}\, , \quad a < 0, \quad c < 0,$}}} \\
\multicolumn{2}{|c|}{{\color{violet} $\bm{c_2 < \gamma_2 < {\color{red}c} < c_0 < \gamma_1 < c_1 < 0$}}} \\
\hline
\hline
\multicolumn{2}{|c|}{\begin{minipage}{45em}
\scriptsize
\begin{center}
\vskip.2cm
\underline{{\color{red}{\it Analysis based on solving quadratic equations only --- continued from previous page}}}
\end{center}
{\bf \underline{\color{red}{Top-right pane}}}
\vskip.2cm
{\bf (i)}  If $- d < c \tau_H + H$, then there are either no real roots, or there are two positive roots greater than $-(d+H)/c$,  or there are two positive roots greater than $-(d+H)/c$ together with two positive roots greater than $-d/c$ and smaller than $\tau_h$ (the latter appear when $c \tau_H + H \ge -\delta_3$).
\vskip.1cm
{\bf (ii)} If $c \tau_H + H \le - d < c \tau_h + h$, then there is a positive root $x_2$ such that $-(d+h)/c < x_2 \le -(d+H)/c$ and a positive root $x_1 > \tau_H$. If $-\delta_3 > c \tau_h + h > c \tau_H + H$, there are no other roots. If $c \tau_h + h >  c \tau_H + H  \ge -\delta_3$, there are two more positive roots $x_{3,4}$ greater than $-d/c$ and smaller than $\tau_h$. If $c \tau_h + h \ge -\delta_3 >  c \tau_H + H$, the roots $x_{3,4}$ may or may not be there.
\vskip.1cm
{\bf (iii)} If $ c \tau_h + h \le - d < 0 $, then, if $-\delta_3 \le c \sigma_h + h$, there will be a positive root $x_1 > \tau_H$, two positive roots $x_{2,3}$ such that $\tau_h \le x_{2,3} < -(d+H)/c$, and a positive root $x_4$ such that $-d/c < x_4 < \sigma_h$. If however, $c \sigma_h + h < -\delta_3$, then the roots $x_3$ and $x_4$ may or may not be there. In summary: there is one positive root $x_1 > \tau_H$ and either one or three positive roots greater than $-d/c$ and smaller than $-(d+H)/c$.
\vskip.1cm
{\bf (iv)} If $0 \le -d$ (pictured), then there is a non-positive root greater than or equal to $-d/c$, either zero or two positive roots greater than $\tau_h$ and smaller than $-(d+H)/c$, and a positive root greater than $\tau_H$. \\
\hrule
\vskip.3cm
{\bf \underline{\color{red}{Bottom-right pane}}}
\vskip.2cm
{\bf (i)} If $- d < 0$, then there are either no real roots, or there are two positive roots greater than $\tau_H$, or there are four positive roots: two greater than $\tau_H$ and two greater than $-d/c$ and smaller than $\tau_h$.
\vskip.1cm
{\bf (ii)} If $0 \le - d < c \tau_H + H$, then there is one non-positive root greater than $-d/c$, a positive root $x_3 < \tau_h$, and two positive roots $x_{1,2} > -(d+H)/c$.
\vskip.1cm
{\bf (iii)} If $c \tau_H + H \le - d < c \tau_h + h$, then there is one negative root $x_4 > -d/c$, a positive root $x_3 < \tau_h$, a positive root $x_2$ such that $-(d+h)/c < x_2 \le -(d+H)/c$ and another positive root $x_1 > \tau_H$.
\vskip.1cm
{\bf (iv)} If $c \tau_h + h \le -d$ (pictured), then there is a negative root greater than $-d/c$, a positive root greater than $\tau_H$ and either zero or two positive roots greater than $\tau_h$ and smaller than or equal to $-(d+H)/c$. \\
\end{minipage}}\\
\hline
\end{tabular}
\end{center}


\begin{center}
\begin{tabular}{|c||c|}
\hline
\multicolumn{2}{|c|}{\footnotesize \bf Figure 4.4}  \\
\hline
\hline
\multicolumn{2}{|c|}{{\color{violet} $\bm{\frac{9}{8}\,\frac{a^2}{4} < b \le \frac{4}{3}\,\frac{a^2}{4}\, , \quad a < 0, \quad c < 0,$}}} \\
\multicolumn{2}{|c|}{{\color{violet} $\bm{c_2 <  \gamma_2 < c_0 < {\color{red}c} < \gamma_1 < c_1 < 0$}}}  \\
\hline
\hline
\includegraphics[width=80mm]{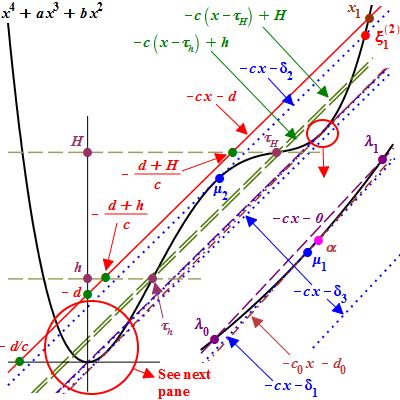} & \includegraphics[width=80mm]{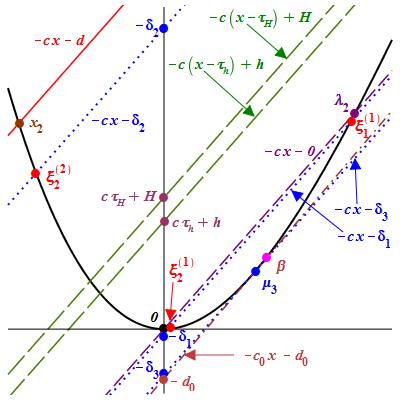} \\
\hline
\multicolumn{2}{|c|}{
\begin{minipage}{45em}
\scriptsize
\vskip.5cm
\begin{center}
\underline{{\color{red}{\it Notes}}}
\end{center}
\vskip-.3cm
As $c_2 < c < c_1$, the quartic has three stationary points $\mu_i$ and the number of real roots can be either 0, or 2, or 4. There are three tangents $- c x - \delta_i$ to $x^4 + ax^3 + bx^2$ --- each at $\mu_i$.
\vskip.1cm
As $\gamma_2 < c < \gamma_1$, the straight line $-cx$ intersects $x^4 + ax^3 + bx^2$ at three points: $\lambda_{0,1,2}$. Thus $-\delta_1 < 0$ and $-\delta_2 > 0$. Also: $-\delta_3 < 0$.
\vskip.1cm
The straight line joining the two points of curvature change ($\tau_h$ and $\tau_H$) has the same slope $-c_0 = -(1/2)a(b-a^2/4) > 0$ as that of the straight line $ - c_0 x - d_0$, tangent to $x^4 + ax^3 + bx^2$ at two points: $\alpha$ and $\beta = -\alpha - a/2$, where $\alpha$ and $\beta$ are the simultaneous double roots of the varied quartic $x^4 + ax^3 + bx^2 + c_0 x + d_0$ [with $d_0 = (1/4)(b-a^2/4)^2 > 0$].
\vskip.1cm
At $c = c_0$, one has  $-\delta_1 = -\delta_3$ and $c \tau_h + h = c \tau_H + H$.
\vskip.1cm
As $c_0 < c < 0$, obviously $-\delta_3 < -\delta_1 < 0 < c \tau_h + h < c \tau_H + H < -\delta_2$.
\vskip.1cm
Consideration of whether $- d \le c \tau_h + h$ or $-d > c \tau_h + h$ and also whether $- d \le c \tau_H + H$ or $-d > c \tau_H + H$ would yield sharper bounds in the analysis based on solving cubic equations. \\
\hrule
\vskip.3cm
\begin{center}
\underline{{\color{red}{\it Analysis based on solving cubic equations}}}
\end{center}
\vskip-.2cm
{\bf (i)} If $-d < -\delta_3$, then there are no real roots.
\vskip.1cm
{\bf (ii)} If $-\delta_3 \le -d < -\delta_1$, then there is one positive root $x_2$ such that max$\{ \xi^{(1)}_2, -d/c \} < x_2 \le \mu_3$  and another positive root $x_1$ such that $\mu_3 \le x_1 < \xi^{(1)}_1$.
\vskip.1cm
{\bf (iii)} If $-\delta_1 \le -d < 0$, then there are four positive roots: $x_4$ such that $- d/c < x_4 \le \xi^{(1)}_2, \,\,$ $x_3$ such that $\xi^{(1)}_1 \le x_3 < \lambda_2, \,\,$ $x_2$ such that $\lambda_0 < x_2 \le \mu_1$, and $x_1$ such that $\mu_1 \le x_1 < \lambda_1$.
\vskip.1cm
{\bf (iv)}  If $0 \le -d < -\delta_2$, then there is one non-positive root $x_4 \ge $ max $\{ \xi^{(2)}_2, -d/c \}$ and tree positive roots: $x_3$ such that $\lambda_2 \le x_3 < \mu_2, \,\,$ $x_2$ such that $\mu_2 < x_2 \le \lambda_0$, and $x_1$ such that $\lambda_1 \le x_1 < \xi^{(2)}_1$.
\vskip.1cm
{\bf (v)}  If $-\delta_2 \le - d $ (pictured), then there is a negative root $x_2$ such that $-d/c < x_2 \le \xi^{(2)}_2$ and a positive root $x_1 \ge \xi^{(2)}_1$. \\
\hrule
\vskip.3cm
\begin{center}
\underline{{\color{red}{\it Analysis based on solving quadratic equations only}}}
\end{center}
\vskip-.2cm
{\bf (i)} If $- d < 0$, then there are either no real roots, or there are two positive roots greater than $-d/c$ and smaller than $\tau_h$, or there are four positive roots, two of which greater than $-d/c$ and smaller than $\tau_h$ and the other two --- greater than $-(d+H)/c$.
\vskip.1cm
{\bf (ii)} If $0 \le - d < c \tau_h + h$, then there is one non-positive root $x_4 \ge -d/c$, and three positive roots: $x_3 < \tau_h$ and $x_{1,2} >  -(d+H)/c$.
\vskip.1cm
{\bf (iii)} If $c \tau_h + h \le - d < c \tau_H + H$, then there is one negative root greater than $-d/c$, a positive root greater than or equal to $\tau_h$ and smaller than $\tau_H$, and two positive roots greater than $-(d+H)/c$.
\vskip.1cm
{\bf (iv)} If $c \tau_H + H \le -d$ (pictured), then there is a negative root $x_2 > -d/c$, a positive root $x_1 > \tau_H$ and either zero or two positive roots greater than $\tau_h$ and smaller than or equal to $-(d+H)/c$. \\
\end{minipage}} \\
\hline
\end{tabular}
\end{center}


\begin{center}
\begin{tabular}{|c||c|}
\hline
{\footnotesize \bf Figure 4.5} &  {\footnotesize \bf Figure 4.6} \\
\hline
\hline
{\color{violet} $\bm{\frac{9}{8}\,\frac{a^2}{4} < b \le \frac{4}{3}\,\frac{a^2}{4}\, , \quad a < 0, \quad c < 0,$}} &
{\color{violet} $\bm{\frac{9}{8}\,\frac{a^2}{4} < b \le \frac{4}{3}\,\frac{a^2}{4}\, , \quad a < 0, \quad c < 0,$}} \\
{\color{violet} $\bm{c_2 <  \gamma_2 < c_0 < \gamma_1 < {\color{red}c} < c_1 < 0$}} &
{\color{violet} $\bm{c_2 < \gamma_2  < c_0 < \gamma_1 < c_1 < {\color{red}c} < 0$}} \\
\hline
\hline
\includegraphics[width=80mm]{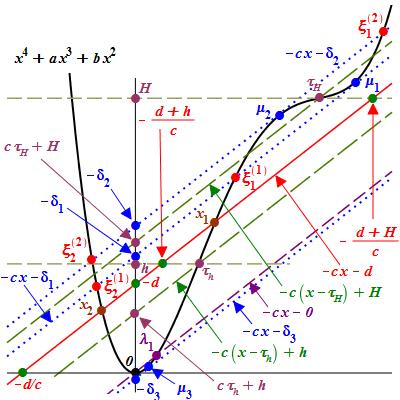} & \includegraphics[width=80mm]{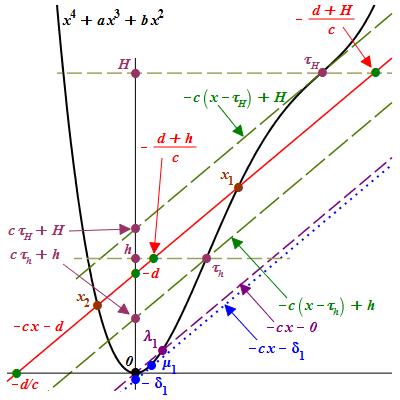} \\
\hline
\multicolumn{1}{|c||}{\begin{minipage}{22em}
\scriptsize
\vskip.2cm
\begin{center}
\underline{{\color{red}{\it Notes}}}
\end{center}
\vskip-.35cm
As $c_2 < c < c_1$, the quartic has three stationary points $\mu_i$ and the number of real roots can be either 0, or 2, or 4. There are three tangents $- c x - \delta_i$ to $x^4 + ax^3 + bx^2$ --- each at $\mu_i$.
\vskip.1cm
As $\gamma_2 < \gamma_1 < c$, the straight line $-cx$ intersects $x^4 + ax^3 + bx^2$ at one point only ($\lambda_1$).  Thus $-\delta_1 > 0$ and $-\delta_2 > 0$. Also: $-\delta_3 < 0$.
\vskip.1cm
One can have either the pictured $-\delta_3 < 0 < c \tau_h + h < -\delta_1 < c \tau_H + H < -\delta_2$ (when $c$ is closer to $c_1$; at $c = c_1$, the two stationary points $\mu_2$ and $\mu_1$ coalesce at $\tau_H$) or $-\delta_3 < 0 < -\delta_1 \le c \tau_h + h <  c \tau_H + H < -\delta_2$ (when $c$ is closer to $\gamma_1$).
\vskip.1cm
Consideration of whether $- d \le c \tau_h + h$ or $-d > c \tau_h + h$ and also whether $- d \le c \tau_H + H$ or $-d > c \tau_H + H$ would yield sharper bounds in the analysis based on solving cubic equations. \\
\hrule
\vskip.2cm
\hskip1.54cm \underline{{\color{red}{\it Analysis based on solving cubic equations}}} \\
\vskip-.2cm
{\bf (i)} If $-d < -\delta_3$, then there are no real roots.
\vskip.1cm
{\bf (ii)} If $-\delta_3 \le -d < 0$, then there is one positive root $x_2$ such that $-d/c < x_2 \le \mu_3$ and another positive root $x_1$ such that $\mu_3 \le x_1 < \lambda_1$.
\vskip.1cm
{\bf (iii)} If $0 \le -d < -\delta_1$ (pictured), then there is one non-positive root $x_2 \ge $ max$\{-d/c, \xi^{(1)}_2 \}$ and one positive root $x_1$ such that $\lambda_1 \le x_1 < \xi^{(1)}_1.$
\vskip.1cm
{\bf (iv)}  If $-\delta_1 \le -d < -\delta_2$, then there is one negative root $x_4$ such that max $\{ \xi^{(2)}_2, -d/c \} \le x_4 < \xi^{(1)}_2$ and tree positive roots: $x_3$ such that $\xi^{(1)}_1 \le x_3 < \mu_2, \,\,$ $x_2$ such that $\mu_2 < x_2 \le \mu_1$, and $x_1$ such that $\mu_1 \le x_1 < \xi^{(2)}_1$.
\vskip.1cm
{\bf (v)}  If $-\delta_2 \le - d $, then there is a negative root $x_2$ such that $-d/c < x_2 \le \xi^{(2)}_2$ and a positive root $x_1 \ge \xi^{(2)}_1$.
\vskip.1cm
\hrule
\vskip.2cm
\hskip.88cm \underline{{\color{red}{\it Analysis based on solving quadratic equations only}}} \\
\vskip-.2cm
{\bf (i)} If $- d < 0$, then there are either no real roots or there are two positive roots greater than $-d/c$ and smaller than $\tau_h$.
\vskip.1cm
{\bf (ii)} If $0 \le - d < c \tau_h + h$, then there is one non-positive root $x_2 \ge -d/c$, one positive root smaller than $\tau_h$ and either zero or two positive roots greater than $-(d+H)/c$ (the latter appear when $-\delta_1 \le c \tau_h + h$).
\vskip.1cm
{\bf (iii)} If $c \tau_h + h \le - d < c \tau_H + H$ (pictured), then there is one negative root greater than $-d/c$, a positive root greater than or equal to $\tau_h$ and smaller than $\tau_H$, and either zero or two positive roots greater than $-(d+H)/c$ (the latter are always present if $-\delta_1 \le c \tau_h + h$).
\vskip.1cm
{\bf (iv)} If $c \tau_H + H \le -d$, then there is a negative root $x_2 > -d/c$, a positive root $x_1 > \tau_H$ and either zero or two positive roots greater than $\tau_h$ and smaller than or equal to $-(d+H)/c$. \\
\end{minipage}}
& \multicolumn{1}{|c|}{\begin{minipage}{22em}
\scriptsize
\vskip-5.5cm
\begin{center}
\underline{{\color{red}{\it Notes}}}
\end{center}
\vskip-.3cm
As $c_2 < c_1 < c$, the quartic has a single stationary point $\mu_1$ and the number of real roots can be either 0 or 2. There is only one tangent $- c x - \delta_1$ to $x^4 + ax^3 + bx^2$ --- the one at $\mu_1$.
\vskip.1cm
As $\gamma_2 < \gamma_1 < c$, the straight line $-cx$ intersects $x^4 + ax^3 + bx^2$ at one point only ($\lambda_1$).
\vskip.1cm
Obviously, $c \tau_H + H > c \tau_h + h > 0 > -\delta_1$.
\vskip.1cm
Consideration of whether $- d \le c \tau_h + h$ or $-d > c \tau_h + h$ and also whether $- d \le c \tau_H + H$ or $-d > c \tau_H + H$ would yield sharper bounds in the analysis based on solving cubic equations. \\
\hrule
\vskip.2cm
\hskip1.54cm \underline{{\color{red}{\it Analysis based on solving cubic equations}}} \\
\vskip-.1cm
{\bf (i)} If $-d < -\delta_1$, then there are no real roots.
\vskip.1cm
{\bf (ii)} If $-\delta_1 \le -d < 0$, then there is one positive root $x_2$ such that $-d/c < x_2 \le \mu_1$ and another positive root $x_1$ such that $\mu_1 \le x_1 < \lambda_1$.
\vskip.1cm
{\bf (iii)} If $0 \le -d$ (pictured), then there is one non-positive root $x_2 \ge -d/c$ and one positive root $x_1 \ge \lambda_1.$
\vskip.1cm
\hrule
\vskip.2cm
\hskip.88cm \underline{{\color{red}{\it Analysis based on solving quadratic equations only}}} \\
\vskip-.1cm
{\bf (i)} If $- d < 0$, then there are either no real roots or there are two positive roots greater than $-d/c$ and smaller than $\tau_h$.
\vskip.1cm
{\bf (ii)} If $0 \le - d < c \tau_h + h$, then there is one non-positive root $x_2 \ge -d/c$ and one positive root smaller than $\tau_h$.
\vskip.1cm
{\bf (iii)} If $c \tau_h + h \le - d < c \tau_H + H$ (pictured), then there is one negative root greater than $-d/c$ and one positive root greater than or equal to $\tau_h$ and smaller than $\tau_H$.
\vskip.1cm
{\bf (iv)} If $c \tau_H + H \le -d$, then there is one negative root $x_2 > -d/c$ and one positive root $x_1 \ge \tau_H$.
\end{minipage}}  \\
\hline
\end{tabular}
\end{center}


\begin{center}
\begin{tabular}{|c||c|}
\hline
{\footnotesize \bf Figure 4.7} &  {\footnotesize \bf Figure 4.8} \\
\hline
\hline
{\color{violet} $\bm{\frac{9}{8}\,\frac{a^2}{4} < b \le \frac{4}{3}\,\frac{a^2}{4}\, , \quad a < 0, \quad c > 0,$}} &
{\color{violet} $\bm{\frac{9}{8}\,\frac{a^2}{4} < b \le \frac{4}{3}\,\frac{a^2}{4}\, , \quad a > 0, \quad c < 0,$}} \\
{\color{violet} $\bm{c_2 < \gamma_2  < c_0 < \gamma_1 < c_1 < 0 < {\color{red}c}$}} &
{\color{violet} $\bm{{\color{red}c} < 0 < c_2 < \gamma_2  < c_0 < \gamma_1 < c_1 $}} \\
\hline
\hline
\includegraphics[width=80mm]{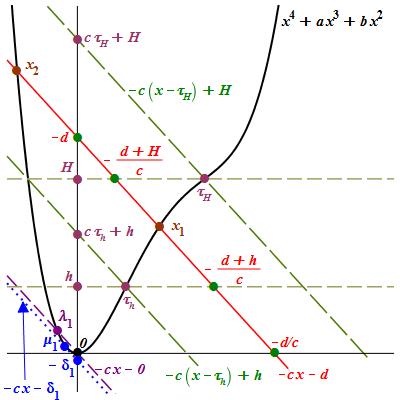} & \includegraphics[width=80mm]{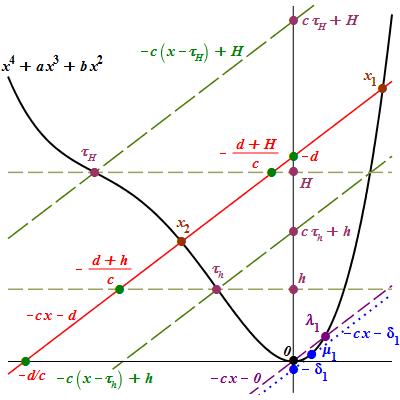} \\
\hline
\multicolumn{1}{|c||}{\begin{minipage}{22em}
\scriptsize
\vskip.3cm
\begin{center}
\underline{{\color{red}{\it Notes}}}
\end{center}
\vskip-.2cm
As $c_2 < c_1 < c$, the quartic has a single stationary point $\mu_1$ and the number of real roots can be either 0 or 2. There is only one tangent $- c x - \delta_1$ to $x^4 + ax^3 + bx^2$ --- the one at $\mu_1$.
\vskip.1cm
As $\gamma_2 < \gamma_1 < c$, the straight line $-cx$ intersects $x^4 + ax^3 + bx^2$ at one point only ($\lambda_1$).
\vskip.1cm
Obviously, $c \tau_H + H > c \tau_h + h > 0 > -\delta_1$.
\vskip.1cm
Consideration of whether $- d \le c \tau_h + h$ or $-d > c \tau_h + h$ and also whether $- d \le c \tau_H + H$ or $-d > c \tau_H + H$ would yield sharper bounds in the analysis based on solving cubic equations. \\
\hrule
\vskip.3cm
\hskip1.54cm \underline{{\color{red}{\it Analysis based on solving cubic equations}}}
\vskip.3cm
{\bf (i)} If $-d < -\delta_1$, then there are no real roots.
\vskip.1cm
{\bf (ii)} If $-\delta_1 \le -d < 0$, then there is one negative root $x_2$ such that $\lambda_1 < x_2 \le \mu_1$ and another negative root $x_1$ such that $\mu_1 \le x_1 < -d/c$.
\vskip.1cm
{\bf (iii)}  If $0 \le - d$ (pictured), then there is a negative root $x_2 \le \lambda_1$ and a non-negative root $x_1 \le -d/c$. \\
\hrule
\vskip.3cm
\hskip.88cm \underline{{\color{red}{\it Analysis based on solving quadratic equations only}}}
\vskip.3cm
{\bf (i)} If $- d < 0$, then there are either no real roots or there are two negative roots smaller than $-d/c$.
\vskip.1cm
{\bf (ii)} If $0 \le - d < c \tau_h + h$, then there is a non-positive root $x_2$ and a non-negative root $x_1$ such that max$\{0, -(d+h)/c \} \le x_1 < $ min$\{ -d/c, \tau_h \}$.
\vskip.1cm
{\bf (iii)} If $c \tau_h + h \le - d < c \tau_H + H$ (pictured), then there is one negative root $x_2$ and one positive root $x_1$ such that max$\{\tau_h, -(d+H)/c \} \le x_1 < $ min$\{ -(d+h)/c, \tau_H \}$.
\vskip.1cm
{\bf (iv)} If $c \tau_H + H \le -d$, then there is a negative root $x_2$ and a positive root $x_1$ such that  $\tau_H \le x_1 < -(d+H)/c$. \\
\end{minipage}}
& \multicolumn{1}{|c|}{\begin{minipage}{22em}
\scriptsize
\vskip.3cm
\begin{center}
\underline{{\color{red}{\it Notes}}}
\end{center}
\vskip-.2cm
As $c < c_2 < c_1$, the quartic has a single stationary point $\mu_1$ and the number of real roots can be either 0 or 2. There is only one tangent $- c x - \delta_1$ to $x^4 + ax^3 + bx^2$ --- the one at $\mu_1$.
\vskip.1cm
As $c < \gamma_2 < \gamma_1$, the straight line $-cx$ intersects $x^4 + ax^3 + bx^2$ at one point only ($\lambda_1$).
\vskip.1cm
Obviously, $c \tau_H + H > c \tau_h + h > 0 > -\delta_1$.
\vskip.1cm
Consideration of whether $- d \le c \tau_h + h$ or $-d > c \tau_h + h$ and also whether $- d \le c \tau_H + H$ or $-d > c \tau_H + H$ would yield sharper bounds in the analysis based on solving cubic equations. \\
\hrule
\vskip.3cm
\hskip1.54cm \underline{{\color{red}{\it Analysis based on solving cubic equations}}}
\vskip.3cm
{\bf (i)} If $-d < -\delta_1$, then there are no real roots.
\vskip.1cm
{\bf (ii)} If $-\delta_1 \le -d < 0$, then there is one positive root $x_2$ such that $-d/c < x_2 \le \mu_1$ and another positive root $x_1$ such that $\mu_1 \le x_1 < \lambda_1$.
\vskip.1cm
{\bf (iii)}  If $0 \le - d$ (pictured), then there is a non-positive root $x_2 \ge -d/c$ and a positive root $x_1 \ge \lambda_1$. \\
\hrule
\vskip.3cm
\hskip.88cm \underline{{\color{red}{\it Analysis based on solving quadratic equations only}}}
\vskip.3cm
{\bf (i)} If $- d < 0$, then there are either no real roots or there are two positive roots greater than $-d/c$.
\vskip.1cm
{\bf (ii)} If $0 \le - d < c \tau_h + h$, then there is a non-negative root $x_1$ and a non-positive root $x_2$ such that max$\{ -d/c, \tau_h \} < x_2 \le $ min$\{0, -(d+h)/c \}$.
\vskip.1cm
{\bf (iii)} If $c \tau_h + h \le - d < c \tau_H + H$ (pictured), then there is one positive root $x_1$ and one negative root $x_2$ such that max$\{ -(d+h)/c, \tau_H \} < x_2 \le $ min$\{\tau_h, -(d+H)/c \}$.
\vskip.1cm
{\bf (iv)} If $c \tau_H + H \le -d$, then there is a positive root $x_1$ and a negative root $x_2$ such that  $-(d+H)/c <  x_2 \le \tau_H$. \\
\end{minipage}}  \\
\hline
\end{tabular}
\end{center}


\begin{center}
\begin{tabular}{|c||c|}
\hline
{\footnotesize \bf Figure 4.9} &  {\footnotesize \bf Figure 4.10} \\
\hline
\hline
{\color{violet} $\bm{\frac{9}{8}\,\frac{a^2}{4} < b \le \frac{4}{3}\,\frac{a^2}{4}\, , \quad a > 0, \quad c > 0,$}} &
{\color{violet} $\bm{\frac{9}{8}\,\frac{a^2}{4} < b \le \frac{4}{3}\,\frac{a^2}{4}\, , \quad a > 0, \quad c > 0,$}} \\
{\color{violet} $\bm{0 < {\color{red}c} < c_2 <  \gamma_2 < c_0  < \gamma_1 < c_1$}} & {\color{violet} $\bm{0 < c_2 < {\color{red}c} <  \gamma_2 < c_0  < \gamma_1 < c_1$}} \\
\hline
\hline
\includegraphics[width=80mm]{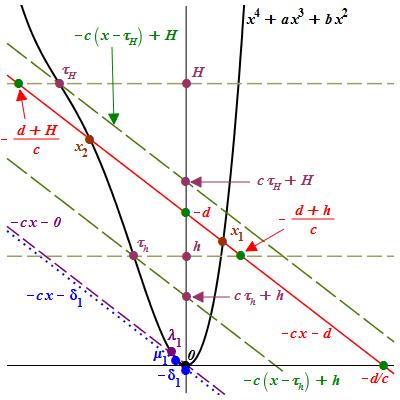} & \includegraphics[width=80mm]{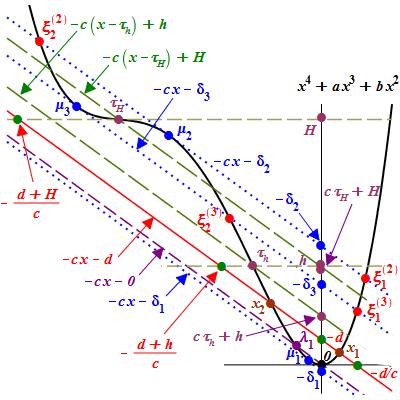} \\
\hline
\multicolumn{1}{|c||}{\begin{minipage}{22em}
\scriptsize
\vskip-5.6cm
\begin{center}
\underline{{\color{red}{\it Notes}}}
\end{center}
\vskip-.3cm
As $c < c_2 < c_1$, the quartic has a single stationary point $\mu_1$ and the number of real roots can be either 0 or 2. There is only one tangent $- c x - \delta_1$ to $x^4 + ax^3 + bx^2$ --- the one at $\mu_1$.
\vskip.1cm
As $c < \gamma_2 < \gamma_1$, the straight line $-cx$ intersects $x^4 + ax^3 + bx^2$ at one point only ($\lambda_1$).
\vskip.1cm
Obviously, $c \tau_H + H > c \tau_h + h > 0 > -\delta_1$.
\vskip.1cm
Consideration of whether $- d \le c \tau_h + h$ or $-d > c \tau_h + h$ and also whether $- d \le c \tau_H + H$ or $-d > c \tau_H + H$ would yield sharper bounds in the analysis based on solving cubic equations. \\
\hrule
\vskip.1cm
\hskip1.54cm \underline{{\color{red}{\it Analysis based on solving cubic equations}}} \\
\vskip-.1cm
{\bf (i)} If $-d < -\delta_1$, then there are no real roots.
\vskip.1cm
{\bf (ii)} If $-\delta_1 \le -d < 0$, then there is one negative root $x_2$ such that $-\lambda_1 < x_2 \le \mu_1$ and another negative root $x_1$ such that $\mu_1 \le x_1 < -d/c$.
\vskip.1cm
{\bf (iii)} If $0 \le -d$ (pictured), then there is one negative root $x_2 \le \lambda_1$ and one non-negative root $x_1 \le -d/c$.
\vskip.1cm
\hrule
\vskip.1cm
\hskip.88cm \underline{{\color{red}{\it Analysis based on solving quadratic equations only}}} \\
\vskip-.1cm
{\bf (i)} If $- d < 0$, then there are either no real roots or there are two negative roots smaller than $-d/c$ and greater than $\tau_h$.
\vskip.1cm
{\bf (ii)} If $0 \le - d < c \tau_h + h$, then there is one non-negative root $x_1 \le -d/c$ and one negative root greater than $\tau_h$.
\vskip.1cm
{\bf (iii)} If $c \tau_h + h \le - d < c \tau_H + H$ (pictured), then there is one positive root smaller than $-d/c$ and one negative root smaller than or equal to $\tau_h$ and greater than $\tau_H$.
\vskip.1cm
{\bf (iv)} If $c \tau_H + H \le -d$, then there is one positive root $x_1 < -d/c$ and one negative root $x_2 \le \tau_H$.
\end{minipage}}
& \multicolumn{1}{|c|}{\begin{minipage}{22em}
\scriptsize
\vskip.2cm
\begin{center}
\underline{{\color{red}{\it Notes}}}
\end{center}
\vskip-.3cm
As $c_2 < c < c_1$, the quartic has three stationary points $\mu_i$ and the number of real roots can be either 0, or 2, or 4. There are three tangents $- c x - \delta_i$ to $x^4 + ax^3 + bx^2$ --- each at $\mu_i$.
\vskip.1cm
As $c < \gamma_2 < \gamma_1$, the straight line $-cx$ intersects $x^4 + ax^3 + bx^2$ at one point only ($\lambda_1$).  Thus $-\delta_2 > 0$ and $-\delta_3 > 0$. Also: $-\delta_1 < 0$.
\vskip.1cm
One can have either the pictured $-\delta_1 < 0 < c \tau_h + h < -\delta_3 < c \tau_H + H < -\delta_2$ (when $c$ is closer to $c_2$; at $c = c_2$, the two stationary points $\mu_3$ and $\mu_2$ coalesce at $\tau_H$) or $-\delta_1 < 0 < -\delta_3 \le c \tau_h + h <  c \tau_H + H < -\delta_2$ (when $c$ is closer to $\gamma_2$).
\vskip.1cm
Consideration of whether $- d \le c \tau_h + h$ or $-d > c \tau_h + h$ and also whether $- d \le c \tau_H + H$ or $-d > c \tau_H + H$ would yield sharper bounds in the analysis based on solving cubic equations. \\
\hrule
\vskip.2cm
\hskip1.54cm \underline{{\color{red}{\it Analysis based on solving cubic equations}}} \\
\vskip-.2cm
{\bf (i)} If $-d < -\delta_1$, then there are no real roots.
\vskip.1cm
{\bf (ii)} If $-\delta_1 \le -d < 0$, then there is one negative root $x_2$ such that $-\lambda_1 < x_2 \le \mu_1$ and another negative root $x_1$ such that $\mu_1 \le x_1 < -d/c$.
\vskip.1cm
{\bf (iii)} If $0 \le -d < -\delta_3$ (pictured), then there is one negative root $x_2$ such that $\xi^{(3)}_2 < x_1 \le \lambda_1$ and one non-negative root $x_1 \le $ min$\{-d/c, \xi^{(3)}_1 \}$.
\vskip.1cm
{\bf (iv)}  If $-\delta_3 \le -d < -\delta_2$, then there is a negative root $x_4$ such that $\xi^{(2)}_2 \le x_4 \le \mu_3$, another negative root
$x_3$ such that $\mu_3 \le x_3 < \mu_2$, a third negative root $x_2$ such that $\mu_2 < x_2 \le \xi^{(3)}_2$  and one positive root $x_1$ such that
$\xi^{(3)}_1 \le x_1 < $ min$\{ \xi^{(2)}_1, -d/c \}$.
\vskip.1cm
{\bf (v)}  If $-\delta_2 \le - d $, then there is a positive root $x_1$ such that $\xi^{(2)}_1 \le x_1 < -d/c $ and a negative root $x_2 \le \xi^{(2)}_2$.
\vskip.1cm
\hrule
\vskip.2cm
\hskip.88cm \underline{{\color{red}{\it Analysis based on solving quadratic equations only}}} \\
\vskip-.2cm
{\bf (i)} If $- d < 0$, then there are either no real roots or there are two negative roots smaller than $-d/c$ and greater than $\tau_h$.
\vskip.1cm
{\bf (ii)} If $0 \le - d < c \tau_h + h$ (pictured), then there is one non-negative root $x_1 \le -d/c$, one negative root greater than $\tau_h$ and either zero or two negative roots smaller than $-(d+H)/c$ (the latter appear when $-\delta_3 \le c \tau_h + h$).
\vskip.1cm
{\bf (iii)} If $c \tau_h + h \le - d < c \tau_H + H$, then there is one positive root smaller than $-d/c$, a negative root smaller than or equal to $\tau_h$ and greater than $\tau_H$, and either zero or two negative roots smaller than $-(d+H)/c$ (the latter are always present if $-\delta_3 \le c \tau_h + h$).
\vskip.1cm
{\bf (iv)} If $c \tau_H + H \le -d$, then there is a positive root $x_1 < -d/c$, a negative root $x_2 < \tau_H$ and either zero or two negative roots smaller than $\tau_h$ and greater than or equal to $-(d+H)/c$. \\
\end{minipage}}  \\
\hline
\end{tabular}
\end{center}


\begin{center}
\begin{tabular}{|c||c|}
\hline
\multicolumn{2}{|c|}{\footnotesize \bf Figure 4.11}  \\
\hline
\hline
\multicolumn{2}{|c|}{{\color{violet} $\bm{\frac{9}{8}\,\frac{a^2}{4} < b \le \frac{4}{3}\,\frac{a^2}{4}\, , \quad a > 0, \quad c > 0,$}}} \\
\multicolumn{2}{|c|}{{\color{violet} $\bm{0 < c_2 <  \gamma_2 <  {\color{red}c} < c_0 < \gamma_1 < c_1$}}}  \\
\hline
\hline
\includegraphics[width=80mm]{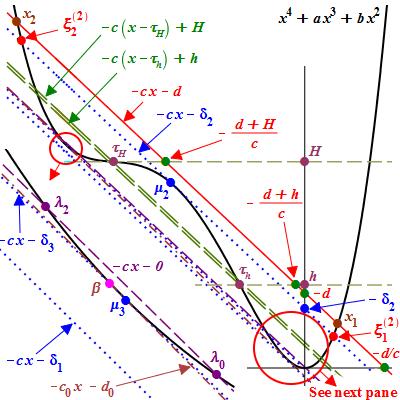} & \includegraphics[width=80mm]{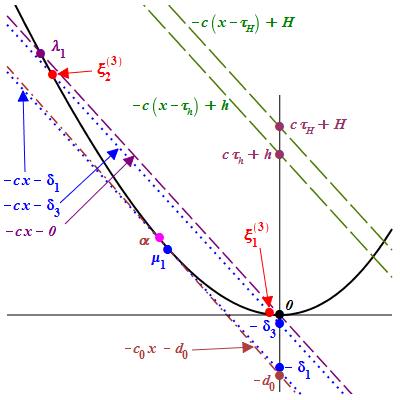} \\
\hline
\multicolumn{2}{|c|}{
\begin{minipage}{45em}
\scriptsize
\vskip.5cm
\begin{center}
\underline{{\color{red}{\it Notes}}}
\end{center}
\vskip-.3cm
As $c_2 < c < c_1$, the quartic has three stationary points $\mu_i$ and the number of real roots can be either 0, or 2, or 4. There are three tangents $- c x - \delta_i$ to $x^4 + ax^3 + bx^2$ --- each at $\mu_i$.
\vskip.1cm
As $\gamma_2 < c < \gamma_1$, the straight line $-cx$ intersects $x^4 + ax^3 + bx^2$ at three points: $\lambda_{0,1,2}$. Thus $-\delta_3 < 0$ and $-\delta_2 > 0$. Also: $-\delta_1 < 0$.
\vskip.1cm
The straight line joining the two points of curvature change ($\tau_h$ and $\tau_H$) has the same slope $-c_0 = -(1/2)a(b-a^2/4) > 0$ as that of the straight line $ - c_0 x - d_0$, tangent to $x^4 + ax^3 + bx^2$ at two points: $\alpha$ and $\beta = -\alpha - a/2$, where $\alpha$ and $\beta$ are the simultaneous double roots of the varied quartic $x^4 + ax^3 + bx^2 + c_0 x + d_0$ [with $d_0 = (1/4)(b-a^2/4)^2 > 0$].
\vskip.1cm
At $c = c_0$, one has  $-\delta_1 = -\delta_3$ and $c \tau_h + h = c \tau_H + H$.
\vskip.1cm
As $0 < c < c_0$, obviously $-\delta_1 < -\delta_3 < 0 < c \tau_h + h < c \tau_H + H < -\delta_2$.
\vskip.1cm
Consideration of whether $- d \le c \tau_h + h$ or $-d > c \tau_h + h$ and also whether $- d \le c \tau_H + H$ or $-d > c \tau_H + H$ would yield sharper bounds in the analysis based on solving cubic equations. \\
\hrule
\vskip.3cm
\begin{center}
\underline{{\color{red}{\it Analysis based on solving cubic equations}}}
\end{center}
\vskip-.2cm
{\bf (i)} If $-d < -\delta_1$, then there are no real roots.
\vskip.1cm
{\bf (ii)} If $-\delta_1 \le -d < -\delta_3$, then there is one negative root $x_2$ such that $\xi^{(3)}_2 < x_2 \le \mu_1$ and another negative root $x_1$ such that $\mu_1 \le x_1 < $ min$\{ \xi^{(3)}_1, -d/c \}$.
\vskip.1cm
{\bf (iii)} If $-\delta_3 \le -d < 0$, then there are four negative roots: $x_4$ such that $\lambda_2 < x_4 \le \mu_3, \,\, x_3$ such that $\mu_3 \le x_3 < \lambda_0, \,\, x_2$ such that $\lambda_1 < x_2 \le \xi^{(3)}_2$, and $x_1$ such that $\xi^{(3)}_1 \le x_1 < -d/c$.
\vskip.1cm
{\bf (iv)}  If $0 \le -d < -\delta_2$, then there is a negative root $x_4$ such that $\xi^{(2)}_2 < x_4 \le \lambda_2$, another negative roor $x_3$ such that $\lambda_0 \le x_3 < \mu_2$, a third negative roor $x_2$ such that $\mu_2 < x_2 \le \lambda_1$ and a positive root $x_1 \le $ min$\{ \xi^{(2)}_1, -d/c \}$.
\vskip.1cm
{\bf (v)}  If $-\delta_2 \le - d $ (pictured), then there is a positive root $x_1$ such that $\xi^{(2)}_1 \le x_1 < -d/c$ and a negative root $x_2 \le \xi^{(2)}_2$. \\
\hrule
\vskip.3cm
\begin{center}
\underline{{\color{red}{\it Analysis based on solving quadratic equations only}}}
\end{center}
\vskip-.2cm
{\bf (i)} If $- d < 0$, then there are either no real roots, or there are two negative roots smaller than $-d/c$ and greater than $\tau_h$, or there are four negative roots, two of which smaller than $-d/c$ and greater than $\tau_h$ and the other two --- smaller than $-(d+H)/c$.
\vskip.1cm
{\bf (ii)} If $0 \le - d < c \tau_h + h$, then there is one non-negative root $x_1 \le -d/c$, and three negative roots: $x_2 > \tau_h$ and $x_{3,4} <  -(d+H)/c$.
\vskip.1cm
{\bf (iii)} If $c \tau_h + h \le - d < c \tau_H + H$, then there is one positive root smaller than $-d/c$, a negative root smaller than or equal to $\tau_h$ and greater than $\tau_H$, and two negative roots smaller than $-(d+H)/c$.
\vskip.1cm
{\bf (iv)} If $c \tau_H + H \le -d$ (pictured), then there is a positive root $x_1 < -d/c$, a negative root $x_2 < \tau_H$ and either zero or two negative roots smaller than $\tau_h$ and greater than or equal to $-(d+H)/c$. \\
\end{minipage}} \\
\hline
\end{tabular}
\end{center}


\begin{center}
\begin{tabular}{|c|c|}
\hline
\multicolumn{2}{|c|}{\footnotesize \bf Figure 4.12}  \\
\multicolumn{2}{|c|}{\footnotesize \bf (continues on next page)}  \\
\hline
\hline
\multicolumn{2}{|c|}{{\color{violet} $\bm{\frac{9}{8}\,\frac{a^2}{4} < b \le \frac{4}{3}\,\frac{a^2}{4}\, , \quad a > 0, \quad c > 0,$}}} \\
\multicolumn{2}{|c|}{{\color{violet} $\bm{0 < c_2  <  \gamma_2 < c_0 < {\color{red}c} < \gamma_1 < c_1$}}} \\
\hline
\hline
\includegraphics[width=80mm,height=80mm]{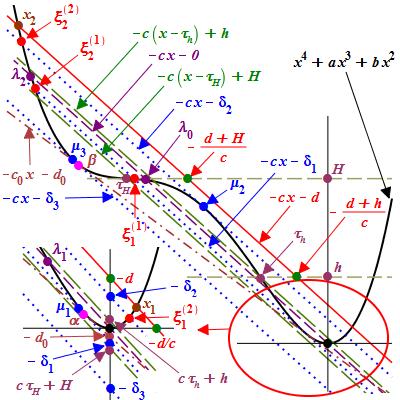} & \includegraphics[width=80mm,height=80mm]{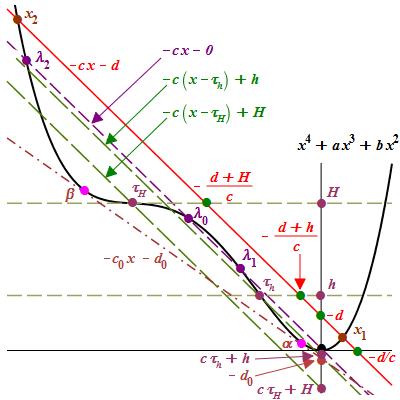} \\
\hline
& \\
& \\
\multicolumn{1}{|c|}{\begin{minipage}{22em}
\scriptsize
\begin{center}
\vskip-8.2cm
\underline{{\color{red}{\it Notes (apply to all panes)}}}
\end{center}
\vskip-.4cm
As $c_2 < c < c_1$, the quartic has three stationary points $\mu_i$ and the number of real roots can be either 0, or 2, or 4. There are three tangents $- c x - \delta_i$ to $x^4 + ax^3 + bx^2$ --- each at $\mu_i$ (only shown on the top-left pane).
\vskip.1cm
As $\gamma_2 < c < \gamma_1$, the straight line $-cx$ intersects $x^4 + ax^3 + bx^2$ at three points: $\lambda_{0,1,2}$. Thus $-\delta_2 > 0$ and $-\delta_3 < 0$. Also: $-\delta_1 < 0$.
\vskip.1cm
The straight line joining the two points of curvature change ($\tau_h$ and $\tau_H$) has the same slope $-c_0 = -(1/2)a(b-a^2/4) < 0$ as that of the straight line $ - c_0 x - d_0$, tangent to $x^4 + ax^3 + bx^2$ at two points: $\alpha$ and $\beta = -\alpha - a/2$, where $\alpha$ and $\beta$ are the simultaneous double roots of the varied quartic $x^4 + ax^3 + bx^2 + c_0 x + d_0$ [with $d_0 = (1/4)(b-a^2/4)^2 > 0$].
\vskip.1cm
As $c > c_0$, one has $- \delta_3 < -\delta_1$ and also $c \tau_H + H < c \tau_h + h$.
\vskip.1cm
\underline{Top-left pane}: one can have $-\delta_3 < c \tau_H + H < -\delta_1 < 0 < c \tau_h + h < -\delta_2$ (pictured) or $-\delta_3 < -\delta_1 \le c \tau_H + H < 0 < c \tau_h + h < -\delta_2$.
\vskip.1cm
\underline{Top-right pane} ($c$ closer to $\gamma_1$): one can have $-\delta_3 < c \tau_H + H < c \tau_h + h < -\delta_1 < 0 <  -\delta_2$, or $-\delta_3 < c \tau_H + H < -\delta_1 < c \tau_h + h < 0 < -\delta_2$, or $-\delta_3 < -\delta_1 < c \tau_H + H  < c \tau_h + h < 0 < -\delta_2$ (when $\tau_h$ and $\tau_H$ are close).
\vskip.1cm
\underline{Bottom-right pane} ($c$ closer to $c_0$): one can only have $-\delta_3 < -\delta_1 < 0 <  c \tau_H + H < c \tau_h + h < -\delta_2$.
\vskip.1cm
Consideration of whether $- d \le c \tau_h + h$ or $-d > c \tau_h + h$ and also whether $- d \le c \tau_H + H$ or $-d > c \tau_H + H$ would yield sharper bounds in the analysis based on solving cubic equations.
\\
\end{minipage}}
 & \includegraphics[width=80mm,height=80mm]{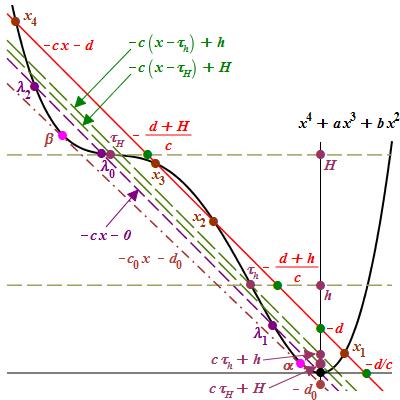} \\
\hline
\multicolumn{2}{|c|}{\begin{minipage}{45em}
\scriptsize
\vskip.1cm
\begin{center}
\underline{{\color{red}{\it Analysis based on solving cubic equations (applies to all panes)}}}
\end{center}
\vskip-.4cm
{\bf (i)} If $-d < -\delta_3$, then there are no real roots.
\vskip.1cm
{\bf (ii)} If $-\delta_3 \le -d < -\delta_1$, then there is one neative root $x_2$ such that $x_2$ such that $\xi^{(1)}_2  < x_2 \le \mu_3$ and another negative root $x_1$ such that $\mu_3 \le x_1 < \xi^{(1)}_1$.
\vskip-.05cm
{\bf (iii)} If $-\delta_1 \le -d < 0$, then there are four negative roots: $x_4$ such that $\lambda_2 < x_4 \le \xi^{(1)}_2, \,\,$ $x_3$ such that $\xi^{(1)}_1 \le x_3 < \lambda_0, \,\,$ $x_2$ such that $\lambda_1 < x_2 \le \mu_1$, and $x_1$ such that $\mu_1 \le x_1 < -d/c$.
\vskip.01cm
{\bf (iv)}  If $0 \le -d < -\delta_2$, then there is a negative root $x_4$ such that $\xi^{(2)}_2 < x_4 \le \lambda_2$, another negative root $x_3$ such that $\lambda_0 \le x_3 < \mu_2$, a third negative root $x_2$ such that $\mu_2 < x_2 \le \lambda_1$ and a non-negative root $x_1 \le $ min$\{ \xi^{(2)}_1, -d/c \}$.
\vskip.01cm
{\bf (v)}  If $-\delta_2 \le - d $ (pictured), then there is a positive root $x_1$ such that $\xi^{(2)}_1 \le x_1 < -d/c$ and a negative root $x_2 \le \xi^{(2)}_2$.
\\
\hrule
\begin{center}
\vskip-.3cm
\underline{{\color{red}{\it Analysis based on solving quadratic equations only (continues on next page)}}}
\end{center}
\vskip-.4cm
{\bf \underline{\color{red}{Top-left pane}}}
\vskip.1cm

{\bf (i)} If $- d < c \tau_H + H$, then there are either no real roots, or there are two negative roots smaller than $-(d+H)/c$,  or there are two negative roots smaller than $-(d+H)/c$ together with two negative roots greater than $\tau_h$ and smaller than $-d/c$ (the latter appear when $c \tau_H + H \ge -\delta_1$).
\vskip.1cm
{\bf (ii)} If $c \tau_H + H \le - d < 0$, then there is one negative root $x_1$ such that $-(d+H)/c \le x_1 < -(d+h)/c$, another negative root $x_2 < \tau_H$ and either zero or two negative roots smaller than $\tau_h$ and greater than $-d/c$ (the latter are always present if $c \tau_H + H \ge -\delta_1$, while for the pictured $c \tau_H + H < -\delta_1$ they may or may not be there).
\vskip.1cm
{\bf (iii)} If $0 \le - d < c \tau_h + h$, then there is one non-negative root $x_1 \le -d/c$, a negative root $x_2 > \tau_h$, another negative root $x_3$ such that $-(d+H)/c < x_3 < -(d+h)/c$, and a third negative root $x_4 < \tau_H$.
\vskip.1cm
{\bf (iv)} If $c \tau_h + h \le -d$ (pictured), then there is one positive root $x_1 < -d/c$, a negative root $x_2 < \tau_H$ and either zero or two negative roots smaller than $\tau_h$ and greater than or equal to $-(d+H)/c$. \\
\vskip-.45cm
$\phantom{e}$
\end{minipage}}\\
\hline
\end{tabular}
\end{center}


\begin{center}
\begin{tabular}{|c|c|}
\hline
\multicolumn{2}{|c|}{\footnotesize \bf Figure 4.12}  \\
\multicolumn{2}{|c|}{\footnotesize \bf (continued from previous page)}  \\
\hline
\hline
\multicolumn{2}{|c|}{{\color{violet} $\bm{\frac{9}{8}\,\frac{a^2}{4} < b \le \frac{4}{3}\,\frac{a^2}{4}\, , \quad a > 0, \quad c > 0,$}}} \\
\multicolumn{2}{|c|}{{\color{violet} $\bm{0 < c_2  <  \gamma_2 < c_0 < {\color{red}c} < \gamma_1 < c_1$}}} \\
\hline
\hline
\multicolumn{2}{|c|}{\begin{minipage}{45em}
\scriptsize
\begin{center}
\vskip.2cm
\underline{{\color{red}{\it Analysis based on solving quadratic equations only --- continued from previous page}}}
\end{center}
{\bf \underline{\color{red}{Top-right pane}}}
\vskip.2cm
{\bf (i)}  If $- d < c \tau_H + H$, then there are either no real roots, or there are two negative roots smaller than $-(d+H)/c$, or there are two negative roots smaller than $-(d+H)/c$ together with two negative roots greater than $\tau_h$ and smaller than $-d/c$ (the latter appear when $c \tau_H + H \ge -\delta_1$).
\vskip.1cm
{\bf (ii)} If $c \tau_H + H \le - d < c \tau_h + h$, then there is one negative root $x_1$ such that $-(d+H)/c \le x_1 < -(d+h)/c$ and another negative root $x_2 < \tau_H$. If $-\delta_1 > c \tau_h + h > c \tau_H + H$, there are no other roots. If $c \tau_h + h >  c \tau_H + H  \ge -\delta_1$, there are two more positive roots greater than $-d/c$ and smaller than $\tau_h$. If $c \tau_h + h \ge -\delta_1 >  c \tau_H + H$, these two roots may or may not be there.
\vskip.1cm
{\bf (iii)} If $ c \tau_h + h \le - d < 0 $, then, if $-\delta_1 \le c \sigma_h + h$, there will be a negative root $x_4 < \tau_H$, two negative roots $x_{2,3}$ such that $-(d+H)/c < x_{2,3} \le \tau_h$, and a negative root $x_1$ such that $\sigma_h < x_1 < -d/c$. If however, $c \sigma_h + h < -\delta_1$, then the roots $x_1$ and $x_2$ may or may not be there. In summary: there is one negative root smaller than $\tau_H$ and either one or three negative roots smaller than $-d/c$ and greater than $-(d+H)/c$.
\vskip.1cm
{\bf (iv)} If $0 \le -d$ (pictured), then there is a non-negative root smaller than or equal to $-d/c$, either zero or two negative roots smaller than $\tau_h$ and greater than $-(d+H)/c$, and a negative root smaller than $\tau_H$. \\
\hrule
\vskip.3cm
{\bf \underline{\color{red}{Bottom-right pane}}}
\vskip.2cm
{\bf (i)} If $- d < 0$, then there are either no real roots, or there are two negative roots smaller than $\tau_H$, or there are four negative roots: two smaller than $\tau_H$ and two smaller than $-d/c$ and greater than $\tau_h$.
\vskip.1cm
{\bf (ii)} If $0 \le - d < c \tau_H + H$, then there is one non-negative root smaller than $-d/c$, a neative root greater than  $\tau_h$, and two negative roots smaller than $-(d+H)/c$.
\vskip.1cm
{\bf (iii)} If $c \tau_H + H \le - d < c \tau_h + h$, then there is one positive root $x_1 < -d/c$, a negative root $x_2 > \tau_h$, a negative root $x_3$ such that $-(d+H)/c \le x_3 < -(d+h)/c$ and another negative root $x_4 < \tau_H$.
\vskip.1cm
{\bf (iv)} If $c \tau_h + h \le -d$ (pictured), then there is a positive root smaller than $-d/c$, a negative root smaller than $\tau_H$ and either zero or two negative roots smaller than $\tau_h$ and greater than or equal to $-(d+H)/c$. \\
\end{minipage}}\\
\hline
\end{tabular}
\end{center}


\begin{center}
\begin{tabular}{|c||c|}
\hline
{\footnotesize \bf Figure 4.13} &  {\footnotesize \bf Figure 4.14} \\
\hline
\hline
{\color{violet} $\bm{\frac{9}{8}\,\frac{a^2}{4} < b \le \frac{4}{3}\,\frac{a^2}{4}\, , \quad a > 0, \quad c > 0,$}} &
{\color{violet} $\bm{\frac{9}{8}\,\frac{a^2}{4} < b \le \frac{4}{3}\,\frac{a^2}{4}\, , \quad a > 0, \quad c > 0,$}} \\
{\color{violet} $\bm{0 < c_2 < \gamma_2 < c_0 < \gamma_1 < {\color{red}c} < c_1$}} & {\color{violet} $\bm{0 < c_2 < \gamma_2 < c_0 < \gamma_1 < c_1 < {\color{red}c}$}} \\
\hline
\hline
\includegraphics[width=80mm]{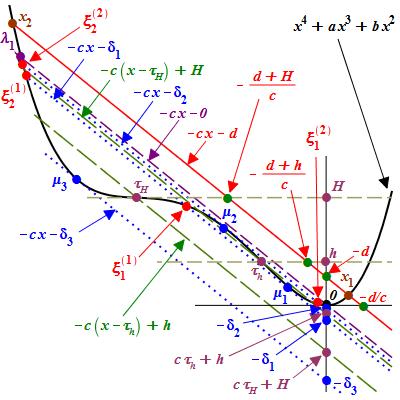} & \includegraphics[width=80mm]{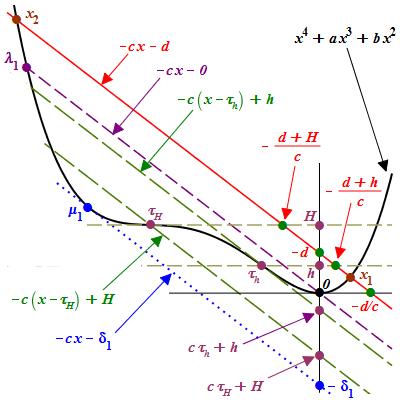} \\
\hline
\multicolumn{1}{|c||}{\begin{minipage}{22em}
\scriptsize
\vskip.2cm
\hskip3.89cm \underline{{\color{red}{\it Notes}}}
\vskip.1cm
As $c_2 < c < c_1$, the quartic has three stationary points $\mu_i$ and the number of real roots can be either 0, or 2, or 4. There are three tangents $- c x - \delta_i$ to $x^4 + ax^3 + bx^2$ --- each at $\mu_i$.
\vskip.1cm
As $\gamma_2 < \gamma_1 < c$, the straight line $-cx$ intersects $x^4 + ax^3 + bx^2$ at one point only ($\lambda_1$).  Thus $-\delta_2 < 0$ and $-\delta_3 < 0$. Also:  $-\delta_1 < 0$.
\vskip.1cm
The straight line joining the two points of curvature change ($\tau_h$ and $\tau_H$) has the same slope $-c_0 = -(1/2)a(b-a^2/4) < 0$ as that of the straight line $ - c_0 x - d_0$, tangent to $x^4 + ax^3 + bx^2$ at two points: $\alpha$ and $\beta = -\alpha - a/2$, where $\alpha$ and $\beta$ are the simultaneous double roots of the varied quartic $x^4 + ax^3 + bx^2 + c_0 x + d_0$ [with $d_0 = (1/4)(b-a^2/4)^2 > 0$].
\vskip.1cm
As $c > c_0 > 0$, one has $-\delta_3 < c \sigma_H + H < -\delta_1 < c \sigma_h + h < -\delta_2 < 0$.
\vskip.1cm
Consideration of whether $- d \le c \tau_h + h$ or $-d > c \tau_h + h$ and also whether $- d \le c \tau_H + H$ or $-d > c \tau_H + H$ would yield sharper bounds in the analysis based on solving cubic equations. \\
\hrule
\vskip.2cm
\hskip1.54cm \underline{{\color{red}{\it Analysis based on solving cubic equations}}}
\vskip.1cm
{\bf (i)} If $-d < -\delta_3$, then there are no real roots.
\vskip.1cm
{\bf (ii)} If $-\delta_3 \le -d < -\delta_1$, then there is one negative root $x_1$ such that $\mu_3 \le x_1 < $ min$\{ \xi^{(1)}_1, -d/c \}$ and another negative root $x_2$ such that $\xi^{(1)}_2 < x_1 \le \mu_3$.
\vskip.1cm
{\bf (iii)}  If $-\delta_1 \le -d < -\delta_2$, then there are four negative roots: $x_1$ such that $\mu_1 \le x_1 < $ min$\{ \xi^{(2)}_1, -d/c \}, \,\,$ $x_2$ such that $\mu_2 < x_2 \le \mu_1, \,\,$ $x_3$ such that $\xi^{(1)}_1 \le x_3 < \mu_2$, and $x_4$ such that $\xi^{(2)}_2 < x_4 \le  \xi^{(1)}_2$.
\vskip.1cm
{\bf (iv)} If $-\delta_2 \le -d < 0$, then there is a negative root $x_1$ such that $\xi^{(2)}_1 \le x_1 < -d/c$ and another negative root $x_2$ such that $\lambda_1 < x_2 \le \xi^{(2)}_2$.
\vskip.1cm
{\bf (v)}  If $0 \le - d$ (pictured), then there is a non-negative root $x_1 \le -d/c$ and a negative root $x_2 \le \lambda_1$. \\
\hrule
\vskip.2cm
\hskip.88cm \underline{{\color{red}{\it Analysis based on solving quadratic equations only}}}
\vskip.1cm
{\bf (i)} If $- d < c \tau_H + H$, then there are either no real roots, or there are two negative roots smaller than $-(d+H)/c$.
\vskip.1cm
{\bf (ii)} If $c \tau_H + H \le - d < c \tau_h + h$, then there is a negative root smaller than $\tau_H$, a negative root greater than or equal to  $-(d+H)/c$ and smaller than $-(d+h)/c$, together with either zero or two negative roots  greater than $\tau_h$ and smaller than $-d/c$.
\vskip.1cm
{\bf (iii)} If $c \tau_h + h \le - d < 0$, then there is a negative root greater than or equal to $-(d+h)/c$ and smaller than $-d/c$, a negative root smaller than $\tau_H$, and either zero or two negative roots smaller than or equal to $\tau_h$ and greater than $-(d+H)/c$.
\vskip.1cm
{\bf (iv)} If $0 \le -d$ (pictured), then there is one non-negtive root $x_1 \le -d/c$ and one negative root $x_2 < \tau_H$.  \\
\end{minipage}}
& \multicolumn{1}{|c|}{\begin{minipage}{22em}
\scriptsize
\vskip-4.35cm
\hskip3.89cm \underline{{\color{red}{\it Notes}}}
\vskip.1cm
As $c_2 < c_1 < c$, the quartic has a single stationary point $\mu_1$ and the number of real roots can be either 0 or 2. There is only one tangent $- c x - \delta_1$ to $x^4 + ax^3 + bx^2$ --- the one at $\mu_1$.
\vskip.1cm
As $\gamma_2 < \gamma_1 < c$, the straight line $-cx$ intersects $x^4 + ax^3 + bx^2$ at one point only ($\lambda_1$).
\vskip.1cm
The straight line joining the two points of curvature change ($\tau_h$ and $\tau_H$) has the same slope $-c_0 = -(1/2)a(b-a^2/4) < 0$ as that of the straight line $ - c_0 x - d_0$, tangent to $x^4 + ax^3 + bx^2$ at two points: $\alpha$ and $\beta = -\alpha - a/2$, where $\alpha$ and $\beta$ are the simultaneous double roots of the varied quartic $x^4 + ax^3 + bx^2 + c_0 x + d_0$ [with $d_0 = (1/4)(b-a^2/4)^2 > 0$].
\vskip.1cm
As $c > c_0 > 0$, one has $-\delta_1 < c \sigma_H + H < c \sigma_h + h < 0$.
\vskip.1cmConsideration of whether $- d \le c \tau_h + h$ or $-d > c \tau_h + h$ and also whether $- d \le c \tau_H + H$ or $-d > c \tau_H + H$ would yield sharper bounds in the analysis based on solving cubic equations. \\
\hrule
\vskip.2cm
\hskip1.54cm \underline{{\color{red}{\it Analysis based on solving cubic equations}}}
\vskip.1cm
{\bf (i)} If $-d < -\delta_1$, then there are no real roots.
\vskip.1cm
{\bf (ii)} If $-\delta_1 \le -d < 0$, then there is one negative root $x_1$ such that $\mu_1 \le x_2 < -d/c$ and another negative root $x_2$ such that $\lambda_1 < x_2 \le \mu_1$.
\vskip.1cm
{\bf (iii)}  If $0 \le - d$ (pictured), then there is a non-negative root $x_1 \le -d/c$ and a negative root $x_2 \le \lambda_1$. \\
\hrule
\vskip.2cm
\hskip.88cm \underline{{\color{red}{\it Analysis based on solving quadratic equations only}}}
\vskip.1cm
{\bf (i)} If $- d < c \tau_H + H$, then there are either no real roots or there are two negative roots smaller than $-(d+H)/c$.
\vskip.1cm
{\bf (ii)} If $c \tau_H + H \le - d < c \tau_h + h$, then there is a negative roots $x_1$ such that $-(d+H)/c \le x_1 < -(d+h)/c$ and another negative root $x_2 < \tau_H$.
\vskip.1cm
{\bf (iii)} If $c \tau_h + h \le - d < 0$, then there are two negative roots: $x_1$ such that $-(d+h)/c \le x_1 < -d/c$ and $x_2 < \tau_H$.
\vskip.1cm
{\bf (iv)} If $0 \le -d$ (pictured), then there is one non-negtive root $x_1 \le -d/c$ and one negative root $x_2 < \tau_H$.  \\
\end{minipage}}  \\
\hline
\end{tabular}
\end{center}


\begin{center}
\begin{tabular}{|c||c|}
\hline
{\footnotesize \bf Figure 5.1} &  {\footnotesize \bf Figure 5.2} \\
\hline
\hline
{\color{violet} $\bm{\frac{4}{3}\,\frac{a^2}{4} < b \le \frac{3}{2}\,\frac{a^2}{4}\, , \quad a < 0, \quad c < 0,$}} &
{\color{violet} $\bm{\frac{4}{3}\,\frac{a^2}{4} < b \le \frac{3}{2}\,\frac{a^2}{4}\, , \quad a < 0, \quad c < 0,$}} \\
{\color{violet} $\bm{{\color{red}c} < c_2 < c_0 < c_1 < 0$}} & {\color{violet} $\bm{c_2 < {\color{red}c} < c_0 < c_1  < 0$}} \\
\hline
\hline
\includegraphics[width=80mm,height=74mm]{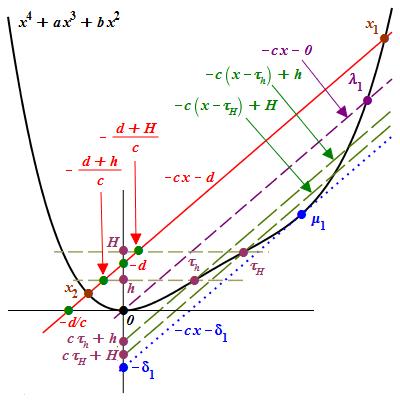} & \includegraphics[width=80mm,height=74mm]{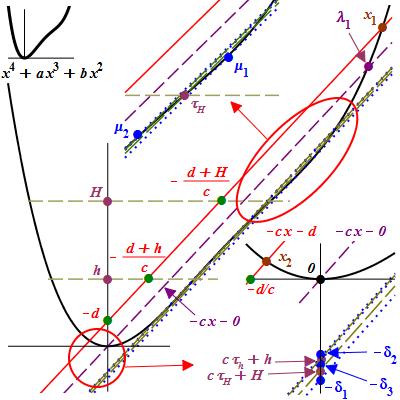} \\
\hline
\multicolumn{1}{|c||}{\begin{minipage}{22em}
\scriptsize
\vskip-5.1cm
\hskip3.89cm \underline{{\color{red}{\it Notes}}}
\vskip.1cm
As $c < c_2 < c_1$, the quartic has a single stationary point $\mu_1$ and the number of real roots can be either 0 or 2. There is only one tangent $- c x - \delta_1$ to $x^4 + ax^3 + bx^2$ --- the one at $\mu_1$.
\vskip.1cm
The straight line $-cx$ intersects $x^4 + ax^3 + bx^2$ at one point only ($\lambda_1$).
\vskip.1cm
The straight line joining the two points of curvature change ($\tau_h$ and $\tau_H$) has the same slope $-c_0 = -(1/2)a(b-a^2/4) > 0$ as that of the straight line $ - c_0 x - d_0$, tangent to $x^4 + ax^3 + bx^2$ at two points: $\alpha$ and $\beta = -\alpha - a/2$, where $\alpha$ and $\beta$ are the simultaneous double roots of the varied quartic $x^4 + ax^3 + bx^2 + c_0 x + d_0$ [with $d_0 = (1/4)(b-a^2/4)^2 > 0$].
\vskip.1cm
As $c < c_0 < 0$, one has $-\delta_1 < c \sigma_H + H < c \sigma_h + h < 0$.
\vskip.1cm
Consideration of whether $- d \le c \tau_h + h$ or $-d > c \tau_h + h$ and also whether $- d \le c \tau_H + H$ or $-d > c \tau_H + H$ would yield sharper bounds in the analysis based on solving cubic equations.
\\
\hrule
\vskip.1cm
\hskip1.54cm \underline{{\color{red}{\it Analysis based on solving cubic equations}}}
\vskip.1cm
{\bf (i)} If $-d < -\delta_1$, then there are no real roots.
\vskip.1cm
{\bf (ii)} If $-\delta_1 \le -d < 0$, then there is one positive root $x_2$ such that $-d/c < x_2 \le \mu_1$ and another positive root $x_1$ such that $\mu_1 \le x_2 < \lambda_1$.
\vskip.1cm
{\bf (iii)}  If $0 \le - d$ (pictured), then there is a non-positive root $x_2 \ge -d/c$ and a positive root $x_1 \ge \lambda_1$. \\
\hrule
\vskip.1cm
\hskip.88cm \underline{{\color{red}{\it Analysis based on solving quadratic equations only}}}
\vskip.1cm
{\bf (i)} If $-d < c \tau_H + H$, then there are either no real roots or there are two positive roots $x_{1,2} > -(d+H)/c$.
\vskip.1cm
{\bf (ii)} If $c \tau_H + H \le - d < c \tau_h + h$, then there are two positive roots: $x_2$ such that $-(d+h)/c < x_2 \le -(d+H)/c$ and $x_1 > \tau_H$.
\vskip.1cm
{\bf (iii)} If $c \tau_h + h \le - d < 0$, then there are two positive roots: $x_2$ such that $-d/c < x_2 \le -(d+h)/c$ and $x_1 > \tau_H$.
\vskip.1cm
{\bf (iv)} If $0 \le -d$ (pictured), then there is one non-positive root $x_2 \ge -d/c$ and one positive root $x_1 > \tau_H$. \\
\end{minipage}}
& \multicolumn{1}{|c|}{\begin{minipage}{22em}
\scriptsize
\vskip.05cm
\hskip3.89cm \underline{{\color{red}{\it Notes}}}
\vskip.05cm
As $c_2 < c < c_1$, the quartic has three stationary points $\mu_i$ and the number of real roots can be either 0, or 2, or 4. There are three tangents $- c x - \delta_i$ to $x^4 + ax^3 + bx^2$ --- each at $\mu_i$.
\vskip.1cm
The straight line $-cx$ intersects $x^4 + ax^3 + bx^2$ at one point only ($\lambda_1$). Thus $-\delta_1 < 0$ and $-\delta_2 < 0$, and $-\delta_3 < 0$.
\vskip.1cm
The straight line joining the two points of curvature change ($\tau_h$ and $\tau_H$) has the same slope $-c_0 = -(1/2)a(b-a^2/4) > 0$ as that of the straight line $ - c_0 x - d_0$, tangent to $x^4 + ax^3 + bx^2$ at two points: $\alpha$ and $\beta = -\alpha - a/2$, where $\alpha$ and $\beta$ are the simultaneous double roots of the varied quartic $x^4 + ax^3 + bx^2 + c_0 x + d_0$ [with $d_0 = (1/4)(b-a^2/4)^2 > 0$].
\vskip.1cm
As $c < c_0 < 0$, one has $-\delta_1 < -\delta_3$ and also $c \sigma_H + H < c \sigma_h + h$.
\vskip.1cm
One can have either the pictured $-\delta_1 < c \sigma_H + H < -\delta_3 < c \sigma_h + h < -\delta_2 < 0$ (when $c$ is closer to $c_2$ where $\mu_2$ and $\mu_3$ coalesce) or $-\delta_1 < -\delta_3 < c \sigma_H + H < c \sigma_h + h < -\delta_2 < 0$ (when $c$ is closer to $c_0$ where $c \tau_h + h$ and $c \tau_h + H$ swap around).
\vskip.1cm
Consideration of whether $- d \le c \tau_h + h$ or $-d > c \tau_h + h$ and also whether $- d \le c \tau_H + H$ or $-d > c \tau_H + H$ would yield sharper bounds in the analysis based on solving cubic equations.
\vskip.1cm
\phantom{e}
\vskip-.3cm
\hrule
\vskip.05cm
\hskip1.54cm \underline{{\color{red}{\it Analysis based on solving cubic equations}}}
\vskip.05cm
{\bf (i)} If $-d < -\delta_1$, then there are no real roots.
\vskip.1cm
{\bf (ii)} If $-\delta_1 \le -d < -\delta_3$, then there is one positive root $x_2$ such that max$\{ \xi^{(3)}_2, -d/c \} < x_2 \le \mu_1$ and another positive root $x_1$ such that $\mu_1 \le x_1 < \xi^{(3)}_1$.
\vskip.1cm
{\bf (iii)}  If $-\delta_3 \le -d < -\delta_2$, then there are four positive roots: $x_4$ such that max$\{ \xi^{(2)}_2, -d/c \} < x_4 \le \mu_3, \,\,$ $x_3$ such that $\mu_3 \le x_3 < \mu_2, \,\,$ $x_2$ such that $\mu_2 < x_2 \le \xi^{(3)}_2$, and $x_1$ such that $\xi^{(3)}_1 \le x_1 < \xi^{(2)}_1$. \vskip.1cm
{\bf (iv)} If $-\delta_2 \le -d < 0$, then there is a positive root $x_2$ such that $-d/c < x_2 \le \xi^{(2)}_2$ and another positive root $x_1$ such that $\xi^{(2)}_1 \le x_1 < \lambda_1$.
\vskip.1cm
{\bf (v)}  If $0 \le - d$ (pictured), then there is a non-positive root $x_2 \ge -d/c$ and a positive root $x_1 \ge \lambda_1$.
\vskip.1cm
\phantom{e}
\vskip-.3cm
\hrule
\vskip.05cm
\hskip.88cm \underline{{\color{red}{\it Analysis based on solving quadratic equations only}}}
\vskip.05cm
{\bf (i)} If $-d < c \tau_H + H$, then there are either no real roots, or there are two positive roots $x_{1,2} > -(d+H)/c$, or there are two positive roots greater than $-(d+H)/c$ together with two positive roots greater than $-d/c$ and smaller than $\tau_h$ (the latter appear when $c \tau_H + H > -\delta_3$).
\vskip.1cm
{\bf (ii)} If $c \tau_H + H \le - d < c \tau_h + h$, then there are two positive roots: $x_2$ such that $-(d+h)/c < x_2 \le -(d+H)/c$ and $x_1 > \tau_H$, together with either zero or two positive roots $x_{3,4}$ such that $-d/c < x_{3,4} < \tau_h$ ($x_{3,4}$ are always present if $c \tau_H + H \ge -\delta_3$).
\vskip.1cm
{\bf (iii)} If $c \tau_h + h \le - d < 0$, then there are two positive roots: $x_2$ such that $-d/c < x_2 < -(d+h)/c$ and $x_1 > \tau_H$, together with either zero or two positive roots greater than or equal to $\tau_h$ and smaller than $-(d+H)/c$.
\vskip.1cm
{\bf (iv)} If $0 \le -d$ (pictured), then there is one non-positive root $x_2 \ge -d/c$ and one positive root $x_1 > \tau_H$.
\\
\end{minipage}}  \\
\hline
\end{tabular}
\end{center}



\begin{center}
\begin{tabular}{|c|c|}
\hline
\multicolumn{2}{|c|}{\footnotesize \bf Figure 5.4}  \\
\multicolumn{2}{|c|}{\footnotesize \bf (Figure 5.3 is on next page)}  \\
\hline
\hline
\multicolumn{2}{|c|}{{\color{violet} $\bm{\frac{4}{3}\,\frac{a^2}{4} < b \le \frac{3}{2}\,\frac{a^2}{4}\, , \quad a < 0, \quad c < 0,$}}} \\
\multicolumn{2}{|c|}{{\color{violet} $\bm{c_2 < c_0 < c_1 < {\color{red}c} < 0$}}} \\
\hline
\hline
\includegraphics[width=80mm,height=80mm]{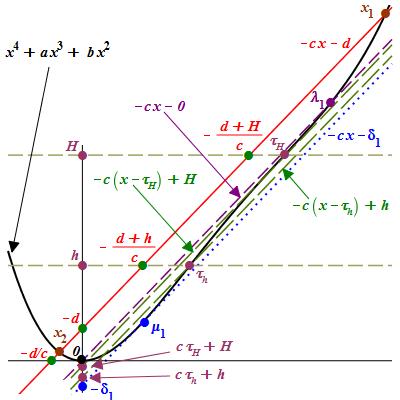} & \includegraphics[width=80mm,height=80mm]{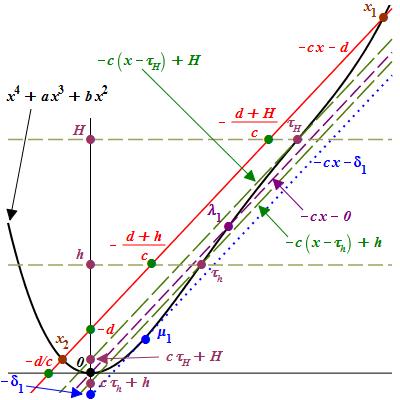} \\
\hline
\multicolumn{1}{|c|}{\begin{minipage}{22em}
\scriptsize
\begin{center}
\vskip-7.67cm
\underline{{\color{red}{\it Notes (apply to all panes)}}}
\end{center}
\vskip-.1cm
As $c_2 < c_1 < c$, the quartic has a single stationary point $\mu_1$ and the number of real roots can be either 0 or 2. There is only one tangent $- c x - \delta_1$ to $x^4 + ax^3 + bx^2$ --- the one at $\mu_1$.
\vskip.1cm
The straight line $-cx$ intersects $x^4 + ax^3 + bx^2$ at one point only ($\lambda_1$).
\vskip.1cm
The straight line joining the two points of curvature change ($\tau_h$ and $\tau_H$) has the same slope $-c_0 = -(1/2)a(b-a^2/4) > 0$ as that of the straight line $ - c_0 x - d_0$, tangent to $x^4 + ax^3 + bx^2$ at two points: $\alpha$ and $\beta = -\alpha - a/2$, where $\alpha$ and $\beta$ are the simultaneous double roots of the varied quartic $x^4 + ax^3 + bx^2 + c_0 x + d_0$ [with $d_0 = (1/4)(b-a^2/4)^2 > 0$].
\vskip.1cm
As $c_0 < c_1 < c < 0$, one can have $-\delta_1 < c \tau_h + h  < c \tau_H + H < 0$ (top-left pane, when $c$ is closer to $c_1$), or $-\delta_1 < c \tau_h + h < 0 < c \tau_H + H$ (top-right pane), or $-\delta_1 < 0 < c \tau_h + h < c \tau_H + H$ (bottom-right pane, when $c$ is closer to 0).
\vskip.1cm
Consideration of whether $- d \le c \tau_h + h$ or $-d > c \tau_h + h$ and also whether $- d \le c \tau_H + H$ or $-d > c \tau_H + H$ would yield sharper bounds in the analysis based on solving cubic equations.
\end{minipage}}
 & \includegraphics[width=80mm,height=80mm]{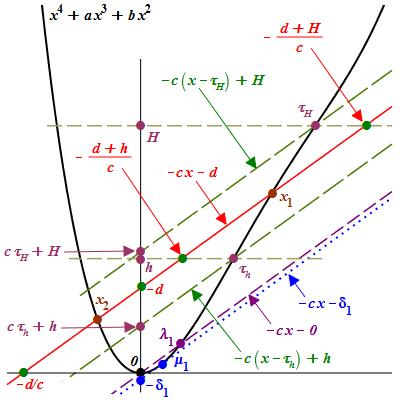} \\
\hline
\multicolumn{2}{|c|}{\begin{minipage}{45em}
\scriptsize
\vskip.1cm
\begin{center}
\underline{{\color{red}{\it Analysis based on solving cubic equations (applies to all panes)}}}
\end{center}
\vskip-.3cm
{\bf (i)} If $-d < -\delta_1$, then there are no real roots.
\vskip.1cm
{\bf (ii)} If $-\delta_1 \le -d < 0$, then there is one positive root $x_2$ such that $-d/c < x_2 \le \mu_1$ and another positive root $x_1$ such that $\mu_1 \le x_1 < \lambda_1$.
\vskip.1cm
\vskip.1cm
{\bf (iii)} If $0 \le -d$ (pictured), then there is one non-positive root $x_2 \ge -d/c$ and one positive root $x_1 \ge \lambda_1.$
\vskip.1cm
\hrule
\begin{center}
\vskip-.3cm
\underline{{\color{red}{\it Analysis based on solving quadratic equations only (clockwise from the top-left pane, with four possibile situations per pane)}}}
\end{center}
\vskip-.3cm
{\bf (i)} If $- d < c \tau_h + h$, then there are either no real roots or there are two positive roots greater than $-d/c$ and smaller than $\tau_h$.
\vskip.1cm
{\bf (ii)}  If $c \tau_h + h \le - d < c \tau_H + H$, then there is a positive root $x_2$ such that $-d/c < x_2 < -(d+h)/c$ and another positive root greater than or equal to $\tau_h$ and smaller than $\tau_H$.
\vskip.1cm
{\bf (iii)} If $c \tau_H + H \le - d < 0$, then there is a positive root greater than $-d/c$ and smaller than $-(d+h)/c$ and another positive root greater than or equal to $\tau_H$.
\vskip.1cm
{\bf (iv)} If $0 \le -d$ (pictured), then there is one non-positive root $x_2 \ge -d/c$ and one positive root $x_1 > \tau_H$.
\\
\hrule
\vskip.1cm
{\bf (v)} If $- d < c \tau_h + h$, then there are either no real roots or there are two positive roots greater than $-d/c$ and smaller than $\tau_h$.
\vskip.1cm
{\bf (vi)} If $c \tau_h + h \le - d < 0$, then there is a positive root $x_2$ such that $-d/c < x_2 < -(d+h)/c$ and another positive root greater than or equal to $\tau_h$ and smaller than $\tau_H$.
\vskip.1cm
{\bf (vii)} If $0 \le - d < c \tau_H + H$, then there is a non-positive root $x_2 \ge -d/c$ and a positive root $x_1$ such that $\tau_h < x_1 < \tau_H$.
\vskip.1cm
{\bf (viii)} If $c \tau_H + H \le -d$ (pictured), then there is one negative root $x_2 > -d/c$ and one positive root $x_1 \ge \tau_H$.
\\
\hrule
\vskip.1cm
{\bf (ix)} If $- d < 0$, then there are either no real roots or there are two positive roots greater than $-d/c$ and smaller than $\tau_h$.
\vskip.1cm
{\bf (x)} If $0 \le - d < c \tau_h + h$, then there is one non-positive root $x_2 \ge -d/c$ and one positive root smaller than $\tau_h$.
\vskip.1cm
{\bf (xi)} If $c \tau_h + h \le - d < c \tau_H + H$ (pictured), then there is one negative root greater than $-d/c$ and one positive root greater than or equal to $\tau_h$ and smaller than $\tau_H$.
\vskip.1cm
{\bf (xii)} If $c \tau_H + H \le -d$, then there is one negative root $x_2 > -d/c$ and one positive root $x_1 \ge \tau_H$.
\\
\end{minipage}}\\
\hline
\end{tabular}
\end{center}


\begin{center}
\begin{tabular}{|c||c|}
\hline
{\footnotesize \bf Figure 5.3} &  {\footnotesize \bf Figure 5.5} \\
\hline
\hline
{\color{violet} $\bm{\frac{4}{3}\,\frac{a^2}{4} < b \le \frac{3}{2}\,\frac{a^2}{4}\, , \quad a < 0, \quad c < 0,$}} &
{\color{violet} $\bm{\frac{4}{3}\,\frac{a^2}{4} < b \le \frac{3}{2}\,\frac{a^2}{4}\, , \quad a < 0, \quad c > 0,$}} \\
{\color{violet} $\bm{c_2 < c_0 < {\color{red}c} < c_1 < 0$}} & {\color{violet} $\bm{c_2 < c_0 < c_1 < 0 < {\color{red}c}$}} \\
\hline
\hline
\includegraphics[width=80mm]{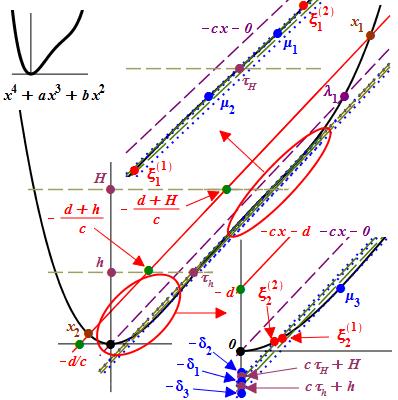} & \includegraphics[width=80mm]{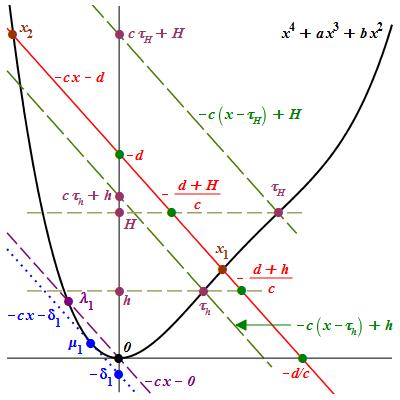} \\
\hline
\multicolumn{1}{|c||}{\begin{minipage}{22em}
\scriptsize
\vskip.05cm
\hskip3.89cm \underline{{\color{red}{\it Notes}}}
\vskip.05cm
As $c_2 < c < c_1$, the quartic has three stationary points $\mu_i$ and the number of real roots can be either 0, or 2, or 4. There are three tangents $- c x - \delta_i$ to $x^4 + ax^3 + bx^2$ --- each at $\mu_i$.
\vskip.1cm
The straight line $-cx$ intersects $x^4 + ax^3 + bx^2$ at one point only ($\lambda_1$). Thus $-\delta_1 < 0$ and $-\delta_2 < 0$, and $-\delta_3 < 0$.
\vskip.1cm
The straight line joining the two points of curvature change ($\tau_h$ and $\tau_H$) has the same slope $-c_0 = -(1/2)a(b-a^2/4) > 0$ as that of the straight line $ - c_0 x - d_0$, tangent to $x^4 + ax^3 + bx^2$ at two points: $\alpha$ and $\beta = -\alpha - a/2$, where $\alpha$ and $\beta$ are the simultaneous double roots of the varied quartic $x^4 + ax^3 + bx^2 + c_0 x + d_0$ [with $d_0 = (1/4)(b-a^2/4)^2 > 0$].
\vskip.1cm
As $c_0 < c < 0$, one has $-\delta_3 < -\delta_1$  and also $c \sigma_h + h < c \sigma_H + H$.
\vskip.1cm
One can have either the pictured $-\delta_3 < c \sigma_h + h < -\delta_1 < c \sigma_H + H < -\delta_2 < 0$ (when $c$ is closer to $c_1$ where $\mu_2$ and $\mu_1$ coalesce) or $-\delta_3 < -\delta_1 < c \sigma_h + h < c \sigma_H + H < -\delta_2 < 0$ (when $c$ is closer to $c_0$ where $c \tau_h + h$ and $c \tau_h + H$ swap around).
\vskip.1cm
Consideration of whether $- d \le c \tau_h + h$ or $-d > c \tau_h + h$ and also whether $- d \le c \tau_H + H$ or $-d > c \tau_H + H$ would yield sharper bounds in the analysis based on solving cubic equations.
\vskip.1cm
\phantom{e}
\vskip-.3cm
\hrule
\vskip.05cm
\hskip1.54cm \underline{{\color{red}{\it Analysis based on solving cubic equations}}} \\
\vskip-.3cm
{\bf (i)} If $-d < -\delta_3$, then there are no real roots.
\vskip.1cm
{\bf (ii)} If $-\delta_3 \le -d < -\delta_1$, then there is one positive root $x_2$ such that max$\{ \xi^{(1)}_2, -d/c \} < x_2 \le \mu_3$ and another positive root $x_1$ such that $\mu_3 \le x_1 < \xi^{(1)}_1$.
\vskip.05cm
{\bf (iii)}  If $-\delta_1 \le -d < -\delta_2$, then there are four positive roots: $x_4$ such that max$\{ \xi^{(2)}_2, -d/c \} < x_4 \le \xi^{(1)}_2, \,\,$ $x_3$ such that $\xi^{(1)}_1 \le x_3 < \mu_1, \,\,$ $x_2$ such that $\mu_2 < x_2 \le \mu_1$, and $x_1$ such that $\mu_1 \le x_1 < \xi^{(2)}_1$. \vskip.1cm
{\bf (iv)} If $-\delta_2 \le -d < 0$, then there is a positive root $x_2$ such that $-d/c < x_2 \le \xi^{(2)}_2$ and another positive root $x_1$ such that $\xi^{(2)}_1 \le x_1 < \lambda_1$.
\vskip.1cm
{\bf (v)}  If $0 \le - d$ (pictured), then there is a non-positive root $x_2 \ge -d/c$ and a positive root $x_1 \ge \lambda_1$.
\vskip.1cm
\hrule
\vskip.05cm
\hskip.88cm \underline{{\color{red}{\it Analysis based on solving quadratic equations only}}}
\vskip.01cm
{\bf (i)} If $-d < c \tau_h + h$, then there are either no real roots, or there are two positive roots greater than $-d/c$ and smaller than $\tau_h$, or there are two positive roots  greater than $-d/c$ and smaller than $\tau_h$, together with two positive roots greater than $-(d+H)/c$ (the latter appear when $c \tau_h + h > -\delta_1$).
\vskip.05cm
{\bf (ii)} If $c \tau_h + h \le - d < c \tau_H + H$, then there is a positive root $x_2$ such that $-d/c < x_{2} < -(d+h)/c$, a positive root $x_1$ such that $\tau_h \le x_1 < \tau_H$, and either zero or two positive roots greater than $-(d+H)/c$ (the latter are always present of $c\tau_h + h \ge -\delta_1$).
\vskip.1cm
{\bf (iii)} If $c \tau_H + H \le - d < 0$, then there is a positive root $x_2$ such that $-d/c < x_2 < -(d+h)/c$, a positive root $x_1 > \tau_H$, and either zero or two positive roots smaller than or equal to $-(d+H)/c$ and greater than $\tau_h$.
\vskip.1cm
{\bf (iv)} If $0 \le -d$ (pictured), then there is one non-positive root $x_2 \ge -d/c$ and one positive root $x_1 > \tau_H$.
\\
\end{minipage}}
& \multicolumn{1}{|c|}{\begin{minipage}{22em}
\scriptsize
\vskip-6.5cm
\begin{center}
\underline{{\color{red}{\it Notes}}}
\end{center}
\vskip-.3cm
As $c_2 < c_1 < c$, the quartic has a single stationary point $\mu_1$ and the number of real roots can be either 0 or 2. There is only one tangent $- c x - \delta_1$ to $x^4 + ax^3 + bx^2$ --- the one at $\mu_1$.
\vskip.1cm
The straight line $-cx$ intersects $x^4 + ax^3 + bx^2$ at one point only ($\lambda_1$).
\vskip.1cm
Obviously, $c \tau_H + H > c \tau_h + h > 0 > -\delta_1$.
\vskip.1cm
Consideration of whether $- d \le c \tau_h + h$ or $-d > c \tau_h + h$ and also whether $- d \le c \tau_H + H$ or $-d > c \tau_H + H$ would yield sharper bounds in the analysis based on solving cubic equations. \\
\hrule
\vskip.2cm
\hskip1.54cm \underline{{\color{red}{\it Analysis based on solving cubic equations}}}
\vskip.2cm
{\bf (i)} If $-d < -\delta_1$, then there are no real roots.
\vskip.1cm
{\bf (ii)} If $-\delta_1 \le -d < 0$, then there is one negative root $x_2$ such that $\lambda_1 < x_2 \le \mu_1$ and another negative root $x_1$ such that $\mu_1 \le x_1 < -d/c$.
\vskip.1cm
{\bf (iii)}  If $0 \le - d$ (pictured), then there is a negative root $x_2 \le \lambda_1$ and a non-negative root $x_1 \le -d/c$. \\
\hrule
\vskip.2cm
\hskip.88cm \underline{{\color{red}{\it Analysis based on solving quadratic equations only}}}
\vskip.2cm
{\bf (i)} If $- d < 0$, then there are either no real roots or there are two negative roots smaller than $-d/c$.
\vskip.1cm
{\bf (ii)} If $0 \le - d < c \tau_h + h$, then there is a non-positive root $x_2$ and a non-negative root $x_1$ such that max$\{0, -(d+h)/c \} \le x_1 < $ min$\{ -d/c, \tau_h \}$.
\vskip.1cm
{\bf (iii)} If $c \tau_h + h \le - d < c \tau_H + H$ (pictured), then there is one negative root $x_2$ and one positive root $x_1$ such that max$\{\tau_h, -(d+H)/c \} \le x_1 < $ min$\{ -(d+h)/c, \tau_H \}$.
\vskip.1cm
{\bf (iv)} If $c \tau_H + H \le -d$, then there is a negative root $x_2$ and a positive root $x_1$ such that  $\tau_H \le x_1 < -(d+H)/c$. \\
\end{minipage}}  \\
\hline
\end{tabular}
\end{center}


\begin{center}
\begin{tabular}{|c|c|}
\hline
\multicolumn{2}{|c|}{\footnotesize \bf Figure 5.7}  \\
\multicolumn{2}{|c|}{\footnotesize \bf (Figure 5.6 is on next page)}  \\
\hline
\hline
\multicolumn{2}{|c|}{{\color{violet} $\bm{\frac{4}{3}\,\frac{a^2}{4} < b \le \frac{3}{2}\,\frac{a^2}{4}\, , \quad a > 0, \quad c > 0,$}}} \\
\multicolumn{2}{|c|}{{\color{violet} $\bm{0 < {\color{red}c} < c_2 < c_0 < c_1}$}} \\
\hline
\hline
\includegraphics[width=80mm,height=80mm]{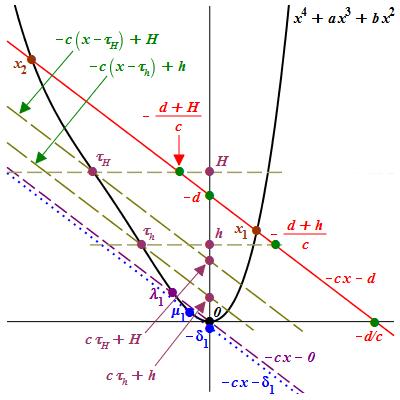} & \includegraphics[width=80mm,height=80mm]{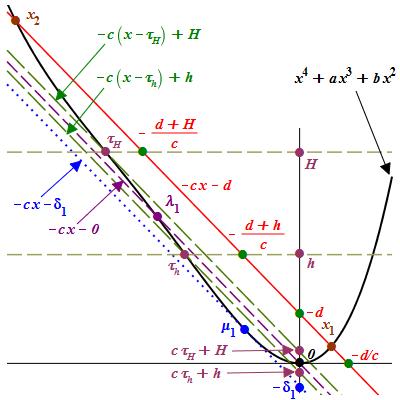} \\
\hline
\multicolumn{1}{|c|}{\begin{minipage}{22em}
\scriptsize
\begin{center}
\vskip-7.67cm
\underline{{\color{red}{\it Notes (apply to all panes)}}}
\end{center}
\vskip-.1cm
As $c < c_2 < c_1$, the quartic has a single stationary point $\mu_1$ and the number of real roots can be either 0 or 2. There is only one tangent $- c x - \delta_1$ to $x^4 + ax^3 + bx^2$ --- the one at $\mu_1$.
\vskip.1cm
The straight line $-cx$ intersects $x^4 + ax^3 + bx^2$ at one point only ($\lambda_1$).
\vskip.1cm
The straight line joining the two points of curvature change ($\tau_h$ and $\tau_H$) has the same slope $-c_0 = -(1/2)a(b-a^2/4) < 0$ as that of the straight line $ - c_0 x - d_0$, tangent to $x^4 + ax^3 + bx^2$ at two points: $\alpha$ and $\beta = -\alpha - a/2$, where $\alpha$ and $\beta$ are the simultaneous double roots of the varied quartic $x^4 + ax^3 + bx^2 + c_0 x + d_0$ [with $d_0 = (1/4)(b-a^2/4)^2 > 0$].
\vskip.1cm
As $\, 0 < c < c_2 < c_0$, one can have $-\delta_1 < 0 < c \tau_h + h  < c \tau_H + H$ (top-left pane, when $c$ is closer to 0), or $-\delta_1 < c \tau_h + h < 0 < c \tau_H + H$ (top-right pane), or $-\delta_1 < c \tau_h + h < c \tau_H + H < 0$ (bottom-right pane, when $c$ is closer to $c_2$).
\vskip.1cm
Consideration of whether $- d \le c \tau_h + h$ or $-d > c \tau_h + h$ and also whether $- d \le c \tau_H + H$ or $-d > c \tau_H + H$ would yield sharper bounds in the analysis based on solving cubic equations.
\end{minipage}}
 & \includegraphics[width=80mm,height=80mm]{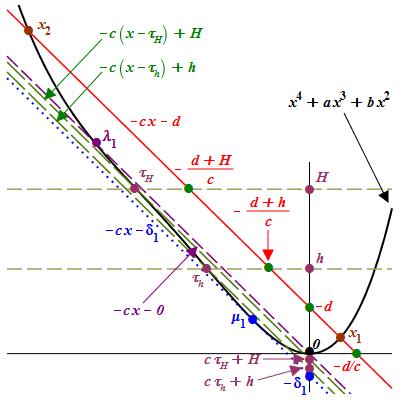} \\
\hline
\multicolumn{2}{|c|}{\begin{minipage}{45em}
\scriptsize
\vskip.1cm
\begin{center}
\underline{{\color{red}{\it Analysis based on solving cubic equations (applies to all panes)}}}
\end{center}
\vskip-.3cm
{\bf (i)} If $-d < -\delta_1$, then there are no real roots.
\vskip.1cm
{\bf (ii)} If $-\delta_1 \le -d < 0$, then there is one negative root $x_1$ such that $\mu_1 \le x_1 < -d/c$ and another negative root $x_2$ such that $\lambda_1 < x_2 \le \mu_1$.
\vskip.1cm
{\bf (iii)} If $0 \le -d$ (pictured), then there is one non-negative root $x_1 \le -d/c$ and one negative root $x_2 \le \lambda_1.$
\vskip.1cm
\phantom{e}
\vskip-.3cm
\hrule
\begin{center}
\vskip-.3cm
\underline{{\color{red}{\it Analysis based on solving quadratic equations only (clockwise from the top-left pane, with four possibile situations per pane)}}}
\end{center}
\vskip-.3cm
{\bf (i)} If $- d < 0$, then there are either no real roots or there are two negative roots smaller than $-d/c$ and greater than $\tau_h$.
\vskip.1cm
{\bf (ii)}  If $0 \le - d < c \tau_h + h$, then there is one non-negative root $x_1 \le -d/c$ and one negative root greater than $\tau_h$.
\vskip.1cm
{\bf (iii)} If $c \tau_h + h \le - d < c \tau_H + H$, then there is one positive root smaller than $-d/c$ and one negative root smaller than or equal to $\tau_h$ and greater than $\tau_H$.
\vskip.1cm
{\bf (iv)} If $c \tau_H + H \le -d$ (pictured), then there is one positive root $x_1 < -d/c$ and one negative root $x_2 \le \tau_H$.
\vskip.1cm
\phantom{e}
\vskip-.3cm
\hrule
\vskip.2cm
{\bf (v)} If $- d < c \tau_h + h$, then there are either no real roots or there are two negative roots smaller than $-d/c$ and greater than $\tau_h$.
\vskip.1cm
{\bf (vi)} If $c \tau_h + h \le - d < 0$, then there is a negative root $x_1$ such that $-(d+h)/c < x_1 < -d/c$ and another negative root smaller than or equal to $\tau_h$ and greater than $\tau_H$.
\vskip.1cm
{\bf (vii)} If $0 \le - d < c \tau_H + H$, then there is a non-negative root $x_1 \le -d/c$ and a negative root $x_2$ such that $\tau_H < x_2 < \tau_h$.
\vskip.1cm
{\bf (viii)} If $c \tau_H + H \le -d$ (pictured), then there is one positive root $x_1 < -d/c$ and one negative root $x_2 \le \tau_H$.
\vskip.1cm
\phantom{e}
\vskip-.3cm
\hrule
\vskip.2cm
{\bf (ix)} If $- d < c \tau_h + h$, then there are either no real roots or there are two negative roots smaller than $-d/c$ and greater than $\tau_h$.
\vskip.1cm
{\bf (x)}  If $c \tau_h + h \le - d < c \tau_H + H$, then there is a negative root $x_1$ such that $-(d+h)/c < x_1 < -d/c$ and another negative root smaller than or equal to $\tau_h$ and greater than $\tau_H$.
\vskip.1cm
{\bf (xi)} If $c \tau_H + H \le - d < 0$, then there is a negative root smaller than $-d/c$ and greater than $-(d+h)/c$ and another negative root smaller than or equal to $\tau_H$.
\vskip.1cm
{\bf (xii)} If $0 \le -d$ (pictured), then there is one non-negative root $x_1 \le -d/c$ and one negative root $x_2 < \tau_H$.
\vskip.1cm
\phantom{e}
\vskip-.3cm
\end{minipage}}\\
\hline
\end{tabular}
\end{center}


\begin{center}
\begin{tabular}{|c||c|}
\hline
{\footnotesize \bf Figure 5.6} &  {\footnotesize \bf Figure 5.8} \\
\hline
\hline
{\color{violet} $\bm{\frac{4}{3}\,\frac{a^2}{4} < b \le \frac{3}{2}\,\frac{a^2}{4}\, , \quad a > 0, \quad c < 0,$}} &
{\color{violet} $\bm{\frac{4}{3}\,\frac{a^2}{4} < b \le \frac{3}{2}\,\frac{a^2}{4}\, , \quad a > 0, \quad c > 0,$}} \\
{\color{violet} $\bm{{\color{red}c} < 0 < c_2  < c_0 < c_1$}} & {\color{violet} $\bm{0 < c_2 < {\color{red}c} < c_0 < c_1$}} \\
\hline
\hline
\includegraphics[width=80mm]{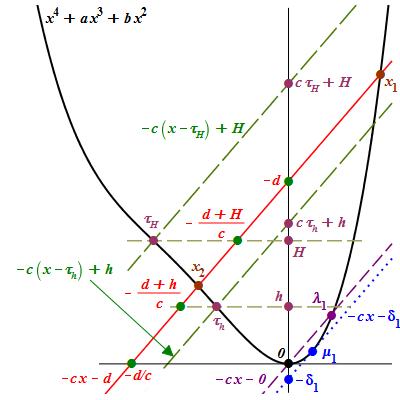} & \includegraphics[width=80mm]{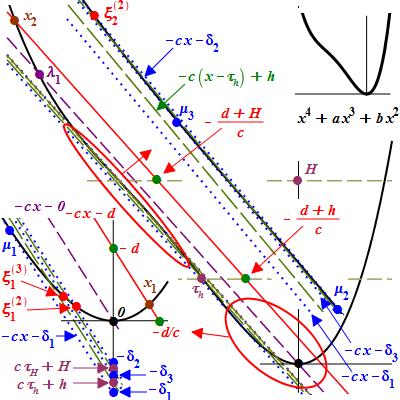} \\
\hline
\multicolumn{1}{|c||}{\begin{minipage}{22em}
\scriptsize
\vskip-7.15cm
\begin{center}
\underline{{\color{red}{\it Notes}}}
\end{center}
\vskip-.4cm
As $c < c_2 < c_1$, the quartic has a single stationary point $\mu_1$ and the number of real roots can be either 0 or 2. There is only one tangent $- c x - \delta_1$ to $x^4 + ax^3 + bx^2$ --- the one at $\mu_1$.
\vskip.1cm
The straight line $-cx$ intersects $x^4 + ax^3 + bx^2$ at one point only ($\lambda_1$).
\vskip.1cm
Obviously, $c \tau_H + H > c \tau_h + h > 0 > -\delta_1$.
\vskip.1cm
Consideration of whether $- d \le c \tau_h + h$ or $-d > c \tau_h + h$ and also whether $- d \le c \tau_H + H$ or $-d > c \tau_H + H$ would yield sharper bounds in the analysis based on solving cubic equations. \\
\hrule
\vskip.1cm
\hskip1.54cm \underline{{\color{red}{\it Analysis based on solving cubic equations}}}
\vskip.1cm
{\bf (i)} If $-d < -\delta_1$, then there are no real roots.
\vskip.1cm
{\bf (ii)} If $-\delta_1 \le -d < 0$, then there is one positive root $x_2$ such that $-d/c < x_2 \le \mu_1$ and another positive root $x_1$ such that $\mu_1 \le x_1 < \lambda_1$.
\vskip.1cm
{\bf (iii)}  If $0 \le - d$ (pictured), then there is a non-positive root $x_2 \ge -d/c$ and a positive root $x_1 \ge \lambda_1$. \\
\hrule
\vskip.1cm
\hskip.88cm \underline{{\color{red}{\it Analysis based on solving quadratic equations only}}}
\vskip.1cm
{\bf (i)} If $- d < 0$, then there are either no real roots or there are two positive roots greater than $-d/c$.
\vskip.1cm
{\bf (ii)} If $0 \le - d < c \tau_h + h$, then there is a non-negative root $x_1$ and a non-positive root $x_2$ such that max$\{ -d/c, \tau_h \} < x_2 \le $ min$\{0, -(d+h)/c \}$.
\vskip.1cm
{\bf (iii)} If $c \tau_h + h \le - d < c \tau_H + H$ (pictured), then there is one positive root $x_1$ and one negative root $x_2$ such that max$\{ -(d+h)/c, \tau_H \} < x_2 \le $ min$\{\tau_h, -(d+H)/c \}$.
\vskip.1cm
{\bf (iv)} If $c \tau_H + H \le -d$, then there is a positive root $x_1$ and a negative root $x_2$ such that  $-(d+H)/c <  x_2 \le \tau_H$. \\
\end{minipage}}
& \multicolumn{1}{|c|}{\begin{minipage}{22em}
\scriptsize
\vskip.05cm
\hskip3.89cm \underline{{\color{red}{\it Notes}}}
\vskip-.05cm
As $c_2 < c < c_1$, the quartic has three stationary points $\mu_i$ and the number of real roots can be either 0, or 2, or 4. There are three tangents $- c x - \delta_i$ to $x^4 + ax^3 + bx^2$ --- each at $\mu_i$.
\vskip.1cm
The straight line $-cx$ intersects $x^4 + ax^3 + bx^2$ at one point only ($\lambda_1$). Thus $-\delta_1 < 0$ and $-\delta_2 < 0$, and $-\delta_3 < 0$.
\vskip.1cm
The straight line joining the two points of curvature change ($\tau_h$ and $\tau_H$) has the same slope $-c_0 = -(1/2)a(b-a^2/4) < 0$ as that of the straight line $ - c_0 x - d_0$, tangent to $x^4 + ax^3 + bx^2$ at two points: $\alpha$ and $\beta = -\alpha - a/2$, where $\alpha$ and $\beta$ are the simultaneous double roots of the varied quartic $x^4 + ax^3 + bx^2 + c_0 x + d_0$ [with $d_0 = (1/4)(b-a^2/4)^2 > 0$].
\vskip.1cm
As $c_0 > c > 0$, one has $-\delta_1 < -\delta_3$  and also $c \sigma_h + h < c \sigma_H + H$.
\vskip.1cm
One can have either the pictured $-\delta_1 < c \sigma_h + h < -\delta_3 < c \sigma_H + H < -\delta_2 < 0$ (when $c$ is closer to $c_2$ where $\mu_2$ and $\mu_3$ coalesce) or $-\delta_1 < -\delta_3 < c \sigma_h + h < c \sigma_H + H < -\delta_2 < 0$ (when $c$ is closer to $c_0$ where $c \tau_h + h$ and $c \tau_h + H$ swap around).
\vskip.1cm
Consideration of whether $- d \le c \tau_h + h$ or $-d > c \tau_h + h$ and also whether $- d \le c \tau_H + H$ or $-d > c \tau_H + H$ would yield sharper bounds in the analysis based on solving cubic equations.
\vskip.1cm
\phantom{e}
\vskip-.3cm
\hrule
\vskip.05cm
\hskip1.54cm \underline{{\color{red}{\it Analysis based on solving cubic equations}}}
\vskip.01cm
{\bf (i)} If $-d < -\delta_1$, then there are no real roots.
\vskip.1cm
{\bf (ii)} If $-\delta_1 \le -d < -\delta_3$, then there is one negative root $x_1$ such that $\mu_1 \le x_1 < $ min$\{ \xi^{(3)}_1, -d/c \}$ and another negative root $x_2$ such that $\xi^{(3)}_2 < x_2 \le \mu_1$.
\vskip.1cm
{\bf (iii)} If $-\delta_3 \le -d < -\delta_2$, then there is a negative root $x_1$ such that $\xi^{(3)}_1 \le x_1 < $ min$\{-d/c, \xi^{(2)}_1 \}$, another negative root $x_2$ such that $\mu_2 < x_2 \le \xi^{(3)}_2$, a third negative root $x_3$ such that $\mu_3 \le x_3 < \mu_2$, and a fourth negative root $x_4$ such that $\xi^{(2)}_2 < x_1 \le \mu_3.$
\vskip.1cm
{\bf (iv)}  If $-\delta_2 \le -d < 0$, then there is a negative root $x_1$ such that $\xi^{(2)}_1 \le x_1 < -d/c$ and another negative root $x_2$ such that $\lambda_1 < x_2 \le \xi^{(2)}_2$.
\vskip.1cm
{\bf (v)}  If $0 \le - d$ (pictured), then there is a non-negative root $x_1 \le -d/c$ and a negative root $x_2 \le \lambda_1$.
\vskip.1cm
\hrule
\vskip.05cm
\hskip.88cm \underline{{\color{red}{\it Analysis based on solving quadratic equations only}}}
\vskip.01cm
{\bf (i)} If $-d < c \tau_h + h$, then there are either no real roots, or there are two negative roots smaller than $-d/c$ and greater than $\tau_h$, or there are two negative roots smaller than $-d/c$ and greater than $\tau_h$, together with two negative roots smaller than $-(d+H)/c$ (the latter appear when $c \tau_h + h > -\delta_3$).
\vskip.05cm
{\bf (ii)} If $c \tau_h + h \le - d < c \tau_H + H$, then there is a negative root $x_1$ such that $-(d+h)/c < x_1 < -d/c$, a negative root $x_2$ such that $\tau_H < x_2 \le \tau_h$, and either zero or two negative roots smaller than $-(d+H)/c$ (the latter are always present of $c\tau_h + h \ge -\delta_3$).
\vskip.1cm
{\bf (iii)} If $c \tau_H + H \le - d < 0$, then there is a negative root $x_1$ such that $-(d+h)/c < x_1 < -d/c$, a negative root $x_2 < \tau_H$, and either zero or two negative roots greater than or equal to $-(d+H)/c$ and smaller than $\tau_h$.
\vskip.1cm
{\bf (iv)} If $0 \le -d$ (pictured), then there is one non-negative root $x_1 \le -d/c$ and one negative root $x_2 > \tau_H$.
\vskip.1cm
\end{minipage}}  \\
\hline
\end{tabular}
\end{center}


\begin{center}
\begin{tabular}{|c||c|}
\hline
{\footnotesize \bf Figure 5.9} &  {\footnotesize \bf Figure 5.10} \\
\hline
\hline
{\color{violet} $\bm{\frac{4}{3}\,\frac{a^2}{4} < b \le \frac{3}{2}\,\frac{a^2}{4}\, , \quad a > 0, \quad c > 0,$}} &
{\color{violet} $\bm{\frac{4}{3}\,\frac{a^2}{4} < b \le \frac{3}{2}\,\frac{a^2}{4}\, , \quad a > 0, \quad c > 0,$}} \\
{\color{violet} $\bm{0 < c_2 < c_0 < {\color{red}c} < c_1$}} & {\color{violet} $\bm{0 < c_2 < c_0 < c_1 < {\color{red}c}$}} \\
\hline
\hline
\includegraphics[width=80mm,height=74mm]{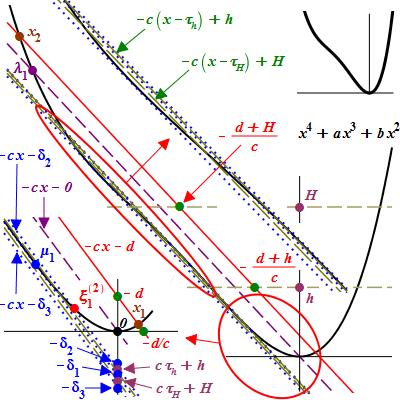} & \includegraphics[width=80mm,height=74mm]{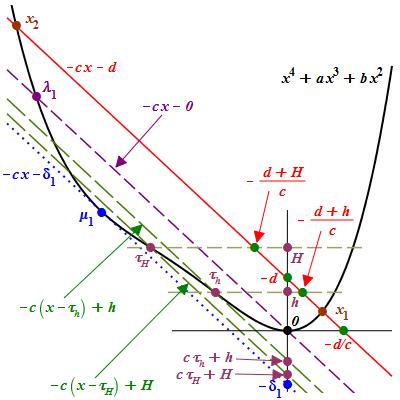} \\
\hline
\multicolumn{1}{|c||}{\begin{minipage}{22em}
\scriptsize
\vskip.1cm
\hskip3.89cm \underline{{\color{red}{\it Notes}}}
\vskip.05cm
As $c_2 < c < c_1$, the quartic has three stationary points $\mu_i$ and the number of real roots can be either 0, or 2, or 4. There are three tangents $- c x - \delta_i$ to $x^4 + ax^3 + bx^2$ --- each at $\mu_i$.
\vskip.1cm
The straight line $-cx$ intersects $x^4 + ax^3 + bx^2$ at one point only ($\lambda_1$). Thus $-\delta_1 < 0$ and $-\delta_2 < 0$, and $-\delta_3 < 0$.
\vskip.1cm
The straight line joining the two points of curvature change ($\tau_h$ and $\tau_H$) has the same slope $-c_0 = -(1/2)a(b-a^2/4) < 0$ as that of the straight line $ - c_0 x - d_0$, tangent to $x^4 + ax^3 + bx^2$ at two points: $\alpha$ and $\beta = -\alpha - a/2$, where $\alpha$ and $\beta$ are the simultaneous double roots of the varied quartic $x^4 + ax^3 + bx^2 + c_0 x + d_0$ [with $d_0 = (1/4)(b-a^2/4)^2 > 0$].
\vskip.1cm
As $c > c_0 > 0$, one has $-\delta_3 < -\delta_1$ and also $c \sigma_H + H < c \sigma_h + h$.
\vskip.1cm
One can have either the pictured $-\delta_3 < c \sigma_H + H < -\delta_1 < c \sigma_h + h < -\delta_2 < 0$ (when $c$ is closer to $c_1$ where $\mu_1$ and $\mu_2$ coalesce) or $-\delta_3 < -\delta_1 < c \sigma_H + H < c \sigma_h + h < -\delta_2 < 0$ (when $c$ is closer to $c_0$ where $c \tau_h + h$ and $c \tau_h + H$ swap around).
\vskip.1cm
Consideration of whether $- d \le c \tau_h + h$ or $-d > c \tau_h + h$ and also whether $- d \le c \tau_H + H$ or $-d > c \tau_H + H$ would yield sharper bounds in the analysis based on solving cubic equations.
\vskip.1cm
\phantom{e}
\vskip-.3cm
\hrule
\vskip.1cm
\hskip1.54cm \underline{{\color{red}{\it Analysis based on solving cubic equations}}}
\vskip.05cm
{\bf (i)} If $-d < -\delta_3$, then there are no real roots.
\vskip.1cm
{\bf (ii)} If $-\delta_3 \le -d < -\delta_1$, then there is one negative root $x_1$ such that $\mu_3 \le x_1 < $ min$\{ \xi^{(1)}_1, -d/c \}$ and another negative root $x_2$ such that $\xi^{(1)}_2 < x_2 \le \mu_3$.
\vskip.1cm
{\bf (iii)} If $-\delta_1 \le -d < -\delta_2$, then there is a negative root $x_1$ such that $\mu_1 \le x_1 < $ min$\{-d/c, \xi^{(2)}_1 \}$, another negative root $x_2$ such that $\mu_2 < x_2 \le \mu_1$, a third negative root $x_3$ such that $\xi^{(1)}_1 \le x_3 < \mu_2$, and a fourth negative root $x_4$ such that $\xi^{(2)}_2 < x_1 \le \xi^{(1)}_2.$
\vskip.1cm
{\bf (iv)}  If $-\delta_2 \le -d < 0$, then there is a negative root $x_1$ such that $\xi^{(2)}_1 \le x_1 < -d/c$ and another negative root $x_2$ such that $\lambda_1 < x_2 \le \xi^{(2)}_2$.
\vskip.1cm
{\bf (v)}  If $0 \le - d$ (pictured), then there is a non-negative root $x_1 \le -d/c$ and a negative root $x_2 \le \lambda_1$.
\vskip.1cm
\phantom{e}
\vskip-.3cm
\hrule
\vskip.1cm
\hskip.88cm \underline{{\color{red}{\it Analysis based on solving quadratic equations only}}}
\vskip.05cm
{\bf (i)} If $- d < c \tau_H + H$, then there are either no real roots, or there are two negative roots smaller than $-(d+H)/c$, or there are two negative roots smaller than $-(d+H)/c$\, together with two negative roots greater than $\tau_h$ and smaller than $-d/c$ (the latter appear when $c \tau_H +H  > -\delta_1$).
\vskip.1cm
{\bf (ii)} If $c \tau_H + H \le - d < c \tau_h + h$, then there are two negative roots: $x_1$ such that $-(d+H)/c \le x_1 < -(d+h)/c$, and $x_2 < \tau_H$ with either zero or two negative roots $x_{3,4}$ such that $\tau_h < x_{3,4} < -d/c$ ($x_{3,4}$ are always present if $c \tau_H + H \ge -\delta_1$).
\vskip.1cm
{\bf (iii)} If $c \tau_h + h \le - d < 0$, then there is a negative root $x_1$ such that $-(d+h)/c \le x_1 < -d/c$, another negative root $x_2 < \tau_H$, and either zero or two negative roots smaller than or equal to $\tau_h$ and greater than $-(d+H)/c$.
\vskip.1cm
{\bf (iv)} If $0 \le -d$ (pictured), then there is a non-negative root $x_1 \le -d/c$ and a negative root $x_2 < \tau_H$.
\\
\end{minipage}}
& \multicolumn{1}{|c|}{\begin{minipage}{22em}
\scriptsize
\vskip-6.2cm
\hskip3.89cm \underline{{\color{red}{\it Notes}}}
\vskip.1cm
As $c_2 < c_1 < c$, the quartic has a single stationary point $\mu_1$ and the number of real roots can be either 0 or 2. There is only one tangent $- c x - \delta_1$ to $x^4 + ax^3 + bx^2$ --- the one at $\mu_1$.
\vskip.1cm
The straight line $-cx$ intersects $x^4 + ax^3 + bx^2$ at one point only ($\lambda_1$).
\vskip.1cm
The straight line joining the two points of curvature change ($\tau_h$ and $\tau_H$) has the same slope $-c_0 = -(1/2)a(b-a^2/4) < 0$ as that of the straight line $ - c_0 x - d_0$, tangent to $x^4 + ax^3 + bx^2$ at two points: $\alpha$ and $\beta = -\alpha - a/2$, where $\alpha$ and $\beta$ are the simultaneous double roots of the varied quartic $x^4 + ax^3 + bx^2 + c_0 x + d_0$ [with $d_0 = (1/4)(b-a^2/4)^2 > 0$].
\vskip.1cm
As $c > c_0 > 0$, one has $-\delta_1 < c \sigma_H + H < c \sigma_h + h < 0$.
\vskip.1cm
Consideration of whether $- d \le c \tau_h + h$ or $-d > c \tau_h + h$ and also whether $- d \le c \tau_H + H$ or $-d > c \tau_H + H$ would yield sharper bounds in the analysis based on solving cubic equations.
\\
\hrule
\vskip.1cm
\hskip1.54cm \underline{{\color{red}{\it Analysis based on solving cubic equations}}}
\vskip.1cm
{\bf (i)} If $-d < -\delta_1$, then there are no real roots.
\vskip.1cm
{\bf (ii)} If $-\delta_1 \le -d < 0$, then there is one negative root $x_1$ such that $\mu_1 \le x_1 < -d/c$ and another negative root $x_2$ such that $\lambda_1 < x_2 \le \mu_1$.
\vskip.1cm
{\bf (iii)}  If $0 \le - d$ (pictured), then there is a non-negative root $x_1 \le -d/c$ and a negative root $x_2 \le \lambda_1$. \\
\hrule
\vskip.1cm
\hskip.88cm \underline{{\color{red}{\it Analysis based on solving quadratic equations only}}}
\vskip.1cm
{\bf (i)} If $- d < c \tau_H + H$, then there are either no real roots, or there are two negative roots $x_{1,2} < -(d+H)/c$.
\vskip.1cm
{\bf (ii)} If $c \tau_H + H \le - d < c \tau_h + h$, then there are two negative roots: $x_1$ such that $-(d+H)/c \le x_1 < -(d+h)/c$ and $x_2 <  \tau_H$.
\vskip.1cm
{\bf (iii)} If $c \tau_h + h \le - d < 0$, then there are two negative roots: $x_1$ such that $-(d+h)/c \le x_1 < -d/c$ and $x_2 < \tau_H$.
\vskip.1cm
{\bf (iv)} If $0 \le -d$ (pictured), then there is one non-negative root $x_1 \le -d/c$ and one negative root $x_1 < \tau_H$.
\\
\end{minipage}}  \\
\hline
\end{tabular}
\end{center}


\begin{center}
\begin{tabular}{|c|c|}
\hline
\multicolumn{2}{|c|}{\footnotesize \bf Figure 6.1}  \\
\hline
\hline
\multicolumn{2}{|c|}{{\color{violet} $\bm{b > \frac{3}{2}\,\frac{a^2}{4}\, , \quad a < 0, \quad c < 0$}}} \\
\hline
\hline
\includegraphics[width=80mm,height=80mm]{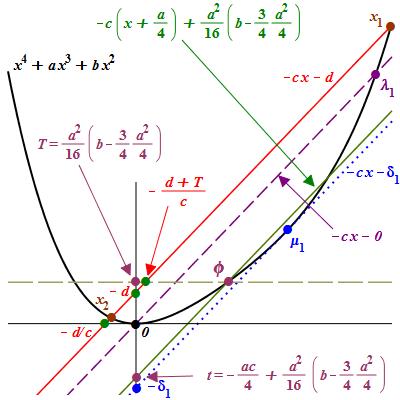} & \includegraphics[width=80mm,height=80mm]{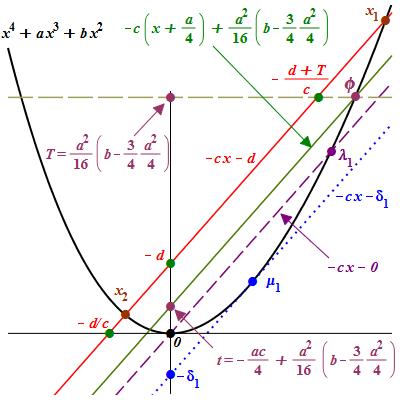} \\
\hline
\multicolumn{1}{|c|}{\begin{minipage}{22em}
\scriptsize
\vskip-7.3cm
\begin{center}
\underline{{\color{red}{\it Notes (apply to all panes)}}}
\end{center}
\vskip-.3cm
The quartic has a single stationary point $\mu_{1}$ and the number of real roots can be either 0 or 2. There is only one tangent $- c x - \delta_1$ to $x^4 + ax^3 + bx^2$ --- the one at $\mu_1$. The intersection point of this tangent with the ordinate is $-\delta_1$ and this is always negative.
\vskip.1cm
The straight line $-cx$ intersects $x^4 + ax^3 + bx^2$ only at $\lambda_1 > \mu_1$.
\vskip.1cm
The third derivative of $x^4 + ax^3 + bx^2$ vanishes at $\phi = -a/4$ --- the only non-zero ``marker" point that can be identifed on the graph of $x^4 + ax^3 + bx^2$.
\vskip.1cm
The ``separator" straight line with equation $-c(x+a/4) + (a^2/16)[b - (3/4)a^2/4]$ passes through the ``marker" point  $\phi$ and intersects the ordinate at point $t$ which is always greater than or equal to $-\delta_1$. This point $t$ is zero if $c = \zeta \equiv (a/4)[b - (3/4)a^2/4] = a^3/64 + c_0/2$, where $c_0 = (1/2)a(b - a^2/4)$. One has $\zeta < 0$, when $a < 0$. Also: $\zeta - c_0 = -(a/4)[b - (5/4)a^2/4]$ and thus $\zeta > c_0$ for $a < 0$ and $b > (3/2)a^2/4$.
\vskip.1cm
One has $\mu_1 = \phi$ if $c = c_0 = (1/2)a(b - a^2/4)$. If $a < 0$ and $b > (3/2) a^2/4$, then $c_0 < 0$. Hence, at $c = c_0$, one has $t = -\delta_1 = -(a^2/16) [b - (5/4)a^2/4] < 0$, since $b > (3/2) a^2/4$.
\vskip.1cm
If $c < c_0 < 0$, one has $\mu_1 > \phi$ and thus $-\delta_1 < t < 0$ --- pictured on the top-left pane.
\vskip.1cm
If $c_0 < c < 0$, one has $\mu_1 < \phi$. There are two possibilities in this
\end{minipage}}
 & \includegraphics[width=80mm,height=80mm]{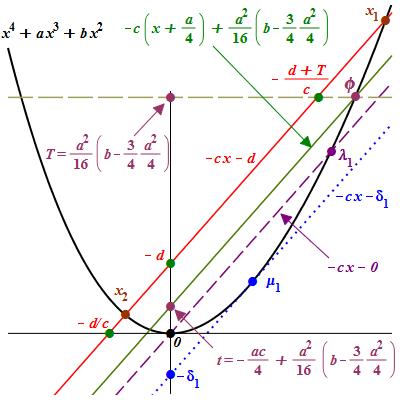} \\
\cline{2-2}
\multicolumn{2}{|c|}{\begin{minipage}{45em}
\scriptsize
\vskip.1cm
case: either the pictured on the top-right pane $-\delta_1 < t \le 0$, which occurs for $c_0 < c \le \zeta < 0$, or the pictured on the bottom-right pane $-\delta_1 < 0 < t$, which occurs for $c_0 < \zeta < c < 0$.
\vskip.1cm
Consideration of whether $- d \le t$ or $-d > t$ would yield sharper bounds in the analysis based on solving cubic equations.
\\
\hrule
\begin{center}
\vskip-.2cm
\underline{{\color{red}{\it Analysis based on solving cubic equations (applies to all panes)}}}
\end{center}
\vskip-.2cm
{\bf (i)} If $-d < -\delta_1$, then there are no real roots.
\vskip.1cm
{\bf (ii)} If $-\delta_1 \le -d < 0$, then there is one positive root $x_2$ such that $-d/c < x_2 \le \mu_1$ and another positive root $x_1$ such that $\mu_1 \le x_1 < \lambda_1$.
\vskip.1cm
{\bf (iii)}  If $0 \le - d$ (pictured), then there is a non-positive root $x_2 \ge -d/c$ and a positive root $x_1 \ge \lambda_1$.
\\
\hrule
\begin{center}
\vskip-.2cm
\underline{{\color{red}{\it Analysis based on solving quadratic equations only}}}
\end{center}
\vskip-.1cm
\underline{{\color{black} Top-left pane: $c < c_0 < \zeta < 0$, hence $\mu_1 > \phi$ and $-\delta_1 < t < 0$}}
\vskip.2cm
{\bf (i)} If $- d < t$, then there are either no real roots, or there are two positive roots $x_{1,2} > -(d+T)/c$.
\vskip.1cm
{\bf (ii)} If $t \le -d < 0$, then there are two positive roots: $x_2$ such that $-d/c < x_2 \le -(d+T)/c$ and $x_1 > \phi = -a/4$.
\vskip.1cm
{\bf (iii)} If $0 \le -d$ (pictured), then there is one non-positive root $x_2 \ge -d/c$ and one positive root $x_1 > \phi = -a/4$.
\vskip.3cm
\underline{{\color{black} Top-right pane: $c_0 < c \le \zeta < 0$, hence $\mu_1 < \phi$ and $-\delta_1 < t \le 0$}}
\vskip.2cm
{\bf (i)} If $- d < t$, then there are either no real roots, or there are two positive roots $x_{1,2}$ such that $-d/c < x_{1,2} < \phi = -a/4$.
\vskip.1cm
{\bf (ii)} If $t \le -d < 0$, then there is a positive root $x_2$ such that $-d/c < x_2 \le -(d+T)/c$ and a positive root $x_1 \ge \phi = -a/4$.
\vskip.1cm
{\bf (iii)} If $0 \le -d$ (pictured), then there is one non-positive root $x_2 \ge -d/c$  and one positive root $x_1 > \phi = -a/4$.
\vskip.3cm
\underline{{\color{black} Bottom-right pane: $c_0 < \zeta < c < 0$, hence $\mu_1 < \phi$ and $-\delta_1 < 0 < t$}}
\vskip.2cm
{\bf (i)} If $- d < 0$, then there are either no real roots, or there are two positive roots $x_{1,2}$ such that $-d/c < x_{1,2} < \phi = -a/4$.
\vskip.1cm
{\bf (ii)} If $0 \le -d < t$, then there is one non-positive root $x_2 \ge -d/c$ and a positive root $x_1 < \phi = -a/4$.
\vskip.1cm
{\bf (iii)} If $t \le -d$ (pictured), then there is one negative root $x_2 > -d/c$  and one positive root $x_1 \ge \phi = -a/4$. \\
\end{minipage}}\\
\hline
\end{tabular}
\end{center}


\begin{center}
\begin{tabular}{|c||c|}
\hline
{\footnotesize \bf Figure 6.2} &  {\footnotesize \bf Figure 6.3} \\
\hline
\hline
{\color{violet} $\bm{b > \frac{3}{2}\,\frac{a^2}{4}\, , \quad a < 0, \quad c > 0$}} &
{\color{violet} $\bm{b > \frac{3}{2}\,\frac{a^2}{4}\, , \quad a > 0, \quad c < 0$}} \\
\hline
\hline
\includegraphics[width=80mm, height=80mm]{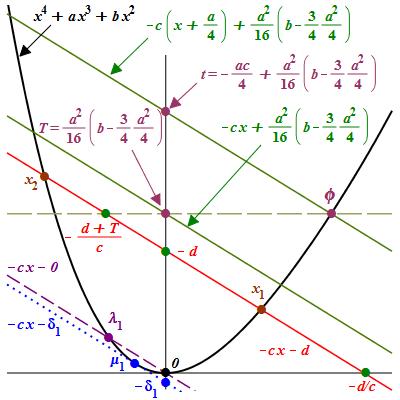} & \includegraphics[width=80mm, height=80mm]{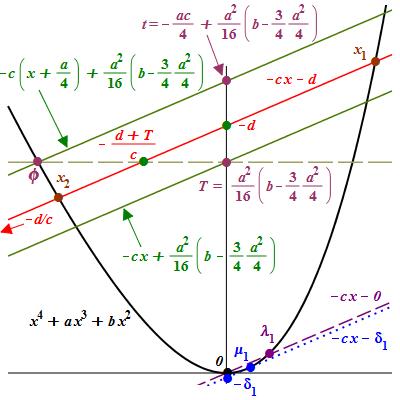} \\
\hline
\multicolumn{1}{|c||}{\begin{minipage}{22em}
\scriptsize
\vskip.2cm
\hskip3.89cm \underline{{\color{red}{\it Notes}}}
\vskip.2cm
The quartic has a single stationary point $\mu_{1}$ and the number of real roots can be either 0 or 2. There is only one tangent $- c x - \delta_1$ to $x^4 + ax^3 + bx^2$ --- the one at $\mu_1$. The intersection point of this tangent with the ordinate is $-\delta_1$ and this is always negative.
\vskip.1cm
The straight line $-cx$ intersects $x^4 + ax^3 + bx^2$ only at $\lambda_1 < \mu_1$.
\vskip.1cm
The third derivative of $x^4 + ax^3 + bx^2$ vanishes at $\phi = -a/4$ --- the only non-zero ``marker" point that can be identifed on the graph of $x^4 + ax^3 + bx^2$.
\vskip.1cm
The ``separator" straight line with equation $-c(x+a/4) + (a^2/16)[b - (3/4)a^2/4]$ passes through the ``marker" point  $\phi$ and intersects the ordinate at point $t > 0$.
\vskip.1cm
The additional ``separator" straight line with equation $- cx + (a^2/16) [b - (3/4)a^2/4]$ provides sharper bounds for the analysis. It intersects the ordinate at point $T = t + ac/4 < t$, since $a < 0$ and $c > 0$.
\vskip.1cm
Thus $-\delta_1 < 0 < T < t$.
\vskip.1cm
Consideration of whether $- d \le t$ or $-d > t$  and also whether $- d \le T$ or $-d > T$ would yield sharper bounds in the analysis based on solving cubic equations. \\
\hrule
\vskip.2cm
\hskip1.54cm \underline{{\color{red}{\it Analysis based on solving cubic equations}}}
\vskip.2cm
{\bf (i)} If $-d < -\delta_1$, then there are no real roots.
\vskip.1cm
{\bf (ii)} If $-\delta_1 \le -d < 0$, then there is one negative root $x_1$ such that $\mu_1 \le x_1 < -d/c$ and another negative root $x_2$ such that $\lambda_1 < x_2 \le \mu_1$.
\vskip.1cm
{\bf (iii)}  If $0 \le - d$ (pictured), then there is a non-negative root $x_1 \le -d/c$ and a negative root $x_2 \le \lambda_1$.   \\
\hrule
\end{minipage}}
&
\multicolumn{1}{|c|}{\begin{minipage}{22em}
\scriptsize
\vskip.2cm
\hskip3.89cm \underline{{\color{red}{\it Notes}}}
\vskip.2cm
The quartic has a single stationary point $\mu_{1}$ and the number of real roots can be either 0 or 2. There is only one tangent $- c x - \delta_1$ to $x^4 + ax^3 + bx^2$ --- the one at $\mu_1$. The intersection point of this tangent with the ordinate is $-\delta_1$ and this is always negative.
\vskip.1cm
The straight line $-cx$ intersects $x^4 + ax^3 + bx^2$ only at $\lambda_1 > \mu_1$.
\vskip.1cm
The third derivative of $x^4 + ax^3 + bx^2$ vanishes at $\phi = -a/4$ --- the only non-zero ``marker" point that can be identifed on the graph of $x^4 + ax^3 + bx^2$.
\vskip.1cm
The ``separator" straight line with equation $-c(x+a/4) + (a^2/16)[b - (3/4)a^2/4]$ passes through the ``marker" point  $\phi$ and intersects the ordinate at point $t > 0$.
\vskip.1cm
The additional ``separator" straight line with equation $- cx + (a^2/16) [b - (3/4)a^2/4]$ provides sharper bounds for the analysis. It intersects the ordinate at point $T = t + ac/4 < t$, since $a > 0$ and $c < 0$.
\vskip.1cm
Thus $-\delta_1 < 0 < T < t$.
\vskip.1cm
Consideration of whether $- d \le t$ or $-d > t$  and also whether $- d \le T$ or $-d > T$ would yield sharper bounds in the analysis based on solving cubic equations. \\
\hrule
\vskip.2cm
\hskip1.54cm \underline{{\color{red}{\it Analysis based on solving cubic equations}}}
\vskip.2cm
{\bf (i)} If $-d < -\delta_1$, then there are no real roots.
\vskip.1cm
{\bf (ii)} If $-\delta_1 \le -d < 0$, then there is one positive root $x_2$ such that $-d/c < x_2 \le \mu_1$ and another positive root $x_1$ such that $\mu_1 \le x_1 < \lambda_1$.
\vskip.1cm
{\bf (iii)}  If $0 \le - d$ (pictured), then there is a non-positive root $x_2 \ge -d/c$ and a positive root $x_1 \ge \lambda_1$.  \\
\hrule
\end{minipage}}
\\
\multicolumn{1}{|c||}{\begin{minipage}{22em}
\scriptsize
\vskip.2cm
\hskip.88cm \underline{{\color{red}{\it Analysis based on solving quadratic equations only}}}
\vskip.2cm
{\bf (i)} If $- d < 0$, then there are either no real roots or there are two negative roots $x_{1,2} < -d/c$.
\vskip.1cm
{\bf (ii)} If $0 \le -d < T$ (pictured), then there is one non-negative root $x_1 \le -d/c$ and one negative root $x_2$.
\vskip.1cm
{\bf (iii)} If $T \le -d < t$, then there is one positive root $x_1$ such that $-(d+T)/c < x_1 < \phi = -a/4$ and one negative root $x_2$.
\vskip.1cm
{\bf (iv)} If $t \le -d $, then there is one positive root $x_1$ such that $\phi = -a/4 \le x_1 < -(d+T)/c $ and one negative root $x_2$. \\
\end{minipage}}
&
\multicolumn{1}{|c|}{\begin{minipage}{22em}
\scriptsize
\vskip.2cm
\hskip.88cm \underline{{\color{red}{\it Analysis based on solving quadratic equations only}}}
\vskip.2cm
{\bf (i)} If $- d < 0$, then there are either no real roots or there are two positive roots $x_{1,2} > -d/c$.
\vskip.1cm
{\bf (ii)} If $0 \le -d < T$, then there is one non-positive root $x_2 \ge -d/c$ and one positive root $x_1$.
\vskip.1cm
{\bf (iii)} If $T \le -d < t$ (pictured), then there is one negative root $x_2$ such that $\phi = -a/4 < x_2 < -(d+T)/c$ and one positive root $x_1$.
\vskip.1cm
{\bf (iv)} If $t \le -d $, then there is one negative root $x_2$ such that $-(d+T)/c <  x_2 \le \phi = -a/4$ and one positive root $x_1$.\\
\end{minipage}}\\
\hline
\end{tabular}
\end{center}


\begin{center}
\begin{tabular}{|c|c|}
\hline
\multicolumn{2}{|c|}{\footnotesize \bf Figure 6.4}  \\
\hline
\hline
\multicolumn{2}{|c|}{{\color{violet} $\bm{b > \frac{3}{2}\,\frac{a^2}{4}\, , \quad a > 0, \quad c > 0$}}} \\
\hline
\hline
\includegraphics[width=80mm,height=80mm]{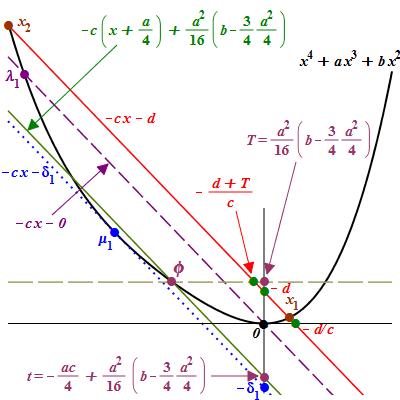} & \includegraphics[width=80mm,height=80mm]{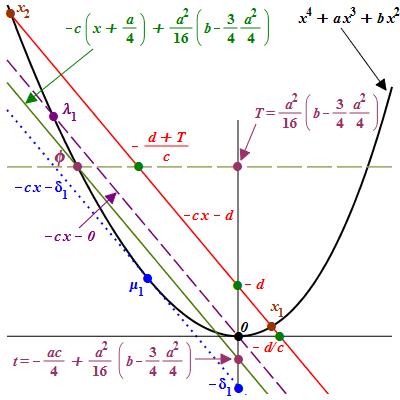} \\
\hline
\multicolumn{1}{|c|}{\begin{minipage}{22em}
\scriptsize
\vskip-7.3cm
\begin{center}
\underline{{\color{red}{\it Notes (apply to all panes)}}}
\end{center}
\vskip-.3cm
The quartic has a single stationary point $\mu_{1}$ and the number of real roots can be either 0 or 2. There is only one tangent $- c x - \delta_1$ to $x^4 + ax^3 + bx^2$ --- the one at $\mu_1$. The intersection point of this tangent with the ordinate is $-\delta_1$ and this is always negative.
\vskip.1cm
The straight line $-cx$ intersects $x^4 + ax^3 + bx^2$ only at $\lambda_1 > \mu_1$.
\vskip.1cm
The third derivative of $x^4 + ax^3 + bx^2$ vanishes at $\phi = -a/4$ --- the only non-zero ``marker" point that can be identifed on the graph of $x^4 + ax^3 + bx^2$.
\vskip.1cm
The ``separator" straight line with equation $-c(x+a/4) + (a^2/16)[b - (3/4)a^2/4]$ passes through the ``marker" point  $\phi$ and intersects the ordinate at point $t$ which is always greater than or equal to $-\delta_1$. This point $t$ is zero if $c = \zeta \equiv (a/4)[b - (3/4)a^2/4] = a^3/64 + c_0/2$, where $c_0 = (1/2)a(b - a^2/4)$. One has $\zeta > 0$, when $a > 0$. Also: $\zeta - c_0 = -(a/4)[b - (5/4)a^2/4]$ and thus $\zeta < c_0$ for $a > 0$ and $b > (3/2)a^2/4$.
\vskip.1cm
One has $\mu_1 = \phi$ if $c = c_0 = (1/2)a(b - a^2/4)$. If $a > 0$ and $b > (3/2) a^2/4$, then $c_0 > 0$. Hence, at $c = c_0$, one has $t = -\delta_1 = -(a^2/16) [b - (5/4)a^2/4] < 0$, since $b > (3/2) a^2/4$.
\vskip.1cm
If $c > c_0 > 0$, one has $\mu_1 < \phi$ and thus $-\delta_1 < t < 0$ --- pictured on the top-left pane.
\vskip.1cm
If $c_0 > c > 0$, one has $\mu_1 > \phi$. There are two possibilities in this
\end{minipage}}
 & \includegraphics[width=80mm,height=80mm]{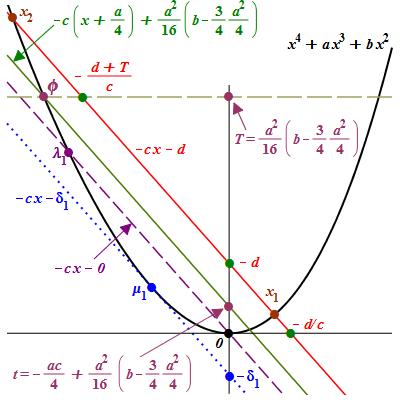} \\
\cline{2-2}
\multicolumn{2}{|c|}{\begin{minipage}{45em}
\scriptsize
\vskip.1cm
case: either the pictured on the top-right pane $-\delta_1 < t \le 0$, which occurs for $c_0 > c \ge \zeta > 0$, or the pictured on the bottom-right pane $-\delta_1 < 0 < t$, which occurs for $c_0  > \zeta > c > 0$.
\vskip.1cm
Consideration of whether $- d \le t$ or $-d > t$ would yield sharper bounds in the analysis based on solving cubic equations.
\\
\hrule
\begin{center}
\vskip-.2cm
\underline{{\color{red}{\it Analysis based on solving cubic equations (applies to all panes)}}}
\end{center}
\vskip-.2cm
{\bf (i)} If $-d < -\delta_1$, then there are no real roots.
\vskip.1cm
{\bf (ii)} If $-\delta_1 \le -d < 0$, then there is one negative root $x_1$ such that $\mu_1 \le x_1 < -d/c$ and another negative root $x_2$ such that $\lambda_1 < x_2 \le \mu_1$.
\vskip.1cm
{\bf (iii)}  If $0 \le - d$ (pictured), then there is a non-negative root $x_1 \le -d/c$ and a negative root $x_2 \le \lambda_1$.
\\
\hrule
\begin{center}
\vskip-.2cm
\underline{{\color{red}{\it Analysis based on solving quadratic equations only}}}
\end{center}
\vskip-.2cm
\underline{{\color{black} Top-left pane: $c > c_0 > \zeta > 0$, hence $\mu_1 < \phi$ and $-\delta_1 < t < 0$}}
\vskip.2cm
{\bf (i)} If $- d < t$, then there are either no real roots, or there are two negative roots $x_{1,2} < -(d+T)/c$.
\vskip.1cm
{\bf (ii)} If $t \le -d < 0$, then there are two negative roots: $x_1$ such that $-(d+T)/c \le x_1 < -d/c$ and $x_2 < \phi = -a/4$.
\vskip.1cm
{\bf (iii)} If $0 \le -d$ (pictured), then there is one non-negative root $x_1 \le -d/c$ and one negative root $x_2 < \phi = -a/4$.
\vskip.3cm
\underline{{\color{black} Top-right pane: $c_0 > c \ge \zeta > 0$, hence $\mu_1 > \phi$ and $-\delta_1 < t \le 0$}}
\vskip.2cm
{\bf (i)} If $- d < t$, then there are either no real roots, or there are two negative roots $x_{1,2}$ such that $\phi = -a/4 < x_{1,2} < -d/c $.
\vskip.1cm
{\bf (ii)} If $t \le -d < 0$, then there is a negative root $x_1$ such that $-(d+T)/c \le x_1 < -d/c $ and a negative root $x_2 \le \phi = -a/4$.
\vskip.1cm
{\bf (iii)} If $0 \le -d$ (pictured), then there is one non-negative root $x_1 \le -d/c$  and one negative root $x_2 < \phi = -a/4$.
\vskip.3cm
\underline{{\color{black} Bottom-right pane: $c_0 > \zeta > c > 0$, hence $\mu_1 > \phi$ and $-\delta_1 < 0 < t$}}
\vskip.2cm
{\bf (i)} If $- d < 0$, then there are either no real roots, or there are two negative roots $x_{1,2}$ such that $\phi = -a/4 < x_{1,2} < -d/c $.
\vskip.1cm
{\bf (ii)} If $0 \le -d < t$, then there is one non-negative root $x_1 \le -d/c$ and a negative root $x_2 > \phi = -a/4$.
\vskip.1cm
{\bf (iii)} If $t \le -d$ (pictured), then there is one positive root $x_1 < -d/c$  and one negative root $x_2 \le \phi = -a/4$.
\\
\end{minipage}}\\
\hline
\end{tabular}
\end{center}


\newpage

\restoregeometry

\section*{Acknowledgements}
It is a pleasure to thank Elena Tonkova for useful discussions and Milena E. Mihaylova for the help with the figures.

\end{document}